\documentclass{amsart}
\usepackage[utf8]{inputenc}
\usepackage{geometry}

\usepackage{amsmath}
\usepackage{amsfonts}
\usepackage{amssymb}
\usepackage{amsthm}
\usepackage{mathrsfs}
\usepackage{fixdif}
\usepackage{ytableau}
\usepackage{subcaption}
\usepackage{bbm}

\usepackage[table]{xcolor}
\usepackage{tikz}
\usetikzlibrary{cd}
\usetikzlibrary{calc}
\definecolor{blue}{HTML}{0072B2}
\definecolor{orange}{HTML}{E69F00}
\definecolor{green}{HTML}{009E73}
\definecolor{purple}{HTML}{CC79A7}

\usepackage[hidelinks]{hyperref}

\usepackage[backend=biber,sorting=nyt]{biblatex}
\bibliography{main}
\renewbibmacro{in:}{}
\DeclareFieldFormat{pages}{#1}
\DeclareFieldFormat[article,unpublished,inproceedings]{title}{#1}

\theoremstyle{plain}
\newtheorem{theorem}{Theorem}[section]
\newtheorem{conjecture}[theorem]{Conjecture}

\newtheorem{proposition}[theorem]{Proposition}
\newtheorem{corollary}[theorem]{Corollary}
\newtheorem{lemma}[theorem]{Lemma}

\theoremstyle{definition}
\newtheorem{definition}[theorem]{Definition}

\newtheorem{example}[theorem]{Example}
\newtheorem{counterexample}[theorem]{Counterexample}

\theoremstyle{remark}
\newtheorem*{remark}{Remark}
\newtheorem*{warning}{Warning}

\newcommand{\bbE}{\mathbb{E}}

\newcommand{\bbQ}{\mathbb{Q}}
\newcommand{\bbR}{\mathbb{R}}

\newcommand{\bbZ}{\mathbb{Z}}

\newcommand{\calE}{\mathcal{E}}

\newcommand{\calP}{\mathcal{P}}

\newcommand{\calS}{\mathcal{S}}
\newcommand{\calT}{\mathcal{T}}

\newcommand{\Z}{\bbZ}
\newcommand{\Q}{\bbQ}
\newcommand{\R}{\bbR}

\newcommand{\E}{\bbE}

\DeclareMathOperator{\len}{len}
\newcommand{\simto}{\xrightarrow{\sim}}

\newcommand{\id}{\mathrm{id}}

\DeclareMathOperator{\typ}{typ}
\newcommand{\brackbinom}{\genfrac{[}{]}{0pt}{}}

\newcommand{\lparallel}{{\backslash\mkern-5mu\backslash}}

\newsavebox{\stairbox}
\savebox{\stairbox}{\tikz[baseline=(stairwayanchor.base)]{
    \node (stairwayanchor) {\quad}; 
    \draw[black] 
        ($(stairwayanchor.south west) + (0.0em,0.25em)$)
        -- ++(0,0.40em)
        -- ++(0.40em,0) -- ++(0,0.40em)
        -- ++(0.40em,0)
        -- ++(0,-0.80em)
        -- cycle
    ;
}}
\newcommand{\stair}{\mathbin{\usebox{\stairbox}}}

\usepackage[automake,toc,sort=use]{glossaries-extra}
\setglossarystyle{long3col}

\newglossary*{symbols}{List of Symbols}
\newglossary*{terms}{List of Terms}
\makeglossaries

\newglossaryentry{symmetric-group}{
    type=symbols,
    name={$S_n$},
    description={The symmetric group on $n$ letters}
}
\newglossaryentry{bracket-n}{
    type=symbols,
    name={$[n]$},
    description={The set $\{1, \dots, n\}$}
}
\newglossaryentry{N-sigma}{
    type=symbols,
    name={$N_\sigma(\pi)$},
    description={The number of pattern occurrences of $\sigma$ in $\pi$}
}
\newglossaryentry{lambda-bracket-n}{
    type=symbols,
    name={$\lambda[n]$},
    description={Partition obtained from $\lambda$ by adding part of size $n-\lvert\lambda\rvert$}
}
\newglossaryentry{chi-lambda}{
    type=symbols,
    name={$\chi^\lambda$},
    description={Irreducible symmetric group character for the partition $\lambda$}
}
\newglossaryentry{chi-rho-lambda}{
    type=symbols,
    name={$\chi_\rho^\lambda$},
    description={Value of $\chi^\lambda$ on the conjugacy class of $\rho$}
}
\newglossaryentry{a-sigma-lambda}{
    type=symbols,
    name={$a_{\sigma_1,\dots,\sigma_d}^\lambda$},
    description={``Permutation pattern character polynomials'' defined by Gaetz and Pierson}
}
\newglossaryentry{id-k}{
    type=symbols,
    name={$\id_k$},
    description={Identity permutation on $k$ letters}
}
\newglossaryentry{m-j}{
    type=symbols,
    name={$m_j$},
    description={Function $\bigsqcup_{n\ge0} S_n \to \Z_{\ge0}$ counting $j$-cycles in its input}
}
\newglossaryentry{m-j-rho}{
    type=symbols,
    name={$m_j(\rho)$},
    description={Number of parts of size $j$ in a partition $\rho$}
}
\newglossaryentry{M-rho-sigma}{
    type=symbols,
    name={$M_\rho^{\sigma_1,\dots,\sigma_d}$},
    description={Mixed moment of $N_{\sigma_1}, \dots, N_{\sigma_d}$ on the conjugacy class of $\rho$}
}
\newglossaryentry{C-rho}{
    type=symbols,
    name={$C_\rho$},
    description={Conjugacy class in $S_{\lvert\rho\rvert}$ corresponding to cycle type $\rho$}
}
\newglossaryentry{abs-rho}{
    type=symbols,
    name={$\lvert\rho\rvert$},
    description={Size of a partition $\rho$}
}
\newglossaryentry{A-sigma-lambda}{
    type=symbols,
    name={$A_{\sigma_1,\dots,\sigma_d}^\lambda$},
    description={Mixed moment of $N_{\sigma_1}, \dots, N_{\sigma_d}$ with $\chi^\lambda$}
}
\newglossaryentry{typ-pi}{
    type=symbols,
    name={$\typ(\pi)$},
    description={Cycle type of a permutation $\pi$}
}
\newglossaryentry{typ-alpha}{
    type=symbols,
    name={$\typ(\alpha)$},
    description={Type of a stitch $\alpha$}
}
\newglossaryentry{lambda-i}{
    type=symbols,
    name={$\lambda_i$},
    description={The $i$th part of a partition $\lambda$}
}
\newglossaryentry{theta-f-i}{
    type=symbols,
    name={$\theta_f^{(i)}$},
    description={Coefficients appearing in Theorems~\ref{main-thm-1} and \ref{E-J-rho-main-thm}}
}
\newglossaryentry{theta-f-i-ell}{
    type=symbols,
    name={$\theta_f^{i,\ell}$},
    description={Coefficients appearing in Theorem~\ref{main-thm-combo}}
}
\newglossaryentry{a-f-k}{
    type=symbols,
    name={$a_f^{(k)}$},
    description={Generalization of the polynomials $a_{\id_k}^\lambda$, appearing in Theorem~\ref{main-thm-combo}}
}
\newglossaryentry{b-f-j}{
    type=symbols,
    name={$b_f^{(j)}$},
    description={Polynomials appearing in Theorems~\ref{main-thm-2} and \ref{main-thm-combo}}
}
\newglossaryentry{E-f-n-k}{
    type=symbols,
    name={$E_f(n,k)$},
    description={Expected value of $f(\pi) \cdot N_{\id_k}(\pi)$ for $\pi \in S_n$}
}
\newglossaryentry{Delta-Q}{
    type=symbols,
    name={$\Delta_\Q$},
    description={Ring of character polynomials $\Q[m_1,m_2,\dots]$}
}
\newglossaryentry{Delta-Q-le-d}{
    type=symbols,
    name={$\Delta_\Q^{\le d}$},
    description={Subspace of character polynomials of degree at most $d$}
}
\newglossaryentry{Delta-Q-d}{
    type=symbols,
    name={$\Delta_\Q^d$},
    description={Subspace of character polynomials which are pure of degree $d$}
}
\newglossaryentry{I-rho}{
    type=symbols,
    name={$I_\rho$},
    description={Binomial basis for $\Delta_\Q$}
}
\newglossaryentry{X-lam}{
    type=symbols,
    name={$X^\lambda$},
    description={Character polynomial basis for $\Delta_\Q$}
}
\newglossaryentry{J-rho}{
    type=symbols,
    name={$J_\rho$},
    description={Basis for $\Delta_\Q$ obtained by orthogonalizing the $I_\rho$}
}
\newglossaryentry{Y-lam}{
    type=symbols,
    name={$Y^\lambda$},
    description={Basis for $\Delta_\Q$ used in the construction of the $X^\lambda$}
}
\newglossaryentry{nth-moment}{
    type=symbols,
    name={$\langle f\rangle_n$},
    description={The $n$th moment of a character polynomial $f$, $\frac{1}{n!}\sum_{\pi \in S_n} f(\pi)$}
}
\newglossaryentry{stable-moment}{
    type=symbols,
    name={$\langle f\rangle$},
    description={The stable moment of a character polynomial $f$, $\lim_{n \to \infty} \langle f\rangle_n$}
}
\newglossaryentry{inner-product}{
    type=symbols,
    name={$\langle f, g\rangle$},
    description={The inner product of character polynomials $f$ and $g$}
}
\newglossaryentry{hall-inner-product}{
    type=symbols,
    name={$\langle \phi, \psi\rangle$},
    description={The Hall inner product of symmetric functions $\phi$ and $\psi$}
}
\newglossaryentry{z-rho}{
    type=symbols,
    name={$z_\rho$},
    description={Size of the centralizer of a permutation of cycle type $\rho$}
}
\newglossaryentry{z-tau}{
    type=symbols,
    name={$z_\tau$},
    description={Size of the stabilizer of a stitch of type $\tau$}
}
\newglossaryentry{calS}{
    type=symbols,
    name={$\calS$},
    description={Set of all stitches}
}
\newglossaryentry{calS-n}{
    type=symbols,
    name={$\calS_n$},
    description={Set of all stitches with $n$ columns}
}
\newglossaryentry{calS-n-k}{
    type=symbols,
    name={$\calS_{n,k}$},
    description={Set of all stitches with $n$ columns and $k$ wires}
}
\newglossaryentry{calS-parallel}{
    type=symbols,
    name={$\calS^\parallel$},
    description={Set of all parallel stitches}
}
\newglossaryentry{calS-n-parallel}{
    type=symbols,
    name={$\calS_n^\parallel$},
    description={Set of all parallel stitches with $n$ columns}
}
\newglossaryentry{calS-n-k-parallel}{
    type=symbols,
    name={$\calS_{n,k}^\parallel$},
    description={Set of all parallel stitches with $n$ columns and $k$ wires}
}
\newglossaryentry{calS-tau}{
    type=symbols,
    name={$\calS(\tau)$},
    description={Set of all stitches of type $\tau$}
}
\newglossaryentry{calS-parallel-tau}{
    type=symbols,
    name={$\calS^\parallel(\tau)$},
    description={Set of all parallel stitches of type $\tau$}
}
\newglossaryentry{calS-bracket-r}{
    type=symbols,
    name={$\calS[r]$},
    description={Set of all $r$-colored stitches}
}
\newglossaryentry{calS-parallel-bracket-r}{
    type=symbols,
    name={$\calS^\parallel[r]$},
    description={Set of all parallel $r$-colored stitches}
}
\newglossaryentry{calS-lparallel}{
    type=symbols,
    name={$\calS^\lparallel$},
    description={Set of all left-leaning stitches}
}
\newglossaryentry{calS-n-lparallel}{
    type=symbols,
    name={$\calS_n^\lparallel$},
    description={Set of all left-leaning stitches with $n$ columns}
}
\newglossaryentry{calS-n-k-lparallel}{
    type=symbols,
    name={$\calS_{n,k}^\lparallel$},
    description={Set of all left-leaning stitches with $n$ columns and $k$ wires}
}
\newglossaryentry{calS-lparallel-tau}{
    type=symbols,
    name={$\calS^\lparallel(\tau)$},
    description={Set of all left-leaning stitches of type $\tau$}
}
\newglossaryentry{calS-n-k-parallel-bracket-r}{
    type=symbols,
    name={$\calS_{n,k}^\parallel[r]$},
    description={Set of $r$-colored parallel stitches with $n$ columns and $k$ wires}
}
\newglossaryentry{tilde-calS-nu-lparallel-tau}{
    type=symbols,
    name={$\tilde{\calS}_\nu^\lparallel(\tau)$},
    description={Set of left-leaning stitches of type $\tau$ containing exactly $\nu$ irreducible components}
}
\newglossaryentry{hat-calS-nu-lparallel-tau}{
    type=symbols,
    name={$\hat{\calS}_\nu^\lparallel(\tau)$},
    description={Set of tuples $(\alpha_1, \dots, \alpha_\nu)$ of nonempty left-leaning stitches where $\alpha_1 \oplus \dots \oplus \alpha_\nu$ has type $\tau$}
}
\newglossaryentry{calS-nu-lparallel-tau}{
    type=symbols,
    name={$\calS_\nu^\lparallel(\tau)$},
    description={Set of tuples $(\alpha_1, \dots, \alpha_\nu)$ of (possibly empty) left-leaning stitches where $\alpha_1 \oplus \dots \oplus \alpha_\nu$ has type $\tau$}
}
\newglossaryentry{tilde-calS-d-parallel-tau}{
    type=symbols,
    name={$\tilde{\calS}_d^\parallel(\tau)$},
    description={Set of parallel stitches of type $\tau$ with exactly $d$ non-poles in their irreducible decompositions}
}
\newglossaryentry{calT}{
    type=symbols,
    name={$\calT$},
    description={Set of all stitch types}
}
\newglossaryentry{calT-n}{
    type=symbols,
    name={$\calT_n$},
    description={Set of all stitch types with $n$ columns}
}
\newglossaryentry{calT-n-k}{
    type=symbols,
    name={$\calT_{n,k}$},
    description={Set of all stitch types with $n$ columns and $k$ wires}
}
\newglossaryentry{calT-parallel}{
    type=symbols,
    name={$\calT^\parallel$},
    description={Set of all parallelizable stitch types}
}
\newglossaryentry{calT-n-parallel}{
    type=symbols,
    name={$\calT_n^\parallel$},
    description={Set of all parallelizable stitch types with $n$ columns}
}
\newglossaryentry{calT-n-k-parallel}{
    type=symbols,
    name={$\calT_{n,k}^\parallel$},
    description={Set of all parallelizable stitch types with $n$ columns and $k$ wires}
}
\newglossaryentry{calT-lparallel}{
    type=symbols,
    name={$\calT^\lparallel$},
    description={Set of all cycle-free stitch types}
}
\newglossaryentry{calT-n-lparallel}{
    type=symbols,
    name={$\calT_n^\lparallel$},
    description={Set of all cycle-free stitch types with $n$ columns}
}
\newglossaryentry{calT-n-k-lparallel}{
    type=symbols,
    name={$\calT_{n,k}^\lparallel$},
    description={Set of all cycle-free stitch types with $n$ columns and $k$ wires}
}
\newglossaryentry{calT-ell-m-t-nu-parallel}{
    type=symbols,
    name={$\calT_{\ell,m,t,\nu}^\parallel$},
    description={Set of all parallelizable stitch types with $\ell$ proper chains, $m$ gaps, $t$ jumps, and $\nu$ poles}
}
\newglossaryentry{G-alpha}{
    type=symbols,
    name={$G_\alpha$},
    description={The directed graph corresponding to a stitch $\alpha$}
}
\newglossaryentry{len}{
    type=symbols,
    name={$\len(\rho)$},
    description={The number of parts of a partition $\rho$}
}
\newglossaryentry{oplus}{
    type=symbols,
    name={$\alpha_1 \oplus \alpha_2$},
    description={Direct sum of two stitches}
}
\newglossaryentry{bar-alpha}{
    type=symbols,
    name={$\bar{\alpha}$},
    description={Reduction of a stitch $\alpha$}
}
\newglossaryentry{bar-tau}{
    type=symbols,
    name={$\bar{\tau}$},
    description={Reduction of a stitch type $\tau$}
}
\newglossaryentry{stitch-subseteq}{
    type=symbols,
    name={$\alpha_1 \subseteq \alpha_2$},
    description={The stitch $\alpha_1$ is obtained by removing wires from $\alpha_2$}
}
\newglossaryentry{stitch-binom}{
    type=symbols,
    name={$\binom{\tau}{\omega}$},
    description={Number of ways to obtain a stitch of type $\omega$ by removing wires from a stitch of type $\tau$}
}
\newglossaryentry{stitch-brackbinom}{
    type=symbols,
    name={$\brackbinom{\tau}{\omega}$},
    description={Number of ways to obtain a stitch of type $\tau$ by adding wires from a stitch of type $\omega$}
}
\newglossaryentry{stitch-ell-m-t-nu-binom}{
    type=symbols,
    name={$\binom{\rho}{\ell,m,t,\nu}$},
    description={Number of ways to obtain a stitch with $0$ proper cycles, $\ell$ proper chains, $m$ gaps, $t$ jumps, and $\nu$ poles by removing wires from a permutation of cycle type $\rho$}
}
\newglossaryentry{calE-f}{
    type=symbols,
    name={$\calE_f$},
    description={Generating function in $x$ and $y$ built from the $E_f(n, k)$}
}
\newglossaryentry{alpha-rvert-i}{
    type=symbols,
    name={$\alpha\rvert_i$},
    description={The $i$th color restriction of a colored stitch $\alpha$}
}
\newglossaryentry{partial-x-i}{
    type=symbols,
    name={$\partial_x^i$},
    description={The $i$th partial derivative (or antiderivative if $i$ is negative) with respect to $x$}
}
\newglossaryentry{h-tau-j-k}{
    type=symbols,
    name={$h_\tau(j,k)$},
    description={Enumerates certain families of $2$-colored stitches}
}
\newglossaryentry{H-tau}{
    type=symbols,
    name={$H_\tau$},
    description={Generating functions appearing in Lemma~\ref{E-from-H-tau-lem}}
}
\newglossaryentry{h-ell-m-t-nu-j-k}{
    type=symbols,
    name={$h_{\ell,m,t,\nu}(j,k)$},
    description={Enumerates certain families of $2$-colored stitches}
}
\newglossaryentry{H-ell-m-t-nu}{
    type=symbols,
    name={$H_{\ell,m,t,\nu}$},
    description={Generating functions appearing in Corollary~\ref{E-from-H-ell-m-t-nu-cor}}
}
\newglossaryentry{tilde-H-ell-m-t-nu}{
    type=symbols,
    name={$\tilde{H}_{\ell,m,t,\nu}$},
    description={Shared value of $z_\tau H_\tau$ for $\tau \in \calT_{\ell,m,t,\nu}^\parallel$}
}
\newglossaryentry{f-tilde-lparallel}{
    type=symbols,
    name={$\tilde{f}_\lparallel(j_1,t_1;\dots;j_r,t_r)$},
    description={Number of integral $r$-colored left-leaning stitches $\alpha$ such that $\alpha\rvert_i$ has $j_i$ chains and $t_i$ jumps}
}
\newglossaryentry{F-tilde-lparallel}{
    type=symbols,
    name={$\tilde{F}_\lparallel$},
    description={Generating function built from $\tilde{f}_\lparallel(j_1,t_1;\dots;j_r,t_r)$}
}
\newglossaryentry{H}{
    type=symbols,
    name={$H$},
    description={Generating function in variables $L$, $T$, $N$, $x$, and $y$ built from the $H_{\ell,m,t,\nu}$}
}
\newglossaryentry{Q}{
    type=symbols,
    name={$Q$},
    description={Formal power series $\sqrt{(1-x-y)^2-4xy}$}
}
\newglossaryentry{U}{
    type=symbols,
    name={$U$},
    description={Formal power series $\frac{1+x-y+Q}{2Q}$}
}
\newglossaryentry{V}{
    type=symbols,
    name={$V$},
    description={Formal power series $\frac{1-x+y+Q}{2Q}$}
}
\newglossaryentry{W}{
    type=symbols,
    name={$W$},
    description={Formal power series $\frac{1-x-y+Q}{2Q}$}
}
\newglossaryentry{deg-x}{
    type=symbols,
    name={$\deg_x(p)$},
    description={Degree of a multivariate polynomial $p$ with respect to a variable $x$}
}
\newglossaryentry{indicator}{
    type=symbols,
    name={$\mathbbm{1}_P$},
    description={Indicator function for the condition $P$}
}
\newglossaryentry{p-ell-m-t-nu-i}{
    type=symbols,
    name={$p_{\ell,m,t,\nu}^{(i)}$},
    description={Polynomials appearing in Theorem~\ref{H-from-p-thm}}
}
\newglossaryentry{hat-p-ell-m-t-nu-i}{
    type=symbols,
    name={$\hat{p}_{\ell,m,t,\nu}^{(i)}$},
    description={Renormalized version of $p_{\ell,m,t,\nu}^{(i)}$}
}
\newglossaryentry{delta-f-j}{
    type=symbols,
    name={$\delta_f^{(j)}$},
    description={Polynomials appearing in Theorem~\ref{E-from-delta-thm}}
}
\newglossaryentry{hat-varepsilon-f-i}{
    type=symbols,
    name={$\hat{\varepsilon}_f^{(i)}$},
    description={Polynomials appearing in Lemma~\ref{E-from-epsilon-hat-lem}}
}
\newglossaryentry{q-a-b-r-j-k}{
    type=symbols,
    name={$q_{a,b}^{(r)}(j,k)$},
    description={Series coefficients of $U^a V^b Q^{-1-r}$}
}
\newglossaryentry{varepsilon-f}{
    type=symbols,
    name={$\varepsilon_f$},
    description={Polynomials appearing in Theorem~\ref{eps-properties-thm}}
}
\newglossaryentry{stitch-derivative}{
    type=symbols,
    name={$\partial\alpha$},
    description={Derivative of a stitch $\alpha$}
}
\newglossaryentry{stitch-type-derivative}{
    type=symbols,
    name={$\partial\tau$},
    description={Derivative of a stitch type $\tau$}
}
\newglossaryentry{stitch-splitting}{
    type=symbols,
    name={$\hat{\alpha}$},
    description={The splitting of a stitch $\alpha$}
}
\newglossaryentry{corner-set}{
    type=symbols,
    name={$\Gamma(S)$},
    description={The corner set of a staircase diagram $S$}
}
\newglossaryentry{stair}{
    type=symbols,
    name={$\stair(\alpha)$},
    description={The staircase diagram of a stitch $\alpha$}
}
\newglossaryentry{falling-factorial}{
    type=symbols,
    name={$(n)_k$},
    description={The falling factorial $\frac{n!}{(n-k)!}$}
}
\newglossaryentry{graph-deg}{
    type=symbols,
    name={$\deg_G(i)$},
    description={Degree of vertex $i$ in a graph $G$}
}
\newglossaryentry{total-deg}{
    type=symbols,
    name={$\deg(p)$},
    description={Total degree of a possibly multivariate polynomial $p$}
}
\newglossaryentry{character-polynomial-deg}{
    type=symbols,
    name={$\deg(f)$},
    description={Degree of a character polynomial $f$}
}
\newglossaryentry{2F1}{
    type=symbols,
    name={${}_2F_1\left(a,b;c;z\right)$},
    description={Gauss's hypergeometric function}
}
\newglossaryentry{f-lparallel}{
    type=symbols,
    name={$f_\lparallel(j_1,k_1;\dots;j_r,k_r)$},
    description={Number of $r$-colored left-leaning stitches $\alpha$ such that $\alpha\rvert_i$ has $j_i$ chains and $k_i$ wires}
}
\newglossaryentry{f-lparallel-nu}{
    type=symbols,
    name={$f_\lparallel^{(\nu)}(j_1,k_1;\dots;j_r,k_r)$},
    description={Number of tuples of (possibly empty) $r$-colored left-leaning stitches $(\alpha_1,\dots,\alpha_\nu)$ such that $(\alpha_1 \oplus \dots \oplus \alpha_\nu)\rvert_i$ has $j_i$ chains and $k_i$ wires}
}
\newglossaryentry{F-lparallel}{
    type=symbols,
    name={$F_\lparallel$},
    description={Generating function for $f_\lparallel(j_1,k_1;\dots;j_r,k_r)$}
}
\newglossaryentry{prec}{
    type=symbols,
    name={$\lambda_1 \prec \lambda_2$},
    description={$\lambda_1$ can be obtained by removing at most one box in each column from $\lambda_2$}
}
\newglossaryentry{prec-prime}{
    type=symbols,
    name={$\lambda_1 \prec' \lambda_2$},
    description={$\lambda_1$ can be obtained by removing at most one box in each row from $\lambda_2$}
}
\newglossaryentry{binomial}{
    type=symbols,
    name={$\binom{x}{k}$},
    description={The generalized binomial coefficient $\frac{1}{k!} \cdot \prod_{i=0}^{k-1}(x-i)$}
}

\newglossaryentry{pattern-occurrence}{
    type=terms,
    name={occurrence of $\sigma$ in $\pi$},
    description={Subsequence of $\pi$ with same ordering as $\sigma$}
}
\newglossaryentry{pattern-containment}{
    type=terms,
    name={$\pi$ contains the pattern $\sigma$},
    description={$N_\sigma(\pi) > 0$}
}
\newglossaryentry{pattern-avoidance}{
    type=terms,
    name={$\pi$ avoids the pattern $\sigma$},
    description={$N_\sigma(\pi) = 0$}
}
\newglossaryentry{stitch}{
    type=terms,
    name={stitch},
    description={A partial permutation visualized as two rows of points connected by wires}
}
\newglossaryentry{character-polynomial}{
    type=terms,
    name={character polynomial},
    description={An element of the ring $\Delta_\Q = \Q[m_1,m_2,\dots]$}
}
\newglossaryentry{pure-character-polynomial}{
    type=terms,
    name={$f$ is pure of degree $d$},
    description={The character polynomial $f$ is an element of the subspace $\Delta_\Q^d \subseteq \Delta_\Q$}
}
\newglossaryentry{standard-form}{
    type=terms,
    name={standard form of a stitch},
    description={Unique way of writing a stitch as $([n], W)$ for some $W \subseteq [n] \times [n]$}
}
\newglossaryentry{wire-diagram}{
    type=terms,
    name={wire diagram of a stitch},
    description={Way of visualizing a stitch as two rows of points connected by wires}
}
\newglossaryentry{column}{
    type=terms,
    name={columns of a stitch},
    description={Pairs of lower and upper points in a wire diagram}
}
\newglossaryentry{wire}{
    type=terms,
    name={wires of a stitch},
    description={Lines in a wire diagram}
}
\newglossaryentry{proper-chain}{
    type=terms,
    name={proper chain of a stitch},
    description={Chain of length greater than $1$}
}
\newglossaryentry{proper-cycle}{
    type=terms,
    name={proper cycle of a stitch},
    description={Cycle of length greater than $1$}
}
\newglossaryentry{chain}{
    type=terms,
    name={chain of length $j$ in a stitch},
    description={Chain in the digraph of a stitch with $j$ vertices}
}
\newglossaryentry{cycle}{
    type=terms,
    name={cycle of length $j$ in a stitch},
    description={Cycle in the digraph of a stitch with $j$ vertices}
}
\newglossaryentry{parallel-stitch}{
    type=terms,
    name={$\alpha$ is parallel},
    description={No two wires in the wire diagram of $\alpha$ intersect}
}
\newglossaryentry{parallelizable-stitch-type}{
    type=terms,
    name={$\tau$ is parallelizable},
    description={The stitch type $\tau$ contains no proper cycles}
}
\newglossaryentry{stitch-type}{
    type=terms,
    name={stitch type},
    description={A pair of partitions which specifies a conjugacy class of stitches}
}
\newglossaryentry{stitch-coloring}{
    type=terms,
    name={$r$-coloring of $\alpha$},
    description={A way of coloring the columns of $\alpha$ so wires connect columns with the same color}
}
\newglossaryentry{irreducible-stitch}{
    type=terms,
    name={$\alpha$ is irreducible},
    description={The stitch $\alpha$ cannot be written as a direct sum of smaller stitches}
}
\newglossaryentry{left-leaning-stitch}{
    type=terms,
    name={$\alpha$ is left- (resp.\@ right-) leaning},
    description={The wires $(i, j)$ of $\alpha$ satisfy $i > j$ (resp.\@ $i < j$)}
}
\newglossaryentry{gap}{
    type=terms,
    name={gaps of a stitch},
    description={Columns with no wires connected to them}
}
\newglossaryentry{reduction}{
    type=terms,
    name={reduction of $\alpha$},
    description={The stitch $\bar{\alpha}$ obtained by removing all gaps from $\alpha$}
}
\newglossaryentry{reduced}{
    type=terms,
    name={$\alpha$ is reduced},
    description={The stitch $\alpha$ has no gaps}
}
\newglossaryentry{integral}{
    type=terms,
    name={$\alpha$ is integral},
    description={The stitch $\alpha$ is irreducible and reduced}
}
\newglossaryentry{pole}{
    type=terms,
    name={poles of a stitch},
    description={Columns with a single vertical wire connected to the top and the bottom}
}
\newglossaryentry{jump}{
    type=terms,
    name={jumps of a stitch},
    description={Non-pole columns with wires connected to both the top and the bottom}
}
\newglossaryentry{stitch-derivative-term}{
    type=terms,
    name={derivative of $\alpha$},
    description={Stitch obtained from $\alpha$ by reinterpreting wires as columns}
}
\newglossaryentry{stitch-antiderivative}{
    type=terms,
    name={antiderivative of $\alpha$},
    description={A stitch whose derivative is $\alpha$}
}
\newglossaryentry{tie}{
    type=terms,
    name={tie of $\alpha$},
    description={Either a jump or a pole}
}
\newglossaryentry{splitting}{
    type=terms,
    name={splitting of a stitch},
    description={Stitch obtained by replacing all ties with join sites}
}
\newglossaryentry{split}{
    type=terms,
    name={$\alpha$ is split},
    description={$\alpha$ has no ties}
}
\newglossaryentry{join-site}{
    type=terms,
    name={join site of $\alpha$},
    description={A pair of adjacent columns where the left column only has a wire connected to its lower vertex and the right column only has a wire connected to its upper vertex}
}
\newglossaryentry{staircase-diagram}{
    type=terms,
    name={staircase diagram},
    description={A subset $S \subseteq [n]^2$ such that if a given box is contained in $S$, then all boxes down and to the right of it are also contained in $S$}
}
\newglossaryentry{corner-set-term}{
    type=terms,
    name={corner of a staircase diagram},
    description={A box with no boxes above or to the left of it}
}
\newglossaryentry{stitch-extension}{
    type=terms,
    name={extension of $\alpha$},
    description={A stitch $\beta$ with the same number of columns as $\alpha$ satisfying $\alpha \subseteq \beta$}
}
\newglossaryentry{substitch}{
    type=terms,
    name={substitch of $\alpha$},
    description={A stitch obtained from $\alpha$ by removing some of its wires while retaining the same columns}
}
\newglossaryentry{empty-stitch}{
    type=terms,
    name={empty stitch},
    description={The stitch $(\varnothing, \varnothing)$ with no wires or columns}
}

\title{Explicit Formulas for Permutation Pattern Character Polynomials}
\author{Jonas Iskander}
\date{\today}

\begin{document}

\begin{abstract}
    Given permutations $\pi \in \glslink{symmetric-group}{S_n}$ and $\sigma \in \glslink{symmetric-group}{S_k}$, let $\glslink{N-sigma}{N_\sigma(\pi)}$ denote the number of \glslink{pattern-occurrence}{occurrences of $\sigma$ in $\pi$}. While pattern avoidance and the distribution of \glslink{pattern-occurrence}{pattern occurrences} in permutations have been extensively studied, their interactions with the group structure on $\glslink{symmetric-group}{S_n}$ are still poorly understood. Gaetz and Ryba showed that the expected value of $\glslink{chi-lambda}{\chi^{\glslink{lambda-bracket-n}{\lambda[n]}}}(\pi)\glslink{N-sigma}{N_\sigma(\pi)}$ for $\pi \in \glslink{symmetric-group}{S_n}$ is given by a polynomial $\glslink{a-sigma-lambda}{a_\sigma^\lambda}(n)$. More recently, Gaetz and Pierson derived explicit formulas for $\glslink{a-sigma-lambda}{a_{\glslink{id-k}{\id_k}}^\lambda}(n)$ when $\glslink{abs-rho}{\lvert\lambda\rvert} \le 2$, which led them to conjecture that the polynomials $\glslink{a-sigma-lambda}{a_{\glslink{id-k}{\id_k}}^\lambda}(n)$ are real-rooted and nonnegative for $n \ge k$. We show that for all partitions $\lambda$, the polynomials $\glslink{a-sigma-lambda}{a_{\glslink{id-k}{\id_k}}^\lambda}(n)$ admit explicit closed forms in $n$ and $k$. These formulas allow us to exhibit counterexamples to Gaetz and Pierson's real-rootedness conjecture as well as to prove special cases of their nonnegativity conjecture. Our results imply that the expected value of $f \cdot \glslink{N-sigma}{N_{\glslink{id-k}{\id_k}}}$ on $\glslink{symmetric-group}{S_n}$ admits a closed form whenever $f$ is a permutation statistic expressible as a polynomial in the functions $\glslink{m-j}{m_j} \colon \bigsqcup_{n \ge 0} \glslink{symmetric-group}{S_n} \to \Z$ which count $j$-cycles in their inputs.
\end{abstract}

\maketitle

\section{Introduction}

Let $\pi \in \glslink{symmetric-group}{S_n}$ and $\sigma \in \glslink{symmetric-group}{S_k}$ be permutations. An \glslink{pattern-occurrence}{\textit{occurrence of the pattern $\sigma$ in $\pi$}} is a sequence $1 \le i_1 < \dots < i_k \le n$ such that for all $1 \le j, j' \le k$, we have $\pi(i_j) \le \pi(i_{j'})$ if and only if $\sigma(j) \le \sigma(j')$. We denote the number of such \glslink{pattern-occurrence}{occurrences} by $\glslink{N-sigma}{N_\sigma(\pi)}$---if $\glslink{N-sigma}{N_\sigma(\pi)} > 0$, we say that \glslink{pattern-containment}{$\pi$ \textit{contains} the pattern $\sigma$}; otherwise, we say that \glslink{pattern-avoidance}{$\pi$ \textit{avoids} $\sigma$}.

It is natural to ask how \glslink{pattern-occurrence}{pattern occurrences} interact with the group structure on $\glslink{symmetric-group}{S_n}$. Several authors have initiated the study of such interactions by investigating the relationships between \glslink{pattern-occurrence}{pattern occurrences} and conjugacy classes in the symmetric group. More precisely, for any list of patterns $\sigma_1,\dots,\sigma_d$ and any partition $\rho$, consider the mixed moment \begin{equation}
    \glslink{M-rho-sigma}{M_\rho^{\sigma_1,\dots,\sigma_d}} := \frac{1}{\lvert \glslink{C-rho}{C_\rho}\rvert} \sum_{\pi \in \glslink{C-rho}{C_\rho}} \glslink{N-sigma}{N_{\sigma_1}(\pi)} \cdots \glslink{N-sigma}{N_{\sigma_d}(\pi)}, \label{M-def}
\end{equation} where $\glslink{C-rho}{C_\rho} \subseteq \glslink{symmetric-group}{S_{\glslink{abs-rho}{\lvert\rho\rvert}}}$ is the conjugacy class consisting of permutations with cycle type $\rho$. In 2013--2014, Hultman and Gill investigated the values $\glslink{M-rho-sigma}{M_\rho^\sigma}$ for $\sigma \in \glslink{symmetric-group}{S_2}$ or $\glslink{symmetric-group}{S_3}$ \cite{hultman,gill}. More recently, in 2021, Gaetz and Ryba showed that when $\sigma_1 = \cdots = \sigma_d \in \glslink{symmetric-group}{S_k}$, the moment $\glslink{M-rho-sigma}{M_\rho^{\sigma_1,\dots,\sigma_d}}$ is given by a polynomial in $n, \glslink{m-j}{m_1}, \dots, \glslink{m-j}{m_{dk}}$, where $n := \glslink{abs-rho}{\lvert\rho\rvert}$ and $\glslink{m-j}{m_j}$ denotes the number of parts of size $j$ in $\rho$ \cite{gaetz-ryba}. Soon after, in 2022, Gaetz and Pierson removed the condition that $\sigma_1=\cdots=\sigma_d$ and showed that for any collection of patterns $\sigma_i \in \glslink{symmetric-group}{S_{k_i}}$, $\glslink{M-rho-sigma}{M_\rho^{\sigma_1,\dots,\sigma_d}}$ is given by a polynomial in $n$, $\glslink{m-j}{m_1}, \dots, \glslink{m-j}{m_{k_1+\dots+k_d}}$ \cite{gaetz-pierson}. In early 2023, Levet, Liu, Loth, Stucky, Sundaram, and Yin proved similar results for so-called ``permutation constraint statistics,'' which further generalize these mixed moments, and they applied their methods to compute the expected values over fixed conjugacy classes of ``weighted inversion statistics'' \cite{levet-et-al}.

In addition to studying the statistics of \glslink{pattern-occurrence}{pattern occurrences} restricted to individual conjugacy classes, Gaetz, Pierson, and Ryba showed that it can be fruitful to study \glslink{pattern-occurrence}{pattern occurrences} weighted by symmetric group characters, and they proved several nontrivial results in this direction. In analogy with \eqref{M-def}, for any list of patterns $\sigma_1, \dots, \sigma_d$ and any partition $\lambda$, define \begin{equation*}
    \glslink{A-sigma-lambda}{A_{\sigma_1,\dots,\sigma_d}^\lambda} := \frac{1}{\lvert\lambda\rvert!}\sum_{\pi \in \glslink{symmetric-group}{S_{\glslink{abs-rho}{\lvert\lambda\rvert}}}} \glslink{chi-lambda}{\chi^\lambda}(\pi) \cdot \glslink{N-sigma}{N_{\sigma_1}(\pi)} \cdots \glslink{N-sigma}{N_{\sigma_d}(\pi)},
\end{equation*} where $\glslink{chi-lambda}{\chi^\lambda}$ denotes the irreducible symmetric group character corresponding to $\lambda$. In addition, given a partition $\lambda = (\glslink{lambda-i}{\lambda_1}, \dots, \glslink{lambda-i}{\lambda_p})$ and any integer $n \ge \glslink{abs-rho}{\lvert\lambda\rvert} + \glslink{lambda-i}{\lambda_1}$, let $\glslink{lambda-bracket-n}{\lambda[n]}$ denote the partition $(n-\glslink{abs-rho}{\lvert\lambda\rvert}, \glslink{lambda-i}{\lambda_1}, \dots, \glslink{lambda-i}{\lambda_p}) \vdash n$. Then we have the following theorem of Gaetz and Pierson, which generalizes Theorem~1.1(b) from \cite{gaetz-ryba}:

\begin{theorem}[Gaetz--Pierson \cite{gaetz-pierson}, Theorem~1.3] \label{gp-poly-thm}
    For any $\sigma_1 \in \glslink{symmetric-group}{S_{k_1}}, \dots, \sigma_d \in \glslink{symmetric-group}{S_{k_d}}$, there exists a polynomial $\glslink{a-sigma-lambda}{a_{\sigma_1, \dots, \sigma_d}^\lambda}(n)$ in $n$ of degree at most $k_1+\dots+k_d-\glslink{abs-rho}{\lvert\lambda\rvert}$ such that \begin{equation}
        \glslink{A-sigma-lambda}{A_{\sigma_1,\dots,\sigma_d}^{\glslink{lambda-bracket-n}{\lambda[n]}}} = \frac{1}{n!} \sum_{\pi \in \glslink{symmetric-group}{S_n}} \glslink{chi-lambda}{\chi^{\glslink{lambda-bracket-n}{\lambda[n]}}}(\pi) \cdot \glslink{N-sigma}{N_{\sigma_1}(\pi)} \cdots \glslink{N-sigma}{N_{\sigma_d}(\pi)} = \glslink{a-sigma-lambda}{a_{\sigma_1,\dots,\sigma_d}^\lambda}(n)
    \end{equation} for all $n \ge k_1+\cdots+k_d+\glslink{abs-rho}{\lvert\lambda\rvert}$. When $\glslink{abs-rho}{\lvert\lambda\rvert} > k_1+\dots+k_d$, we interpret the statement that $\glslink{a-sigma-lambda}{a_{\sigma_1, \dots, \sigma_d}^\lambda}$ has degree at most $k_1+\dots+k_d-\glslink{abs-rho}{\lvert\lambda\rvert}$ to mean that $\glslink{a-sigma-lambda}{a_{\sigma_1, \dots, \sigma_d}^\lambda}$ is the zero polynomial.
\end{theorem}

\begin{remark}
    While the quantity $\glslink{A-sigma-lambda}{A_{\sigma_1,\dots,\sigma_d}^{\glslink{lambda-bracket-n}{\lambda[n]}}}$ is technically only defined for $n \ge \glslink{abs-rho}{\lvert\lambda\rvert}+\glslink{lambda-i}{\lambda_1}$, we will see in Section~\ref{background-on-character-polynomials-section} that $\glslink{chi-lambda}{\chi^{\glslink{lambda-bracket-n}{\lambda[n]}}}$ and $\glslink{A-sigma-lambda}{A_{\sigma_1,\dots,\sigma_d}^{\glslink{lambda-bracket-n}{\lambda[n]}}}$ have a natural interpretation for all $n \ge 0$---see Theorem~\ref{macdonald-thm} and equation \eqref{A-redef}. To minimize distractions, we will state all theorems using these extended definitions for $\glslink{chi-lambda}{\chi^{\glslink{lambda-bracket-n}{\lambda[n]}}}$ and $\glslink{A-sigma-lambda}{A_{\sigma_1,\dots,\sigma_d}^{\glslink{lambda-bracket-n}{\lambda[n]}}}$, allowing us to omit the $n \ge \glslink{abs-rho}{\lvert\lambda\rvert}+\glslink{lambda-i}{\lambda_1}$ condition.
\end{remark}

A natural next question is whether the quantities $\glslink{M-rho-sigma}{M_\rho^\sigma}$ and $\glslink{A-sigma-lambda}{A_\sigma^\lambda}$ exhibit any special properties in the case where $\sigma$ is the identity permutation $\glslink{id-k}{\id_k}$, and \glslink{pattern-occurrence}{occurrences of $\sigma$ in a permutation $\pi$} correspond to length-$k$ increasing subsequences of $\pi$. On this front, Gaetz and Pierson proposed the following two related conjectures \cite[Conjecture~1.4]{gaetz-pierson}:

\begin{conjecture}[Gaetz--Pierson] \label{root-conj}
    For any partition $\lambda$ and all $k \ge 0$, the polynomial $\glslink{a-sigma-lambda}{a_{\glslink{id-k}{\id_k}}^\lambda}(n)$ is real-rooted with all roots less than $k$.
\end{conjecture}

\begin{conjecture}[Gaetz--Pierson] \label{pos-conj}
    For all $n, k \ge 0$ and $\lambda \vdash n$, we have $\glslink{A-sigma-lambda}{A_{\glslink{id-k}{\id_k}}^\lambda} \ge 0$.
\end{conjecture}

\noindent
Conjecture~\ref{pos-conj} is of particular interest, as it would imply that the class function on $\glslink{symmetric-group}{S_n}$ given by $\pi \mapsto \glslink{M-rho-sigma}{M_{\glslink{typ-pi}{\typ(\pi)}}^{\glslink{id-k}{\id_k}}}$, where $\glslink{typ-pi}{\typ(\pi)}$ denotes the cycle type of $\pi$, is a rational multiple of the character of some representation of $\glslink{symmetric-group}{S_n}$.

For $\lambda \in \{\varnothing, (1), (2), (1,1)\}$, Gaetz and Pierson refined Theorem~\ref{gp-poly-thm} and computed the following explicit formulas for $\glslink{a-sigma-lambda}{a_{\glslink{id-k}{\id_k}}^\lambda}(n)$, which allowed them to deduce Conjecture~\ref{root-conj} for these values of $\lambda$. We state their formulas below in an equivalent form better suited for our purposes. In the theorem statement as well as throughout the rest of the paper, we use the generalized binomial coefficient $\glslink{binomial}{\binom{x}{k}} := \frac{1}{k!} \cdot \prod_{i=0}^{k-1} (x-i)$, defined for all integers $k \ge 0$ and $x$ an element of any $\Q$-algebra.

\begin{theorem}[Gaetz--Pierson \cite{gaetz-pierson}, Theorem~2.7] \label{small-lam-formula}
    For all $k \ge 0$ such that the right hand sides are defined, we have \begin{alignat*}{2}
        \glslink{a-sigma-lambda}{a_{\glslink{id-k}{\id_k}}^\varnothing}(n) &= \frac{1}{k!} \glslink{binomial}{\binom{n}{k}}, &\qquad \glslink{a-sigma-lambda}{a_{\glslink{id-k}{\id_k}}^{(1)}}(n) &= \frac{2^{k-1}}{(2k-1)!!} \glslink{binomial}{\binom{n-\frac{1}{2}}{k-1}} - \frac{1}{k!} \glslink{binomial}{\binom{n-1}{k-1}}, \\
        \glslink{a-sigma-lambda}{a_{\glslink{id-k}{\id_k}}^{(2)}}(n) &= -\frac{1}{n \cdot k!} \glslink{binomial}{\binom{n-2}{k-1}} + \frac{1}{2(k-1)!} \glslink{binomial}{\binom{n-1}{k-2}} - 2^{k-2}\frac{(2k-4)n+(2k-1)}{(k-1)n \cdot (2k-1)!!} \glslink{binomial}{\binom{n-\frac{3}{2}}{k-2}}, \span\span \\
        \glslink{a-sigma-lambda}{a_{\glslink{id-k}{\id_k}}^{(1,1)}}(n) &= \frac{1}{k!} \glslink{binomial}{\binom{n-2}{k-2}} + \frac{1}{n \cdot k!} \glslink{binomial}{\binom{n-2}{k-1}} + \frac{1}{2(k-1)!} \glslink{binomial}{\binom{n-1}{k-2}} - 2^{k-2} \frac{2kn-(2k-1)}{(k-1)n \cdot (2k-1)!!} \glslink{binomial}{\binom{n-\frac{3}{2}}{k-2}}. \span\span
    \end{alignat*} Moreover, when $\lambda = (1), (2), (1,1)$, we have $\glslink{A-sigma-lambda}{A_{\glslink{id-k}{\id_k}}^{\glslink{lambda-bracket-n}{\lambda[n]}}} = \glslink{a-sigma-lambda}{a_{\glslink{id-k}{\id_k}}^\lambda}(n)$ for all $n \ge k+\glslink{abs-rho}{\lvert\lambda\rvert}-1$, giving a slight improvement on the prediction of Theorem~\ref{gp-poly-thm}.
\end{theorem}

The main goal of this paper is to prove the following generalization of Theorem~\ref{small-lam-formula} which applies to \textit{all} partitions $\lambda$. In the statement of the theorem as well as the remainder of the paper, we define $x! := \Gamma(x+1)$ for all $x \in \R$, with the convention that $x! := \infty$ whenever $x$ is a negative integer and $\frac{y}{\infty} := 0$ whenever $y$ is a finite real number. (Any expression that would require evaluating $\frac{\infty}{\infty}$ is regarded as undefined.)

\begin{theorem} \label{main-thm-1}
    Let $\lambda$ be a partition of size at least $1$. Then there exist coefficients $\glslink{theta-f-i}{\theta_{\glslink{X-lam}{X^\lambda}}^{(1)}}, \dots, \glslink{theta-f-i}{\theta_{\glslink{X-lam}{X^\lambda}}^{(3\glslink{abs-rho}{\lvert\lambda\rvert})}} \in \Q$ such that for all $k \ge 0$, we have \begin{equation*}
        \frac{n!}{(n-\glslink{abs-rho}{\lvert\lambda\rvert})!} \cdot \glslink{a-sigma-lambda}{a_{\glslink{id-k}{\id_k}}^\lambda}(n) = \sum_{1 \le i \le 3\glslink{abs-rho}{\lvert\lambda\rvert}} \glslink{theta-f-i}{\theta_{\glslink{X-lam}{X^\lambda}}^{(i)}} \cdot \frac{\left(\frac{i}{2}\right)!}{\left(k-\frac{i}{2}\right)!} \glslink{binomial}{\binom{n-\glslink{abs-rho}{\lvert\lambda\rvert}+\frac{i}{2}}{k}}
    \end{equation*} as polynomials in $n$, where we interpret $\frac{n!}{(n-\glslink{abs-rho}{\lvert\lambda\rvert})!}$ as the falling factorial $n \cdot (n-1) \cdots (n-\glslink{abs-rho}{\lvert\lambda\rvert}+1)$. Moreover, we have $\glslink{A-sigma-lambda}{A_{\glslink{id-k}{\id_k}}^{\glslink{lambda-bracket-n}{\lambda[n]}}} = \glslink{a-sigma-lambda}{a_{\glslink{id-k}{\id_k}}^\lambda}(n)$ for all $n \ge k+\glslink{abs-rho}{\lvert\lambda\rvert}-1$, assuming that $(n, k) \ne (\glslink{abs-rho}{\lvert\lambda\rvert}-1, 0)$.
\end{theorem}

\begin{figure}[h]
    \renewcommand{\arraystretch}{1.5}
    \setlength{\tabcolsep}{2pt}
    \begin{tabular}{|c|c|c|c|c|c|c|c|c|c|c|c|c|}
        \hline
        $\lambda$ & $\glslink{theta-f-i}{\theta_{\glslink{X-lam}{X^\lambda}}^{(1)}}$ & $\glslink{theta-f-i}{\theta_{\glslink{X-lam}{X^\lambda}}^{(2)}}$ & $\glslink{theta-f-i}{\theta_{\glslink{X-lam}{X^\lambda}}^{(3)}}$ & $\glslink{theta-f-i}{\theta_{\glslink{X-lam}{X^\lambda}}^{(4)}}$ & $\glslink{theta-f-i}{\theta_{\glslink{X-lam}{X^\lambda}}^{(5)}}$ & $\glslink{theta-f-i}{\theta_{\glslink{X-lam}{X^\lambda}}^{(6)}}$ & $\glslink{theta-f-i}{\theta_{\glslink{X-lam}{X^\lambda}}^{(7)}}$ & $\glslink{theta-f-i}{\theta_{\glslink{X-lam}{X^\lambda}}^{(8)}}$ & $\glslink{theta-f-i}{\theta_{\glslink{X-lam}{X^\lambda}}^{(9)}}$ & $\glslink{theta-f-i}{\theta_{\glslink{X-lam}{X^\lambda}}^{(10)}}$ & $\glslink{theta-f-i}{\theta_{\glslink{X-lam}{X^\lambda}}^{(11)}}$ & $\glslink{theta-f-i}{\theta_{\glslink{X-lam}{X^\lambda}}^{(12)}}$ \\ \hline
        $(1)$ & $\phantom{+}\frac{1}{2}$ & $-1$ & $\phantom{+}\frac{2}{3}$ & & & & & & & & & \\ \hline
        $(2)$ & $\phantom{+}\frac{3}{4}$ & $-1$ & $\phantom{+}0$ & $\phantom{+}\frac{1}{2}$ & $-\frac{4}{15}$ & $\phantom{+}\frac{1}{12}$ & & & & & & \\ \hline
        $(1, 1)$ & $-\frac{1}{4}$ & $\phantom{+}1$ & $-\frac{4}{3}$ & $\phantom{+}1$ & $-\frac{4}{15}$ & $\phantom{+}\frac{1}{12}$ & & & & & & \\ \hline
        $(3)$ & $\phantom{+}\frac{15}{8}$ & $-3$ & $\phantom{+}\frac{15}{16}$ & $\phantom{+}\frac{1}{2}$ & $-\frac{1}{4}$ & $-\frac{1}{12}$ & $\phantom{+}\frac{3}{35}$ & $-\frac{1}{48}$ & $\phantom{+}\frac{4}{945}$ & & & \\ \hline
        $(2, 1)$ & $-\frac{3}{8}$ & $\phantom{+}2$ & $-\frac{21}{8}$ & $\phantom{+}1$ & $\phantom{+}\frac{13}{30}$ & $-\frac{1}{2}$ & $\phantom{+}\frac{26}{105}$ & $-\frac{1}{24}$ & $\phantom{+}\frac{8}{945}$ & & & \\ \hline
        $(1, 1, 1)$ & $\phantom{+}0$ & $-1$ & $\phantom{+}\frac{39}{16}$ & $-\frac{5}{2}$ & $\phantom{+}\frac{89}{60}$ & $-\frac{7}{12}$ & $\phantom{+}\frac{17}{105}$ & $-\frac{1}{48}$ & $\phantom{+}\frac{4}{945}$ & & & \\ \hline
        $(4)$ & $\phantom{+}\frac{105}{16}$ & $-12$ & $\phantom{+}\frac{175}{32}$ & $\phantom{+}1$ & $-\frac{5}{4}$ & $\phantom{+}\frac{1}{6}$ & $\phantom{+}\frac{2}{35}$ & $\phantom{+}0$ & $-\frac{4}{315}$ & $\phantom{+}\frac{1}{180}$ & $-\frac{8}{10395}$ & $\phantom{+}\frac{1}{8640}$ \\ \hline
        $(3, 1)$ & $-\frac{15}{16}$ & $\phantom{+}6$ & $-\frac{295}{32}$ & $\phantom{+}5$ & $-\frac{1}{4}$ & $-\frac{2}{3}$ & $\phantom{+}\frac{4}{35}$ & $\phantom{+}\frac{1}{8}$ & $-\frac{4}{63}$ & $\phantom{+}\frac{1}{48}$ & $-\frac{8}{3465}$ & $\phantom{+}\frac{1}{2880}$ \\ \hline
        $(2, 2)$ & $\phantom{+}0$ & $\phantom{+}0$ & $\phantom{+}\frac{21}{16}$ & $-3$ & $\phantom{+}\frac{13}{5}$ & $-1$ & $\phantom{+}\frac{2}{35}$ & $\phantom{+}\frac{1}{8}$ & $-\frac{16}{315}$ & $\phantom{+}\frac{11}{720}$ & $-\frac{16}{10395}$ & $\phantom{+}\frac{1}{4320}$ \\ \hline
        $(2, 1, 1)$ & $\phantom{+}0$ & $-2$ & $\phantom{+}\frac{165}{32}$ & $-5$ & $\phantom{+}\frac{23}{12}$ & $\phantom{+}\frac{1}{3}$ & $-\frac{74}{105}$ & $\phantom{+}\frac{3}{8}$ & $-\frac{20}{189}$ & $\phantom{+}\frac{1}{40}$ & $-\frac{8}{3465}$ & $\phantom{+}\frac{1}{2880}$ \\ \hline
        $(1, 1, 1, 1)$ & $\phantom{+}0$ & $\phantom{+}0$ & $-\frac{69}{32}$ & $\phantom{+}5$ & $-\frac{99}{20}$ & $\phantom{+}\frac{17}{6}$ & $-\frac{16}{15}$ & $\phantom{+}\frac{7}{24}$ & $-\frac{52}{945}$ & $\phantom{+}\frac{7}{720}$ & $-\frac{8}{10395}$ & $\phantom{+}\frac{1}{8640}$ \\ \hline
    \end{tabular}
    \caption{The values of the coefficients $\glslink{theta-f-i}{\theta_{\glslink{X-lam}{X^\lambda}}^{(i)}}$ from Theorem~\ref{main-thm-1} for all partitions $\lambda$ of size at most $4$.} \label{theta-X-figure}
\end{figure}

\noindent
In the process of proving Theorem~\ref{main-thm-1}, we will also prove the following theorem, which provides explicit formulas for the values $\glslink{A-sigma-lambda}{A_{\glslink{id-k}{\id_k}}^{\glslink{lambda-bracket-n}{\lambda[n]}}}$ when $k \le n < k+\glslink{abs-rho}{\lvert\lambda\rvert}-1$ (see Figure~\ref{main-thms-fig}).

\begin{theorem} \label{main-thm-2}
    Let $\lambda$ be a partition of size at least one, and let $j \ge 0$. Then there exists a polynomial $\glslink{b-f-j}{b_{\glslink{X-lam}{X^\lambda}}^{(j)}} \in \Q[n]$ of degree at most $\glslink{abs-rho}{\lvert\lambda\rvert}+2j$ such that\begin{equation*}
        \glslink{A-sigma-lambda}{A_{\glslink{id-k}{\id_{n-j}}}^{\glslink{lambda-bracket-n}{\lambda[n]}}} = \frac{\glslink{b-f-j}{b_{\glslink{X-lam}{X^\lambda}}^{(j)}}(n)}{n!}
    \end{equation*} for all $n \ge \max(j, \glslink{abs-rho}{\lvert\lambda\rvert}-1)$.
\end{theorem}

\begin{figure}[b]
    \begin{tikzpicture}[scale=0.8]
        \foreach \x in {0,...,7}
            \draw[gray] (\x, 0) -- (\x, 7);
        \foreach \y in {0,...,7}
            \draw[gray] (0, \y) -- (7, \y);
        \foreach \x in {0,...,6}
            \draw[black] (\x+0.5, 0) node[below] {$\x$};
        \foreach \y in {0,...,6}
            \draw[black] (0, \y+0.5) node[left] {$\y$};
        \draw[black, thick, ->] (0, 0) -- (8, 0) node[right] {$n-k$};
        \draw[black, thick, ->] (0, 0) -- (0, 8) node[above] {$k$};

        \draw[black, thick] (0,3) -- (1,3) -- (1,2) -- (2,2) -- (2,1) -- (4,1) -- (4,0);
        \draw[black, thick] (3,1) -- (3,7);

        \foreach \x in {0,...,2} {
            \draw[orange, thick, <->] (\x+0.5, 3-\x+0.5) -- (\x+0.5, 6.5);
            \draw[orange] (\x+0.5, 4.5) node[fill=white] {$\scriptstyle \glslink{b-f-j}{b_{\glslink{X-lam}{X^\lambda}}^{(\x)}}$};
        }
        \draw[orange] (1.5, 5.5) node[fill=white] {Theorem~\ref{main-thm-2}};
        \fill[orange, opacity=0.1] (0,3) -- (1,3) -- (1,2) -- (2,2) -- (2,1) -- (3,1) -- (3,7) -- (0,7);
        \fill[purple, opacity=0.1] (0,3) -- (1,3) -- (1,2) -- (2,2) -- (2,1) -- (4,1) -- (4,0) -- (0,0);
        \draw[blue] (5, 5.5) node[fill=white] {Theorem~\ref{main-thm-1}};
        \fill[blue, opacity=0.1] (4,0) -- (7,0) -- (7,7) -- (3,7) -- (3,1) -- (4,1);
    \end{tikzpicture}
    \caption{A visual representation of how Theorems~\ref{main-thm-1} and \ref{main-thm-2} can be combined to obtain formulas for $\glslink{A-sigma-lambda}{A_{\glslink{id-k}{\id_k}}^{\glslink{lambda-bracket-n}{\lambda[n]}}}$ for various values of $n$ and $k$. Here, we assume $\lambda$ is a partition of size $4$, so Theorem~\ref{main-thm-1} can be applied to compute $\glslink{A-sigma-lambda}{A_{\glslink{id-k}{\id_k}}^{\glslink{lambda-bracket-n}{\lambda[n]}}}$ whenever $n-k \ge 3$ and $(n, k) \ne (3, 0)$ (the blue region). The polynomials $\glslink{b-f-j}{b_{\glslink{X-lam}{X^\lambda}}^{(0)}}$, $\glslink{b-f-j}{b_{\glslink{X-lam}{X^\lambda}}^{(1)}}$, and $\glslink{b-f-j}{b_{\glslink{X-lam}{X^\lambda}}^{(2)}}$ from Theorem~\ref{main-thm-2} are used to fill in most of the remaining ``gap'' where Theorem~\ref{main-thm-1} does not apply (the orange region). The purple region that is left over represents a finite number of values $\glslink{A-sigma-lambda}{A_{\glslink{id-k}{\id_k}}^{\glslink{lambda-bracket-n}{\lambda[n]}}}$ that can be computed separately.} \label{main-thms-fig}
\end{figure}

\noindent
As a final core result, we will show that the ``leading terms'' of the formulas in Theorems~\ref{main-thm-1} and \ref{main-thm-2} admit explicit expressions. Below and throughout the rest of the paper, we adopt the standard convention that $(-1)!! := 1$.

\begin{theorem} \label{main-thm-3}
    Let $\lambda$ be a partition of size at least one, and let $\glslink{theta-f-i}{\theta_{\glslink{X-lam}{X^\lambda}}^{(i)}}$ and $\glslink{b-f-j}{b_{\glslink{X-lam}{X^\lambda}}^{(j)}}$ be as in Theorems~\ref{main-thm-1} and \ref{main-thm-2}, respectively. \begin{itemize}
        \item[(i)] The last two coefficients $\glslink{theta-f-i}{\theta_{\glslink{X-lam}{X^\lambda}}^{(i)}}$ are given by \begin{equation*}
            \glslink{theta-f-i}{\theta_{\glslink{X-lam}{X^\lambda}}^{(3\glslink{abs-rho}{\lvert\lambda\rvert})}} = \frac{2^{\glslink{abs-rho}{\lvert\lambda\rvert}}\glslink{chi-lambda}{\chi^\lambda}(\glslink{id-k}{\id_{\glslink{abs-rho}{\lvert\lambda\rvert}}})}{(\glslink{abs-rho}{\lvert\lambda\rvert}-1)!!(3\glslink{abs-rho}{\lvert\lambda\rvert})!!} \qquad \text{and} \qquad \glslink{theta-f-i}{\theta_{\glslink{X-lam}{X^\lambda}}^{(3\glslink{abs-rho}{\lvert\lambda\rvert}-1)}} = -\frac{2^{\glslink{abs-rho}{\lvert\lambda\rvert}}\glslink{chi-lambda}{\chi^\lambda}(\glslink{id-k}{\id_{\glslink{abs-rho}{\lvert\lambda\rvert}}})}{(\glslink{abs-rho}{\lvert\lambda\rvert}-2)!!(3\glslink{abs-rho}{\lvert\lambda\rvert}-1)!!}.
        \end{equation*}
        \item[(ii)] For all $j \ge 0$, the coefficients of $n^{\glslink{abs-rho}{\lvert\lambda\rvert}+2j}$ and $n^{\glslink{abs-rho}{\lvert\lambda\rvert}+2j-1}$ in $\glslink{b-f-j}{b_{\glslink{X-lam}{X^\lambda}}^{(j)}}(n)$ are \begin{equation*}
            \frac{2^j\glslink{chi-lambda}{\chi^\lambda}(\glslink{id-k}{\id_{\glslink{abs-rho}{\lvert\lambda\rvert}}})}{(\glslink{abs-rho}{\lvert\lambda\rvert}-1)!!(\glslink{abs-rho}{\lvert\lambda\rvert}+2j)!!}
        \end{equation*} and \begin{equation*}
            -2^j\glslink{chi-lambda}{\chi^\lambda}(\glslink{id-k}{\id_{\glslink{abs-rho}{\lvert\lambda\rvert}}}) \left(\frac{\glslink{abs-rho}{\lvert\lambda\rvert}^2+(3j-1)\glslink{abs-rho}{\lvert\lambda\rvert}+2j(j-1)}{2(\glslink{abs-rho}{\lvert\lambda\rvert}-1)!!(\glslink{abs-rho}{\lvert\lambda\rvert}+2j)!!} + \frac{1}{(\glslink{abs-rho}{\lvert\lambda\rvert}-2)!!(\glslink{abs-rho}{\lvert\lambda\rvert}+2j-1)!!}\right),
        \end{equation*} respectively.
    \end{itemize}
\end{theorem}

We apply Theorems~\ref{main-thm-1}, \ref{main-thm-2}, and \ref{main-thm-3} in a variety of directions, starting with the conjectures of Gaetz and Pierson. With regard to Conjecture~\ref{root-conj}, Theorem~\ref{main-thm-1} vastly increases the ranges of values of $\lambda$ and $k$ for which it is feasible to compute the polynomials $\glslink{a-sigma-lambda}{a_{\glslink{id-k}{\id_k}}^\lambda}(n)$, which in particular gives rise to the following counterexample (verified using exact computation in the code accompanying this paper; see the code availability statement below):

\begin{counterexample} \label{roots-counterexample}
    The polynomial $\glslink{a-sigma-lambda}{a_{\glslink{id-k}{\id_{14}}}^{(2,1)}}(n)$ is given by \begin{align*}
        &\tfrac{768381017}{226851446883715670016000000}n^{11} - \tfrac{3475721021}{15123429792247711334400000}n^{10} + \tfrac{155888476969}{22685144688371567001600000}n^9 \\
        &- \tfrac{20018079127}{168038108802752348160000}n^8 + \tfrac{292128669229}{220458160236847104000000}n^7 - \tfrac{7116427733081}{720163323440367206400000}n^6 \\
        &+ \tfrac{29910904675037}{596977491799251763200000}n^5 - \tfrac{6480061684471}{37808574480619278336000}n^4 + \tfrac{2720816654300431}{7089107715116114688000000}n^3 \\
        &- \tfrac{166778615642761}{315071454005160652800000}n^2 + \tfrac{449551268741}{1125255192875573760000}n - \tfrac{222663967}{1786119353770752000},
    \end{align*} and has degree $11$ but only $9$ (simple) real roots.
\end{counterexample}

\noindent
At the same time, Theorem~\ref{main-thm-1} allows us to prove the following weakened form of Conjecture~\ref{root-conj}, which shows that the polynomials $\glslink{a-sigma-lambda}{a_{\glslink{id-k}{\id_k}}^\lambda}$ are ``close'' to being real-rooted.

\begin{corollary} \label{real-roots-cor}
    Let $\lambda$ be a partition of size at least $1$, and let $k \ge \glslink{abs-rho}{\lvert\lambda\rvert}$. Assume that the polynomial $\glslink{a-sigma-lambda}{a_{\glslink{id-k}{\id_k}}^\lambda}(n)$ is not identically zero. Then $\glslink{a-sigma-lambda}{a_{\glslink{id-k}{\id_k}}^\lambda}(n)$ has at most $\glslink{abs-rho}{\lvert\lambda\rvert}-1$ conjugate pairs of non-real roots.
\end{corollary}

\noindent
As for Conjecture~\ref{pos-conj}, Theorems~\ref{main-thm-1}, \ref{main-thm-2}, and \ref{main-thm-3} together imply that the coefficients $\glslink{A-sigma-lambda}{A_{\glslink{id-k}{\id_k}}^{\glslink{lambda-bracket-n}{\lambda[n]}}}$ are ``eventually'' nonnegative in the following sense:

\begin{corollary} \label{eventually-pos-cor}
    For any partition $\lambda$, there exists an effectively computable integer $t_\lambda \ge 0$ such that $\glslink{A-sigma-lambda}{A_{\glslink{id-k}{\id_k}}^{\glslink{lambda-bracket-n}{\lambda[n]}}} \ge 0$ for all $n \ge k \ge t_\lambda$.
\end{corollary}

\noindent
Indeed, by implementing the algorithm behind Corollary~\ref{eventually-pos-cor} on a computer and explicitly calculating the values $t_\lambda$ for small values of $\lambda$, we obtain the following result.

\begin{proposition} \label{nonnegativity-prop}
    For all partitions $\lambda$ of size at most $10$, we have $\glslink{A-sigma-lambda}{A_{\glslink{id-k}{\id_k}}^{\glslink{lambda-bracket-n}{\lambda[n]}}} \ge 0$ for all $n \ge k \ge 0$ with $n \ge \glslink{abs-rho}{\lvert\lambda\rvert}+\glslink{lambda-i}{\lambda_1}$.
\end{proposition}

Beyond the conjectures of Gaetz and Pierson, we use Theorem~\ref{main-thm-1} to deduce the following ``local rule'' satisfied by the polynomials $\glslink{a-sigma-lambda}{a_{\glslink{id-k}{\id_k}}^\lambda}(n)$.

\begin{corollary} \label{local-rule-cor}
    Let $\lambda$ be any partition. Then for all $k \ge 0$, we have \begin{equation*}
        (n+1)\glslink{a-sigma-lambda}{a_{\glslink{id-k}{\id_{k+1}}}^\lambda}(n+1) = \glslink{a-sigma-lambda}{a_{\glslink{id-k}{\id_k}}^\lambda}(n) + (n+k+2-\glslink{abs-rho}{\lvert\lambda\rvert}) \glslink{a-sigma-lambda}{a_{\glslink{id-k}{\id_{k+1}}}^\lambda}(n)
    \end{equation*} as polynomials in $n$.
\end{corollary}

\noindent
We also provide several applications that are purely combinatorial. For instance, we use Theorems~\ref{main-thm-1}--\ref{main-thm-3} to show that the covariances of $\glslink{N-sigma}{N_{\glslink{id-k}{\id_k}}(\pi)}$ with certain families of class functions $\bigsqcup_{n\ge0} \glslink{symmetric-group}{S_n} \to \Q$ admit explicit formulas.

\begin{theorem} \label{main-thm-combo}
    For all $j \ge 1$, let $\glslink{m-j}{m_j} \colon \bigsqcup_{n\ge0} \glslink{symmetric-group}{S_n} \to \Z_{\ge0}$ denote the function which measures the number of $j$-cycles in its input, and let $f \in \Q[\glslink{m-j}{m_1},\glslink{m-j}{m_2},\dots]$ be an arbitrary nonconstant polynomial. Let $d$ be the total degree of $f$, where we treat each variable $\glslink{m-j}{m_j}$ as having degree $j$. For each $n \ge k \ge 0$, define \begin{equation*}
        \glslink{E-f-n-k}{E_f(n, k)} := \frac{1}{n!} \sum_{\pi \in \glslink{symmetric-group}{S_n}} f(\pi) \cdot \glslink{N-sigma}{N_{\glslink{id-k}{\id_k}}(\pi)}.
    \end{equation*} Then: \begin{itemize}
        \item[(i)] There exist coefficients $\glslink{theta-f-i-ell}{\theta_f^{i,\ell}}$ such that \begin{equation*}
            \glslink{E-f-n-k}{E_f(n, k)} = \glslink{theta-f-i-ell}{\theta_f^{0,0}} \cdot \frac{1}{k!} \glslink{binomial}{\binom{n}{k}} +  \sum_{\substack{1 \le \ell \le d \\ 1 \le i \le 3\ell}} \glslink{theta-f-i-ell}{\theta_f^{i,\ell}} \cdot \frac{(n-\ell)!}{n!} \cdot \frac{(\frac{i}{2})!}{(k-\frac{i}{2})!} \glslink{binomial}{\binom{n-\ell+\frac{i}{2}}{k}}
        \end{equation*} for all $k \ge 0$ and $n \ge k+d-1$ with $(n, k) \ne (d-1, 0)$.
        \item[(ii)] For fixed $k \ge 0$, there exists a polynomial $\glslink{a-f-k}{a_f^{(k)}} \in \Q[n]$ of degree at most $k$ such that \begin{equation*}
            \glslink{E-f-n-k}{E_f(n, k)} = \glslink{a-f-k}{a_f^{(k)}}(n)
        \end{equation*} for all $n \ge k+d-1$ with $(n, k) \ne (d-1, 0)$.
        \item[(iii)] For all $j \ge 0$, there exists a polynomial $\glslink{b-f-j}{b_f^{(j)}} \in \Q[n]$ of degree at most $d+2j$ such that \begin{equation*}
            \glslink{E-f-n-k}{E_f(n, n-j)} = \frac{\glslink{b-f-j}{b_f^{(j)}}(n)}{n!}
        \end{equation*} for all $n \ge \max(j, d-1)$.
    \end{itemize}
\end{theorem}

\noindent
As an immediate corollary, we deduce that, somewhat counterintuitively, the information of $\glslink{E-f-n-k}{E_f(n,k)}$ for all $n \ge k \ge 0$ is not enough to recover the polynomial $f$.

\begin{corollary} \label{kernel-cor}
    Adopting the notation of Theorem~\ref{main-thm-combo}, there exist nonzero polynomials $f \in \Q[\glslink{m-j}{m_1},\glslink{m-j}{m_2},\dots]$ of arbitrarily large degree such that $\glslink{E-f-n-k}{E_f(n, k)} = 0$ for all $n \ge k \ge 0$. In particular, there exists such a polynomial with degree at most $24$.
\end{corollary}

\noindent
While explicit examples of such polynomials $f$ remain out of reach for the methods of this paper, our methods do make it possible to write down examples of polynomials $f$ that ``come close'' to satisfying the condition of Corollary~\ref{kernel-cor}.

\begin{example} \label{E-f-n-k-zero-example}
    Take $f = 23\glslink{m-j}{m}_2 - 37\glslink{m-j}{m}_2^2 + 16\glslink{m-j}{m}_2^3 - 2\glslink{m-j}{m}_2^4 + 9\glslink{m-j}{m}_4 - 6\glslink{m-j}{m}_4^2 - 3\glslink{m-j}{m}_6 + 6\glslink{m-j}{m}_2\glslink{m-j}{m}_6$. Then $\glslink{E-f-n-k}{E_f(n, k)} = 0$ for all $n \ge k \ge 0$ with $n-k \ge 7$.
\end{example}

\noindent
As a final combinatorial result, we prove an approximate formula for the conditional expectation of $\glslink{N-sigma}{N_{\glslink{id-k}{\id_k}}(\pi)}$ given certain ``local group theoretic information'' about a random permutation $\pi$.

\begin{theorem} \label{conditional-expectation-thm}
    Fix a partition $\rho$ and a real number $0 < \epsilon < \frac{1}{2}$, and for each $j \ge 1$, let $\glslink{m-j-rho}{m_j(\rho)}$ denote the number of parts of size $j$ in $\rho$. Let $n \ge \max(\glslink{abs-rho}{\lvert\rho\rvert}, 1)$ and $\epsilon n \le k \le (1-\epsilon)n$, and set $t := \frac{k}{n} \in [\epsilon, 1-\epsilon]$. Suppose a permutation $\pi \in \glslink{symmetric-group}{S_n}$ and a subset $R \subseteq [n]$ of size $\lvert\rho\rvert$ are selected uniformly at random. Then \begin{multline*}
        \frac{\E[\glslink{N-sigma}{N_{\glslink{id-k}{\id_k}}(\pi)} \mid \text{$\pi(R) = R$ and $\pi\rvert_R$ has cycle type $\rho$}]}{\E[\glslink{N-sigma}{N_{\glslink{id-k}{\id_k}}(\pi)}]} \\
        = \left(\frac{\glslink{m-j-rho}{m_1(\rho)}}{2}\right)! \cdot \left(\frac{nt^3}{1-t}\right)^{\frac{\glslink{m-j-rho}{m_1(\rho)}}{2}} \left(\prod_{j\ge2} (1-t^j)^{\glslink{m-j-rho}{m_j(\rho)}}\right) + O\left(n^{\frac{\glslink{m-j-rho}{m_1(\rho)}-1}{2}}\right),
    \end{multline*} where the implicit constant of the big $O$ depends only on $\rho$ and $\epsilon$.
\end{theorem}

\noindent
Informally, Theorem~\ref{conditional-expectation-thm} suggests that ``permutations with lots of fixed points tend to have more long increasing subsequences,'' whereas ``permutations with lots of small cycles of size at least $2$ tend to have fewer long increasing subsequences.''

The outline of this paper is as follows. In Section~\ref{background-on-character-polynomials-section}, we review background on the ring of \glslink{character-polynomial}{character polynomials} that will be used throughout the paper and show how Theorem~\ref{main-thm-combo} follows as a simple consequence of Theorems~\ref{main-thm-1} and \ref{main-thm-2}. In Section~\ref{stitches-section}, we introduce the notion of \glslink{stitch}{\textit{stitches}}---partial permutations which we choose to visualize in a particular way---along with a variety of key definitions and properties surrounding them. In Section~\ref{proof-overview-section}, we use the background from Sections~\ref{background-on-character-polynomials-section} and \ref{stitches-section} to give a broad overview of the structure of the proof of Theorems~\ref{main-thm-1}, \ref{main-thm-2} and \ref{main-thm-3}. In Section~\ref{counting-parallel-stitches-section}, we prove a variety of explicit formulas for enumerating certain families of \glslink{stitch}{stitches}. In Section~\ref{generating-functions-for-colored-stitches-section}, we then apply the results of Section~\ref{counting-parallel-stitches-section} to study generating functions for certain families of \glslink{stitch-coloring}{colored} \glslink{stitch}{stitches}. In Section~\ref{explicit-formulas-for-H-section}, we use the results of Section~\ref{generating-functions-for-colored-stitches-section} to obtain explicit formulas for a certain family of generating functions which are more immediately relevant to the specific task of the paper. In Section~\ref{explicit-formulas-for-E-section}, we use the explicit formulas derived in Section~\ref{explicit-formulas-for-H-section} to prove Theorems~\ref{main-thm-1}, \ref{main-thm-2}, and \ref{main-thm-3}, along with several slight variations of these theorems. Finally, in Section~\ref{applications-section}, we give proofs of the various corollaries and applications listed in the introduction.

\subsection*{Code availability}

Several of the claims of this paper rely on computer-assisted computations. A SageMath notebook implementing the definitions and formulas in this paper is available at \url{https://github.com/JonasIskander/character-polynomials}. In particular, the notebook provides reproducible verifications of the computer-assisted claims in Counterexample~\ref{roots-counterexample}, Proposition~\ref{nonnegativity-prop}, and Example~\ref{E-f-n-k-zero-example}, as well as the computations underlying Figures~\ref{theta-X-figure} and \ref{theta-J-figure} and numerical checks of the main theorems.

\section*{Acknowledgements}

This work was done as part of the University of Minnesota Duluth REU with support from Jane Street Capital, the NSA (grant number H98230-22-1-0015), the NSF (grant number DMS-2052036), and Harvard University. The author thanks Joseph Gallian for organizing the REU, as well as Amanda Burcroff, Noah Kravitz, Mitchell Lee, Yelena Mandelshtam, Maya Sankar and Katherine Tung for their helpful advice and suggestions. In particular, the proof of Lemma~\ref{cayley-lem} is based on ideas of Mitchell Lee.

During the research and writing of this paper, both SageMath and Mathematica were used for exploratory computations, pattern-finding, and algebraic manipulations. Various LLMs including Claude, ChatGPT, and Gemini were consulted throughout the writing of this paper, primarily to search for references for obscure propositions or identities or help with organization and proofreading, but also occasionally to assist with the proof process. In particular, the proofs of Propositions~\ref{stitch-unique-factorization-prop} and \ref{irreducible-parallel-classification-prop} are based on suggestions from ChatGPT, and parts of the proof of Proposition~\ref{J-prop} given in Appendix~\ref{J-appendix} are based on conversations with Claude.

\section{Background on Character Polynomials} \label{background-on-character-polynomials-section}

In this section, we recount a number of (mostly) standard definitions and propositions relating to the theory of \glslink{character-polynomial}{character polynomials} that we will use throughout the paper. We loosely follow the exposition in Sections~1--2 of \cite{narayanan-et-al}, but adopt different notation and terminology better suited to the context of this paper.

We start with the following definition, which allows us to more easily refer to the polynomials appearing in Theorem~\ref{main-thm-combo}.

\begin{definition}[cf. Section 1 of \cite{narayanan-et-al}] \label{m-j-def}
    For all $j \ge 1$, define the \textit{$j$th cycle-counting function} $\glslink{m-j}{m_j} \colon \bigsqcup_{n\ge0} \glslink{symmetric-group}{S_n} \to \Z_{\ge0}$ by \begin{equation*}
        \glslink{m-j}{m_j}(\pi) := \text{the number of $j$-cycles in $\pi$}.
    \end{equation*} A \glslink{character-polynomial}{\textit{character polynomial}} (over $\Q$) is a function $f \colon \bigsqcup_{n\ge0} \glslink{symmetric-group}{S_n} \to \Q$ which can be expressed as a polynomial in the cycle-counting functions $\glslink{m-j}{m_j}$. We denote the \glslink{character-polynomial}{\textit{ring of character polynomials}} (over $\Q$) by $\glslink{Delta-Q}{\Delta_\Q}$. We equip this ring with a degree map $\glslink{character-polynomial-deg}{\deg} \colon \glslink{Delta-Q}{\Delta_\Q} \to \{-\infty\} \cup \Z_{\ge0}$ determined by the property that $\glslink{character-polynomial-deg}{\deg(\glslink{m-j}{m_j})} := j$ for all $j$.
\end{definition}

\begin{remark}
    The authors of \cite{narayanan-et-al} use variables $X_j$ in place of $\glslink{m-j}{m_j}$ and $P$ in place of $\glslink{Delta-Q}{\Delta_\Q}$. In addition, they equip $P$ with a graded structure in which the variable $X_j$ is homogeneous of degree $j$. However, to avoid confusion with a different notion of ``homogeneity'' we will introduce later in this section, we choose to only endow $\glslink{Delta-Q}{\Delta_\Q}$ with the structure of a \textit{filtered} ring, so the variables $\glslink{m-j}{m_j}$ are not considered homogeneous.
\end{remark}

Like the ring of symmetric functions, the \glslink{character-polynomial}{ring of character polynomials} admits a number of natural bases indexed by integer partitions. For the purposes of this paper, three main bases will prove useful, which we will denote by $\glslink{I-rho}{I_\rho}$, $\glslink{X-lam}{X^\lambda}$, and $\glslink{J-rho}{J_\rho}$. The first of these is defined below.

\begin{definition}[cf.\@ Definition 2.2 of \cite{narayanan-et-al}]
    Let $\rho$ be a partition, and for each $j \ge 1$, let $\glslink{m-j-rho}{m_j(\rho)}$ denote the number of parts of size $j$ in $\rho$. We define \begin{equation*}
        \glslink{I-rho}{I_\rho} := \prod_{j\ge1} \glslink{binomial}{\binom{\glslink{m-j}{m_j}}{\glslink{m-j-rho}{m_j(\rho)}}} \in \glslink{Delta-Q}{\Delta_\Q},
    \end{equation*} where the $\glslink{m-j}{m_j}$ on the top of the binomial refers to the function $\glslink{m-j}{m_j} \in \glslink{Delta-Q}{\Delta_\Q}$ from Definition~\ref{m-j-def}.
\end{definition}

\noindent
Combinatorially, for any permutation $\pi \in \glslink{symmetric-group}{S_n}$, $\glslink{I-rho}{I_\rho}(\pi)$ counts subsets $R \subseteq \glslink{bracket-n}{[n]}$ of size $\glslink{abs-rho}{\lvert\rho\rvert}$ for which $\pi(R) = R$ and $\pi\rvert_R$ has cycle type $\rho$. This characterization makes the $\glslink{I-rho}{I_\rho}$ the go-to basis when one wishes to study \glslink{character-polynomial}{character polynomials} using enumerative methods.

In the same way that the $\glslink{I-rho}{I_\rho}$ provide a natural combinatorial interpretation for elements of $\glslink{Delta-Q}{\Delta_\Q}$, the next key basis for our purposes, $\glslink{X-lam}{X^\lambda}$, gives a natural interpretation of elements of $\glslink{Delta-Q}{\Delta_\Q}$ in terms of symmetric group characters. Specifically, we have the following theorem.

\begin{theorem}[Macdonald {\cite[pp.~123--124]{macdonald}}] \label{macdonald-thm}
    For each partition $\lambda$, there exists a unique polynomial $\glslink{X-lam}{X^\lambda} \in \glslink{Delta-Q}{\Delta_\Q}$ such that for all $n \ge \glslink{abs-rho}{\lvert\lambda\rvert} + \glslink{lambda-i}{\lambda_1}$ and $\pi \in \glslink{symmetric-group}{S_n}$, we have \begin{equation*}
        \glslink{chi-lambda}{\chi^{\glslink{lambda-bracket-n}{\lambda[n]}}}(\pi) = \glslink{X-lam}{X^\lambda}(\pi).
    \end{equation*} Explicitly, setting $\glslink{Y-lam}{Y^\lambda} := \sum_{\rho \vdash \glslink{abs-rho}{\lvert\lambda\rvert}} \glslink{chi-rho-lambda}{\chi_\rho^\lambda} \cdot \glslink{I-rho}{I_\rho}$ for all $\lambda$, one has \begin{equation*}
        \glslink{X-lam}{X^\lambda} = \sum_{\substack{\glslink{prec-prime}{\mu \prec' \lambda}}} (-1)^{\glslink{abs-rho}{\lvert\lambda\rvert}-\glslink{abs-rho}{\lvert\mu\rvert}} \cdot \glslink{Y-lam}{Y^\mu},
    \end{equation*} where $\glslink{prec-prime}{\mu \prec' \lambda}$ indicates that $\mu$ may be obtained from $\lambda$ by removing at most one box in each row.
\end{theorem}

\noindent
In light of the above theorem, it is natural to extend the definition of $\glslink{A-sigma-lambda}{A_{\sigma_1,\dots,\sigma_d}^{\glslink{lambda-bracket-n}{\lambda[n]}}}$ from the introduction to the range $0 \le n < \glslink{abs-rho}{\lvert\lambda\rvert}+\glslink{lambda-i}{\lambda_1}$ by setting \begin{equation}
    \glslink{A-sigma-lambda}{A_{\sigma_1,\dots,\sigma_d}^{\glslink{lambda-bracket-n}{\lambda[n]}}} := \frac{1}{n!} \sum_{\pi \in \glslink{symmetric-group}{S_n}} \glslink{X-lam}{X^\lambda}(\pi) \cdot \glslink{N-sigma}{N_{\sigma_1}(\pi)} \cdots \glslink{N-sigma}{N_{\sigma_d}(\pi)}. \label{A-redef}
\end{equation} In particular, in the notation of Theorem~\ref{main-thm-combo}, we see that $\glslink{A-sigma-lambda}{A_{\glslink{id-k}{\id_k}}^{\glslink{lambda-bracket-n}{\lambda[n]}}}$ is precisely the quantity $\glslink{E-f-n-k}{E_f(n, k)}$ in the case where $f = \glslink{X-lam}{X^\lambda}$.

Before introducing the final key basis $\glslink{J-rho}{J_\rho}$, which the author is not aware of having been studied previously, we must review some standard structural properties of the ring $\glslink{Delta-Q}{\Delta_\Q}$. We start by recalling the notion of \textit{moments} of a \glslink{character-polynomial}{character polynomial}.

\begin{definition} [cf.\@ Definition 2.1 of \cite{narayanan-et-al}]
    Let $f \in \glslink{Delta-Q}{\Delta_\Q}$, and let $n \ge 0$. We define the \textit{$n$th moment} of $f$ to be the expected value of $f$ on $\glslink{symmetric-group}{S_n}$: \begin{equation*}
        \glslink{nth-moment}{\langle f\rangle_n} := \frac{1}{n!} \sum_{\pi \in \glslink{symmetric-group}{S_n}} f(\pi).
    \end{equation*}
\end{definition}

\noindent
A central result in the theory of \glslink{character-polynomial}{character polynomials} is the fact that moments eventually stabilize:

\begin{proposition}[cf.\@ Theorem 2.3 and Corollary 2.4 of \cite{narayanan-et-al}] \label{I-rho-expectation-prop}
    For any partition $\rho$ and any $n \ge 0$, we have \begin{equation*}
        \glslink{nth-moment}{\langle \glslink{I-rho}{I_\rho}\rangle_n} = \begin{cases}
            \frac{1}{\glslink{z-rho}{z_\rho}} & \text{if $n \ge \glslink{abs-rho}{\lvert\rho\rvert}$,} \\
            0 & \text{otherwise,}
        \end{cases}
    \end{equation*} where $\glslink{z-rho}{z_\rho}$ denotes the size of the centralizer of a permutation of cycle type $\rho$. Consequently, for any polynomial $f \in \glslink{Delta-Q}{\Delta_\Q}$, the moments $\glslink{nth-moment}{\langle f\rangle_n}$ achieve a single constant value for $n \ge \glslink{character-polynomial-deg}{\deg(f)}$.
\end{proposition}

\begin{definition}[cf.\@ Definition 2.5 of \cite{narayanan-et-al}]
    Let $f \in \glslink{Delta-Q}{\Delta_\Q}$. We define the \textit{stable moment} of $f$ to be the eventual constant value of the $n$th moments: \begin{equation*}
        \glslink{stable-moment}{\langle f\rangle} := \lim_{n\to\infty} \glslink{nth-moment}{\langle f\rangle_n}.
    \end{equation*}
\end{definition}

The existence of stable moments gives rise to a natural inner product structure on $\glslink{Delta-Q}{\Delta_\Q}$.

\begin{definition}
    We equip $\glslink{Delta-Q}{\Delta_\Q}$ with an inner product $\glslink{inner-product}{\langle \cdot, \cdot\rangle}$ defined by $\glslink{inner-product}{\langle f, g\rangle} := \glslink{stable-moment}{\langle fg\rangle}$.
\end{definition}

\noindent
In light of this additional structure, $\glslink{Delta-Q}{\Delta_\Q}$ admits a canonical decomposition into orthogonal subspaces.

\begin{definition} \label{pure-def}
    For all $d \ge 0$, let $\glslink{Delta-Q-le-d}{\Delta_\Q^{\le d}} := \{f \in \glslink{Delta-Q}{\Delta_\Q} : \glslink{character-polynomial-deg}{\deg(f)} \le d\}$. We define the \glslink{pure-character-polynomial}{\textit{space of pure \glslink{character-polynomial}{character polynomials} of degree $d$}}, denoted $\glslink{Delta-Q-d}{\Delta_\Q^d}$, to be the orthogonal complement of $\glslink{Delta-Q-le-d}{\Delta_\Q^{\le(d-1)}}$ inside $\glslink{Delta-Q-le-d}{\Delta_\Q^{\le d}}$. Given an element $f \in \glslink{Delta-Q}{\Delta_\Q}$, we say that $f$ is \glslink{pure-character-polynomial}{\textit{pure of degree $d$}} if $f \in \glslink{Delta-Q-d}{\Delta_\Q^d}$. In particular, in this terminology, the $0$ \glslink{character-polynomial}{character polynomial} is \glslink{pure-character-polynomial}{pure} of every degree.
\end{definition}

\noindent
With the above definition, we may decompose $\glslink{Delta-Q}{\Delta_\Q}$ as a direct sum of pairwise orthogonal spaces $\bigoplus_{d\ge0} \glslink{Delta-Q-d}{\Delta_\Q^d}$. The following fact, which is an easy corollary to Theorem~\ref{macdonald-thm}, further suggests why both the inner product on $\glslink{Delta-Q}{\Delta_\Q}$ and the subspaces $\glslink{Delta-Q-d}{\Delta_\Q^d}$ are particularly natural to study:

\begin{proposition}
    The $\glslink{X-lam}{X^\lambda}$ form an orthonormal basis of $\glslink{Delta-Q}{\Delta_\Q}$ with respect to the inner product $\glslink{inner-product}{\langle \cdot,\cdot\rangle}$. Moreover, for fixed $d \ge 0$, the $\glslink{X-lam}{X^\lambda}$ ranging over $\lambda \vdash d$ form an orthonormal basis for $\glslink{Delta-Q-d}{\Delta_\Q^d}$.
\end{proposition}

\begin{warning}
    While the decomposition $\glslink{Delta-Q}{\Delta_\Q} = \bigoplus_{d\ge0} \glslink{Delta-Q-d}{\Delta_\Q^d}$ is reminiscent of the decomposition of a graded ring into its homogeneous pieces, this decomposition does not endow $\glslink{Delta-Q}{\Delta_\Q}$ with the structure of a graded ring because the product of an element of $\glslink{Delta-Q-d}{\Delta_\Q^{d_1}}$ with an element of $\glslink{Delta-Q-d}{\Delta_\Q^{d_2}}$ is not necessarily \glslink{pure-character-polynomial}{pure of degree $d_1+d_2$}. For example, one verifies that $\glslink{X-lam}{X^{(1)}} \cdot \glslink{X-lam}{X^{(1)}} = \glslink{X-lam}{X^{(2)}} + \glslink{X-lam}{X^{(1,1)}} + \glslink{X-lam}{X^{(1)}} + \glslink{X-lam}{X^\varnothing} \notin \glslink{Delta-Q-d}{\Delta_\Q^2}$.
\end{warning}

With the above definitions, we may now define the third and final basis for $\glslink{Delta-Q}{\Delta_\Q}$ that we will use.

\begin{definition}
    For each partition $\rho$, we define $\glslink{J-rho}{J_\rho}$ to be the orthogonal projection of $\glslink{I-rho}{I_\rho}$ onto the subspace of \glslink{pure-character-polynomial}{pure} \glslink{character-polynomial}{character polynomials} $\glslink{Delta-Q-d}{\Delta_\Q^{\glslink{abs-rho}{\lvert\rho\rvert}}}$.
\end{definition}

\noindent
In the context of this paper, it will turn out that \glslink{pure-character-polynomial}{pure} \glslink{character-polynomial}{character polynomials} behave particularly nicely in certain formulas. For our purposes, the advantage of the basis $\glslink{J-rho}{J_\rho}$ will come from the fact that it combines the \glslink{pure-character-polynomial}{pureness} of the basis $\glslink{X-lam}{X^\lambda}$ with the explicit combinatorial interpretability of the basis $\glslink{I-rho}{I_\rho}$. For the reader's convenience, we quickly summarize the key properties of the polynomials $\glslink{J-rho}{J_\rho}$ in the proposition below; see Appendix~\ref{J-appendix} for proofs of these properties.

\begin{proposition} \label{J-prop}
    \phantom{.}\begin{itemize}
        \item[(i)] The $\glslink{J-rho}{J_\rho}$ form an orthogonal basis of $\glslink{Delta-Q}{\Delta_\Q}$ with respect to the inner product $\glslink{inner-product}{\langle\cdot,\cdot\rangle}$ satisfying $\glslink{inner-product}{\langle \glslink{J-rho}{J_\rho}, \glslink{J-rho}{J_\rho}\rangle} = \frac{1}{\glslink{z-rho}{z_\rho}}$. Moreover, for fixed $d \ge 0$, the $\glslink{J-rho}{J_\rho}$ ranging over $\rho \vdash d$ form an orthogonal basis for $\glslink{Delta-Q-d}{\Delta_\Q^d}$.
        \item[(ii)] The bases $\glslink{X-lam}{X^\lambda}$ and $\glslink{J-rho}{J_\rho}$ are related by the equations \begin{equation*}
            \glslink{J-rho}{J_\rho} = \frac{1}{\glslink{z-rho}{z_\rho}} \sum_{\lambda \vdash \glslink{abs-rho}{\lvert\rho\rvert}} \glslink{chi-rho-lambda}{\chi_\rho^\lambda} \cdot \glslink{X-lam}{X^\lambda} \qquad \text{and} \qquad \glslink{X-lam}{X^\lambda} = \sum_{\rho \vdash \glslink{abs-rho}{\lvert\lambda\rvert}} \glslink{chi-rho-lambda}{\chi_\rho^\lambda} \cdot \glslink{J-rho}{J_\rho}.
        \end{equation*}
        \item[(iii)] The bases $\glslink{I-rho}{I_\rho}$ and $\glslink{J-rho}{J_\rho}$ are related by the equations \begin{equation*}
            \glslink{J-rho}{J_\rho} = \sum_{\rho' \sqcup \mu = \rho} \frac{(-1)^{\glslink{len}{\len(\mu)}}}{\glslink{z-rho}{z_\mu}} \cdot \glslink{I-rho}{I_{\rho'}} \qquad \text{and} \qquad \glslink{I-rho}{I_\rho} = \sum_{\rho' \sqcup \mu = \rho} \frac{1}{\glslink{z-rho}{z_\mu}} \cdot \glslink{J-rho}{J_{\rho'}},
        \end{equation*} where the sums are over pairs of partitions $\rho'$ and $\mu$ with multiset sum $\rho$.
    \end{itemize}
\end{proposition}

For completeness, we include a short proof of Theorem~\ref{main-thm-combo} assuming Theorems~\ref{main-thm-1} and \ref{main-thm-2}, as this implication follows quickly from the background in this section.

\begin{proof}[Proof of Theorem~\ref{main-thm-combo}]
    When $f = \glslink{X-lam}{X^\lambda}$ for some $\lambda$ of size at least $1$, claims (i-ii) reduce to Theorem~\ref{main-thm-1}, whereas claim (iii) reduces to Theorem~\ref{main-thm-2}. Moreover, it is a standard result that $\glslink{E-f-n-k}{E_{\glslink{X-lam}{X^\varnothing}}(n, k)} = \frac{1}{n!} \sum_{\pi \in \glslink{symmetric-group}{S_n}} \glslink{N-sigma}{N_{\glslink{id-k}{\id_k}}(\pi)} = \frac{1}{k!} \glslink{binomial}{\binom{n}{k}}$ for all $n \ge k \ge 0$. Since the $\glslink{X-lam}{X^\lambda}$ form a basis for $\glslink{Delta-Q}{\Delta_\Q}$ and $\glslink{E-f-n-k}{E_f(n, k)}$ is linear in $f$, the result for arbitrary $f \in \glslink{Delta-Q}{\Delta_\Q}$ follows immediately.
\end{proof}

\section{Stitches} \label{stitches-section}

\subsection{Basic definitions}

Partial permutations will be ubiquitous throughout our arguments, and we will find it useful to adopt a particular way of visualizing them together with a broad vocabulary surrounding them. So as to avoid confusion with any pre-existing notation or terminology surrounding partial permutations, as well as to emphasize our nonstandard visualization of them, we will adopt the term ``\glslink{stitch}{stitch}'' for a partial permutation as it is treated in this paper.

\begin{definition}
    A \glslink{stitch}{\textit{stitch}} is a tuple $\alpha = (C, W)$, where $C$ is a finite totally ordered set and $W \subseteq C \times C$ is a set of pairs such that for each $i \in C$, there is at most one $j \in C$ such that $(i, j) \in W$ and at most one $j' \in C$ such that $(j', i) \in W$. We identify two \glslink{stitch}{stitches} $\alpha = (C, W)$ and $\alpha' = (C', W')$ if $\lvert C\rvert = \lvert C'\rvert$ and the unique order-preserving bijection $C \to C'$ sends $W$ to $W'$; thus, every \glslink{stitch}{stitch} may be written in a unique \glslink{standard-form}{\textit{standard form}} as $(\glslink{bracket-n}{[n]}, W)$ for some $n \ge 0$ and $W \subseteq \glslink{bracket-n}{[n]} \times \glslink{bracket-n}{[n]}$. We visualize a \glslink{stitch}{stitch} $\alpha = (C, W)$ by drawing two rows of $\lvert C\rvert$ points labeled by the elements of $C$, with a line segment connecting the lower point with label $i$ to the upper point with label $j$ for each $(i, j) \in W$ (see Figure~\ref{wire-diagram-example-fig}); we refer to such a picture as the \glslink{wire-diagram}{\textit{wire diagram}} of $\alpha$. We refer to the elements of $C$ as the \glslink{column}{\textit{columns}} of $\alpha$ and the elements of $W$ as the \glslink{wire}{\textit{wires}} of $\alpha$. We denote the set of all \glslink{stitch}{stitches} by $\glslink{calS}{\calS}$, the set of \glslink{stitch}{stitches} with $n$ \glslink{column}{columns} by $\glslink{calS-n}{\calS_n}$, and the set of \glslink{stitch}{stitches} with $n$ \glslink{column}{columns} and $k$ \glslink{wire}{wires} by $\glslink{calS-n-k}{\calS_{n,k}}$. In addition, we treat $\glslink{symmetric-group}{S_n}$ as a subset of $\glslink{calS-n}{\calS_n}$ by identifying each permutation $\pi \in \glslink{symmetric-group}{S_n}$ with the \glslink{stitch}{stitch} $(\glslink{bracket-n}{[n]}, \{(i, \pi(i)) : i \in \glslink{bracket-n}{[n]}\})$.
\end{definition}

\noindent
Because this paper focuses on \glslink{pattern-occurrence}{occurrences} of the pattern $\glslink{id-k}{\id_k}$---that is, increasing subsequences---in permutations, the following definition will also be central.

\begin{definition}
    Let $\alpha = (C, W)$ be a \glslink{stitch}{stitch}. We say that $\alpha$ is \glslink{parallel-stitch}{\textit{parallel}} if for any two distinct \glslink{wire}{wires} $(i_1, j_1), (i_2, j_2) \in W$, we have $i_1 < i_2$ if and only if $j_1 < j_2$---or equivalently, if no two \glslink{wire}{wires} in the \glslink{wire-diagram}{wire diagram} of $\alpha$ intersect (see Figure~\ref{wire-diagram-example-fig}). We denote the set of all \glslink{parallel-stitch}{parallel} \glslink{stitch}{stitches} by $\glslink{calS-parallel}{\calS^\parallel}$, the set of \glslink{parallel-stitch}{parallel} \glslink{stitch}{stitches} with $n$ \glslink{column}{columns} by $\glslink{calS-n-parallel}{\calS_n^\parallel}$, and the set of \glslink{parallel-stitch}{parallel} \glslink{stitch}{stitches} with $n$ \glslink{column}{columns} and $k$ \glslink{wire}{wires} by $\glslink{calS-n-k-parallel}{\calS_{n,k}^\parallel}$.
\end{definition}

\begin{figure}[h!]
    \centering
    \begin{subfigure}{\linewidth}
        \centering
        \begin{tikzpicture}
            \fill[black] (1,0) circle (2pt) node[below=4pt]{1};
            \fill[black] (2,0) circle (2pt) node[below=4pt]{2};
            \fill[black] (3,0) circle (2pt) node[below=4pt]{3};
            \fill[black] (4,0) circle (2pt) node[below=4pt]{4};
            \fill[black] (5,0) circle (2pt) node[below=4pt]{5};
            \fill[black] (6,0) circle (2pt) node[below=4pt]{6};
            \fill[black] (7,0) circle (2pt) node[below=4pt]{7};
            \fill[black] (8,0) circle (2pt) node[below=4pt]{8};
            \fill[black] (9,0) circle (2pt) node[below=4pt]{9};
            \fill[black] (10,0) circle (2pt) node[below=4pt]{10};

            \draw[black,dashed] (1,0) -- (1,2) node[above=4pt]{1};
            \draw[black,dashed] (2,0) -- (2,2) node[above=4pt]{2};
            \draw[black,dashed] (3,0) -- (3,2) node[above=4pt]{3};
            \draw[black,dashed] (4,0) -- (4,2) node[above=4pt]{4};
            \draw[black,dashed] (5,0) -- (5,2) node[above=4pt]{5};
            \draw[black,dashed] (6,0) -- (6,2) node[above=4pt]{6};
            \draw[black,dashed] (7,0) -- (7,2) node[above=4pt]{7};
            \draw[black,dashed] (8,0) -- (8,2) node[above=4pt]{8};
            \draw[black,dashed] (9,0) -- (9,2) node[above=4pt]{9};
            \draw[black,dashed] (10,0) -- (10,2) node[above=4pt]{10};

            \fill[black] (1,2) circle (2pt);
            \fill[black] (2,2) circle (2pt);
            \fill[black] (3,2) circle (2pt);
            \fill[black] (4,2) circle (2pt);
            \fill[black] (5,2) circle (2pt);
            \fill[black] (6,2) circle (2pt);
            \fill[black] (7,2) circle (2pt);
            \fill[black] (8,2) circle (2pt);
            \fill[black] (9,2) circle (2pt);
            \fill[black] (10,2) circle (2pt);

            \draw[black,thick] (4,0) -- (1,2);
            \draw[black,thick] (5,0) -- (2,2);
            \draw[black,thick] (6,0) -- (6,2);
            \draw[black,thick] (9,0) -- (7,2);
            \draw[black,thick] (3,0) -- (5,2);
            \draw[black,thick] (10,0) -- (3,2);
        \end{tikzpicture}
        \caption{The non-\glslink{parallel-stitch}{parallel} \glslink{stitch}{stitch} $(\glslink{bracket-n}{[10]}, \{(3,5), (4,1), (5,2), (6,6), (9,7), (10,3)\})$.}
    \end{subfigure}
    
    \vspace{0.5em}
    
    \begin{subfigure}{\linewidth}
        \centering
        \begin{tikzpicture}
            \fill[black] (1,0) circle (2pt) node[below=4pt]{1};
            \fill[black] (2,0) circle (2pt) node[below=4pt]{2};
            \fill[black] (3,0) circle (2pt) node[below=4pt]{3};
            \fill[black] (4,0) circle (2pt) node[below=4pt]{4};
            \fill[black] (5,0) circle (2pt) node[below=4pt]{5};
            \fill[black] (6,0) circle (2pt) node[below=4pt]{6};
            \fill[black] (7,0) circle (2pt) node[below=4pt]{7};
            \fill[black] (8,0) circle (2pt) node[below=4pt]{8};
            \fill[black] (9,0) circle (2pt) node[below=4pt]{9};
            \fill[black] (10,0) circle (2pt) node[below=4pt]{10};

            \draw[black,dashed] (1,0) -- (1,2) node[above=4pt]{1};
            \draw[black,dashed] (2,0) -- (2,2) node[above=4pt]{2};
            \draw[black,dashed] (3,0) -- (3,2) node[above=4pt]{3};
            \draw[black,dashed] (4,0) -- (4,2) node[above=4pt]{4};
            \draw[black,dashed] (5,0) -- (5,2) node[above=4pt]{5};
            \draw[black,dashed] (6,0) -- (6,2) node[above=4pt]{6};
            \draw[black,dashed] (7,0) -- (7,2) node[above=4pt]{7};
            \draw[black,dashed] (8,0) -- (8,2) node[above=4pt]{8};
            \draw[black,dashed] (9,0) -- (9,2) node[above=4pt]{9};
            \draw[black,dashed] (10,0) -- (10,2) node[above=4pt]{10};

            \fill[black] (1,2) circle (2pt);
            \fill[black] (2,2) circle (2pt);
            \fill[black] (3,2) circle (2pt);
            \fill[black] (4,2) circle (2pt);
            \fill[black] (5,2) circle (2pt);
            \fill[black] (6,2) circle (2pt);
            \fill[black] (7,2) circle (2pt);
            \fill[black] (8,2) circle (2pt);
            \fill[black] (9,2) circle (2pt);
            \fill[black] (10,2) circle (2pt);

            \draw[black,thick] (2,0) -- (1,2);
            \draw[black,thick] (3,0) -- (3,2);
            \draw[black,thick] (4,0) -- (6,2);
            \draw[black,thick] (5,0) -- (7,2);
            \draw[black,thick] (7,0) -- (8,2);
            \draw[black,thick] (9,0) -- (10,2);
        \end{tikzpicture}
        \caption{The \glslink{parallel-stitch}{parallel} \glslink{stitch}{stitch} $(\glslink{bracket-n}{[10]}, \{(2,1), (3,3), (4,6), (5,7), (7,8), (9,10)\})$.}
    \end{subfigure}
    \caption{Two \glslink{stitch}{stitches} with $10$ \glslink{column}{columns} and $6$ \glslink{wire}{wires}.} \label{wire-diagram-example-fig}
\end{figure}

\noindent
Since our work involves reasoning about conjugacy classes in $\glslink{symmetric-group}{S_n}$, it also makes sense that the following definition will be important.

\begin{definition}
    Let $\alpha = (\glslink{bracket-n}{[n]}, W)$ be a \glslink{stitch}{stitch} with $n$ \glslink{column}{columns}, and let $\pi \in \glslink{symmetric-group}{S_n}$ be a permutation. We define the \textit{conjugation of $\alpha$ by $\pi$} to be the \glslink{stitch}{stitch} obtained by permuting the \glslink{column}{columns} of $\alpha$ by $\pi$ without disconnecting any \glslink{wire}{wires}. More formally, we set \begin{equation*}
        \pi\alpha\pi^{-1} := (\glslink{bracket-n}{[n]}, \{(\pi(i),\pi(j)) : (i, j) \in W\}).
    \end{equation*} We confirm that this agrees with the group theoretic meaning of $\pi\alpha\pi^{-1}$ when $\alpha$ is a permutation.
\end{definition}

With the above definition, it is natural to ask whether conjugacy classes of \glslink{stitch}{stitches} admit a simple characterization. Indeed, observe that we may interpret any \glslink{stitch}{stitch} $\alpha = (C, W)$ as a directed graph $\glslink{G-alpha}{G_\alpha}$ with vertex set $C$ in which each vertex has at most one incoming and at most one outgoing edge. Moreover, given two \glslink{stitch}{stitches} $\alpha_1 = (\glslink{bracket-n}{[n]}, W_1)$ and $\alpha_2 = (\glslink{bracket-n}{[n]}, W_2)$, the statement that $\alpha_2 = \pi\alpha_1\pi^{-1}$ for some $\pi \in \glslink{symmetric-group}{S_n}$ is precisely the statement that the digraphs $\glslink{G-alpha}{G_{\alpha_1}}$ and $\glslink{G-alpha}{G_{\alpha_2}}$ are isomorphic. If we consider a \glslink{cycle}{\textit{cycle of length $j \ge 1$}} to be a digraph isomorphic to the graph with vertex set $\glslink{bracket-n}{[j]}$ and edge set $\{(i,i+1) : 1 \le i \le j-1\} \cup \{(j, 1)\}$, and a \glslink{chain}{\textit{chain of length $j \ge 1$}} to be a digraph isomorphic to the graph with vertex set $\glslink{bracket-n}{[j]}$ and edge set $\{(i,i+1) : 1 \le i \le j-1\}$ (note that under this convention, a \glslink{chain}{chain of length $j$} has $j$ vertices and $j-1$ edges), then it is a standard result that any digraph in which each vertex has at most one incoming and one outgoing edge decomposes uniquely as a disjoint union of \glslink{cycle}{cycles} and \glslink{chain}{chains}. It follows that the conjugacy class of a \glslink{stitch}{stitch} $\alpha$ is determined by the numbers of \glslink{cycle}{cycles} and \glslink{chain}{chains} in $\glslink{G-alpha}{G_\alpha}$ of each length. We codify this with the following definitions.

\begin{definition}
    \phantom{.}\begin{itemize}
        \item[(i)] A \glslink{cycle}{\textit{cycle}} (resp.\@ \glslink{chain}{\textit{chain}}) of a \glslink{stitch}{stitch} $\alpha$ of length $j$ is a \glslink{cycle}{cycle} (resp.\@ \glslink{chain}{chain}) of length $j$ in the directed graph $\glslink{G-alpha}{G_\alpha}$ associated to $\alpha$.
        \item[(ii)] A \glslink{stitch-type}{\textit{stitch type}} is a pair $\tau = (\rho, \mu)$, where $\rho$ and $\mu$ are integer partitions called the \textit{cycle partition} and the \textit{chain partition} of $\tau$, respectively.
        \item[(iii)] We say that a \glslink{stitch-type}{stitch type} $\tau = (\rho, \mu)$ \textit{has $\glslink{abs-rho}{\lvert\rho\rvert}+\glslink{abs-rho}{\lvert\mu\rvert}$ \glslink{column}{columns} and $\glslink{abs-rho}{\lvert\rho\rvert}+\glslink{abs-rho}{\lvert\mu\rvert}-\glslink{len}{\len(\mu)}$ \glslink{wire}{wires}}, and we denote the set of all \glslink{stitch-type}{stitch types}, the set of \glslink{stitch-type}{stitch types} with $n$ \glslink{column}{columns}, and the set of \glslink{stitch-type}{stitch types} with $n$ \glslink{column}{columns} and $k$ \glslink{wire}{wires} by $\glslink{calT}{\calT}$, $\glslink{calT-n}{\calT_n}$, and $\glslink{calT-n-k}{\calT_{n,k}}$, respectively. We define a \glslink{cycle}{\textit{cycle}} (resp.\@ \glslink{chain}{\textit{chain}}) of $\tau$ to be a part of $\rho$ (resp.\@ $\mu$) of size $j$.
        \item[(iv)] Given a \glslink{stitch}{stitch} $\alpha$, we define the \glslink{stitch-type}{\textit{type of $\alpha$}} to be the \glslink{stitch-type}{stitch type} $\glslink{typ-alpha}{\typ(\alpha)} = (\rho, \mu)$, where for all $j \ge 1$, $\glslink{m-j-rho}{m_j(\rho)}$ is the number of \glslink{cycle}{cycles} of length $j$ in $\glslink{G-alpha}{G_\alpha}$ and $\glslink{m-j-rho}{m_j(\mu)}$ is the number of \glslink{chain}{chains} of length $j$ in $\glslink{G-alpha}{G_\alpha}$. For each \glslink{stitch-type}{stitch type} $\tau$, we define $\glslink{calS-tau}{\calS(\tau)}$ to be the set of all \glslink{stitch}{stitches} of \glslink{stitch-type}{type} $\tau$ and $\glslink{calS-parallel-tau}{\calS^\parallel(\tau)}$ to be the set of all \glslink{parallel-stitch}{parallel} \glslink{stitch}{stitches} of \glslink{stitch-type}{type} $\tau$. It is an easy exercise that the numbers of \glslink{column}{columns}, \glslink{wire}{wires} and \glslink{cycle}{cycles} and \glslink{chain}{chains} of fixed length of a \glslink{stitch}{stitch} are equal to the corresponding numbers for its \glslink{stitch-type}{type}.
        \item[(v)] We say that a \glslink{stitch-type}{stitch type} $\tau = (\rho, \mu)$ is \glslink{parallelizable-stitch-type}{\textit{parallelizable}} if $\glslink{calS-parallel-tau}{\calS^\parallel(\tau)}$ is nonempty, or equivalently, if $\rho = 1^\nu$ for some $\nu \ge 0$. We denote the set of \glslink{parallelizable-stitch-type}{parallelizable} \glslink{stitch-type}{stitch types} by $\glslink{calT-parallel}{\calT^\parallel}$, the set of \glslink{parallelizable-stitch-type}{parallelizable} \glslink{stitch-type}{stitch types} with $n$ \glslink{column}{columns} by $\glslink{calT-n-parallel}{\calT_n^\parallel}$, and the set of \glslink{parallelizable-stitch-type}{parallelizable} \glslink{stitch-type}{stitch types} with $n$ \glslink{column}{columns} and $k$ \glslink{wire}{wires} by $\glslink{calT-n-k-parallel}{\calT_{n,k}^\parallel}$.
    \end{itemize}
\end{definition}

\noindent
See Figure~\ref{stitch-types-example-fig} for some examples.

\begin{figure}[h!]
    \centering
    \begin{subfigure}{\linewidth}
        \centering
        \begin{tikzpicture}
            \fill[orange] (1,0) circle (2pt) node[below=4pt]{1};
            \fill[black] (2,0) circle (2pt) node[below=4pt]{2};
            \fill[blue] (3,0) circle (2pt) node[below=4pt]{3};
            \fill[orange] (4,0) circle (2pt) node[below=4pt]{4};
            \fill[orange] (5,0) circle (2pt) node[below=4pt]{5};
            \fill[purple] (6,0) circle (2pt) node[below=4pt]{6};
            \fill[orange] (7,0) circle (2pt) node[below=4pt]{7};
            \fill[purple] (8,0) circle (2pt) node[below=4pt]{8};
            \fill[purple] (9,0) circle (2pt) node[below=4pt]{9};
            \fill[blue] (10,0) circle (2pt) node[below=4pt]{10};
            \fill[black] (11,0) circle (2pt) node[below=4pt]{11};
            \fill[black] (12,0) circle (2pt) node[below=4pt]{12};

            \draw[orange,dashed] (1,0) -- (1,2) node[above=4pt]{1};
            \draw[black,dashed] (2,0) -- (2,2) node[above=4pt]{2};
            \draw[blue,dashed] (3,0) -- (3,2) node[above=4pt]{3};
            \draw[orange,dashed] (4,0) -- (4,2) node[above=4pt]{4};
            \draw[orange,dashed] (5,0) -- (5,2) node[above=4pt]{5};
            \draw[purple,dashed] (6,0) -- (6,2) node[above=4pt]{6};
            \draw[orange,dashed] (7,0) -- (7,2) node[above=4pt]{7};
            \draw[purple,dashed] (8,0) -- (8,2) node[above=4pt]{8};
            \draw[purple,dashed] (9,0) -- (9,2) node[above=4pt]{9};
            \draw[blue,dashed] (10,0) -- (10,2) node[above=4pt]{10};
            \draw[black,dashed] (11,0) -- (11,2) node[above=4pt]{11};
            \draw[black,dashed] (12,0) -- (12,2) node[above=4pt]{12};

            \fill[orange] (1,2) circle (2pt);
            \fill[black] (2,2) circle (2pt);
            \fill[blue] (3,2) circle (2pt);
            \fill[orange] (4,2) circle (2pt);
            \fill[orange] (5,2) circle (2pt);
            \fill[purple] (6,2) circle (2pt);
            \fill[orange] (7,2) circle (2pt);
            \fill[purple] (8,2) circle (2pt);
            \fill[purple] (9,2) circle (2pt);
            \fill[blue] (10,2) circle (2pt);
            \fill[black] (11,2) circle (2pt);
            \fill[black] (12,2) circle (2pt);

            \draw[purple,thick] (6,0) -- (8,2);
            \draw[blue,thick] (10,0) -- (3,2);
            \draw[orange,thick] (5,0) -- (1,2);
            \draw[purple,thick] (9,0) -- (6,2);
            \draw[black,thick] (12,0) -- (12,2);
            \draw[orange,thick] (4,0) -- (7,2);
            \draw[orange,thick] (1,0) -- (4,2);
            \draw[orange,thick] (7,0) -- (5,2);
        \end{tikzpicture}
        \caption{A \glslink{stitch}{stitch} of \glslink{stitch-type}{type} $((4, 1), (3,2,1,1))$. The \glslink{cycle}{cycle} of length $4$ is highlighted in orange, while the \glslink{chain}{chains} of lengths $2$ and $3$ are highlighted in blue and purple, respectively.}
    \end{subfigure}
    
    \vspace{1em}
    
    \begin{subfigure}{\linewidth}
        \centering
        \begin{tikzpicture}
            \fill[orange] (1,0) circle (2pt) node[below=4pt]{1};
            \fill[orange] (2,0) circle (2pt) node[below=4pt]{2};
            \fill[blue] (3,0) circle (2pt) node[below=4pt]{3};
            \fill[black] (4,0) circle (2pt) node[below=4pt]{4};
            \fill[orange] (5,0) circle (2pt) node[below=4pt]{5};
            \fill[blue] (6,0) circle (2pt) node[below=4pt]{6};
            \fill[orange] (7,0) circle (2pt) node[below=4pt]{7};
            \fill[black] (8,0) circle (2pt) node[below=4pt]{8};
            \fill[black] (9,0) circle (2pt) node[below=4pt]{9};
            \fill[purple] (10,0) circle (2pt) node[below=4pt]{10};
            \fill[black] (11,0) circle (2pt) node[below=4pt]{11};
            \fill[black] (12,0) circle (2pt) node[below=4pt]{12};
            \fill[black] (13,0) circle (2pt) node[below=4pt]{13};

            \draw[orange,dashed] (1,0) -- (1,2) node[above=4pt]{1};
            \draw[orange,dashed] (2,0) -- (2,2) node[above=4pt]{2};
            \draw[blue,dashed] (3,0) -- (3,2) node[above=4pt]{3};
            \draw[black,dashed] (4,0) -- (4,2) node[above=4pt]{4};
            \draw[orange,dashed] (5,0) -- (5,2) node[above=4pt]{5};
            \draw[blue,dashed] (6,0) -- (6,2) node[above=4pt]{6};
            \draw[orange,dashed] (7,0) -- (7,2) node[above=4pt]{7};
            \draw[black,dashed] (8,0) -- (8,2) node[above=4pt]{8};
            \draw[black,dashed] (9,0) -- (9,2) node[above=4pt]{9};
            \draw[purple,dashed] (10,0) -- (10,2) node[above=4pt]{10};
            \draw[black,dashed] (11,0) -- (11,2) node[above=4pt]{11};
            \draw[black,dashed] (12,0) -- (12,2) node[above=4pt]{12};
            \draw[black,dashed] (13,0) -- (13,2) node[above=4pt]{13};

            \fill[orange] (1,2) circle (2pt);
            \fill[orange] (2,2) circle (2pt);
            \fill[blue] (3,2) circle (2pt);
            \fill[black] (4,2) circle (2pt);
            \fill[orange] (5,2) circle (2pt);
            \fill[blue] (6,2) circle (2pt);
            \fill[orange] (7,2) circle (2pt);
            \fill[black] (8,2) circle (2pt);
            \fill[black] (9,2) circle (2pt);
            \fill[purple] (10,2) circle (2pt);
            \fill[black] (11,2) circle (2pt);
            \fill[black] (12,2) circle (2pt);
            \fill[black] (13,2) circle (2pt);

            \draw[orange,thick] (2,0) -- (1,2);
            \draw[orange,thick] (5,0) -- (2,2);
            \draw[blue,thick] (6,0) -- (3,2);
            \draw[orange,thick] (7,0) -- (5,2);
            \draw[black,thick] (8,0) -- (9,2);
            \draw[purple,thick] (10,0) -- (10,2);
            \draw[black,thick] (13,0) -- (11,2);
        \end{tikzpicture}
        \caption{A \glslink{parallel-stitch}{parallel} \glslink{stitch}{stitch} of \glslink{stitch-type}{type} $((1), (4,2,2,2,1,1))$. The unique \glslink{chain}{chain} of length $4$ and one of the \glslink{chain}{chains} of length $2$ are highlighted in orange and blue, respectively, while the unique $1$-\glslink{cycle}{cycle} is highlighted in purple.}
    \end{subfigure}
    \caption{Cycles and \glslink{chain}{chains} within \glslink{stitch}{stitches}.} \label{stitch-types-example-fig}
\end{figure}

\subsection{Notions of \texorpdfstring{\glslink{stitch}{stitch}}{stitch} decomposition}

It will often be useful to speak of ``decomposing'' \glslink{stitch}{stitches} into smaller parts. We will introduce two main ways of doing this, the first of which is through ``\glslink{stitch-coloring}{stitch colorings}'' as defined below.

\begin{definition}
    Let $\alpha = (C, W)$ be a \glslink{stitch}{stitch}, and let $r \ge 0$. An \glslink{stitch-coloring}{\textit{$r$-coloring}} of $\alpha$ is a function $\gamma \colon C \to \glslink{bracket-n}{[r]}$ labelling the \glslink{column}{columns} of $\alpha$ with numbers from $1$ to $r$ such that for each \glslink{wire}{wire} $(i, j) \in W$, $\gamma(i) = \gamma(j)$. Visually, an \glslink{stitch-coloring}{$r$-coloring} of $\alpha$ is a way of coloring the \glslink{column}{columns} of $\alpha$ with $r$ colors so that the endpoints of each \glslink{wire}{wire} have the same color. In such a visualization, we also color each \glslink{wire}{wire} based on the shared color of its endpoints. See Figure~\ref{colorings-example-fig}(A) for an example.
\end{definition}

\begin{definition}
    An \glslink{stitch-coloring}{\textit{$r$-colored stitch}} is a \glslink{stitch}{stitch} $\alpha$ equipped with a chosen \glslink{stitch-coloring}{$r$-coloring} $\gamma$. Given such an \glslink{stitch-coloring}{$r$-colored} \glslink{stitch}{stitch} $\alpha = (C, W)$, for each $i \in \glslink{bracket-n}{[r]}$, we define its \textit{$i$th restriction} to be the \glslink{stitch}{stitch} \begin{equation*}
        \glslink{alpha-rvert-i}{\alpha\rvert_i} := (\gamma^{-1}(i), W \cap (\gamma^{-1}(i) \times \gamma^{-1}(i)))
    \end{equation*} obtained by ignoring all \glslink{column}{columns} and \glslink{wire}{wires} other than those with color $i$---see Figure~\ref{colorings-example-fig}(B-C) for an example. We denote the set of all (resp.\@ \glslink{parallel-stitch}{parallel}) \glslink{stitch-coloring}{$r$-colored} \glslink{stitch}{stitches} by $\glslink{calS-bracket-r}{\calS[r]}$ (resp.\@ $\glslink{calS-parallel-bracket-r}{\calS^\parallel[r]}$).
\end{definition}

\begin{figure}[h!]
    \centering
    \begin{subfigure}{\linewidth}
        \centering
        \begin{tikzpicture}
            \fill[blue] (1,0) circle (2pt) node[below=4pt]{1};
            \fill[orange] (2,0) circle (2pt) node[below=4pt]{2};
            \fill[orange] (3,0) circle (2pt) node[below=4pt]{3};
            \fill[blue] (4,0) circle (2pt) node[below=4pt]{4};
            \fill[orange] (5,0) circle (2pt) node[below=4pt]{5};
            \fill[blue] (6,0) circle (2pt) node[below=4pt]{6};
            \fill[blue] (7,0) circle (2pt) node[below=4pt]{7};
            \fill[orange] (8,0) circle (2pt) node[below=4pt]{8};
            \fill[orange] (9,0) circle (2pt) node[below=4pt]{9};
            \fill[blue] (10,0) circle (2pt) node[below=4pt]{10};
            \fill[orange] (11,0) circle (2pt) node[below=4pt]{11};
            \fill[orange] (12,0) circle (2pt) node[below=4pt]{12};
            \fill[orange] (13,0) circle (2pt) node[below=4pt]{13};

            \draw[blue,dashed] (1,0) -- (1,2) node[above=4pt]{1};
            \draw[orange,dashed] (2,0) -- (2,2) node[above=4pt]{2};
            \draw[orange,dashed] (3,0) -- (3,2) node[above=4pt]{3};
            \draw[blue,dashed] (4,0) -- (4,2) node[above=4pt]{4};
            \draw[orange,dashed] (5,0) -- (5,2) node[above=4pt]{5};
            \draw[blue,dashed] (6,0) -- (6,2) node[above=4pt]{6};
            \draw[blue,dashed] (7,0) -- (7,2) node[above=4pt]{7};
            \draw[orange,dashed] (8,0) -- (8,2) node[above=4pt]{8};
            \draw[orange,dashed] (9,0) -- (9,2) node[above=4pt]{9};
            \draw[blue,dashed] (10,0) -- (10,2) node[above=4pt]{10};
            \draw[orange,dashed] (11,0) -- (11,2) node[above=4pt]{11};
            \draw[orange,dashed] (12,0) -- (12,2) node[above=4pt]{12};
            \draw[orange,dashed] (13,0) -- (13,2) node[above=4pt]{13};

            \fill[blue] (1,2) circle (2pt);
            \fill[orange] (2,2) circle (2pt);
            \fill[orange] (3,2) circle (2pt);
            \fill[blue] (4,2) circle (2pt);
            \fill[orange] (5,2) circle (2pt);
            \fill[blue] (6,2) circle (2pt);
            \fill[blue] (7,2) circle (2pt);
            \fill[orange] (8,2) circle (2pt);
            \fill[orange] (9,2) circle (2pt);
            \fill[blue] (10,2) circle (2pt);
            \fill[orange] (11,2) circle (2pt);
            \fill[orange] (12,2) circle (2pt);
            \fill[orange] (13,2) circle (2pt);

            \draw[orange,thick] (2,0) -- (3,2);
            \draw[orange,thick] (3,0) -- (5,2);
            \draw[blue,thick] (4,0) -- (6,2);
            \draw[orange,thick] (5,0) -- (9,2);
            \draw[blue,thick] (6,0) -- (10,2);
            \draw[orange,thick] (12,0) -- (11,2);
            \draw[orange,thick] (13,0) -- (12,2);
        \end{tikzpicture}
        \caption{A \glslink{stitch-coloring}{$2$-colored} \glslink{parallel-stitch}{parallel} \glslink{stitch}{stitch} $\alpha$, with color $1$ shown in blue and color $2$ shown in orange.}
    \end{subfigure}

    \vspace{1em}

    \begin{subfigure}{0.35\linewidth}
        \centering
        \begin{tikzpicture}
            \fill[blue] (1,0) circle (2pt) node[below=4pt]{1};
            \fill[blue] (2,0) circle (2pt) node[below=4pt]{2};
            \fill[blue] (3,0) circle (2pt) node[below=4pt]{3};
            \fill[blue] (4,0) circle (2pt) node[below=4pt]{4};
            \fill[blue] (5,0) circle (2pt) node[below=4pt]{5};

            \draw[blue,dashed] (1,0) -- (1,2) node[above=4pt]{1};
            \draw[blue,dashed] (2,0) -- (2,2) node[above=4pt]{2};
            \draw[blue,dashed] (3,0) -- (3,2) node[above=4pt]{3};
            \draw[blue,dashed] (4,0) -- (4,2) node[above=4pt]{4};
            \draw[blue,dashed] (5,0) -- (5,2) node[above=4pt]{5};

            \fill[blue] (1,2) circle (2pt);
            \fill[blue] (2,2) circle (2pt);
            \fill[blue] (3,2) circle (2pt);
            \fill[blue] (4,2) circle (2pt);
            \fill[blue] (5,2) circle (2pt);

            \draw[blue,thick] (2,0) -- (3,2);
            \draw[blue,thick] (3,0) -- (5,2);
        \end{tikzpicture}
        \caption{The restriction $\glslink{alpha-rvert-i}{\alpha\rvert_1}$.}
    \end{subfigure} \hspace{1em} \begin{subfigure}{0.55\linewidth}
        \centering
        \begin{tikzpicture}
            \fill[orange] (1,0) circle (2pt) node[below=4pt]{1};
            \fill[orange] (2,0) circle (2pt) node[below=4pt]{2};
            \fill[orange] (3,0) circle (2pt) node[below=4pt]{3};
            \fill[orange] (4,0) circle (2pt) node[below=4pt]{4};
            \fill[orange] (5,0) circle (2pt) node[below=4pt]{5};
            \fill[orange] (6,0) circle (2pt) node[below=4pt]{6};
            \fill[orange] (7,0) circle (2pt) node[below=4pt]{7};
            \fill[orange] (8,0) circle (2pt) node[below=4pt]{8};

            \draw[orange,dashed] (1,0) -- (1,2) node[above=4pt]{1};
            \draw[orange,dashed] (2,0) -- (2,2) node[above=4pt]{2};
            \draw[orange,dashed] (3,0) -- (3,2) node[above=4pt]{3};
            \draw[orange,dashed] (4,0) -- (4,2) node[above=4pt]{4};
            \draw[orange,dashed] (5,0) -- (5,2) node[above=4pt]{5};
            \draw[orange,dashed] (6,0) -- (6,2) node[above=4pt]{6};
            \draw[orange,dashed] (7,0) -- (7,2) node[above=4pt]{7};
            \draw[orange,dashed] (8,0) -- (8,2) node[above=4pt]{8};

            \fill[orange] (1,2) circle (2pt);
            \fill[orange] (2,2) circle (2pt);
            \fill[orange] (3,2) circle (2pt);
            \fill[orange] (4,2) circle (2pt);
            \fill[orange] (5,2) circle (2pt);
            \fill[orange] (6,2) circle (2pt);
            \fill[orange] (7,2) circle (2pt);
            \fill[orange] (8,2) circle (2pt);

            \draw[orange,thick] (1,0) -- (2,2);
            \draw[orange,thick] (2,0) -- (3,2);
            \draw[orange,thick] (3,0) -- (5,2);
            \draw[orange,thick] (7,0) -- (6,2);
            \draw[orange,thick] (8,0) -- (7,2);
        \end{tikzpicture}
        \caption{The restriction $\glslink{alpha-rvert-i}{\alpha\rvert_2}$.}
    \end{subfigure}
    
    \caption{\glslink{stitch-coloring}{Colored} \glslink{parallel-stitch}{parallel} \glslink{stitch}{stitches} and their restrictions.} \label{colorings-example-fig}
\end{figure}

The second main way we will speak of ``decomposing'' stitches comes from the observation that, much of the time, one can ``divide the \glslink{wire-diagram}{wire diagram} of a \glslink{stitch}{stitch} down the middle'' using strategically placed vertical lines. For example, in Figure~\ref{stitch-types-example-fig}(B), one can place vertical lines between \glslink{column}{columns} $7$ and $8$, \glslink{column}{columns} $9$ and $10$, or \glslink{column}{columns} $10$ and $11$. Similarly, in Figure~\ref{colorings-example-fig}(A), one can place vertical lines between \glslink{column}{columns} $1$ and $2$ or \glslink{column}{columns} $10$ and $11$. We formalize this observation using the notion of the \textit{direct sum} of \glslink{stitch}{stitches}.

\begin{definition}
    Let $\alpha_1 = ([n_1], W_1)$ and $\alpha_2 = ([n_2], W_2)$ be two \glslink{stitch}{stitches} in \glslink{standard-form}{standard form}. We define the \textit{direct sum of $\alpha_1$ and $\alpha_2$} to be the \glslink{stitch}{stitch} \begin{equation*}
        \glslink{oplus}{\alpha_1 \oplus \alpha_2} := ([n_1+n_2], W_1 \sqcup \{(n_1+i, n_1+j) : (i, j) \in W_2\})
    \end{equation*} obtained by placing the \glslink{wire-diagram}{wire diagram} of $\alpha_2$ to the right of the \glslink{wire-diagram}{wire diagram} of $\alpha_1$. We extend the definition of $\glslink{oplus}{\oplus}$ to \glslink{stitch-coloring}{colored} \glslink{stitch}{stitches} in the obvious way.
\end{definition}

It is easy to verify that the operation of direct sum is associative, with the \glslink{empty-stitch}{empty stitch} $(\varnothing, \varnothing)$ acting as identity. The notion of direct sum in turn gives rise to a natural notion of irreducibility and a form of unique factorization for \glslink{stitch}{stitches}.

\begin{definition}
    We say that a \glslink{stitch}{stitch} $\alpha$ is \glslink{irreducible-stitch}{\textit{irreducible}} if (i) $\alpha$ is not the \glslink{empty-stitch}{empty stitch}, and (ii) whenever $\alpha = \glslink{oplus}{\alpha_1 \oplus \alpha_2}$ for \glslink{stitch}{stitches} $\alpha_1$ and $\alpha_2$, one of $\alpha_1$ and $\alpha_2$ is the \glslink{empty-stitch}{empty stitch}.
\end{definition}

\begin{proposition} \label{stitch-unique-factorization-prop}
    Every \glslink{stitch}{stitch} $\alpha$ admits a unique decomposition $\alpha = \glslink{oplus}{\alpha_1 \oplus \cdots \oplus \alpha_r}$ as a direct sum of \glslink{irreducible-stitch}{irreducible} \glslink{stitch}{stitches}. (If $r = 0$, we adopt the convention that $\glslink{oplus}{\alpha_1 \oplus \cdots \oplus \alpha_r}$ is the \glslink{empty-stitch}{empty stitch}.)
\end{proposition}

\begin{proof}
    Write $\alpha = (\glslink{bracket-n}{[n]}, W)$ in \glslink{standard-form}{standard form}. For each $0 \le c \le n$, say that $c$ is a \textit{cut} of $\alpha$ if no \glslink{wire}{wire} of $\alpha$ has one endpoint in $\glslink{bracket-n}{[c]}$ and the other in $\glslink{bracket-n}{[n]} \setminus \glslink{bracket-n}{[c]}$. Equivalently, $c$ is a cut of $\alpha$ if and only if $\alpha = \alpha' \oplus \alpha''$ for some \glslink{stitch}{stitches} $\alpha'$ and $\alpha''$ with $c$ and $n-c$ \glslink{column}{columns}, respectively. Let $0 = c_0 < \dots < c_r = n$ be the complete list of cuts of $\alpha$, and for each $1 \le i \le r$, let $\alpha_i := (\{c_{i-1}+1,\dots,c_i\}, W \cap \{c_{i-1}+1,\dots,c_i\}^2)$ be the \glslink{stitch}{stitch} obtained by restricting $\alpha$ to the \glslink{column}{columns} $\{c_{i-1}+1,\dots,c_i\}$. Then we may write $\alpha = \alpha_1 \oplus \dots \oplus \alpha_r$. Moreover, each $\alpha_i$ is \glslink{irreducible-stitch}{irreducible}: otherwise, $\alpha_i$ would have a cut strictly between $0$ and $c_i-c_{i-1}$, which would result in a cut of $\alpha$ lying strictly between $c_{i-1}$ and $c_i$. This proves existence. For uniqueness, observe that given any other \glslink{irreducible-stitch}{irreducible} decomposition $\alpha = \glslink{oplus}{\beta_1 \oplus \dots \oplus \beta_s}$ where each $\beta_i$ has $n_i$ \glslink{column}{columns}, the values $\sum_{j=1}^i n_j$ for $0 \le i \le s$ must all be cuts of $\alpha$, and \glslink{irreducible-stitch}{irreducibility} of the $\beta_i$ guarantees that there are no other cuts. In other words, we have $s = r$, and for all $0 \le i \le s$, we have $c_i = \sum_{j=1}^i n_j$. In particular, this implies that for each $1 \le i \le r$, $\alpha_i$ and $\beta_i$ are obtained by restricting $\alpha$ to the same set $\{c_{i-1}+1,\dots,c_i\}$ of \glslink{column}{columns}, and so the two \glslink{irreducible-stitch}{irreducible} decompositions are equal.
\end{proof}

\noindent
By examining the form of \glslink{irreducible-stitch}{irreducible} \glslink{parallel-stitch}{\textit{parallel}} \glslink{stitch}{stitches} specifically, we are also led to make the following definition and observation.

\begin{definition}
    Let $\alpha = (C, W)$ be a \glslink{parallel-stitch}{parallel} \glslink{stitch}{stitch}. \begin{itemize}
        \item[(i)] We say that $\alpha$ is \glslink{left-leaning-stitch}{\textit{left-leaning}} if $i > j$ for all $(i, j) \in W$.
        \item[(ii)] We say that $\alpha$ is \glslink{left-leaning-stitch}{\textit{right-leaning}} if $i < j$ for all $(i, j) \in W$.
    \end{itemize} We denote the set of all \glslink{left-leaning-stitch}{left-leaning} \glslink{stitch}{stitches} by $\glslink{calS-lparallel}{\calS^\lparallel}$, the set of \glslink{left-leaning-stitch}{left-leaning} \glslink{stitch}{stitches} with $n$ \glslink{column}{columns} by $\glslink{calS-n-lparallel}{\calS_n^\lparallel}$, and the set of \glslink{left-leaning-stitch}{left-leaning} \glslink{stitch}{stitches} with $n$ \glslink{column}{columns} and $k$ \glslink{wire}{wires} by $\glslink{calS-n-k-lparallel}{\calS_{n,k}^\lparallel}$. In addition, we denote the set of all \glslink{cycle}{cycle}-free \glslink{stitch-type}{stitch types} (resp.\@ with $n$ \glslink{column}{columns}, with $n$ \glslink{column}{columns} and $k$ \glslink{wire}{wires}) by $\glslink{calT-lparallel}{\calT^\lparallel}$ (resp.\@ $\glslink{calT-n-lparallel}{\calT_n^\lparallel}$, $\glslink{calT-n-k-lparallel}{\calT_{n,k}^\lparallel}$), and given a \glslink{stitch-type}{stitch type} $\tau \in \glslink{calT-lparallel}{\calT^\lparallel}$, we denote the set of all \glslink{left-leaning-stitch}{left-leaning} \glslink{stitch}{stitches} with \glslink{stitch-type}{type} $\tau$ by $\glslink{calS-lparallel-tau}{\calS^\lparallel(\tau)}$.
\end{definition}

\begin{proposition} \label{irreducible-parallel-classification-prop}
    Every \glslink{irreducible-stitch}{irreducible} \glslink{parallel-stitch}{parallel} \glslink{stitch}{stitch} is exactly one of the following: \begin{itemize}
        \item[(i)] A \glslink{left-leaning-stitch}{left-leaning} \glslink{stitch}{stitch} with at least one \glslink{wire}{wire};
        \item[(ii)] A \glslink{left-leaning-stitch}{right-leaning} \glslink{stitch}{stitch} with at least one \glslink{wire}{wire};
        \item[(iii)] The \glslink{stitch}{stitch} $(\glslink{bracket-n}{[1]}, \varnothing)$, which we refer to as a \glslink{gap}{``gap''}; or
        \item[(iv)] The \glslink{stitch}{stitch} $(\glslink{bracket-n}{[1]}, \{(1,1)\})$, which we refer to as a \glslink{pole}{``pole''}.
    \end{itemize}
\end{proposition}

\begin{proof}
    Let $\alpha = (\glslink{bracket-n}{[n]}, W)$ be an \glslink{irreducible-stitch}{irreducible} \glslink{parallel-stitch}{parallel} \glslink{stitch}{stitch} written in \glslink{standard-form}{standard form}. If $n = 1$, then $\alpha$ is either $(\glslink{bracket-n}{[1]}, \varnothing)$ or $(\glslink{bracket-n}{[1]}, \{(1,1)\})$, giving cases (iii) and (iv). Thus, suppose that $n \ge 2$. For each $1 \le i \le n-1$, refer to the vertical line between \glslink{column}{columns} $i$ and $i+1$ in the \glslink{wire-diagram}{wire diagram} of $\alpha$ as the \textit{$i$th checkpoint}. Since $\alpha$ is \glslink{irreducible-stitch}{irreducible}, every checkpoint must be crossed by at least one \glslink{wire}{wire}; otherwise, $\alpha$ would have a nontrivial cut in the sense defined in the proof of Proposition~\ref{stitch-unique-factorization-prop}. In particular, since there is at least one such checkpoint, $\alpha$ must have at least one \glslink{wire}{wire}. Now, $\alpha$ cannot have any vertical wires: any vertical wires at a \glslink{column}{column} $i$ would imply that the $(i-1)$st checkpoint (if $i \ge 2$) and $i$th checkpoint (if $i \le n-1$) were not crossed by any \glslink{wire}{wires}. Moreover, all the \glslink{wire}{wires} crossing a given checkpoint must lean in the same direction: otherwise, two \glslink{wire}{wires} would be forced to cross each other. Because of this, we may characterize each checkpoint as either ``left-leaning'' or ``right-leaning'' according to which way all \glslink{wire}{wires} crossing the checkpoint lean. Finally, observe that if for some $1 \le i \le n-2$, the $i$th checkpoint and the $(i+1)$st checkpoint leaned in different directions, then this would also force two \glslink{wire}{wires} to cross. Thus, all checkpoints must lean in the same direction. Since every \glslink{wire}{wire} crosses at least one checkpoint (we have already ruled out vertical wires), we deduce that $\alpha$ is either \glslink{left-leaning-stitch}{left-leaning} or \glslink{left-leaning-stitch}{right-leaning}.
\end{proof}

\noindent
See Figure~\ref{decomposition-example-fig} for an example.

\begin{figure}[h!]
    \centering
    \begin{subfigure}{\linewidth}
        \centering
        \begin{tikzpicture}
            \draw[black,dashed] (1, 0) -- (1, 2); \fill[black] (1, 0) circle (2pt) node[below=4pt]{1}; \fill[black] (1, 2) circle (2pt) node[above=4pt]{1};
            \draw[black,dashed] (2, 0) -- (2, 2); \fill[black] (2, 0) circle (2pt) node[below=4pt]{2}; \fill[black] (2, 2) circle (2pt) node[above=4pt]{2};
            \draw[black,dashed] (3, 0) -- (3, 2); \fill[black] (3, 0) circle (2pt) node[below=4pt]{3}; \fill[black] (3, 2) circle (2pt) node[above=4pt]{3};
            \draw[black,dashed] (4, 0) -- (4, 2); \fill[black] (4, 0) circle (2pt) node[below=4pt]{4}; \fill[black] (4, 2) circle (2pt) node[above=4pt]{4};
            \draw[black,dashed] (5, 0) -- (5, 2); \fill[black] (5, 0) circle (2pt) node[below=4pt]{5}; \fill[black] (5, 2) circle (2pt) node[above=4pt]{5};
            \draw[black,dashed] (6, 0) -- (6, 2); \fill[black] (6, 0) circle (2pt) node[below=4pt]{6}; \fill[black] (6, 2) circle (2pt) node[above=4pt]{6};
            \draw[black,dashed] (7, 0) -- (7, 2); \fill[black] (7, 0) circle (2pt) node[below=4pt]{7}; \fill[black] (7, 2) circle (2pt) node[above=4pt]{7};
            \draw[black,dashed] (8, 0) -- (8, 2); \fill[black] (8, 0) circle (2pt) node[below=4pt]{8}; \fill[black] (8, 2) circle (2pt) node[above=4pt]{8};
            \draw[black,dashed] (9, 0) -- (9, 2); \fill[black] (9, 0) circle (2pt) node[below=4pt]{9}; \fill[black] (9, 2) circle (2pt) node[above=4pt]{9};
            \draw[black,dashed] (10, 0) -- (10, 2); \fill[black] (10, 0) circle (2pt) node[below=4pt]{10}; \fill[black] (10, 2) circle (2pt) node[above=4pt]{10};
            \draw[black,dashed] (11, 0) -- (11, 2); \fill[black] (11, 0) circle (2pt) node[below=4pt]{11}; \fill[black] (11, 2) circle (2pt) node[above=4pt]{11};
            \draw[black,dashed] (12, 0) -- (12, 2); \fill[black] (12, 0) circle (2pt) node[below=4pt]{12}; \fill[black] (12, 2) circle (2pt) node[above=4pt]{12};
            \draw[black,thick] plot [] coordinates {(9, 0)(9, 2)};
            \draw[black,thick] plot [] coordinates {(1, 0)} plot [] coordinates {(3, 0)(1, 2)} plot [] coordinates {(5, 0)(3, 2)};
            \draw[black,thick] plot [] coordinates {(2, 0)} plot [] coordinates {(4, 0)(2, 2)};
            \draw[black,thick] plot [] coordinates {(7, 0)} plot [] coordinates {(6, 0)(7, 2)};
            \draw[black,thick] plot [] coordinates {(8, 0)};
            \draw[black,thick] plot [] coordinates {(10, 0)} plot [] coordinates {(12, 0)(10, 2)};
            \draw[black,thick] plot [] coordinates {(11, 0)};
        \end{tikzpicture}
        \caption*{$\alpha$}
    \end{subfigure}
    
    \vspace{1em}
    
    \begin{subfigure}{0.3\linewidth}
        \centering
        \begin{tikzpicture}
            \draw[black,dashed] (1, 0) -- (1, 2); \fill[black] (1, 0) circle (2pt) node[below=4pt]{1}; \fill[black] (1, 2) circle (2pt) node[above=4pt]{1};
            \draw[black,dashed] (2, 0) -- (2, 2); \fill[black] (2, 0) circle (2pt) node[below=4pt]{2}; \fill[black] (2, 2) circle (2pt) node[above=4pt]{2};
            \draw[black,dashed] (3, 0) -- (3, 2); \fill[black] (3, 0) circle (2pt) node[below=4pt]{3}; \fill[black] (3, 2) circle (2pt) node[above=4pt]{3};
            \draw[black,dashed] (4, 0) -- (4, 2); \fill[black] (4, 0) circle (2pt) node[below=4pt]{4}; \fill[black] (4, 2) circle (2pt) node[above=4pt]{4};
            \draw[black,dashed] (5, 0) -- (5, 2); \fill[black] (5, 0) circle (2pt) node[below=4pt]{5}; \fill[black] (5, 2) circle (2pt) node[above=4pt]{5};
            \draw[black,thick] plot [] coordinates {(1, 0)} plot [] coordinates {(3, 0)(1, 2)} plot [] coordinates {(5, 0)(3, 2)};
            \draw[black,thick] plot [] coordinates {(2, 0)} plot [] coordinates {(4, 0)(2, 2)};
        \end{tikzpicture}
        \caption*{$\alpha_1$}
    \end{subfigure} \hspace{1em} \begin{subfigure}{0.17\linewidth}
        \centering
        \begin{tikzpicture}
            \draw[black,dashed] (6, 0) -- (6, 2); \fill[black] (6, 0) circle (2pt) node[below=4pt]{1}; \fill[black] (6, 2) circle (2pt) node[above=4pt]{1};
            \draw[black,dashed] (7, 0) -- (7, 2); \fill[black] (7, 0) circle (2pt) node[below=4pt]{2}; \fill[black] (7, 2) circle (2pt) node[above=4pt]{2};
            \draw[black,thick] plot [] coordinates {(7, 0)} plot [] coordinates {(6, 0)(7, 2)};
        \end{tikzpicture}
        \caption*{$\alpha_2$}
    \end{subfigure} \begin{subfigure}{0.1\linewidth}
        \centering
        \begin{tikzpicture}
            \draw[black,dashed] (8, 0) -- (8, 2); \fill[black] (8, 0) circle (2pt) node[below=4pt]{1}; \fill[black] (8, 2) circle (2pt) node[above=4pt]{1};
        \end{tikzpicture}
        \caption*{$\alpha_3$}
    \end{subfigure} \begin{subfigure}{0.1\linewidth}
        \centering
        \begin{tikzpicture}
            \draw[black,dashed] (9, 0) -- (9, 2); \fill[black] (9, 0) circle (2pt) node[below=4pt]{1}; \fill[black] (9, 2) circle (2pt) node[above=4pt]{1};
            \draw[black,thick] plot [] coordinates {(9, 0)(9, 2)};
        \end{tikzpicture}
        \caption*{$\alpha_4$}
    \end{subfigure} \begin{subfigure}{0.2\linewidth}
        \centering
        \begin{tikzpicture}
            \draw[black,dashed] (10, 0) -- (10, 2); \fill[black] (10, 0) circle (2pt) node[below=4pt]{1}; \fill[black] (10, 2) circle (2pt) node[above=4pt]{1};
            \draw[black,dashed] (11, 0) -- (11, 2); \fill[black] (11, 0) circle (2pt) node[below=4pt]{2}; \fill[black] (11, 2) circle (2pt) node[above=4pt]{2};
            \draw[black,dashed] (12, 0) -- (12, 2); \fill[black] (12, 0) circle (2pt) node[below=4pt]{3}; \fill[black] (12, 2) circle (2pt) node[above=4pt]{3};
            \draw[black,thick] plot [] coordinates {(10, 0)} plot [] coordinates {(12, 0)(10, 2)};
            \draw[black,thick] plot [] coordinates {(11, 0)};
        \end{tikzpicture}
        \caption*{$\alpha_5$}
    \end{subfigure}
    \caption{A \glslink{parallel-stitch}{parallel} \glslink{stitch}{stitch} $\alpha$, which decomposes into $5$ \glslink{irreducible-stitch}{irreducible} components as $\alpha = \glslink{oplus}{\alpha_1 \oplus \alpha_2 \oplus \alpha_3 \oplus \alpha_4 \oplus \alpha_5}$. Here, $\alpha_1$ and $\alpha_5$ are \glslink{left-leaning-stitch}{left-leaning}, $\alpha_2$ is \glslink{left-leaning-stitch}{right-leaning}, $\alpha_3$ is a \glslink{gap}{gap}, and $\alpha_4$ is a \glslink{pole}{pole}.} \label{decomposition-example-fig}
\end{figure}

Still on the subject of ``reducing large \glslink{stitch}{stitches} in terms of smaller \glslink{stitch}{stitches},'' one notices that \glslink{column}{columns} of \glslink{stitch}{stitches} that are not connected to any \glslink{wire}{wires} can be freely removed without affecting the overall structure of a \glslink{stitch}{stitch}. This motivates the following definition.

\begin{definition}
    A \glslink{gap}{\textit{gap}} of a \glslink{stitch}{stitch} $\alpha = (C, W)$ is a \glslink{column}{column} $i \in C$ such that no element of $W$ is of the form $(i, j)$ or $(j, i)$ for $j \in C$. A \glslink{gap}{\textit{gap}} of a \glslink{stitch-type}{stitch type} $\tau = (\rho, \mu)$ is a part of $\mu$ of size $1$. Equivalently, a \glslink{gap}{gap} of $\alpha$ (resp.\@ $\tau$) is a \glslink{chain}{chain} of $\alpha$ (resp.\@ $\tau$) of length $1$. Given a \glslink{stitch}{stitch} $\alpha = (C, W)$, if $G$ denotes the set of \glslink{gap}{gaps} of $\alpha$, we define the \glslink{reduction}{\textit{reduction}} of $\alpha$ to be the \glslink{stitch}{stitch} \begin{equation*}
        \glslink{bar-alpha}{\bar{\alpha}} := (C \setminus G, W).
    \end{equation*} If $\alpha$ has no \glslink{gap}{gaps} (or equivalently, $\alpha = \glslink{bar-alpha}{\bar{\alpha}}$), we say that $\alpha$ is \glslink{reduced}{\textit{reduced}}. Similarly, given a \glslink{stitch-type}{stitch type} $\tau = (\rho, \mu)$, we define the \glslink{reduction}{\textit{reduction}} of $\tau$ to be the \glslink{stitch-type}{stitch type} $\glslink{bar-tau}{\bar{\tau}}$ obtained by removing all parts of size $1$ from $\mu$, and we say that $\tau$ is \glslink{reduced}{reduced} if $\tau = \glslink{bar-tau}{\bar{\tau}}$. If a \glslink{stitch}{stitch} $\alpha$ is both \glslink{irreducible-stitch}{irreducible} and \glslink{reduced}{reduced}, we say that $\alpha$ is \glslink{integral}{\textit{integral}}.
\end{definition}

\begin{example}
    The \glslink{stitch}{stitch} in Figure~\ref{stitch-types-example-fig}(B) has $2$ \glslink{gap}{gaps} at \glslink{column}{columns} $4$ and $12$. The \glslink{stitch}{stitch} in Figure~\ref{colorings-example-fig}(A) has $3$ \glslink{gap}{gaps} at \glslink{column}{columns} $1$, $7$, and $8$.
\end{example}

\noindent
See Figure~\ref{reduction-example-fig} for an example of a \glslink{stitch}{stitch} and its \glslink{reduction}{reduction}.

\begin{figure}[h!]
    \centering
    \begin{subfigure}{\linewidth}
        \centering
        \begin{tikzpicture}
            \draw[black,dashed] (1, 0) -- (1, 2); \fill[black] (1, 0) circle (2pt) node[below=4pt]{1}; \fill[black] (1, 2) circle (2pt) node[above=4pt]{1};
            \draw[black,dashed] (2, 0) -- (2, 2); \fill[black] (2, 0) circle (2pt) node[below=4pt]{2}; \fill[black] (2, 2) circle (2pt) node[above=4pt]{2};
            \draw[orange,dashed] (3, 0) -- (3, 2); \fill[orange] (3, 0) circle (2pt) node[below=4pt]{3}; \fill[orange] (3, 2) circle (2pt) node[above=4pt]{3};
            \draw[black,dashed] (4, 0) -- (4, 2); \fill[black] (4, 0) circle (2pt) node[below=4pt]{4}; \fill[black] (4, 2) circle (2pt) node[above=4pt]{4};
            \draw[black,dashed] (5, 0) -- (5, 2); \fill[black] (5, 0) circle (2pt) node[below=4pt]{5}; \fill[black] (5, 2) circle (2pt) node[above=4pt]{5};
            \draw[black,dashed] (6, 0) -- (6, 2); \fill[black] (6, 0) circle (2pt) node[below=4pt]{6}; \fill[black] (6, 2) circle (2pt) node[above=4pt]{6};
            \draw[orange,dashed] (7, 0) -- (7, 2); \fill[orange] (7, 0) circle (2pt) node[below=4pt]{7}; \fill[orange] (7, 2) circle (2pt) node[above=4pt]{7};
            \draw[orange,dashed] (8, 0) -- (8, 2); \fill[orange] (8, 0) circle (2pt) node[below=4pt]{8}; \fill[orange] (8, 2) circle (2pt) node[above=4pt]{8};
            \draw[black,dashed] (9, 0) -- (9, 2); \fill[black] (9, 0) circle (2pt) node[below=4pt]{9}; \fill[black] (9, 2) circle (2pt) node[above=4pt]{9};
            \draw[orange,dashed] (10, 0) -- (10, 2); \fill[orange] (10, 0) circle (2pt) node[below=4pt]{10}; \fill[orange] (10, 2) circle (2pt) node[above=4pt]{10};
            \draw[black,dashed] (11, 0) -- (11, 2); \fill[black] (11, 0) circle (2pt) node[below=4pt]{11}; \fill[black] (11, 2) circle (2pt) node[above=4pt]{11};
            \draw[black,dashed] (12, 0) -- (12, 2); \fill[black] (12, 0) circle (2pt) node[below=4pt]{12}; \fill[black] (12, 2) circle (2pt) node[above=4pt]{12};
            \draw[black,thick] plot [] coordinates {(1, 0)(1, 2)};
            \draw[black,thick] plot [] coordinates {(5, 0)(5, 2)};
            \draw[black,thick] plot [] coordinates {(4, 0)} plot [] coordinates {(2, 0)(4, 2)};
            \draw[orange,thick] plot [] coordinates {(3, 0)};
            \draw[black,thick] plot [] coordinates {(9, 0)} plot [] coordinates {(6, 0)(9, 2)};
            \draw[orange,thick] plot [] coordinates {(7, 0)};
            \draw[orange,thick] plot [] coordinates {(8, 0)};
            \draw[orange,thick] plot [] coordinates {(10, 0)};
            \draw[black,thick] plot [] coordinates {(11, 0)} plot [] coordinates {(12, 0)(11, 2)};
        \end{tikzpicture}
        \caption{A \glslink{parallel-stitch}{parallel} \glslink{stitch}{stitch} $\alpha := (\glslink{bracket-n}{[12]}, \{(1,1), (2,4), (5,5), (6,9), (12,11)\})$ of \glslink{stitch-type}{type} $((1,1), (2,2,2,1,1,1,1))$, with \glslink{gap}{gaps} highlighted in orange.}
    \end{subfigure}
    
    \vspace{0.5em}
    
    \begin{subfigure}{\linewidth}
        \centering
        \begin{tikzpicture}
            \draw[black,dashed] (1, 0) -- (1, 2); \fill[black] (1, 0) circle (2pt) node[below=4pt]{1}; \fill[black] (1, 2) circle (2pt) node[above=4pt]{1};
            \draw[black,dashed] (2, 0) -- (2, 2); \fill[black] (2, 0) circle (2pt) node[below=4pt]{2}; \fill[black] (2, 2) circle (2pt) node[above=4pt]{2};
            \draw[black,dashed] (3, 0) -- (3, 2); \fill[black] (3, 0) circle (2pt) node[below=4pt]{3}; \fill[black] (3, 2) circle (2pt) node[above=4pt]{3};
            \draw[black,dashed] (4, 0) -- (4, 2); \fill[black] (4, 0) circle (2pt) node[below=4pt]{4}; \fill[black] (4, 2) circle (2pt) node[above=4pt]{4};
            \draw[black,dashed] (5, 0) -- (5, 2); \fill[black] (5, 0) circle (2pt) node[below=4pt]{5}; \fill[black] (5, 2) circle (2pt) node[above=4pt]{5};
            \draw[black,dashed] (6, 0) -- (6, 2); \fill[black] (6, 0) circle (2pt) node[below=4pt]{6}; \fill[black] (6, 2) circle (2pt) node[above=4pt]{6};
            \draw[black,dashed] (7, 0) -- (7, 2); \fill[black] (7, 0) circle (2pt) node[below=4pt]{7}; \fill[black] (7, 2) circle (2pt) node[above=4pt]{7};
            \draw[black,dashed] (8, 0) -- (8, 2); \fill[black] (8, 0) circle (2pt) node[below=4pt]{8}; \fill[black] (8, 2) circle (2pt) node[above=4pt]{8};
            \draw[black,thick] plot [] coordinates {(1, 0)(1, 2)};
            \draw[black,thick] plot [] coordinates {(4, 0)(4, 2)};
            \draw[black,thick] plot [] coordinates {(3, 0)} plot [] coordinates {(2, 0)(3, 2)};
            \draw[black,thick] plot [] coordinates {(6, 0)} plot [] coordinates {(5, 0)(6, 2)};
            \draw[black,thick] plot [] coordinates {(7, 0)} plot [] coordinates {(8, 0)(7, 2)};
        \end{tikzpicture}
        \caption{The \glslink{reduction}{reduction} of $\alpha$, which has \glslink{stitch-type}{type} $((1,1), (2,2,2))$ and is obtained by removing the \glslink{gap}{gaps} from $\alpha$.}
    \end{subfigure}
    \caption{A \glslink{stitch}{stitch} and its \glslink{reduction}{reduction}.} \label{reduction-example-fig}
\end{figure}

\subsection{Additional \texorpdfstring{\glslink{stitch}{stitch}}{stitch} statistics}

In addition to the quantities associated with \glslink{stitch}{stitches} and \glslink{stitch-type}{stitch types} that we have already introduced, there are several more statistics that will be important for us. The first of these generalizes the quantity traditionally denoted $\glslink{z-rho}{z_\rho}$ that measures the size of the centralizer of a permutation of cycle type $\rho$.

\begin{definition}
    For each \glslink{stitch-type}{stitch type} $\tau = (\rho, \mu)$, we define the \textit{stabilizer size} of $\tau$ to be the quantity \begin{equation*}
        \glslink{z-tau}{z_\tau} := \left(\prod_{j\ge1} j^{\glslink{m-j-rho}{m_j(\rho)}}\glslink{m-j-rho}{m_j(\rho)}!\right)\left(\prod_{j\ge1} \glslink{m-j-rho}{m_j(\mu)}!\right).
    \end{equation*} Equivalently, $\glslink{z-tau}{z_\tau}$ is the number of permutations which stabilize a \glslink{stitch}{stitch} of \glslink{stitch-type}{type} $\tau$ under conjugation. The latter characterization demonstrates that $\lvert\glslink{calS-tau}{\calS(\tau)}\rvert = \frac{n!}{\glslink{z-tau}{z_\tau}}$, where $n$ is the number of \glslink{column}{columns} of $\tau$.
\end{definition}

Our next two \glslink{stitch}{stitch} statistics give us ways of quantifying how frequently \glslink{stitch}{stitches} of one \glslink{stitch-type}{type} are ``contained in'' \glslink{stitch}{stitches} of another \glslink{stitch-type}{type}.

\begin{definition}
    Let $\alpha_1 = (C, W_1)$ and $\alpha_2 = (C, W_2)$ be \glslink{stitch}{stitches} with the same number of \glslink{column}{columns}. We write $\glslink{stitch-subseteq}{\alpha_1 \subseteq \alpha_2}$ and say that $\alpha_1$ is a \glslink{substitch}{\textit{substitch}} of $\alpha_2$ if $W_1 \subseteq W_2$.
\end{definition}

\begin{definition}
    Let $\tau$ and $\omega$ be \glslink{stitch-type}{stitch types} with the same number of \glslink{column}{columns}. We define \begin{equation*}
        \glslink{stitch-binom}{\binom{\tau}{\omega}} := \#\{\beta \in \glslink{calS-tau}{\calS(\omega)} : \glslink{stitch-subseteq}{\beta \subseteq \alpha_0}\} \qquad \text{and} \qquad \glslink{stitch-brackbinom}{\brackbinom{\tau}{\omega}} := \#\{\alpha \in \glslink{calS-tau}{\calS(\tau)} : \glslink{stitch-subseteq}{\beta_0 \subseteq \alpha}\},
    \end{equation*} where $\alpha_0 \in \glslink{calS-tau}{\calS(\tau)}$ and $\beta_0 \in \glslink{calS-tau}{\calS(\omega)}$ are arbitrary representatives (it is not hard to see that neither expression depends on the choice of representative). In other words, $\glslink{stitch-binom}{\binom{\tau}{\omega}}$ measures the number of ways to start with a \glslink{stitch}{stitch} of \glslink{stitch-type}{type} $\tau$ and obtain a \glslink{stitch}{stitch} of \glslink{stitch-type}{type} $\omega$ by removing \glslink{wire}{wires}, whereas $\glslink{stitch-brackbinom}{\brackbinom{\tau}{\omega}}$ measures the number of ways to start with a \glslink{stitch}{stitch} of \glslink{stitch-type}{type} $\omega$ and obtain a \glslink{stitch}{stitch} of \glslink{stitch-type}{type} $\tau$ by adding \glslink{wire}{wires}.
\end{definition}

\noindent
As the following lemma demonstrates, the coefficients $\glslink{stitch-binom}{\binom{\tau}{\omega}}$ and $\glslink{stitch-brackbinom}{\brackbinom{\tau}{\omega}}$ are closely related.

\begin{lemma} \label{binom-relation-lem}
    For any two \glslink{stitch-type}{stitch types} $\tau$ and $\omega$ with the same number of \glslink{column}{columns}, we have \begin{equation*}
        \glslink{z-tau}{z_\omega} \glslink{stitch-binom}{\binom{\tau}{\omega}} = \glslink{z-tau}{z_\tau} \glslink{stitch-brackbinom}{\brackbinom{\tau}{\omega}}.
    \end{equation*}
\end{lemma}

\begin{proof}
    Let $n$ be the shared number of \glslink{column}{columns} of $\tau$ and $\omega$, and let $\alpha_0 \in \glslink{calS-tau}{\calS(\tau)}$ and $\beta_0 \in \glslink{calS-tau}{\calS(\omega)}$ be arbitrary representatives. Then we may write \begin{align*}
        \frac{n!}{\glslink{z-tau}{z_\tau}} \glslink{stitch-binom}{\binom{\tau}{\omega}} &= \frac{n!}{\glslink{z-tau}{z_\tau}} \cdot \#\{\beta \in \glslink{calS-tau}{\calS(\omega)} : \glslink{stitch-subseteq}{\beta \subseteq \alpha_0}\} \\
        &= \#\{(\alpha,\beta) \in \glslink{calS-tau}{\calS(\tau)} \times \glslink{calS-tau}{\calS(\omega)} : \glslink{stitch-subseteq}{\beta \subseteq \alpha}\} \\
        &= \frac{n!}{\glslink{z-tau}{z_\omega}} \cdot \#\{\alpha \in \glslink{calS-tau}{\calS(\tau)} : \glslink{stitch-subseteq}{\beta_0 \subseteq \alpha}\} = \frac{n!}{\glslink{z-tau}{z_\omega}} \glslink{stitch-brackbinom}{\brackbinom{\tau}{\omega}}.
    \end{align*} From here, multiplying by $\frac{\glslink{z-tau}{z_\tau} \glslink{z-tau}{z_\omega}}{n!}$ gives the desired result.
\end{proof}

Finally, we introduce one last supplemental collection of statistics that it will prove worthy to assign names to.

\begin{definition}
    \phantom{.} \begin{itemize}
        \item[(i)] A \glslink{pole}{\textit{pole}} of a \glslink{stitch}{stitch} $\alpha = (C, W)$ is a \glslink{column}{column} $i \in C$ such that $(i, i) \in W$, or equivalently, a vertical line in the \glslink{wire-diagram}{wire diagram} of $\alpha$. A \glslink{pole}{\textit{pole}} of a \glslink{stitch-type}{stitch type} $\tau = (\rho, \mu)$ is a part of size $1$ in $\rho$.
        \item[(ii)] A \glslink{proper-cycle}{\textit{proper cycle}} (resp.\@ \glslink{proper-chain}{\textit{proper chain}}) of a \glslink{stitch}{stitch} $\alpha$ is a \glslink{cycle}{cycle} (resp.\@ \glslink{chain}{chain}) of $\alpha$ of length at least $2$. A \glslink{proper-cycle}{\textit{proper cycle}} (resp.\@ \glslink{proper-chain}{\textit{proper chain}}) of a \glslink{stitch-type}{stitch type} $\tau$ is a \glslink{cycle}{cycle} (resp.\@ \glslink{chain}{chain}) of $\tau$ of length at least $2$.
        \item[(iii)] A \glslink{jump}{\textit{jump}} of a \glslink{stitch}{stitch} $\alpha = (C, W)$ is a \glslink{column}{column} $i \in C$ such that $(i, i) \notin W$, and there exist $j, j' \in C$ such that $(i, j) \in W$ and $(j', i) \in W$. We define the \glslink{jump}{\textit{number of jumps}} of a \glslink{stitch-type}{stitch type} $\tau = (\rho, \mu)$ to be the quantity $\glslink{abs-rho}{\lvert\rho\rvert}-\glslink{m-j-rho}{m_1(\rho)}+\glslink{abs-rho}{\lvert\mu\rvert}-2\glslink{len}{\len(\mu)}+\glslink{m-j-rho}{m_1(\mu)}$.
    \end{itemize}
\end{definition}

\noindent
Given a \glslink{stitch}{stitch} $\alpha$ or a \glslink{stitch-type}{stitch type} $\tau$, we remark that the number of \glslink{column}{columns} $n$, the number of \glslink{wire}{wires} $k$, the numbers of \glslink{cycle}{cycles} and \glslink{chain}{chains} $c$ and $j$, the numbers of \glslink{proper-cycle}{proper cycles} and \glslink{proper-chain}{proper chains} $p$ and $\ell$, the numbers of \glslink{pole}{poles} and \glslink{gap}{gaps} $\nu$ and $m$, and the number of \glslink{jump}{jumps} $t$ satisfy a variety of linear relations. For the reader's convenience, we list the main relations that will be used freely throughout the paper below: \begin{equation*}
    n = j + k,\qquad j = m + \ell,\qquad c = \nu + p,\qquad n = 2\ell+m+t+\nu.
\end{equation*}

\begin{example}
    The \glslink{stitch}{stitch} in Figure~\ref{stitch-types-example-fig}(A) has $12$ \glslink{column}{columns}, $8$ \glslink{wire}{wires}, $2$ \glslink{cycle}{cycles}, $1$ \glslink{proper-cycle}{proper cycle}, $4$ \glslink{chain}{chains}, $2$ \glslink{proper-chain}{proper chains}, $1$ \glslink{pole}{pole}, $2$ \glslink{gap}{gaps}, and $5$ \glslink{jump}{jumps}. The \glslink{stitch}{stitch} in Figure~\ref{stitch-types-example-fig}(B) has $13$ \glslink{column}{columns}, $7$ \glslink{wire}{wires}, $1$ \glslink{cycle}{cycle}, $1$ \glslink{proper-cycle}{proper cycle}, $6$ \glslink{chain}{chains}, $4$ \glslink{proper-chain}{proper chains}, $1$ \glslink{pole}{pole}, $2$ \glslink{gap}{gaps}, and $2$ \glslink{jump}{jumps}.
\end{example}

\section{Proof Overview} \label{proof-overview-section}

The proofs of our main theorems are long and involved, requiring a large number of technical lemmas counting various objects or finding explicit formulas for various generating functions. To give the reader a sense of direction, we use this section to give a bird's eye view of the broad structure and key steps of the proof before delving into the details in later sections.

For the remainder of the paper, adopt the definition of $\glslink{E-f-n-k}{E_f(n,k)}$ from Theorem~\ref{main-thm-combo}. Our main theorems---Theorems~\ref{main-thm-1}, \ref{main-thm-2}, and \ref{main-thm-3}---all have to do with obtaining explicit formulas for $\glslink{A-sigma-lambda}{A_{\glslink{id-k}{\id_k}}^{\glslink{lambda-bracket-n}{\lambda[n]}}} = \glslink{E-f-n-k}{E_{\glslink{X-lam}{X^\lambda}}(n, k)}$. Since the expression $\glslink{E-f-n-k}{E_f(n, k)}$ is linear in $f$, it suffices to understand the quantities $\glslink{E-f-n-k}{E_f(n, k)}$ as $f$ ranges through a basis of $\glslink{Delta-Q}{\Delta_\Q}$ of our choosing. Eventually, we will see that the functions $\glslink{E-f-n-k}{E_{\glslink{J-rho}{J_\rho}}(n, k)}$ behave particularly nicely. However, to start, it will be easiest to focus on the quantities $\glslink{E-f-n-k}{E_{\glslink{I-rho}{I_\rho}}(n, k)}$, which have a simpler combinatorial interpretation.

Beyond changing our focus from the values $\glslink{E-f-n-k}{E_{\glslink{X-lam}{X^\lambda}}(n, k)}$ to the values $\glslink{E-f-n-k}{E_{\glslink{I-rho}{I_\rho}}(n, k)}$, it will prove helpful to focus on certain generating functions obtained from the quantities $\glslink{E-f-n-k}{E_f(n, k)}$, rather than the quantities themselves.

\begin{definition}
    For each $f \in \glslink{Delta-Q}{\Delta_\Q}$, we define a power series $\glslink{calE-f}{\calE_f} \in \Q[[x, y]]$ by \begin{equation*}
        \glslink{calE-f}{\calE_f} := \sum_{j, k \ge 0} \frac{(j+k)!}{j!} \cdot \glslink{E-f-n-k}{E_f(j+k, k)} \cdot x^j y^k.
    \end{equation*}
\end{definition}

\noindent
We will also find it useful to adopt the following notation that can simultaneously be used to refer to certain derivatives and certain antiderivatives of power series.

\begin{definition}
    Let $F = \sum_{j,k\ge0} c_{j,k}\cdot x^j y^k \in \Q[[x, y]]$, and let $i \in \Z$. We define \begin{equation*}
        \glslink{partial-x-i}{\partial_x^i} F := \sum_{\substack{j \ge \max(0,-i) \\ k \ge 0}} \frac{(j+i)!}{j!} c_{j+i,k} \cdot x^j y^k.
    \end{equation*}
\end{definition}

Our launch point for the rest of the argument will be the following lemma, which gives an expression for the generating functions $\glslink{calE-f}{\calE_{\glslink{I-rho}{I_\rho}}}$ in terms of certain generating functions $\glslink{H-tau}{H_\tau}$ that count certain families of \glslink{stitch-coloring}{colored} \glslink{stitch}{stitches}.

\begin{lemma} \label{E-from-H-tau-lem}
    For each \glslink{parallelizable-stitch-type}{parallelizable} \glslink{stitch-type}{stitch type} $\tau$ and each $j, k \ge 0$, let \begin{equation*}
        \glslink{h-tau-j-k}{h_\tau(j, k)} := \#\left\{\alpha \in \glslink{calS-parallel-bracket-r}{\calS^\parallel[2]} \,\middle\vert\, \begin{gathered}
            \text{$\glslink{alpha-rvert-i}{\alpha\rvert_1}$ has \glslink{stitch-type}{type} $\tau$} \\
            \text{$\glslink{alpha-rvert-i}{\alpha\rvert_2}$ has $j$ \glslink{chain}{chains} and $k$ \glslink{wire}{wires}}
        \end{gathered}\right\},
    \end{equation*} and for each \glslink{parallelizable-stitch-type}{parallelizable} \glslink{stitch-type}{stitch type} $\tau$, define a generating function $\glslink{H-tau}{H_\tau} \in \Z[[x, y]]$ by \begin{equation*}
        \glslink{H-tau}{H_\tau} := \sum_{j, k \ge 0} \glslink{h-tau-j-k}{h_\tau(j, k)} \cdot x^j y^k.
    \end{equation*} Then for all partitions $\rho$, we have \begin{equation}
        \glslink{z-rho}{z_\rho} \glslink{calE-f}{\calE_{\glslink{I-rho}{I_\rho}}} = \sum_{\substack{i+d=\glslink{abs-rho}{\lvert\rho\rvert}}} \sum_{\tau \in \glslink{calT-n-k-parallel}{\calT_{\glslink{abs-rho}{\lvert\rho\rvert},d}^\parallel}} \glslink{stitch-binom}{\binom{\rho}{\tau}} \cdot y^d \glslink{partial-x-i}{\partial_x^{-i}}(\glslink{z-tau}{z_\tau} \glslink{H-tau}{H_\tau}). \label{E-from-H-tau}
    \end{equation}
\end{lemma}

\begin{proof}
    Let $j, k \ge 0$, and set $n := j+k$. For a fixed permutation $\pi \in \glslink{symmetric-group}{S_n}$, observe that, in the notation of Section~\ref{stitches-section}, \begin{equation*}
        \glslink{N-sigma}{N_{\glslink{id-k}{\id_k}}(\pi)} = \#\{\alpha \in \glslink{calS-n-k-parallel}{\calS_{n,k}^\parallel} : \glslink{stitch-subseteq}{\alpha \subseteq \pi}\}.
    \end{equation*} Consequently, we may write \begin{align*}
        n!\glslink{E-f-n-k}{E_{\glslink{I-rho}{I_\rho}}(n, k)} &= \sum_{\pi \in \glslink{symmetric-group}{S_n}} \glslink{I-rho}{I_\rho}(\pi) \cdot \glslink{N-sigma}{N_{\glslink{id-k}{\id_k}}(\pi)} = \sum_{\substack{\alpha \in \glslink{calS-n-k-parallel}{\calS_{n,k}^\parallel},\ \pi \in \glslink{symmetric-group}{S_n} \\ \glslink{stitch-subseteq}{\alpha \subseteq \pi}}} \glslink{I-rho}{I_\rho}(\pi).
    \end{align*} Recalling the combinatorial interpretation of $\glslink{I-rho}{I_\rho}(\pi)$, we therefore conclude that \begin{equation*}
        n!\glslink{E-f-n-k}{E_{\glslink{I-rho}{I_\rho}}(n, k)} = \#\left\{(\alpha,\pi,R) \in \glslink{calS-n-k-parallel}{\calS_{n,k}^\parallel} \times \glslink{symmetric-group}{S_n} \times \calP(\glslink{bracket-n}{[n]}) \,\middle|\, \begin{gathered}
            \text{$\glslink{stitch-subseteq}{\alpha \subseteq \pi}$, $\pi(R) = R$,} \\
            \text{$\pi\rvert_R$ has cycle type $\rho$}
        \end{gathered}\right\}.
    \end{equation*} Next, note that for any tuple $(\alpha,\pi,R)$ in the above set, the requirements that $\glslink{stitch-subseteq}{\alpha \subseteq \pi}$ and $\pi(R) = R$ ensure that $R$ gives rise to a \glslink{stitch-coloring}{$2$-coloring} of $\alpha$. Consequently, the right hand side becomes \begin{equation*}
        \#\left\{(\alpha,\pi_1,\pi_2) \in \glslink{calS-n-k-parallel-bracket-r}{\calS_{n,k}^\parallel[2]} \times \glslink{symmetric-group}{S_{\glslink{abs-rho}{\lvert\rho\rvert}}} \times S_{n-\glslink{abs-rho}{\lvert\rho\rvert}} \,\middle|\, \begin{gathered}
            \text{$\glslink{alpha-rvert-i}{\alpha\rvert_1}$ has $\glslink{abs-rho}{\lvert\rho\rvert}$ \glslink{column}{columns}, $\glslink{stitch-subseteq}{\glslink{alpha-rvert-i}{\alpha\rvert_1} \subseteq \pi_1}$, $\glslink{stitch-subseteq}{\glslink{alpha-rvert-i}{\alpha\rvert_2} \subseteq \pi_2}$,} \\
            \text{$\pi_1$ has cycle type $\rho$}
        \end{gathered}\right\},
    \end{equation*} where $\glslink{calS-n-k-parallel-bracket-r}{\calS_{n,k}^\parallel[2]}$ denotes the set of \glslink{stitch-coloring}{$2$-colored} \glslink{parallel-stitch}{parallel} \glslink{stitch}{stitches} with $n$ \glslink{column}{columns} and $k$ \glslink{wire}{wires}. Now, for fixed $\alpha \in \glslink{calS-n-k-parallel-bracket-r}{\calS_{n,k}^\parallel[2]}$ such that $\glslink{alpha-rvert-i}{\alpha\rvert_1}$ has $\glslink{abs-rho}{\lvert\rho\rvert}$ \glslink{column}{columns} and $d$ \glslink{wire}{wires} (so $\glslink{alpha-rvert-i}{\alpha\rvert_2}$ has $k-d$ \glslink{wire}{wires}), if we set $i := \glslink{abs-rho}{\lvert\rho\rvert}-d$, then the number of permutations $\pi_2 \in \glslink{symmetric-group}{S_{n-\glslink{abs-rho}{\lvert\rho\rvert}}}$ satisfying $\glslink{stitch-subseteq}{\glslink{alpha-rvert-i}{\alpha\rvert_2} \subseteq \pi_2}$ is given by $(n-\glslink{abs-rho}{\lvert\rho\rvert}-(k-d))! = (j-i)!$ if $j \ge i$ and $0$ otherwise. Meanwhile, the number of permutations $\pi_1 \in \glslink{symmetric-group}{S_{\glslink{abs-rho}{\lvert\rho\rvert}}}$ with cycle type $\rho$ satisfying $\glslink{stitch-subseteq}{\glslink{alpha-rvert-i}{\alpha\rvert_1} \subseteq \pi_1}$ is by definition the coefficient $\glslink{stitch-brackbinom}{\brackbinom{\rho}{\tau}}$ from Section~\ref{stitches-section}, where $\tau$ is the \glslink{stitch-type}{type} of $\glslink{alpha-rvert-i}{\alpha\rvert_1}$. We thus conclude that \begin{align*}
        n!\glslink{E-f-n-k}{E_{\glslink{I-rho}{I_\rho}}(n,k)} &= \sum_{\substack{i+d=\glslink{abs-rho}{\lvert\rho\rvert} \\ i \le j}} \sum_{\tau \in \glslink{calT-n-k-parallel}{\calT_{\glslink{abs-rho}{\lvert\rho\rvert},d}^\parallel}} \sum_{\substack{\alpha \in \glslink{calS-n-k-parallel-bracket-r}{\calS_{n,k}^\parallel[2]} \\ \text{$\glslink{alpha-rvert-i}{\alpha\rvert_1}$ has \glslink{stitch-type}{type} $\tau$}}} (j-i)! \cdot \glslink{stitch-brackbinom}{\brackbinom{\rho}{\tau}} = \sum_{\substack{i+d=\glslink{abs-rho}{\lvert\rho\rvert} \\ i \le j}} (j-i)! \sum_{\tau \in \glslink{calT-n-k-parallel}{\calT_{\glslink{abs-rho}{\lvert\rho\rvert},d}^\parallel}} \glslink{stitch-brackbinom}{\brackbinom{\rho}{\tau}} \cdot \glslink{h-tau-j-k}{h_\tau(j-i, k-d)},
    \end{align*} where the second step follows using the fact that the number of \glslink{column}{columns} in a \glslink{stitch}{stitch} equals the number of \glslink{wire}{wires} plus the number of \glslink{chain}{chains}. From here, dividing both sides by $j!$, multiplying both sides by $\glslink{z-rho}{z_\rho}$, and applying Lemma~\ref{binom-relation-lem} gives \begin{equation*}
        \glslink{z-rho}{z_\rho} \cdot \frac{(j+k)!}{j!} \cdot \glslink{E-f-n-k}{E_{\glslink{I-rho}{I_\rho}}(j+k, k)} = \sum_{\substack{i+d=\glslink{abs-rho}{\lvert\rho\rvert} \\ i \le j}} \frac{(j-i)!}{j!} \sum_{\tau \in \glslink{calT-n-k-parallel}{\calT_{\glslink{abs-rho}{\lvert\rho\rvert},d}^\parallel}} \glslink{stitch-binom}{\binom{\rho}{\tau}} \cdot \glslink{z-tau}{z_\tau} \glslink{h-tau-j-k}{h_\tau(j-i, k-d)}.
    \end{equation*} Taking the sum of each side multiplied by $x^j y^k$ then gives the desired identity of generating functions.
\end{proof}

\noindent
In light of Lemma~\ref{E-from-H-tau-lem}, the problem of computing the generating functions $\glslink{calE-f}{\calE_{\glslink{I-rho}{I_\rho}}}$ reduces to the problem of computing the generating functions $\glslink{H-tau}{H_\tau}$, and hence the problem of enumerating \glslink{stitch-coloring}{colored} \glslink{stitch}{stitches} satisfying certain properties.

While equation~\eqref{E-from-H-tau} might suggest that we will need to compute each generating function $\glslink{H-tau}{H_\tau}$ independently, it will turn out that there are nontrivial relationships between the $\glslink{H-tau}{H_\tau}$ that allow us to further simplify equation~\eqref{E-from-H-tau}. At the heart of these relationships is the following nontrivial theorem, whose proof is the content of Section~\ref{counting-parallel-stitches-section}.

\begin{theorem} \label{parallel-count-thm}
    Let $\tau$ be a \glslink{parallelizable-stitch-type}{parallelizable} \glslink{stitch-type}{stitch type} with $k$ \glslink{wire}{wires}, $\nu$ \glslink{pole}{poles}, and $\ell$ \glslink{proper-chain}{proper chains}. Then \begin{equation*}
        \frac{\lvert\glslink{calS-parallel-tau}{\calS^\parallel(\tau)}\rvert}{\lvert\glslink{calS-tau}{\calS(\tau)}\rvert} = \begin{cases}
            1 & \text{if $\ell = 0$,} \\
            \displaystyle 2\sum_{i=0}^{\ell-1} \glslink{binomial}{\binom{\ell-1}{i}} \frac{(\nu+1+i)!}{(k+1+i)!} & \text{otherwise.}
        \end{cases}
    \end{equation*} In particular, the proportion of \glslink{stitch}{stitches} of \glslink{stitch-type}{type} $\tau$ which are \glslink{parallel-stitch}{parallel} depends only on the numbers of \glslink{wire}{wires}, \glslink{pole}{poles}, and \glslink{proper-chain}{proper chains} of $\tau$.
\end{theorem}

\noindent
As one might expect, this ``uniformity'' in the frequencies of \glslink{parallel-stitch}{parallel} \glslink{stitch}{stitches} with given types gives rise to a certain uniformity of the functions $\glslink{H-tau}{H_\tau}$ as $\tau$ varies:

\begin{corollary} \label{H-tau-from-H-ell-m-t-nu-cor}
    For each $\ell, m, t, \nu \ge 0$, let \begin{equation*}
        \glslink{h-ell-m-t-nu-j-k}{h_{\ell,m,t,\nu}(j, k)} := \#\left\{\alpha \in \glslink{calS-parallel-bracket-r}{\calS^\parallel[2]} \,\middle\vert\, \begin{gathered}
            \text{$\glslink{alpha-rvert-i}{\alpha\rvert_1}$ has $\ell$ \glslink{proper-chain}{proper chains}, $m$ \glslink{gap}{gaps}, $t$ \glslink{jump}{jumps}, and $\nu$ \glslink{pole}{poles}} \\
            \text{$\glslink{alpha-rvert-i}{\alpha\rvert_2}$ has $j$ \glslink{chain}{chains} and $k$ \glslink{wire}{wires}}
        \end{gathered}\right\}
    \end{equation*} for all $j, k \ge 0$, and define \begin{equation*}
        \glslink{H-ell-m-t-nu}{H_{\ell,m,t,\nu}} := \sum_{j,k\ge0} \glslink{h-ell-m-t-nu-j-k}{h_{\ell,m,t,\nu}(j, k)} \cdot x^j y^k \in \Z[[x, y]].
    \end{equation*} Let $\tau$ be a \glslink{parallelizable-stitch-type}{parallelizable} \glslink{stitch-type}{stitch type} with $\ell$ \glslink{proper-chain}{proper chains}, $m$ \glslink{gap}{gaps}, $t$ \glslink{jump}{jumps}, and $\nu$ \glslink{pole}{poles}. Then \begin{equation*}
        \glslink{z-tau}{z_\tau} \glslink{H-tau}{H_{\tau}} = \frac{\ell!m!\nu!}{\glslink{binomial}{\binom{\ell+t-1}{t}}} \glslink{H-ell-m-t-nu}{H_{\ell,m,t,\nu}}.
    \end{equation*}
\end{corollary}

\noindent
By applying Corollary~\ref{H-tau-from-H-ell-m-t-nu-cor} to equation~\eqref{E-from-H-tau}, we thus obtain a formula for the $\glslink{calE-f}{\calE_{\glslink{I-rho}{I_\rho}}}$ in terms of the narrower family of functions $\glslink{H-ell-m-t-nu}{H_{\ell,m,t,\nu}}$.

\begin{corollary} \label{E-from-H-ell-m-t-nu-cor}
    Let $\rho$ be a partition, and for each $\ell, m, t, \nu \ge 0$, define \begin{equation*}
        \glslink{stitch-ell-m-t-nu-binom}{\binom{\rho}{\ell,m,t,\nu}} := \#\left\{\beta \in \glslink{calS-n}{\calS_{\glslink{abs-rho}{\lvert\rho\rvert}}} \,\middle\vert\, \begin{gathered}
            \text{$\beta$ has $0$ \glslink{proper-cycle}{proper cycles}, $\ell$ \glslink{proper-chain}{proper chains}, $m$ \glslink{gap}{gaps},} \\
            \text{$t$ \glslink{jump}{jumps}, and $\nu$ \glslink{pole}{poles} and satisfies $\glslink{stitch-subseteq}{\beta \subseteq \pi_0}$}
        \end{gathered}\right\},
    \end{equation*} where $\pi_0$ is an arbitrary permutation with cycle type $\rho$. Then we have \begin{equation*}
        \glslink{z-rho}{z_\rho}\glslink{calE-f}{\calE_{\glslink{I-rho}{I_\rho}}} = \sum_{\substack{2\ell+m+t+\nu=\glslink{abs-rho}{\lvert\rho\rvert} \\ \text{$\ell \ge 1$ or $t = 0$}}} \glslink{stitch-ell-m-t-nu-binom}{\binom{\rho}{\ell,m,t,\nu}} \cdot \frac{\ell!m!\nu!}{\glslink{binomial}{\binom{\ell+t-1}{t}}} \cdot y^{\ell+t+\nu} \glslink{partial-x-i}{\partial_x^{-\ell-m}} \glslink{H-ell-m-t-nu}{H_{\ell,m,t,\nu}}.
    \end{equation*}
\end{corollary}

\noindent
An efficient method for computing the coefficients $\glslink{stitch-ell-m-t-nu-binom}{\binom{\rho}{\ell,m,t,\nu}}$, which we use in our computer calculations, is given in Appendix~\ref{stitch-binom-appendix}.

Having reduced the problem to essentially the problem of counting \glslink{stitch-coloring}{$2$-colored} \glslink{parallel-stitch}{parallel} \glslink{stitch}{stitches} whose restrictions are constrained to have certain numbers of \glslink{proper-chain}{proper chains}, \glslink{gap}{gaps}, \glslink{jump}{jumps}, \glslink{pole}{poles}, or \glslink{wire}{wires}, we will use Section~\ref{generating-functions-for-colored-stitches-section} to take a broader look at generating functions for \glslink{stitch-coloring}{colored} \glslink{parallel-stitch}{parallel} \glslink{stitch}{stitches} satisfying various constraints. The most useful takeaway from this section will be the following theorem on the algebraicity of a certain family of generating functions counting \glslink{integral}{integral} \glslink{stitch-coloring}{$r$-colored} \glslink{left-leaning-stitch}{left-leaning} \glslink{stitch}{stitches}.

\begin{theorem} \label{F-tilde-lparallel-algebraicity-thm}
    Let $\glslink{f-tilde-lparallel}{\tilde{f}_\lparallel(j_1,t_1;\dots;j_r,t_r)}$ denote the number of \glslink{integral}{integral} \glslink{stitch-coloring}{$r$-colored} \glslink{left-leaning-stitch}{left-leaning} \glslink{stitch}{stitches} $\alpha$ such that for all $1 \le i \le r$, the $i$th restriction $\glslink{alpha-rvert-i}{\alpha\rvert_i}$ has $j_i$ \glslink{chain}{chains} and $t_i$ \glslink{jump}{jumps}. Define a formal power series $\glslink{F-tilde-lparallel}{\tilde{F}_\lparallel} \in \Z[[x_1,y_1;\dots;x_r,y_r]]$ by \begin{equation*}
        \glslink{F-tilde-lparallel}{\tilde{F}_\lparallel} := \sum_{\substack{j_1,\dots,j_r \ge 0 \\ t_1,\dots,t_r \ge 0}} \glslink{f-tilde-lparallel}{\tilde{f}_\lparallel(j_1,t_1;\dots;j_r,t_r)} \cdot x_1^{j_1} y_1^{t_1} \cdots x_r^{j_r} y_r^{t_r}.
    \end{equation*} Then $\glslink{F-tilde-lparallel}{\tilde{F}_\lparallel}$ satisfies the algebraic equation \begin{equation*}
        \glslink{F-tilde-lparallel}{\tilde{F}_\lparallel} = \sum_{i \in \glslink{bracket-n}{[r]}} \frac{x_i}{1-y_i-\glslink{F-tilde-lparallel}{\tilde{F}_\lparallel}}.
    \end{equation*}
\end{theorem}

\noindent
While the bearing of Theorem~\ref{F-tilde-lparallel-algebraicity-thm} on computing the functions $\glslink{H-ell-m-t-nu}{H_{\ell,m,t,\nu}}$ may not be immediately clear, it turns out that a large family of the functions $\glslink{H-ell-m-t-nu}{H_{\ell,m,t,\nu}}$ can be expressed in terms of the generating function $\glslink{F-tilde-lparallel}{\tilde{F}_\lparallel}$ for \glslink{stitch-coloring}{$2$-colored} stitches:

\begin{lemma} \label{H-from-F-lparallel-tilde-lem}
    Define a single generating function $\glslink{H}{H} \in \Z[[L, T, N; x, y]]$ by \begin{equation*}
        \glslink{H}{H} := \sum_{\substack{\ell,t,\nu \ge 0}} \glslink{H-ell-m-t-nu}{H_{\ell,0,t,\nu}} \cdot L^\ell T^t N^\nu.
    \end{equation*} Then we have \begin{equation*}
        \glslink{H}{H} = \frac{1}{1-x-y-N - 2(1-x)\glslink{F-tilde-lparallel}{\tilde{F}_\lparallel}\left(\frac{L}{(1-x)^2}, \frac{T}{1-x}; \frac{xy}{(1-x)^2}, \frac{y}{1-x}\right)}.
    \end{equation*}
\end{lemma}

\noindent
By applying the algebraic relation from Theorem~\ref{F-tilde-lparallel-algebraicity-thm} to the function $\glslink{H}{H}$, we thus find that $\glslink{H}{H}$ itself satisfies a certain algebraic relation of degree $3$:

\begin{corollary} \label{H-algebraicity-cor}
    The generating function $\glslink{H}{H}$ satisfies the algebraic equation \begin{equation*}
        L = -\frac{(1+(1-x+y-2T+N)\glslink{H}{H})(1+2N \glslink{H}{H}+(N^2-(1-x-y)^2+4xy)\glslink{H}{H}^2)}{4\glslink{H}{H}^2(1+(1-x-y+N)\glslink{H}{H})}.
    \end{equation*}
\end{corollary}

Having established that the generating function $\glslink{H}{H}$ built from the functions $\glslink{H-ell-m-t-nu}{H_{\ell,0,t,\nu}}$ satisfies the relation in Corollary~\ref{H-algebraicity-cor}, we proceed in Section~\ref{explicit-formulas-for-H-section} by applying Lagrange inversion to compute explicit formulas for the functions $\glslink{H-ell-m-t-nu}{H_{\ell,0,t,\nu}}$. In addition, we show that these explicit formulas give rise to explicit formulas for the functions $\glslink{H-ell-m-t-nu}{H_{\ell,m,t,\nu}}$ when $m \ge 1$, as well as for their partial derivatives with respect to $x$. We assemble the results into the following theorem, which shows that all of the explicit formulas in question have a shared form. In both the theorem and the rest of the paper, we find it useful to define the following power series in $x$ and $y$: \begin{gather*}
    \glslink{Q}{Q} := \sqrt{(1-x-y)^2-4xy}, \\
    \glslink{U}{U} := \frac{1+x-y+\glslink{Q}{Q}}{2\glslink{Q}{Q}}, \qquad \glslink{V}{V} := \frac{1-x+y+\glslink{Q}{Q}}{2\glslink{Q}{Q}}, \qquad \glslink{W}{W} := \frac{1-x-y+\glslink{Q}{Q}}{2\glslink{Q}{Q}} = \glslink{U}{U}\glslink{V}{V}\glslink{Q}{Q}.
\end{gather*} (To be clear, we take $Q$ to be the unique square root with constant coefficient $1$.) In addition, given a polynomial $p \in \Q[z, w]$, we use $\glslink{total-deg}{\deg(p)}$, $\glslink{deg-x}{\deg_z(p)}$, and $\glslink{deg-x}{\deg_w(p)}$ to denote the total degree of $p$, the degree of $p$ as an element of $(\Q[w])[z]$, and the degree of $p$ as an element of $(\Q[z])[w]$, respectively.

\begin{theorem} \label{H-from-p-thm}
    Let $\ell, m, t, \nu \ge 0$ and $i \ge -1$, and suppose that either (i) $\ell = t = 0$ and $i \ge 0$, (ii) $\ell \ge 1$ and $i \ge 0$, or (iii) $\ell = t = \nu = 0$, $i = -1$ and $m \ge 1$. Then there exists a polynomial $\glslink{p-ell-m-t-nu-i}{p_{\ell,m,t,\nu}^{(i)}} \in \Q[z,w]$ such that \begin{equation*}
        \glslink{partial-x-i}{\partial_x^i} \glslink{H-ell-m-t-nu}{H_{\ell,m,t,\nu}} = \glslink{W}{W}^{\ell+m} \glslink{V}{V}^{i-\ell-t} \glslink{Q}{Q}^{-1-2\ell-m-t-\nu-i} \glslink{p-ell-m-t-nu-i}{p_{\ell,m,t,\nu}^{(i)}}\left(-\frac{1}{\glslink{W}{W}}, \frac{1}{\glslink{V}{V}}\right) - \frac{\glslink{indicator}{\mathbbm{1}_{i=-1}}}{m},
    \end{equation*} where $\glslink{indicator}{\mathbbm{1}_{i=-1}}$ is defined to be $1$ if $i = -1$ and $0$ otherwise. (If $i \ge 0$ and $m = 0$, we simply omit the last term.) Moreover: \begin{itemize}
        \item[(i)] If $\ell = t = 0$ and $i \ge 0$, then $\glslink{deg-x}{\deg_z(\glslink{p-ell-m-t-nu-i}{p_{\ell,m,t,\nu}^{(i)}})} \le m$ and $\glslink{deg-x}{\deg_w(\glslink{p-ell-m-t-nu-i}{p_{\ell,m,t,\nu}^{(i)}})} \le i$.
        \item[(ii)] If $\ell \ge 1$ and $i \ge 0$, then $\glslink{deg-x}{\deg_z(\glslink{p-ell-m-t-nu-i}{p_{\ell,m,t,\nu}^{(i)}})} \le \ell+m$ and $\glslink{total-deg}{\deg(\glslink{p-ell-m-t-nu-i}{p_{\ell,m,t,\nu}^{(i)}})} \le \ell+m+i-1$.
        \item[(iii)] If $\ell = t = \nu = 0$, $i = -1$, and $m \ge 1$, then $\glslink{total-deg}{\deg(\glslink{p-ell-m-t-nu-i}{p_{\ell,m,t,\nu}^{(i)}})} \le m-1$.
        \item[(iv)] If $\ell = m = t = 0$, then we have \begin{equation*}
            \glslink{p-ell-m-t-nu-i}{p_{\ell,m,t,\nu}^{(i)}}(z, 0) = \frac{2^{2i}\left(\frac{\nu-1}{2}+i\right)!}{\left(\frac{\nu-1}{2}\right)!} \qquad \text{and} \qquad (\partial_w \glslink{p-ell-m-t-nu-i}{p_{\ell,m,t,\nu}^{(i)}})(z, 0) = -\frac{2^{2i-1}i\left(\frac{\nu-1}{2}+i\right)!}{\left(\frac{\nu-1}{2}\right)!}.
        \end{equation*}
    \end{itemize}
\end{theorem}

In Section~\ref{explicit-formulas-for-E-section}, we finally turn our focus to obtaining what essentially amounts to explicit formulas for the generating functions $\glslink{calE-f}{\calE_f}$. We start with the following ``generating function equivalent'' of Theorems~\ref{main-thm-2} and \ref{main-thm-3}(ii), which follows quickly by applying Theorem~\ref{H-from-p-thm} to the equation from Corollary~\ref{E-from-H-ell-m-t-nu-cor}.

\begin{theorem} \label{E-from-delta-thm}
    Let $f \in \glslink{Delta-Q}{\Delta_\Q}$ and $j \ge 0$. Then there exists a Laurent polynomial $\glslink{delta-f-j}{\delta_f^{(j)}} \in \Q[v,\frac{1}{v}]$ such that \begin{equation*}
        \left.\glslink{partial-x-i}{\partial_x^j}\glslink{calE-f}{\calE_f}\right|_{x=0} = \glslink{delta-f-j}{\delta_f^{(j)}}\left(\frac{1}{1-y}\right).
    \end{equation*} Moreover, if $f$ is nonzero and $d := \glslink{character-polynomial-deg}{\deg(f)}$, then: \begin{itemize}
        \item[(i)] The smallest power of $v$ appearing in $\glslink{delta-f-j}{\delta_f^{(j)}}$ is at least $2+j-\max(1,d)$.
        \item[(ii)] The largest power of $v$ appearing in $\glslink{delta-f-j}{\delta_f^{(j)}}$ is at most $1+2j+d$.
        \item[(iii)] The coefficients of $v^{1+2j+d}$ and $v^{2j+d}$ in $\glslink{delta-f-j}{\delta_f^{(j)}}$ are given by \begin{equation*}
            c_0 \cdot \frac{2^{2j}(\frac{d-1}{2}+j)!}{(\frac{d-1}{2})!} \qquad \text{and} \qquad c_1 \cdot \frac{2^{2j}(\frac{d-2}{2}+j)!}{(\frac{d-2}{2})!} - c_0 \cdot \frac{2^{2j-1}(j+2d)(\frac{d-1}{2}+j)!}{(\frac{d-1}{2})!},
        \end{equation*} respectively, where $c_0$ and $c_1$ are, respectively, the coefficients of $\glslink{I-rho}{I_{1^d}}$ and $\glslink{I-rho}{I_{1^{d-1}}}$ (assuming $d \ge 1$) in the $\glslink{I-rho}{I_\rho}$-basis expansion of $f$. If $d = 0$, we simply omit the term involving $c_1$.
    \end{itemize}
\end{theorem}

\noindent
From Theorem~\ref{E-from-delta-thm}, it is short work to prove Theorems~\ref{main-thm-2} and \ref{main-thm-3}(ii) as they are stated in the introduction. Intuitively, Theorem~\ref{E-from-delta-thm} tells us that in order to fully understand the functions $\glslink{calE-f}{\calE_f}$, it suffices to obtain an explicit formula for a sufficiently high partial derivative $\glslink{partial-x-i}{\partial_x^i} \glslink{calE-f}{\calE_f}$.

Towards this end of computing closed forms for partial derivatives of $\glslink{calE-f}{\calE_f}$, our starting point is the following lemma, which follows directly from the combination of Corollary~\ref{E-from-H-ell-m-t-nu-cor} and Theorem~\ref{H-from-p-thm}.

\begin{lemma} \label{E-from-epsilon-hat-lem}
    Let $f \in \glslink{Delta-Q}{\Delta_\Q}$ be a \glslink{character-polynomial}{character polynomial} of degree at least $1$, and let $i \ge \glslink{character-polynomial-deg}{\deg(f)}-1$. Then there exists a polynomial $\glslink{hat-varepsilon-f-i}{\hat{\varepsilon}_f^{(i)}} \in \Q[u, v]$ such that \begin{equation*}
        \glslink{partial-x-i}{\partial_x^i}\glslink{calE-f}{\calE_f} = \glslink{hat-varepsilon-f-i}{\hat{\varepsilon}_f^{(i)}}(\glslink{U}{U}, \glslink{V}{V}).
    \end{equation*}
\end{lemma}

\noindent
While Lemma~\ref{E-from-epsilon-hat-lem} already gives an explicit formula for $\glslink{partial-x-i}{\partial_x^{\glslink{character-polynomial-deg}{\deg(f)}-1}}\glslink{calE-f}{\calE_f}$, it turns out that Theorem~\ref{gp-poly-thm} from the introduction of this paper allows us to significantly simplify these explicit formulas in the case where $f$ is \glslink{pure-character-polynomial}{pure}. The reason for this comes down to the following characterization of the series coefficients of the functions $\glslink{U}{U}^a \glslink{V}{V}^b \glslink{Q}{Q}^{-1-r}$, which gives formulas for the series coefficients of the functions $\glslink{U}{U}^a \glslink{V}{V}^b$ as a special case.

\begin{theorem} \label{q-a-b-r-formula-thm}
    For each $r \in \R$ and $a, b \ge 0$, define \begin{equation*}
        \sum_{j,k\ge0} \glslink{q-a-b-r-j-k}{q_{a,b}^{(r)}(j, k)} \cdot x^j y^k := \glslink{U}{U}^a \glslink{V}{V}^b \glslink{Q}{Q}^{-1-r}.
    \end{equation*} Then for all $a, b, j, k \ge 0$, $\glslink{q-a-b-r-j-k}{q_{a,b}^{(r)}(j,k)}$ is given by a polynomial in $r$, and we have \begin{equation*}
        \glslink{q-a-b-r-j-k}{q_{a,b}^{(r)}(j, k)} = \sum_{\substack{0 \le \alpha \le a \\ 0 \le \beta \le b}} \glslink{binomial}{\binom{a}{\alpha}} \glslink{binomial}{\binom{b}{\beta}} \frac{(r+j+k)!(\frac{r-\alpha+\beta}{2}+a-1)!(\frac{r+\alpha-\beta}{2}+b-1)!(\frac{r+\alpha+\beta}{2}+j+k)!}{2j!k!(r+a+b-1)!(\frac{r+\alpha+\beta}{2}-1)!(\frac{r-\alpha+\beta}{2}+j)!(\frac{r+\alpha-\beta}{2}+k)!}
    \end{equation*} whenever the right hand side is defined.
\end{theorem}

\begin{remark}
    Since we are primarily interested in the case where $r = -1$, it may be tempting to try to state and prove the theorem for only this value of $r$. However, by stating the theorem for general $r$, we allow ourselves to determine $\glslink{q-a-b-r-j-k}{q_{a,b}^{(-1)}(j, k)}$ whenever the formula for $\glslink{q-a-b-r-j-k}{q_{a,b}^{(-1)}(j, k)}$ from Theorem~\ref{q-a-b-r-formula-thm} is ill-defined by taking the limit of $\glslink{q-a-b-r-j-k}{q_{a,b}^{(r)}(j, k)}$ as $r \to -1$.
\end{remark}

\noindent
If we apply Theorem~\ref{q-a-b-r-formula-thm} to the formula from Lemma~\ref{E-from-epsilon-hat-lem}, we quickly find that if any term of the form $\glslink{U}{U}^a \glslink{V}{V}^b$ with $a \ge 1$ appeared in the polynomial $\glslink{hat-varepsilon-f-i}{\hat{\varepsilon}_{\glslink{X-lam}{X^\lambda}}^{(\glslink{abs-rho}{\lvert\lambda\rvert}-1)}}(\glslink{U}{U}, \glslink{V}{V})$, it would force the degrees of the polynomials $\glslink{a-sigma-lambda}{a_{\glslink{id-k}{\id_k}}^\lambda}(n)$ to be greater than $k-\glslink{abs-rho}{\lvert\lambda\rvert}$, contradicting Theorem~\ref{gp-poly-thm}. Thus, for \glslink{pure-character-polynomial}{pure} $f$, we find that $\glslink{partial-x-i}{\partial_x^{\glslink{character-polynomial-deg}{\deg(f)}-1}}\glslink{calE-f}{\calE_f}$ admits a simpler characterization:

\begin{lemma} \label{E-from-epsilon-lem}
    Suppose $f \in \glslink{Delta-Q}{\Delta_\Q}$ is \glslink{pure-character-polynomial}{pure of degree $d \ge 1$}. Then there exists a polynomial $\glslink{varepsilon-f}{\varepsilon_f} \in \Q[v]$ such that \begin{equation*}
        \glslink{partial-x-i}{\partial_x^{d-1}} \glslink{calE-f}{\calE_f} = \glslink{varepsilon-f}{\varepsilon_f}(\glslink{V}{V}).
    \end{equation*} In other words, the polynomial $\glslink{hat-varepsilon-f-i}{\hat{\varepsilon}_f^{(d-1)}} \in \Q[u, v]$ from Lemma~\ref{E-from-epsilon-hat-lem} is constant in $u$.
\end{lemma}

\noindent
Armed with Lemma~\ref{E-from-epsilon-lem}, we finally revisit the computation that produced Lemma~\ref{E-from-epsilon-hat-lem} specifically for $f = \glslink{J-rho}{J_\rho}$, taking into account the fact that we may replace the variable $u$ with any convenient value. This results in the following ``generating function equivalent'' of Theorems~\ref{main-thm-1} and \ref{main-thm-3}(i), from which the statements of these theorems given in the introduction quickly follow.

\begin{theorem} \label{eps-properties-thm}
    Suppose $f \in \glslink{Delta-Q}{\Delta_\Q}$ is \glslink{pure-character-polynomial}{pure of degree $d \ge 1$}. Then the polynomial $\glslink{varepsilon-f}{\varepsilon_f}$ from Lemma~\ref{E-from-epsilon-lem} satisfies the following properties: \begin{itemize}
        \item[(i)] $\glslink{total-deg}{\deg(\glslink{varepsilon-f}{\varepsilon_f})} \le 3d-1$;
        \item[(ii)] $\glslink{varepsilon-f}{\varepsilon_f}(0) = 0$;
        \item[(iii)] $(v-1)^{d-1} \mid \glslink{varepsilon-f}{\varepsilon_f}'(v)$.
        \item[(iv)] The coefficients of $v^{3d-1}$ and $v^{3d-2}$ in $\glslink{varepsilon-f}{\varepsilon_f}$ are given by \begin{equation*}
            c_0 \cdot \frac{2^{2d-2}(\frac{3d-3}{2})!}{(\frac{d-1}{2})!} \qquad \text{and} \qquad -c_0 \cdot 2^{2d-2}\left(\frac{(\frac{3d-1}{2})!}{(\frac{d-1}{2})!} + \frac{(\frac{3d-4}{2})!}{(\frac{d-2}{2})!}\right).
        \end{equation*} respectively, where $c_0$ denotes the coefficient of $\glslink{J-rho}{J_{1^d}}$ in the $\glslink{J-rho}{J_\rho}$-basis expansion of $f$.
    \end{itemize}
\end{theorem}

\section{Counting \texorpdfstring{\glslink{parallel-stitch}{parallel}}{parallel} \texorpdfstring{\glslink{stitch}{stitches}}{stitch} of a given \texorpdfstring{\glslink{stitch-type}{type}}{type}} \label{counting-parallel-stitches-section}

In general, counting sets of \glslink{parallel-stitch}{parallel} \glslink{stitch}{stitches} by \glslink{stitch-type}{type} tends to be more difficult than counting sets of \glslink{stitch}{stitches} without the \glslink{parallel-stitch}{parallel} constraint, as the former usually have far less symmetry. The goal of this section will be to develop tools for counting such sets of \glslink{parallel-stitch}{parallel} \glslink{stitch}{stitches}, eventually culminating in a formula for the number $\lvert\glslink{calS-parallel-tau}{\calS^\parallel(\tau)}\rvert$ of \glslink{parallel-stitch}{parallel} \glslink{stitch}{stitches} of any fixed \glslink{stitch-type}{type} $\tau$. Along the way, we will also obtain formulas for the numbers of \glslink{left-leaning-stitch}{left-leaning} \glslink{stitch}{stitches}, as well as \glslink{left-leaning-stitch}{left-leaning} \glslink{stitch}{stitches} with a fixed number of \glslink{irreducible-stitch}{irreducible} components, with a fixed \glslink{stitch-type}{type} $\tau$.

\subsection{Derivatives of \texorpdfstring{\glslink{parallel-stitch}{parallel}}{parallel} \texorpdfstring{\glslink{stitch}{stitches}}{stitch}}

In Section~\ref{stitches-section}, we developed a number of ways of ``reducing'' or ``decomposing'' larger \glslink{stitch}{stitches} in terms of smaller \glslink{stitch}{stitches}. Using this infrastructure, it is already possible to make a minor simplification in the problems of this section.

\begin{lemma} \label{reduction-lem}
    Let $\beta$ be a \glslink{reduced}{reduced} \glslink{parallel-stitch}{parallel} \glslink{stitch}{stitch} with $n$ \glslink{column}{columns}. Then the number of \glslink{parallel-stitch}{parallel} \glslink{stitch}{stitches} $\alpha$ with \glslink{reduction}{reduction} $\glslink{bar-alpha}{\bar{\alpha}} = \beta$ containing $m$ \glslink{gap}{gaps} is $\glslink{binomial}{\binom{n+m}{m}}$.
\end{lemma}

\begin{proof}
    A \glslink{parallel-stitch}{parallel} \glslink{stitch}{stitch} $\alpha$ with \glslink{reduction}{reduction} $\beta$ containing $m$ \glslink{gap}{gaps} necessarily has $n+m$ \glslink{column}{columns}, $m$ of which are \glslink{gap}{gaps}. It is not hard to see that $\alpha$ is determined by the set of \glslink{column}{columns} which are \glslink{gap}{gaps}. Moreover, if we start with $n+m$ \glslink{column}{columns} and choose any $m$ of them to be \glslink{gap}{gaps}, then by simply drawing the \glslink{wire}{wires} of $\beta$ on the remaining $n$ \glslink{column}{columns}, we obtain a \glslink{parallel-stitch}{parallel} \glslink{stitch}{stitch} with $m$ \glslink{gap}{gaps} and \glslink{reduction}{reduction} $\beta$. We conclude that \glslink{parallel-stitch}{parallel} \glslink{stitch}{stitches} with \glslink{reduction}{reduction} $\beta$ containing $m$ \glslink{gap}{gaps} are in bijection with $m$-element subsets of $\glslink{bracket-n}{[n+m]}$.
\end{proof}

\begin{corollary} \label{reduction-cor}
    \phantom{.}\begin{itemize}
        \item[(i)] For any \glslink{parallelizable-stitch-type}{parallelizable} \glslink{stitch-type}{stitch type} $\tau \in \glslink{calT-parallel}{\calT^\parallel}$ with $n$ \glslink{column}{columns} and $m$ \glslink{gap}{gaps}, we have \begin{equation*}
            \lvert\glslink{calS-parallel-tau}{\calS^\parallel(\tau)}\rvert = \glslink{binomial}{\binom{n}{m}} \cdot \glslink{calS-parallel-tau}{\lvert\calS^\parallel(\glslink{bar-tau}{\bar{\tau}})}\rvert.
        \end{equation*}
        \item[(ii)] For any \glslink{cycle}{cycle}-free \glslink{stitch-type}{stitch type} $\tau \in \glslink{calT-lparallel}{\calT^\lparallel}$ with $n$ \glslink{column}{columns} and $m$ \glslink{gap}{gaps}, we have \begin{equation*}
            \lvert\glslink{calS-lparallel-tau}{\calS^\lparallel(\tau)}\rvert = \glslink{binomial}{\binom{n}{m}} \cdot \lvert\glslink{calS-lparallel-tau}{\calS^\lparallel(\glslink{bar-tau}{\bar{\tau}})}\rvert.
        \end{equation*}
    \end{itemize}
\end{corollary}

\begin{proof}
    Immediate from Lemma~\ref{reduction-lem}.
\end{proof}

While narrowing our focus to \glslink{reduced}{reduced} \glslink{stitch}{stitches} helps us to eliminate some distractions, counting \glslink{reduced}{reduced} \glslink{stitch}{stitches} by \glslink{stitch-type}{type} is still not an easy problem. The key tool that will allow us to start counting sets of \glslink{parallel-stitch}{parallel} \glslink{stitch}{stitches} will be the following definition.

\begin{definition}
    Let $\alpha = (C, W)$ be a \glslink{parallel-stitch}{parallel} \glslink{stitch}{stitch} on $n$ \glslink{column}{columns}, and equip the set $W$ with the total order defined by writing $(i, j) \le (i', j')$ if and only if $i \le i'$ and $j \le j'$. We define the \glslink{stitch-derivative-term}{\textit{derivative}} of $\alpha$ to be the \glslink{parallel-stitch}{parallel} \glslink{stitch}{stitch} \begin{equation*}
        \glslink{stitch-derivative}{\partial\alpha} := (W, \{((i,j), (i',j')) \in W^2 : j = i'\}).
    \end{equation*} In other words, $\glslink{stitch-derivative}{\partial\alpha}$ has one \glslink{column}{column} for each \glslink{wire}{wire} of $\alpha$ and keeps track of when the upper endpoint of one \glslink{wire}{wire} lies in the same \glslink{column}{column} as the lower endpoint of another \glslink{wire}{wire}. We define the \glslink{stitch-derivative-term}{\textit{derivative}} of a \glslink{parallelizable-stitch-type}{parallelizable} \glslink{stitch-type}{stitch type} $\tau = (1^\nu,\mu)$ to be the \glslink{stitch-type}{stitch type} $\glslink{stitch-type-derivative}{\partial\tau}$ obtained from $\tau$ by removing all parts of size $1$ from $\mu$, then decreasing the size of each remaining part of $\mu$ by $1$. (As usual, it is not difficult to confirm that these definitions imply that the \glslink{stitch-type}{type} of the \glslink{stitch-derivative-term}{derivative} of a \glslink{stitch}{stitch} $\alpha$ equals the \glslink{stitch-derivative-term}{derivative} of the \glslink{stitch-type}{type} of $\alpha$.) An \glslink{stitch-antiderivative}{\textit{antiderivative}} of a \glslink{stitch}{stitch} $\alpha$ (resp.\@ \glslink{stitch-type}{stitch type} $\tau$) is a \glslink{stitch}{stitch} (resp.\@ \glslink{stitch-type}{stitch type}) whose \glslink{stitch-derivative-term}{derivative} is $\alpha$ (resp.\@ $\tau$). (See Figure~\ref{derivative-fig}.)
\end{definition}

\begin{figure}[h!]
    \centering
    \begin{subfigure}{\linewidth}
        \centering
        \begin{tikzpicture}
            \draw[black,dashed] (1, 0) -- (1, 2); \fill[black] (1, 0) circle (2pt); \fill[black] (1, 2) circle (2pt);
            \draw[black,dashed] (2, 0) -- (2, 2); \fill[black] (2, 0) circle (2pt); \fill[black] (2, 2) circle (2pt);
            \draw[black,dashed] (3, 0) -- (3, 2); \fill[black] (3, 0) circle (2pt); \fill[black] (3, 2) circle (2pt);
            \draw[black,dashed] (4, 0) -- (4, 2); \fill[black] (4, 0) circle (2pt); \fill[black] (4, 2) circle (2pt);
            \draw[black,dashed] (5, 0) -- (5, 2); \fill[black] (5, 0) circle (2pt); \fill[black] (5, 2) circle (2pt);
            \draw[black,dashed] (6, 0) -- (6, 2); \fill[black] (6, 0) circle (2pt); \fill[black] (6, 2) circle (2pt);
            \draw[black,dashed] (7, 0) -- (7, 2); \fill[black] (7, 0) circle (2pt); \fill[black] (7, 2) circle (2pt);
            \draw[black,dashed] (8, 0) -- (8, 2); \fill[black] (8, 0) circle (2pt); \fill[black] (8, 2) circle (2pt);
            \draw[black,dashed] (9, 0) -- (9, 2); \fill[black] (9, 0) circle (2pt); \fill[black] (9, 2) circle (2pt);
            \draw[black,dashed] (10, 0) -- (10, 2); \fill[black] (10, 0) circle (2pt); \fill[black] (10, 2) circle (2pt);
            \draw[black,dashed] (11, 0) -- (11, 2); \fill[black] (11, 0) circle (2pt); \fill[black] (11, 2) circle (2pt);
            \draw[black,dashed] (12, 0) -- (12, 2); \fill[black] (12, 0) circle (2pt); \fill[black] (12, 2) circle (2pt);
            \draw[black,dashed] (13, 0) -- (13, 2); \fill[black] (13, 0) circle (2pt); \fill[black] (13, 2) circle (2pt);
            \draw[black,dashed] (14, 0) -- (14, 2); \fill[black] (14, 0) circle (2pt); \fill[black] (14, 2) circle (2pt);
            \draw[black,dashed] (15, 0) -- (15, 2); \fill[black] (15, 0) circle (2pt); \fill[black] (15, 2) circle (2pt);
            \draw[black,thick] (1,0) -- (1,2) node[midway, left=2pt] {$1$};
            \draw[black,thick] (8,0) -- (8,2) node[midway, left=2pt] {$6$};
            \draw[black,thick] (5, 0) -- (7, 2) node[midway, left=3pt] {$5$} (3, 0) -- (5, 2) node[midway, left=3pt] {$3$} (2, 0) -- (3, 2) node[midway, left=2pt] {$2$};
            \draw[black,thick] (4, 0) -- (6, 2) node[midway, left=3pt] {$4$};
            \draw[black,thick] (11, 0) -- (9, 2) node[midway, left=3pt] {$7$};
            \draw[black,thick] (13, 0) -- (15, 2) node[midway, left=3pt] {$9$} (12, 0) -- (13, 2) node[midway, left=2pt] {$8$};
        \end{tikzpicture}
        \caption{A \glslink{parallel-stitch}{parallel} \glslink{stitch}{stitch} $\alpha := (\glslink{bracket-n}{[15]}, \{(1,1), (2,3), (3,5), (4,6), (5,7), (8,8), (11,9), (12,13), (13,15)\})$ of \glslink{stitch-type}{type} $((1,1), (4,3,2,2,1,1))$, with \glslink{wire}{wires} labelled in order.}
    \end{subfigure}
    
    \vspace{0.5em}
    
    \begin{subfigure}{\linewidth}
        \centering
        \begin{tikzpicture}
            \draw[black,dashed] (1, 0) -- (1, 2); \fill[black] (1, 0) circle (2pt) node[below=4pt]{1}; \fill[black] (1, 2) circle (2pt) node[above=4pt]{1};
            \draw[black,dashed] (2, 0) -- (2, 2); \fill[black] (2, 0) circle (2pt) node[below=4pt]{2}; \fill[black] (2, 2) circle (2pt) node[above=4pt]{2};
            \draw[black,dashed] (3, 0) -- (3, 2); \fill[black] (3, 0) circle (2pt) node[below=4pt]{3}; \fill[black] (3, 2) circle (2pt) node[above=4pt]{3};
            \draw[black,dashed] (4, 0) -- (4, 2); \fill[black] (4, 0) circle (2pt) node[below=4pt]{4}; \fill[black] (4, 2) circle (2pt) node[above=4pt]{4};
            \draw[black,dashed] (5, 0) -- (5, 2); \fill[black] (5, 0) circle (2pt) node[below=4pt]{5}; \fill[black] (5, 2) circle (2pt) node[above=4pt]{5};
            \draw[black,dashed] (6, 0) -- (6, 2); \fill[black] (6, 0) circle (2pt) node[below=4pt]{6}; \fill[black] (6, 2) circle (2pt) node[above=4pt]{6};
            \draw[black,dashed] (7, 0) -- (7, 2); \fill[black] (7, 0) circle (2pt) node[below=4pt]{7}; \fill[black] (7, 2) circle (2pt) node[above=4pt]{7};
            \draw[black,dashed] (8, 0) -- (8, 2); \fill[black] (8, 0) circle (2pt) node[below=4pt]{8}; \fill[black] (8, 2) circle (2pt) node[above=4pt]{8};
            \draw[black,dashed] (9, 0) -- (9, 2); \fill[black] (9, 0) circle (2pt) node[below=4pt]{9}; \fill[black] (9, 2) circle (2pt) node[above=4pt]{9};
            \draw[black,thick] plot [] coordinates {(1, 0)(1, 2)};
            \draw[black,thick] plot [] coordinates {(6, 0)(6, 2)};
            \draw[black,thick] plot [] coordinates {(2, 2)} plot [] coordinates {(3, 2)(2, 0)} plot [] coordinates {(5, 2)(3, 0)};
            \draw[black,thick] plot [] coordinates {(4, 2)};
            \draw[black,thick] plot [] coordinates {(7, 2)};
            \draw[black,thick] plot [] coordinates {(8, 2)} plot [] coordinates {(9, 2)(8, 0)};
        \end{tikzpicture}
        \caption{The \glslink{stitch-derivative-term}{derivative} of $\alpha$, which has \glslink{stitch-type}{type} $((1,1), (3,2,1,1))$.}
    \end{subfigure}
    \caption{A \glslink{stitch}{stitch} and its \glslink{stitch-derivative-term}{derivative}.} \label{derivative-fig}
\end{figure}

In the spirit of ``reducing large \glslink{parallel-stitch}{parallel} \glslink{stitch}{stitches} to smaller ones,'' it is natural to ask whether the \glslink{reduced}{reduced} \glslink{stitch-antiderivative}{antiderivatives} of a given \glslink{stitch}{stitch} can be enumerated in a systematic way. To this end, the following definitions prove useful.

\begin{definition}
    Let $\alpha$ be a \glslink{reduced}{reduced} \glslink{parallel-stitch}{parallel} \glslink{stitch}{stitch}. A \glslink{tie}{\textit{tie}} of $\alpha$ is a \glslink{column}{column} which has \glslink{wire}{wires} connected to both its upper and lower vertices---that is, either a \glslink{jump}{jump} or a \glslink{pole}{pole}. We say that $\alpha$ is \glslink{split}{\textit{split}} if it has no \glslink{tie}{ties}. A \glslink{join-site}{\textit{join site of $\alpha$}} is a pair of adjacent \glslink{column}{columns} such that the left \glslink{column}{column} only has a \glslink{wire}{wire} connected to its lower vertex and the right \glslink{column}{column} only has a \glslink{wire}{wire} connected to its upper vertex. We define the \glslink{splitting}{\textit{splitting}} of $\alpha$ to be the \glslink{stitch}{stitch} $\glslink{stitch-splitting}{\hat{\alpha}}$ obtained by replacing each \glslink{tie}{tie} by a \glslink{join-site}{join site}, diverting the previously connected \glslink{wire}{wires} appropriately (see Figure~\ref{splitting-fig}).
\end{definition}

\begin{figure}[h!]
    \centering
    \begin{subfigure}{\linewidth}
        \centering
        \begin{tikzpicture}
            \draw[gray,densely dotted] (1, 0)--(9, 0) (1, 2)--(9, 2);
            \draw[black,dashed] (1, 0) -- (1, 2); \fill[black] (1, 0) circle (2pt) node[below=4pt]{1}; \fill[black] (1, 2) circle (2pt) node[above=4pt]{1};
            \draw[orange,dashed] (2, 0) -- (2, 2); \fill[orange] (2, 0) circle (2pt) node[below=4pt]{2}; \fill[orange] (2, 2) circle (2pt) node[above=4pt]{2};
            \draw[black,dashed] (3, 0) -- (3, 2); \fill[black] (3, 0) circle (2pt) node[below=4pt]{3}; \fill[black] (3, 2) circle (2pt) node[above=4pt]{3};
            \draw[blue,dashed] (4, 0) -- (4, 2); \fill[blue] (4, 0) circle (2pt) node[below=4pt]{4}; \fill[blue] (4, 2) circle (2pt) node[above=4pt]{4};
            \draw[black,dashed] (5, 0) -- (5, 2); \fill[black] (5, 0) circle (2pt) node[below=4pt]{5}; \fill[black] (5, 2) circle (2pt) node[above=4pt]{5};
            \draw[black,dashed] (6, 0) -- (6, 2); \fill[black] (6, 0) circle (2pt) node[below=4pt]{6}; \fill[black] (6, 2) circle (2pt) node[above=4pt]{6};
            \draw[purple,dashed] (7, 0) -- (7, 2); \fill[purple] (7, 0) circle (2pt) node[below=4pt]{7}; \fill[purple] (7, 2) circle (2pt) node[above=4pt]{7};
            \draw[black,dashed] (8, 0) -- (8, 2); \fill[black] (8, 0) circle (2pt) node[below=4pt]{8}; \fill[black] (8, 2) circle (2pt) node[above=4pt]{8};
            \draw[black,dashed] (9, 0) -- (9, 2); \fill[black] (9, 0) circle (2pt) node[below=4pt]{9}; \fill[black] (9, 2) circle (2pt) node[above=4pt]{9};
            \draw[black,thick] plot [] coordinates {(4, 0)(4, 2)};
            \draw[black,thick] plot [] coordinates {(3, 0)} plot [] coordinates {(2, 0)(3, 2)} plot [] coordinates {(1, 0)(2, 2)};
            \draw[black,thick] plot [] coordinates {(5, 0)} plot [] coordinates {(7, 0)(5, 2)} plot [] coordinates {(9, 0)(7, 2)};
            \draw[black,thick] plot [] coordinates {(6, 0)} plot [] coordinates {(8, 0)(6, 2)};
        \end{tikzpicture}
        \caption{A \glslink{stitch}{stitch} $\alpha$ with three \glslink{tie}{ties}, marked in orange, blue, and purple.}
    \end{subfigure}
    
    \vspace{0.5em}
    
    \begin{subfigure}{\linewidth}
        \centering
        \begin{tikzpicture}
            \draw[gray,densely dotted] (1, 0)--(12, 0) (1, 2)--(12, 2);
            \draw[black,dashed] (1, 0) -- (1, 2); \fill[black] (1, 0) circle (2pt) node[below=4pt]{1}; \fill[black] (1, 2) circle (2pt) node[above=4pt]{1};
            \draw[orange,dashed] (2, 0) -- (2, 2); \fill[orange] (2, 0) circle (2pt) node[below=4pt]{2}; \fill[orange] (2, 2) circle (2pt) node[above=4pt]{2};
            \draw[orange,dashed] (3, 0) -- (3, 2); \fill[orange] (3, 0) circle (2pt) node[below=4pt]{3}; \fill[orange] (3, 2) circle (2pt) node[above=4pt]{3};
            \draw[black,dashed] (4, 0) -- (4, 2); \fill[black] (4, 0) circle (2pt) node[below=4pt]{4}; \fill[black] (4, 2) circle (2pt) node[above=4pt]{4};
            \draw[blue,dashed] (5, 0) -- (5, 2); \fill[blue] (5, 0) circle (2pt) node[below=4pt]{5}; \fill[blue] (5, 2) circle (2pt) node[above=4pt]{5};
            \draw[blue,dashed] (6, 0) -- (6, 2); \fill[blue] (6, 0) circle (2pt) node[below=4pt]{6}; \fill[blue] (6, 2) circle (2pt) node[above=4pt]{6};
            \draw[black,dashed] (7, 0) -- (7, 2); \fill[black] (7, 0) circle (2pt) node[below=4pt]{7}; \fill[black] (7, 2) circle (2pt) node[above=4pt]{7};
            \draw[black,dashed] (8, 0) -- (8, 2); \fill[black] (8, 0) circle (2pt) node[below=4pt]{8}; \fill[black] (8, 2) circle (2pt) node[above=4pt]{8};
            \draw[purple,dashed] (9, 0) -- (9, 2); \fill[purple] (9, 0) circle (2pt) node[below=4pt]{9}; \fill[purple] (9, 2) circle (2pt) node[above=4pt]{9};
            \draw[purple,dashed] (10, 0) -- (10, 2); \fill[purple] (10, 0) circle (2pt) node[below=4pt]{10}; \fill[purple] (10, 2) circle (2pt) node[above=4pt]{10};
            \draw[black,dashed] (11, 0) -- (11, 2); \fill[black] (11, 0) circle (2pt) node[below=4pt]{11}; \fill[black] (11, 2) circle (2pt) node[above=4pt]{11};
            \draw[black,dashed] (12, 0) -- (12, 2); \fill[black] (12, 0) circle (2pt) node[below=4pt]{12}; \fill[black] (12, 2) circle (2pt) node[above=4pt]{12};
            \draw[black,thick] plot [] coordinates {(3, 0)} plot [] coordinates {(1, 0)(3, 2)};
            \draw[black,thick] plot [] coordinates {(4, 0)} plot [] coordinates {(2, 0)(4, 2)};
            \draw[black,thick] plot [] coordinates {(6, 0)} plot [] coordinates {(5, 0)(6, 2)};
            \draw[black,thick] plot [] coordinates {(7, 0)} plot [] coordinates {(9, 0)(7, 2)};
            \draw[black,thick] plot [] coordinates {(8, 0)} plot [] coordinates {(11, 0)(8, 2)};
            \draw[black,thick] plot [] coordinates {(10, 0)} plot [] coordinates {(12, 0)(10, 2)};
        \end{tikzpicture}
        \caption{The \glslink{splitting}{splitting} $\glslink{stitch-splitting}{\hat{\alpha}}$ of $\alpha$, with \glslink{join-site}{join sites} marked according to which \glslink{tie}{tie} they came from.}
    \end{subfigure}
    \caption{A \glslink{reduced}{reduced} \glslink{stitch}{stitch} and its \glslink{splitting}{splitting}.} \label{splitting-fig}
\end{figure}

\begin{definition}
    An \glslink{staircase-diagram}{\textit{$n \times n$ staircase diagram}} is a subset $S$ of an $n \times n$ grid such that if a given grid cell is contained in $S$, then so is every cell down and to the right of it. More formally, $S \subseteq \glslink{bracket-n}{[n]}^2$ is a \glslink{staircase-diagram}{staircase diagram} if for all $(i_0, j_0) \in S$ and $(i, j) \in \glslink{bracket-n}{[n]}^2$ with $i \ge i_0$ and $j_0 \ge j$, we have $(i, j) \in S$. The \glslink{corner-set-term}{\textit{corner set}} of a \glslink{staircase-diagram}{staircase diagram} $S \subseteq \glslink{bracket-n}{[n]}^2$ is the set \begin{equation*}
        \glslink{corner-set}{\Gamma(S)} := \{(i_0, j_0) \in S : \text{if $(i, j) \in S$ with $i_0 \ge i$ and $j \ge j_0$, then $(i, j) = (i_0, j_0)$}\}.
    \end{equation*}
\end{definition}

\begin{definition}
    Let $\alpha = (\glslink{bracket-n}{[n]}, W)$ be a \glslink{split}{split}, \glslink{reduced}{reduced} \glslink{parallel-stitch}{parallel} \glslink{stitch}{stitch}, and let $\{(\ell_1,u_1),\dots,(\ell_k,u_k)\} = W$ be the \glslink{wire}{wires} of $\alpha$ listed in increasing order. We define the \glslink{staircase-diagram}{\textit{staircase diagram of $\alpha$}} by \begin{equation*}
        \glslink{stair}{\stair(\alpha)} := \{(i,j) \in \glslink{bracket-n}{[k]}^2 : u_i > \ell_j\}.
    \end{equation*} (One easily verifies that this is indeed a \glslink{staircase-diagram}{staircase diagram}.) See Figure~\ref{staircase-bij-fig} for an example.
\end{definition}

\noindent
Equipped with the above terminology, we make the following observation:

\begin{lemma} \label{staircase-bij-lem}
    \phantom{.}\begin{itemize}
        \item[(i)] The map \begin{equation*}
            \glslink{stair}{\stair} \colon \{\text{\glslink{split}{split}, \glslink{reduced}{reduced} \glslink{parallel-stitch}{parallel} \glslink{stitch}{stitches} with $k$ \glslink{wire}{wires}}\} \to \{\text{$k \times k$ \glslink{staircase-diagram}{staircase diagrams}}\}
        \end{equation*} is a bijection.
        \item[(ii)] The above map further restricts to a bijection between \glslink{split}{split}, \glslink{reduced}{reduced} \glslink{left-leaning-stitch}{left-leaning} \glslink{stitch}{stitches} with $k$ \glslink{wire}{wires} and $k \times k$ \glslink{staircase-diagram}{staircase diagrams} $S$ ``lying below the diagonal''---that is, satisfying $j < i$ for all boxes $(i,j) \in S$.
        \item[(iii)] For any \glslink{split}{split}, \glslink{reduced}{reduced} \glslink{parallel-stitch}{parallel} \glslink{stitch}{stitch} $\alpha$, we have a bijection \begin{equation*}
            \{\text{\glslink{join-site}{join sites} of $\alpha$}\} \to \{\text{\glslink{corner-set-term}{corners} of $\glslink{stair}{\stair(\alpha)}$}\}
        \end{equation*} which sends a \glslink{join-site}{join site} to the pair $(i, j)$, where $i$ is the index of the \glslink{wire}{wire} connected to the upper right vertex of the \glslink{join-site}{join site} and $j$ is the index of the \glslink{wire}{wire} connected to the lower left vertex of the \glslink{join-site}{join site} (see Figure~\ref{staircase-bij-fig}).
    \end{itemize}
\end{lemma}

\begin{proof}
    We can reinterpret a \glslink{split}{split}, \glslink{reduced}{reduced} \glslink{parallel-stitch}{parallel} \glslink{stitch}{stitch} with $k$ \glslink{wire}{wires} as a word of length $2k$ in the letters $L$ and $U$ containing exactly $k$ of each letter, where $L$ represents a \glslink{column}{column} with a \glslink{wire}{wire} connected to its lower vertex and $U$ represents a \glslink{column}{column} with a \glslink{wire}{wire} connected to its upper vertex. It is then a standard fact that such words are in bijection with $k \times k$ \glslink{staircase-diagram}{staircase diagrams}: given such a word $w_1 \cdots w_{2k}$, we obtain the boundary of the corresponding \glslink{staircase-diagram}{staircase diagram} by starting at the bottom left corner, going up by $1$ whenever we encounter an $L$, and going right by one whenever we encounter a $U$. Under the word interpretation of \glslink{split}{split}, \glslink{reduced}{reduced} \glslink{parallel-stitch}{parallel} \glslink{stitch}{stitches}, such a \glslink{stitch}{stitch} is \glslink{left-leaning-stitch}{left-leaning} if and only if the number of $U$s to the left of any point in the word is always greater than or equal to the number of $L$s to the left of that point. This is precisely the statement that the corresponding \glslink{staircase-diagram}{staircase diagram} lies below the diagonal line $y = x$. Lastly, under this word interpretation of \glslink{split}{split}, \glslink{reduced}{reduced} \glslink{parallel-stitch}{parallel} \glslink{stitch}{stitches}, \glslink{join-site}{join sites} correspond to occurrences of $LU$, which in turn correspond to the claimed \glslink{corner-set-term}{corners} in the \glslink{staircase-diagram}{staircase diagram}.
\end{proof}

\begin{figure}[h!]
    \centering
    \begin{subfigure}{\linewidth}
        \centering
        \begin{tikzpicture}
            \draw[gray,densely dotted] (1, 0)--(14, 0) (1, 2)--(14, 2);
            \draw[black,dashed] (1, 0) -- (1, 2); \fill[black] (1, 0) circle (2pt); \fill[black] (1, 2) circle (2pt);
            \draw[orange,dashed] (2, 0) -- (2, 2); \fill[orange] (2, 0) circle (2pt); \fill[orange] (2, 2) circle (2pt);
            \draw[orange,dashed] (3, 0) -- (3, 2); \fill[orange] (3, 0) circle (2pt); \fill[orange] (3, 2) circle (2pt);
            \draw[black,dashed] (4, 0) -- (4, 2); \fill[black] (4, 0) circle (2pt); \fill[black] (4, 2) circle (2pt);
            \draw[black,dashed] (5, 0) -- (5, 2); \fill[black] (5, 0) circle (2pt); \fill[black] (5, 2) circle (2pt);
            \draw[blue,dashed] (6, 0) -- (6, 2); \fill[blue] (6, 0) circle (2pt); \fill[blue] (6, 2) circle (2pt);
            \draw[blue,dashed] (7, 0) -- (7, 2); \fill[blue] (7, 0) circle (2pt); \fill[blue] (7, 2) circle (2pt);
            \draw[black,dashed] (8, 0) -- (8, 2); \fill[black] (8, 0) circle (2pt); \fill[black] (8, 2) circle (2pt);
            \draw[purple,dashed] (9, 0) -- (9, 2); \fill[purple] (9, 0) circle (2pt); \fill[purple] (9, 2) circle (2pt);
            \draw[purple,dashed] (10, 0) -- (10, 2); \fill[purple] (10, 0) circle (2pt); \fill[purple] (10, 2) circle (2pt);
            \draw[black,dashed] (11, 0) -- (11, 2); \fill[black] (11, 0) circle (2pt); \fill[black] (11, 2) circle (2pt);
            \draw[green,dashed] (12, 0) -- (12, 2); \fill[green] (12, 0) circle (2pt); \fill[green] (12, 2) circle (2pt);
            \draw[green,dashed] (13, 0) -- (13, 2); \fill[green] (13, 0) circle (2pt); \fill[green] (13, 2) circle (2pt);
            \draw[black,dashed] (14, 0) -- (14, 2); \fill[black] (14, 0) circle (2pt); \fill[black] (14, 2) circle (2pt);
            \draw[black,thick] (1, 0) -- (3, 2) node[midway,left=4pt] {1};
            \draw[black,thick] (2, 0) -- (4, 2) node[midway,left=4pt] {2};
            \draw[black,thick] (6, 0) -- (5, 2) node[midway,left=3pt] {3};
            \draw[black,thick] (8, 0) -- (7, 2) node[midway,left=3pt] {4};
            \draw[black,thick] (9, 0) -- (10, 2) node[midway,left=3pt] {5};
            \draw[black,thick] (11, 0) -- (13, 2) node[midway,left=4pt] {6};
            \draw[black,thick] (12, 0) -- (14, 2) node[midway,left=4pt] {7};
        \end{tikzpicture}
        \caption{A \glslink{split}{split}, \glslink{reduced}{reduced} \glslink{stitch}{stitch} $\alpha$ with \glslink{join-site}{join sites} colored.}
    \end{subfigure}
    
    \vspace{0.5em}
    
    \begin{subfigure}{\linewidth}
        \centering
        \begin{tikzpicture}[scale=0.7]
            \draw[black] (0,0) -- (7,0);
            \draw[black] (0,1) -- (7,1);
            \draw[black] (0,2) -- (7,2);
            \draw[black] (0,3) -- (7,3);
            \draw[black] (0,4) -- (7,4);
            \draw[black] (0,5) -- (7,5);
            \draw[black] (0,6) -- (7,6);
            \draw[black] (0,7) -- (7,7);

            \draw[black] (0,0) -- (0,7);
            \draw[black] (1,0) -- (1,7);
            \draw[black] (2,0) -- (2,7);
            \draw[black] (3,0) -- (3,7);
            \draw[black] (4,0) -- (4,7);
            \draw[black] (5,0) -- (5,7);
            \draw[black] (6,0) -- (6,7);
            \draw[black] (7,0) -- (7,7);

            \fill[black,opacity=0.3] (0,0) -- (0,1) -- (0,2 ) -- (1,2) -- (2,2) -- (3,2) -- (3,3) -- (4,3) -- (4,4) -- (4,5) -- (5,5) -- (5,6) -- (5,7) -- (6,7) -- (7,7) -- (7,0);

            \fill[orange, opacity=0.4] (0,1) -- (0,2) -- (1,2) -- (1,1);
            \draw[black] (0.5,1.5) node {$\scriptstyle (1,2)$};
            \fill[blue, opacity=0.4] (3,2) -- (3,3) -- (4,3) -- (4,2);
            \draw[black] (3.5,2.5) node {$\scriptstyle (4,3)$};
            \fill[purple, opacity=0.4] (4,4) -- (4,5) -- (5,5) -- (5,4);
            \draw[black] (4.5,4.5) node {$\scriptstyle (5,5)$};
            \fill[green, opacity=0.4] (5,6) -- (5,7) -- (6,7) -- (6,6);
            \draw[black] (5.5,6.5) node {$\scriptstyle (6,7)$};

            \draw[black,very thick] (0,0) -- (0,1) -- (0,2 ) -- (1,2) -- (2,2) -- (3,2) -- (3,3) -- (4,3) -- (4,4) -- (4,5) -- (5,5) -- (5,6) -- (5,7) -- (6,7) -- (7,7);
        \end{tikzpicture}
        \caption{The stair diagram of $\alpha$, with \glslink{corner-set-term}{corners} colored according to which \glslink{join-site}{join site} they correspond to.}
    \end{subfigure}
    \caption{An example application of Lemma~\ref{staircase-bij-lem}.} \label{staircase-bij-fig}
\end{figure}

Returning to the task of characterizing \glslink{stitch-antiderivative}{antiderivatives} of \glslink{parallel-stitch}{parallel} \glslink{stitch}{stitches}, adopt the following terminology.

\begin{definition}
    Let $\alpha$ be a \glslink{parallel-stitch}{parallel} \glslink{stitch}{stitch} with $n$ \glslink{column}{columns}. An \glslink{stitch-extension}{\textit{extension of $\alpha$}} is any \glslink{parallel-stitch}{parallel} \glslink{stitch}{stitch} $\beta$ with $n$ \glslink{column}{columns} satisfying $\glslink{stitch-subseteq}{\alpha \subseteq \beta}$.
\end{definition}

\noindent
Then we have the following theorem.

\begin{theorem} \label{antiderivative-thm}
    Let $\beta$ be a \glslink{parallel-stitch}{parallel} \glslink{stitch}{stitch} with $k$ \glslink{column}{columns}, and for each \glslink{reduced}{reduced} \glslink{stitch-antiderivative}{antiderivative} $\alpha$ of $\beta$, define \begin{equation*}
        \phi(\alpha) := (\glslink{bracket-n}{[k]}, \glslink{corner-set}{\Gamma}(\glslink{stair}{\stair}(\glslink{stitch-splitting}{\hat{\alpha}}))).
    \end{equation*} Then the map \begin{equation*}
        \phi \colon \{\text{\glslink{reduced}{reduced} \glslink{stitch-antiderivative}{antiderivatives} of $\beta$}\} \to \{\text{\glslink{stitch-extension}{extensions} of $\beta$}\}
    \end{equation*} is a bijection. Moreover, if $\beta$ is \glslink{left-leaning-stitch}{left-leaning}, then $\phi$ restricts to a bijection \begin{equation*}
        \phi \colon \{\text{\glslink{left-leaning-stitch}{left-leaning} \glslink{reduced}{reduced} \glslink{stitch-antiderivative}{antiderivatives} of $\beta$}\} \simto \{\text{\glslink{left-leaning-stitch}{left-leaning} \glslink{stitch-extension}{extensions} of $\beta$}\}.
    \end{equation*}
\end{theorem}

\begin{proof}
    The fact that $\glslink{stitch-subseteq}{\beta \subseteq \phi(\alpha)}$ for any \glslink{reduced}{reduced} \glslink{stitch-antiderivative}{antiderivative} $\alpha$ of $\beta$ follows from Lemma~\ref{staircase-bij-lem}(iii) together with the fact that $\glslink{stitch-splitting}{\hat{\alpha}}$ has a \glslink{join-site}{join site} corresponding to each \glslink{wire}{wire} of $\beta$. To see that $\phi$ is a bijection, it suffices to check that each step in the sequence of maps \begin{align*}
        \{\text{\glslink{reduced}{reduced} \glslink{stitch-antiderivative}{antiderivatives} of $\beta$}\} &\xrightarrow{\alpha \mapsto \glslink{stitch-splitting}{\hat{\alpha}}} \left\{\begin{gathered}
            \text{\glslink{split}{split}, \glslink{reduced}{reduced} \glslink{stitch}{stitches} $\alpha_0$ with $k$ \glslink{wire}{wires} such that} \\
            \text{for each $(i,j) \in \beta$, $\alpha_0$ has a \glslink{join-site}{join site} corresponding} \\
            \text{to $(i,j)$ under Lemma~\ref{staircase-bij-lem}(iii)}
        \end{gathered}\right\} \\
        &\xrightarrow{\glslink{stair}{\stair}} \{\text{$k \times k$ \glslink{staircase-diagram}{staircase diagrams} with \glslink{corner-set-term}{corners} at every \glslink{wire}{wire} of $\beta$}\} \\
        &\xrightarrow{S \mapsto (\glslink{bracket-n}{[k]}, \glslink{corner-set}{\Gamma(S)})} \{\text{\glslink{stitch-extension}{extensions} of $\beta$}\}
    \end{align*} is bijective. It is clear that the first step is bijective, and the bijectivity of the second step follows from Lemma~\ref{staircase-bij-lem}. The bijectivity of the last step follows from the standard fact that a \glslink{staircase-diagram}{staircase diagram} can be uniquely reconstructed from its \glslink{corner-set-term}{corner set}. If $\beta$ is \glslink{left-leaning-stitch}{left-leaning}, then similar arguments show that each of the maps \begin{align*}
        \left\{\begin{gathered}
            \text{\glslink{reduced}{reduced}, \glslink{left-leaning-stitch}{left-leaning}} \\
            \text{\glslink{stitch-antiderivative}{antiderivatives} of $\beta$}
        \end{gathered}\right\} &\xrightarrow{\alpha \mapsto \glslink{stitch-splitting}{\hat{\alpha}}} \left\{\begin{gathered}
            \text{\glslink{split}{split}, \glslink{reduced}{reduced}, \glslink{left-leaning-stitch}{left-leaning} \glslink{stitch}{stitches} $\alpha_0$ with $k$ \glslink{wire}{wires} such that} \\
            \text{for each $(i,j) \in \beta$, $\alpha_0$ has a \glslink{join-site}{join site} corresponding} \\
            \text{to $(i,j)$ under Lemma~\ref{staircase-bij-lem}(iii)}
        \end{gathered}\right\} \\
        &\xrightarrow{\glslink{stair}{\stair}} \left\{\begin{gathered}
            \text{$k \times k$ \glslink{staircase-diagram}{staircase diagrams} lying below the diagonal} \\
            \text{with \glslink{corner-set-term}{corners} at every \glslink{wire}{wire} of $\beta$}
        \end{gathered}\right\} \\
        &\xrightarrow{S \mapsto (\glslink{bracket-n}{[k]}, \glslink{corner-set}{\Gamma(S)})} \{\text{\glslink{left-leaning-stitch}{left-leaning} \glslink{stitch-extension}{extensions} of $\beta$}\}
    \end{align*} is a bijection.
\end{proof}

\noindent
See Figure~\ref{antiderivative-bijection-fig} for an illustration of the procedure used for computing $\phi(\alpha)$ for certain \glslink{stitch-antiderivative}{antiderivatives} $\alpha$ of a given \glslink{stitch}{stitch} $\beta$.

\begin{figure}
    \centering
    \begin{subfigure}{0.49\linewidth}
        \centering
        \begin{tikzpicture}[scale=0.5]
            \draw[orange,dashed] (1, 0) -- (1, 2); \fill[orange] (1, 0) circle (2pt) node[below=4pt]{1}; \fill[orange] (1, 2) circle (2pt) node[above=4pt]{1};
            \draw[black,dashed] (2, 0) -- (2, 2); \fill[black] (2, 0) circle (2pt) node[below=4pt]{2}; \fill[black] (2, 2) circle (2pt) node[above=4pt]{2};
            \draw[blue,dashed] (3, 0) -- (3, 2); \fill[blue] (3, 0) circle (2pt) node[below=4pt]{3}; \fill[blue] (3, 2) circle (2pt) node[above=4pt]{3};
            \draw[blue,dashed] (4, 0) -- (4, 2); \fill[blue] (4, 0) circle (2pt) node[below=4pt]{4}; \fill[blue] (4, 2) circle (2pt) node[above=4pt]{4};
            \draw[black,dashed] (5, 0) -- (5, 2); \fill[black] (5, 0) circle (2pt) node[below=4pt]{5}; \fill[black] (5, 2) circle (2pt) node[above=4pt]{5};
            \draw[black,dashed] (6, 0) -- (6, 2); \fill[black] (6, 0) circle (2pt) node[below=4pt]{6}; \fill[black] (6, 2) circle (2pt) node[above=4pt]{6};
            \draw[orange,dashed] (7, 0) -- (7, 2); \fill[orange] (7, 0) circle (2pt) node[below=4pt]{7}; \fill[orange] (7, 2) circle (2pt) node[above=4pt]{7};
            \draw[blue,dashed] (8, 0) -- (8, 2); \fill[blue] (8, 0) circle (2pt) node[below=4pt]{8}; \fill[blue] (8, 2) circle (2pt) node[above=4pt]{8};
            \draw[blue,dashed] (9, 0) -- (9, 2); \fill[blue] (9, 0) circle (2pt) node[below=4pt]{9}; \fill[blue] (9, 2) circle (2pt) node[above=4pt]{9};
            \draw[black,dashed] (10, 0) -- (10, 2); \fill[black] (10, 0) circle (2pt) node[below=4pt]{10}; \fill[black] (10, 2) circle (2pt) node[above=4pt]{10};
            \draw[orange,dashed] (11, 0) -- (11, 2); \fill[orange] (11, 0) circle (2pt) node[below=4pt]{11}; \fill[orange] (11, 2) circle (2pt) node[above=4pt]{11};
            \draw[black,dashed] (12, 0) -- (12, 2); \fill[black] (12, 0) circle (2pt) node[below=4pt]{12}; \fill[black] (12, 2) circle (2pt) node[above=4pt]{12};
            \draw[black,dashed] (13, 0) -- (13, 2); \fill[black] (13, 0) circle (2pt) node[below=4pt]{13}; \fill[black] (13, 2) circle (2pt) node[above=4pt]{13};
            \draw[black,thick] plot [] coordinates {(1, 0)(1, 2)};
            \draw[black,thick] plot [] coordinates {(11, 0)(11, 2)};
            \draw[black,thick] plot [] coordinates {(2, 0)} plot [] coordinates {(3, 0)(2, 2)};
            \draw[black,thick] plot [] coordinates {(4, 0)} plot [] coordinates {(5, 0)(4, 2)};
            \draw[black,thick] plot [] coordinates {(7, 0)} plot [] coordinates {(6, 0)(7, 2)};
            \draw[black,thick] plot [] coordinates {(10, 0)} plot [] coordinates {(8, 0)(10, 2)};
            \draw[black,thick] plot [] coordinates {(7, 0)} plot [] coordinates {(6, 0)(7, 2)};
            \draw[black,thick] plot [] coordinates {(12, 0)} plot [] coordinates {(13, 0)(12, 2)};
            \draw[black,thick] plot [] coordinates {(7, 0)(9, 2)};
        \end{tikzpicture}

        \begin{tikzcd}
            {} \arrow[d, "\alpha \mapsto \glslink{stitch-splitting}{\hat{\alpha}}"] \\
            {}
        \end{tikzcd}

        \scalebox{0.85}{\begin{tikzpicture}[scale=0.5]
            \draw[orange,dashed] (1, 0) -- (1, 2); \fill[orange] (1, 0) circle (2pt) node[below=4pt]{1}; \fill[orange] (1, 2) circle (2pt) node[above=4pt]{1};
            \draw[orange,dashed] (2, 0) -- (2, 2); \fill[orange] (2, 0) circle (2pt) node[below=4pt]{2}; \fill[orange] (2, 2) circle (2pt) node[above=4pt]{2};
            \draw[black,dashed] (3, 0) -- (3, 2); \fill[black] (3, 0) circle (2pt) node[below=4pt]{3}; \fill[black] (3, 2) circle (2pt) node[above=4pt]{3};
            \draw[blue,dashed] (4, 0) -- (4, 2); \fill[blue] (4, 0) circle (2pt) node[below=4pt]{4}; \fill[blue] (4, 2) circle (2pt) node[above=4pt]{4};
            \draw[blue,dashed] (5, 0) -- (5, 2); \fill[blue] (5, 0) circle (2pt) node[below=4pt]{5}; \fill[blue] (5, 2) circle (2pt) node[above=4pt]{5};
            \draw[black,dashed] (6, 0) -- (6, 2); \fill[black] (6, 0) circle (2pt) node[below=4pt]{6}; \fill[black] (6, 2) circle (2pt) node[above=4pt]{6};
            \draw[black,dashed] (7, 0) -- (7, 2); \fill[black] (7, 0) circle (2pt) node[below=4pt]{7}; \fill[black] (7, 2) circle (2pt) node[above=4pt]{7};
            \draw[orange,dashed] (8, 0) -- (8, 2); \fill[orange] (8, 0) circle (2pt) node[below=4pt]{8}; \fill[orange] (8, 2) circle (2pt) node[above=4pt]{8};
            \draw[orange,dashed] (9, 0) -- (9, 2); \fill[orange] (9, 0) circle (2pt) node[below=4pt]{9}; \fill[orange] (9, 2) circle (2pt) node[above=4pt]{9};
            \draw[blue,dashed] (10, 0) -- (10, 2); \fill[blue] (10, 0) circle (2pt) node[below=4pt]{10}; \fill[blue] (10, 2) circle (2pt) node[above=4pt]{10};
            \draw[blue,dashed] (11, 0) -- (11, 2); \fill[blue] (11, 0) circle (2pt) node[below=4pt]{11}; \fill[blue] (11, 2) circle (2pt) node[above=4pt]{11};
            \draw[black,dashed] (12, 0) -- (12, 2); \fill[black] (12, 0) circle (2pt) node[below=4pt]{12}; \fill[black] (12, 2) circle (2pt) node[above=4pt]{12};
            \draw[orange,dashed] (13, 0) -- (13, 2); \fill[orange] (13, 0) circle (2pt) node[below=4pt]{13}; \fill[orange] (13, 2) circle (2pt) node[above=4pt]{13};
            \draw[orange,dashed] (14, 0) -- (14, 2); \fill[orange] (14, 0) circle (2pt) node[below=4pt]{14}; \fill[orange] (14, 2) circle (2pt) node[above=4pt]{14};
            \draw[black,dashed] (15, 0) -- (15, 2); \fill[black] (15, 0) circle (2pt) node[below=4pt]{15}; \fill[black] (15, 2) circle (2pt) node[above=4pt]{15};
            \draw[black,dashed] (16, 0) -- (16, 2); \fill[black] (16, 0) circle (2pt) node[below=4pt]{16}; \fill[black] (16, 2) circle (2pt) node[above=4pt]{16};
            \draw[black,thick] plot [] coordinates {(2, 0)} plot [] coordinates {(1, 0)(2, 2)};
            \draw[black,thick] plot [] coordinates {(3, 0)} plot [] coordinates {(4, 0)(3, 2)};
            \draw[black,thick] plot [] coordinates {(5, 0)} plot [] coordinates {(6, 0)(5, 2)};
            \draw[black,thick] plot [] coordinates {(9, 0)} plot [] coordinates {(7, 0)(9, 2)};
            \draw[black,thick] plot [] coordinates {(11, 0)} plot [] coordinates {(8, 0)(11, 2)};
            \draw[black,thick] plot [] coordinates {(12, 0)} plot [] coordinates {(10, 0)(12, 2)};
            \draw[black,thick] plot [] coordinates {(14, 0)} plot [] coordinates {(13, 0)(14, 2)};
            \draw[black,thick] plot [] coordinates {(15, 0)} plot [] coordinates {(16, 0)(15, 2)};
        \end{tikzpicture}}

        \begin{tikzcd}
            {} \arrow[d, "\glslink{stair}{\stair}"] \\
            {}
        \end{tikzcd}

        \begin{tikzpicture}[scale=0.7]
            \draw[black] (0,0) -- (8,0);
            \draw[black] (0,1) -- (8,1);
            \draw[black] (0,2) -- (8,2);
            \draw[black] (0,3) -- (8,3);
            \draw[black] (0,4) -- (8,4);
            \draw[black] (0,5) -- (8,5);
            \draw[black] (0,6) -- (8,6);
            \draw[black] (0,7) -- (8,7);
            \draw[black] (0,8) -- (8,8);

            \draw[black] (0,0) -- (0,8);
            \draw[black] (1,0) -- (1,8);
            \draw[black] (2,0) -- (2,8);
            \draw[black] (3,0) -- (3,8);
            \draw[black] (4,0) -- (4,8);
            \draw[black] (5,0) -- (5,8);
            \draw[black] (6,0) -- (6,8);
            \draw[black] (7,0) -- (7,8);
            \draw[black] (8,0) -- (8,8);

            \fill[black,opacity=0.3] (0,0) -- (0,1) -- (2,1) -- (2,2) -- (3,2) -- (3,5) -- (4,5) -- (4,6) -- (6,6) -- (6,7) -- (8,7) -- (8,0);

            \fill[orange, opacity=0.4] (0,0) -- (0,1) -- (1,1) -- (1,0);
            \draw[black] (0.5,0.5) node {$\scriptstyle (1,1)$};
            \fill[orange, opacity=0.4] (3,4) -- (3,5) -- (4,5) -- (4,4);
            \draw[black] (3.5,4.5) node {$\scriptstyle (4,5)$};
            \fill[orange, opacity=0.4] (6,6) -- (6,7) -- (7,7) -- (7,6);
            \draw[black] (6.5,6.5) node {$\scriptstyle (7,7)$};

            \fill[blue, opacity=0.4] (2,1) -- (2,2) -- (3,2) -- (3,1);
            \draw[black] (2.5,1.5) node {$\scriptstyle (3,2)$};
            \fill[blue, opacity=0.4] (4,5) -- (4,6) -- (5,6) -- (5,5);
            \draw[black] (4.5,5.5) node {$\scriptstyle (5,6)$};

            \draw[black,very thick] (0,0) -- (0,1) node[text=orange,midway,shift={(-0.2,0.1)}] {$\scriptstyle 1$} -- (1,1) node[text=orange,midway,shift={(-0.1,0.2)}] {$\scriptstyle 2$} -- (2,1) node[midway,shift={(-0.1,0.2)}] {$\scriptstyle 3$} -- (2,2) node[text=blue,midway,shift={(-0.2,0.1)}] {$\scriptstyle 4$} -- (3,2) node[text=blue,midway,shift={(-0.1,0.2)}] {$\scriptstyle 5$} -- (3,3) node[midway,shift={(-0.2,0.1)}] {$\scriptstyle 6$} -- (3,4) node[midway,shift={(-0.2,0.1)}] {$\scriptstyle 7$} -- (3,5) node[text=orange,midway,shift={(-0.2,0.1)}] {$\scriptstyle 8$} -- (4,5) node[text=orange,midway,shift={(-0.1, 0.2)}] {$\scriptstyle 9$} -- (4,6) node[text=blue,midway,shift={(-0.2,0.1)}] {$\scriptstyle 10$} -- (5,6) node[text=blue,midway,shift={(-0.1, 0.2)}] {$\scriptstyle 11$} -- (6,6) node[midway,shift={(-0.1, 0.2)}] {$\scriptstyle 12$} -- (6,7) node[text=orange,midway,shift={(-0.2,0.1)}] {$\scriptstyle 13$} -- (7,7) node[text=orange,midway,shift={(-0.1, 0.2)}] {$\scriptstyle 14$} -- (8,7) node[midway,shift={(-0.1, 0.2)}] {$\scriptstyle 15$} -- (8,8) node[midway,shift={(-0.2,0.1)}] {$\scriptstyle 16$};
        \end{tikzpicture}

        \begin{tikzcd}
            {} \arrow[d, "{S \mapsto (\glslink{bracket-n}{[8]}, \glslink{corner-set}{\Gamma(S)})}"] \\
            {}
        \end{tikzcd}

        \begin{tikzpicture}[scale=0.5]
            \draw[black,dashed] (1, 0) -- (1, 2); \fill[black] (1, 0) circle (2pt) node[below=4pt]{1}; \fill[black] (1, 2) circle (2pt) node[above=4pt]{1};
            \draw[black,dashed] (2, 0) -- (2, 2); \fill[black] (2, 0) circle (2pt) node[below=4pt]{2}; \fill[black] (2, 2) circle (2pt) node[above=4pt]{2};
            \draw[black,dashed] (3, 0) -- (3, 2); \fill[black] (3, 0) circle (2pt) node[below=4pt]{3}; \fill[black] (3, 2) circle (2pt) node[above=4pt]{3};
            \draw[black,dashed] (4, 0) -- (4, 2); \fill[black] (4, 0) circle (2pt) node[below=4pt]{4}; \fill[black] (4, 2) circle (2pt) node[above=4pt]{4};
            \draw[black,dashed] (5, 0) -- (5, 2); \fill[black] (5, 0) circle (2pt) node[below=4pt]{5}; \fill[black] (5, 2) circle (2pt) node[above=4pt]{5};
            \draw[black,dashed] (6, 0) -- (6, 2); \fill[black] (6, 0) circle (2pt) node[below=4pt]{6}; \fill[black] (6, 2) circle (2pt) node[above=4pt]{6};
            \draw[black,dashed] (7, 0) -- (7, 2); \fill[black] (7, 0) circle (2pt) node[below=4pt]{7}; \fill[black] (7, 2) circle (2pt) node[above=4pt]{7};
            \draw[black,dashed] (8, 0) -- (8, 2); \fill[black] (8, 0) circle (2pt) node[below=4pt]{8}; \fill[black] (8, 2) circle (2pt) node[above=4pt]{8};
            \draw[orange,thick] plot [] coordinates {(1, 0)(1, 2)};
            \draw[blue,thick] plot [] coordinates {(3, 0)(2, 2)};
            \draw[orange,thick] plot [] coordinates {(7, 0)(7, 2)};
            \draw[blue,thick] plot [] coordinates {(2, 0)};
            \draw[orange,thick] plot [] coordinates {(4, 0)(5, 2)};
            \draw[blue,thick] plot [] coordinates {(5, 0)(6, 2)};
            \draw[blue,thick] plot [] coordinates {(8, 0)};
        \end{tikzpicture}
    \end{subfigure} \begin{subfigure}{0.49\linewidth}
        \centering
        \begin{tikzpicture}[scale=0.5]
            \draw[orange,dashed] (1, 0) -- (1, 2); \fill[orange] (1, 0) circle (2pt) node[below=4pt]{1}; \fill[orange] (1, 2) circle (2pt) node[above=4pt]{1};
            \draw[black,dashed] (2, 0) -- (2, 2); \fill[black] (2, 0) circle (2pt) node[below=4pt]{2}; \fill[black] (2, 2) circle (2pt) node[above=4pt]{2};
            \draw[black,dashed] (3, 0) -- (3, 2); \fill[black] (3, 0) circle (2pt) node[below=4pt]{3}; \fill[black] (3, 2) circle (2pt) node[above=4pt]{3};
            \draw[blue,dashed] (4, 0) -- (4, 2); \fill[blue] (4, 0) circle (2pt) node[below=4pt]{4}; \fill[blue] (4, 2) circle (2pt) node[above=4pt]{4};
            \draw[blue,dashed] (5, 0) -- (5, 2); \fill[blue] (5, 0) circle (2pt) node[below=4pt]{5}; \fill[blue] (5, 2) circle (2pt) node[above=4pt]{5};
            \draw[black,dashed] (6, 0) -- (6, 2); \fill[black] (6, 0) circle (2pt) node[below=4pt]{6}; \fill[black] (6, 2) circle (2pt) node[above=4pt]{6};
            \draw[orange,dashed] (7, 0) -- (7, 2); \fill[orange] (7, 0) circle (2pt) node[below=4pt]{7}; \fill[orange] (7, 2) circle (2pt) node[above=4pt]{7};
            \draw[black,dashed] (8, 0) -- (8, 2); \fill[black] (8, 0) circle (2pt) node[below=4pt]{8}; \fill[black] (8, 2) circle (2pt) node[above=4pt]{8};
            \draw[black,dashed] (9, 0) -- (9, 2); \fill[black] (9, 0) circle (2pt) node[below=4pt]{9}; \fill[black] (9, 2) circle (2pt) node[above=4pt]{9};
            \draw[black,dashed] (10, 0) -- (10, 2); \fill[black] (10, 0) circle (2pt) node[below=4pt]{10}; \fill[black] (10, 2) circle (2pt) node[above=4pt]{10};
            \draw[orange,dashed] (11, 0) -- (11, 2); \fill[orange] (11, 0) circle (2pt) node[below=4pt]{11}; \fill[orange] (11, 2) circle (2pt) node[above=4pt]{11};
            \draw[blue,dashed] (12, 0) -- (12, 2); \fill[blue] (12, 0) circle (2pt) node[below=4pt]{12}; \fill[blue] (12, 2) circle (2pt) node[above=4pt]{12};
            \draw[blue,dashed] (13, 0) -- (13, 2); \fill[blue] (13, 0) circle (2pt) node[below=4pt]{13}; \fill[blue] (13, 2) circle (2pt) node[above=4pt]{13};
            \draw[black,thick] plot [] coordinates {(1, 0)(1, 2)};
            \draw[black,thick] plot [] coordinates {(11, 0)(11, 2)};
            \draw[black,thick] plot [] coordinates {(5, 0)} plot [] coordinates {(2, 0)(5, 2)};
            \draw[black,thick] plot [] coordinates {(6, 0)} plot [] coordinates {(3, 0)(6, 2)};
            \draw[black,thick] plot [] coordinates {(8, 0)} plot [] coordinates {(7, 0)(8, 2)} plot [] coordinates {(4, 0)(7, 2)};
            \draw[black,thick] plot [] coordinates {(9, 0)} plot [] coordinates {(10, 0)(9, 2)};
            \draw[black,thick] plot [] coordinates {(13, 0)} plot [] coordinates {(12, 0)(13, 2)};
        \end{tikzpicture}

        \begin{tikzcd}
            {} \arrow[d, "\alpha \mapsto \glslink{stitch-splitting}{\hat{\alpha}}"] \\
            {}
        \end{tikzcd}

        \scalebox{0.85}{\begin{tikzpicture}[scale=0.5]
            \draw[orange,dashed] (1, 0) -- (1, 2); \fill[orange] (1, 0) circle (2pt) node[below=4pt]{1}; \fill[orange] (1, 2) circle (2pt) node[above=4pt]{1};
            \draw[orange,dashed] (2, 0) -- (2, 2); \fill[orange] (2, 0) circle (2pt) node[below=4pt]{2}; \fill[orange] (2, 2) circle (2pt) node[above=4pt]{2};
            \draw[black,dashed] (3, 0) -- (3, 2); \fill[black] (3, 0) circle (2pt) node[below=4pt]{3}; \fill[black] (3, 2) circle (2pt) node[above=4pt]{3};
            \draw[black,dashed] (4, 0) -- (4, 2); \fill[black] (4, 0) circle (2pt) node[below=4pt]{4}; \fill[black] (4, 2) circle (2pt) node[above=4pt]{4};
            \draw[blue,dashed] (5, 0) -- (5, 2); \fill[blue] (5, 0) circle (2pt) node[below=4pt]{5}; \fill[blue] (5, 2) circle (2pt) node[above=4pt]{5};
            \draw[blue,dashed] (6, 0) -- (6, 2); \fill[blue] (6, 0) circle (2pt) node[below=4pt]{6}; \fill[blue] (6, 2) circle (2pt) node[above=4pt]{6};
            \draw[black,dashed] (7, 0) -- (7, 2); \fill[black] (7, 0) circle (2pt) node[below=4pt]{7}; \fill[black] (7, 2) circle (2pt) node[above=4pt]{7};
            \draw[orange,dashed] (8, 0) -- (8, 2); \fill[orange] (8, 0) circle (2pt) node[below=4pt]{8}; \fill[orange] (8, 2) circle (2pt) node[above=4pt]{8};
            \draw[orange,dashed] (9, 0) -- (9, 2); \fill[orange] (9, 0) circle (2pt) node[below=4pt]{9}; \fill[orange] (9, 2) circle (2pt) node[above=4pt]{9};
            \draw[black,dashed] (10, 0) -- (10, 2); \fill[black] (10, 0) circle (2pt) node[below=4pt]{10}; \fill[black] (10, 2) circle (2pt) node[above=4pt]{10};
            \draw[black,dashed] (11, 0) -- (11, 2); \fill[black] (11, 0) circle (2pt) node[below=4pt]{11}; \fill[black] (11, 2) circle (2pt) node[above=4pt]{11};
            \draw[black,dashed] (12, 0) -- (12, 2); \fill[black] (12, 0) circle (2pt) node[below=4pt]{12}; \fill[black] (12, 2) circle (2pt) node[above=4pt]{12};
            \draw[orange,dashed] (13, 0) -- (13, 2); \fill[orange] (13, 0) circle (2pt) node[below=4pt]{13}; \fill[orange] (13, 2) circle (2pt) node[above=4pt]{13};
            \draw[orange,dashed] (14, 0) -- (14, 2); \fill[orange] (14, 0) circle (2pt) node[below=4pt]{14}; \fill[orange] (14, 2) circle (2pt) node[above=4pt]{14};
            \draw[blue,dashed] (15, 0) -- (15, 2); \fill[blue] (15, 0) circle (2pt) node[below=4pt]{15}; \fill[blue] (15, 2) circle (2pt) node[above=4pt]{15};
            \draw[blue,dashed] (16, 0) -- (16, 2); \fill[blue] (16, 0) circle (2pt) node[below=4pt]{16}; \fill[blue] (16, 2) circle (2pt) node[above=4pt]{16};
            \draw[black,thick] plot [] coordinates {(2, 0)} plot [] coordinates {(1, 0)(2, 2)};
            \draw[black,thick] plot [] coordinates {(6, 0)} plot [] coordinates {(3, 0)(6, 2)};
            \draw[black,thick] plot [] coordinates {(7, 0)} plot [] coordinates {(4, 0)(7, 2)};
            \draw[black,thick] plot [] coordinates {(9, 0)} plot [] coordinates {(5, 0)(9, 2)};
            \draw[black,thick] plot [] coordinates {(10, 0)} plot [] coordinates {(8, 0)(10, 2)};
            \draw[black,thick] plot [] coordinates {(11, 0)} plot [] coordinates {(12, 0)(11, 2)};
            \draw[black,thick] plot [] coordinates {(14, 0)} plot [] coordinates {(13, 0)(14, 2)};
            \draw[black,thick] plot [] coordinates {(16, 0)} plot [] coordinates {(15, 0)(16, 2)};
        \end{tikzpicture}}

        \begin{tikzcd}
            {} \arrow[d, "\glslink{stair}{\stair}"] \\
            {}
        \end{tikzcd}

        \vspace{-10pt}
        \begin{tikzpicture}[scale=0.7]
            \draw[black] (0,0) -- (8,0);
            \draw[black] (0,1) -- (8,1);
            \draw[black] (0,2) -- (8,2);
            \draw[black] (0,3) -- (8,3);
            \draw[black] (0,4) -- (8,4);
            \draw[black] (0,5) -- (8,5);
            \draw[black] (0,6) -- (8,6);
            \draw[black] (0,7) -- (8,7);
            \draw[black] (0,8) -- (8,8);

            \draw[black] (0,0) -- (0,8);
            \draw[black] (1,0) -- (1,8);
            \draw[black] (2,0) -- (2,8);
            \draw[black] (3,0) -- (3,8);
            \draw[black] (4,0) -- (4,8);
            \draw[black] (5,0) -- (5,8);
            \draw[black] (6,0) -- (6,8);
            \draw[black] (7,0) -- (7,8);
            \draw[black] (8,0) -- (8,8);

            \fill[black,opacity=0.3] (0,0) -- (0,1) -- (1,1) -- (1,4) -- (3,4) -- (3,5) -- (6,5) -- (6,7) -- (7,7) -- (7,8) -- (8,8) -- (8,0);

            \fill[orange, opacity=0.4] (0,0) -- (0,1) -- (1,1) -- (1,0);
            \draw[black] (0.5,0.5) node {$\scriptstyle (1,1)$};
            \fill[orange, opacity=0.4] (3,4) -- (3,5) -- (4,5) -- (4,4);
            \draw[black] (3.5,4.5) node {$\scriptstyle (4,5)$};
            \fill[orange, opacity=0.4] (6,6) -- (6,7) -- (7,7) -- (7,6);
            \draw[black] (6.5,6.5) node {$\scriptstyle (7,7)$};

            \fill[blue, opacity=0.4] (1,3) -- (1,4) -- (2,4) -- (2,3);
            \draw[black] (1.5,3.5) node {$\scriptstyle (2,4)$};
            \fill[blue, opacity=0.4] (7,7) -- (7,8) -- (8,8) -- (8,7);
            \draw[black] (7.5,7.5) node {$\scriptstyle (7,7)$};

            \draw[black,very thick] (0,0) -- (0,1) node[text=orange,midway,shift={(-0.2,0.1)}] {$\scriptstyle 1$} -- (1,1) node[text=orange,midway,shift={(-0.1,0.2)}] {$\scriptstyle 2$} -- (1,2) node[midway,shift={(-0.2,0.1)}] {$\scriptstyle 3$} -- (1,3) node[midway,shift={(-0.2,0.1)}] {$\scriptstyle 4$} -- (1,4) node[text=blue,midway,shift={(-0.2,0.1)}] {$\scriptstyle 5$} -- (2,4) node[text=blue,midway,shift={(-0.1,0.2)}] {$\scriptstyle 6$} -- (3,4) node[midway,shift={(-0.1,0.2)}] {$\scriptstyle 7$} -- (3,5) node[text=orange,midway,shift={(-0.2,0.1)}] {$\scriptstyle 8$} -- (4,5) node[text=orange,midway,shift={(-0.1, 0.2)}] {$\scriptstyle 9$} -- (5,5) node[midway,shift={(-0.1,0.2)}] {$\scriptstyle 10$} -- (6,5) node[midway,shift={(-0.1, 0.2)}] {$\scriptstyle 11$} -- (6,6) node[midway,shift={(-0.2, 0.1)}] {$\scriptstyle 12$} -- (6,7) node[text=orange,midway,shift={(-0.2,0.1)}] {$\scriptstyle 13$} -- (7,7) node[text=orange,midway,shift={(-0.1, 0.2)}] {$\scriptstyle 14$} -- (7,8) node[text=blue,midway,shift={(-0.2, 0.1)}] {$\scriptstyle 15$} -- (8,8) node[text=blue,midway,shift={(-0.1,0.2)}] {$\scriptstyle 16$};
        \end{tikzpicture}

        \begin{tikzcd}
            {} \arrow[d, "{S \mapsto (\glslink{bracket-n}{[8]}, \glslink{corner-set}{\Gamma(S)})}"] \\
            {}
        \end{tikzcd}

        \begin{tikzpicture}[scale=0.5]
            \draw[black,dashed] (1, 0) -- (1, 2); \fill[black] (1, 0) circle (2pt) node[below=4pt]{1}; \fill[black] (1, 2) circle (2pt) node[above=4pt]{1};
            \draw[black,dashed] (2, 0) -- (2, 2); \fill[black] (2, 0) circle (2pt) node[below=4pt]{2}; \fill[black] (2, 2) circle (2pt) node[above=4pt]{2};
            \draw[black,dashed] (3, 0) -- (3, 2); \fill[black] (3, 0) circle (2pt) node[below=4pt]{3}; \fill[black] (3, 2) circle (2pt) node[above=4pt]{3};
            \draw[black,dashed] (4, 0) -- (4, 2); \fill[black] (4, 0) circle (2pt) node[below=4pt]{4}; \fill[black] (4, 2) circle (2pt) node[above=4pt]{4};
            \draw[black,dashed] (5, 0) -- (5, 2); \fill[black] (5, 0) circle (2pt) node[below=4pt]{5}; \fill[black] (5, 2) circle (2pt) node[above=4pt]{5};
            \draw[black,dashed] (6, 0) -- (6, 2); \fill[black] (6, 0) circle (2pt) node[below=4pt]{6}; \fill[black] (6, 2) circle (2pt) node[above=4pt]{6};
            \draw[black,dashed] (7, 0) -- (7, 2); \fill[black] (7, 0) circle (2pt) node[below=4pt]{7}; \fill[black] (7, 2) circle (2pt) node[above=4pt]{7};
            \draw[black,dashed] (8, 0) -- (8, 2); \fill[black] (8, 0) circle (2pt) node[below=4pt]{8}; \fill[black] (8, 2) circle (2pt) node[above=4pt]{8};
            \draw[orange,thick] plot [] coordinates {(1, 0)(1, 2)};
            \draw[blue,thick] plot [] coordinates {(2, 0)(4, 2)};
            \draw[orange,thick] plot [] coordinates {(7, 0)(7, 2)};
            \draw[blue,thick] plot [] coordinates {(2, 0)};
            \draw[orange,thick] plot [] coordinates {(4, 0)(5, 2)};
            \draw[blue,thick] plot [] coordinates {(8, 0)(8, 2)};
            \draw[blue,thick] plot [] coordinates {(8, 0)};
        \end{tikzpicture}
    \end{subfigure}

    \caption{An illustration of the bijection from Theorem~\ref{antiderivative-thm} applied to two \glslink{stitch-antiderivative}{antiderivatives} of the \glslink{stitch}{stitch} $\beta = (\glslink{bracket-n}{[8]}, \{(1,1), (4,5), (7,7)\})$, yielding two different \glslink{stitch-extension}{extensions} of $\beta$.} \label{antiderivative-bijection-fig}
\end{figure}

\subsection{Counting \texorpdfstring{\glslink{left-leaning-stitch}{left-leaning}}{left-leaning} \texorpdfstring{\glslink{stitch}{stitches}}{stitches}} \label{counting-left-leaning-subsection}

While our main goal for this section will be a formula for the number $\lvert\glslink{calS-parallel-tau}{\calS^\parallel(\tau)}\rvert$ of arbitrary \glslink{parallel-stitch}{parallel} \glslink{stitch}{stitches} of a given \glslink{stitch-type}{type} $\tau$, it will prove useful to start by focusing on \glslink{left-leaning-stitch}{left-leaning} \glslink{stitch}{stitches}. Indeed, we can already combine Lemma~\ref{reduction-lem} and Theorem~\ref{antiderivative-thm} to obtain the following recursive formula for the number of \glslink{left-leaning-stitch}{left-leaning} \glslink{stitch}{stitches} of a given \glslink{stitch-type}{type}.

\begin{lemma} \label{recursion-lem}
    The number $\lvert\glslink{calS-lparallel-tau}{\calS^\lparallel(\tau)}\rvert$ of \glslink{left-leaning-stitch}{left-leaning} \glslink{stitch}{stitches} of a given \glslink{cycle}{cycle}-free \glslink{stitch-type}{stitch type} $\tau \in \glslink{calT-lparallel}{\calT^\lparallel}$ may be computed recursively using the following rules: \begin{itemize}
        \item[(i)] We have \begin{equation*}
            \lvert\glslink{calS-lparallel-tau}{\calS^\lparallel((\varnothing,\varnothing))}\rvert = 1.
        \end{equation*}
        \item[(ii)] For any \glslink{cycle}{cycle}-free \glslink{stitch-type}{stitch type} $\tau$ on $n$ \glslink{column}{columns} containing $m \ge 1$ \glslink{gap}{gaps}, we have \begin{equation*}
            \lvert\glslink{calS-lparallel-tau}{\calS^\lparallel(\tau)}\rvert = \glslink{binomial}{\binom{n}{m}} \cdot \lvert\glslink{calS-lparallel-tau}{\calS^\lparallel(\glslink{bar-tau}{\bar{\tau}})}\rvert.
        \end{equation*}
        \item[(iii)] For any \glslink{reduced}{reduced}, \glslink{cycle}{cycle}-free \glslink{stitch-type}{stitch type} $\tau$ with $k$ \glslink{wire}{wires}, we have \begin{equation*}
            \lvert\glslink{calS-lparallel-tau}{\calS^\lparallel(\tau)}\rvert = \sum_{\omega \in \glslink{calT-n-lparallel}{\calT_k^\lparallel}} \glslink{stitch-binom}{\binom{\omega}{\glslink{stitch-type-derivative}{\partial\tau}}} \cdot \lvert\glslink{calS-lparallel-tau}{\calS^\lparallel(\omega)}\rvert,
        \end{equation*} where $\glslink{calT-n-lparallel}{\calT_k^\lparallel}$ denotes the set of \glslink{cycle}{cycle}-free \glslink{stitch-type}{stitch types} with $k$ \glslink{column}{columns}.
    \end{itemize}
\end{lemma}

\begin{remark}
    As the proof suggests, Lemma~\ref{recursion-lem} holds with only slight modifications if we replace each $\glslink{calS-lparallel-tau}{\calS^\lparallel(\tau)}$ with $\glslink{calS-parallel-tau}{\calS^\parallel(\tau)}$. However, we will end up computing $\lvert\glslink{calS-parallel-tau}{\calS^\parallel(\tau)}\rvert$ using a different method, so we state only the recursive formula for $\lvert\glslink{calS-lparallel-tau}{\calS^\lparallel(\tau)}\rvert$.
\end{remark}

\begin{proof}
    Part (i) is immediate, and part (ii) is just the statement of Corollary~\ref{reduction-cor}. For part (iii), note first that because $\tau$ is \glslink{reduced}{reduced}, a \glslink{left-leaning-stitch}{left-leaning} \glslink{stitch}{stitch} is contained in $\glslink{calS-lparallel-tau}{\calS^\lparallel(\tau)}$ if and only if it is \glslink{reduced}{reduced} and its \glslink{stitch-derivative-term}{derivative} is contained in $\glslink{calS-lparallel-tau}{\calS^\lparallel(\glslink{stitch-type-derivative}{\partial\tau})}$. Applying Theorem~\ref{antiderivative-thm}, we conclude that \begin{align*}
        \lvert\glslink{calS-lparallel-tau}{\calS^\lparallel(\tau)}\rvert &= \sum_{\beta \in \glslink{calS-lparallel-tau}{\calS^\lparallel(\glslink{stitch-type-derivative}{\partial\tau})}} \#\{\text{\glslink{left-leaning-stitch}{left-leaning}, \glslink{reduced}{reduced} \glslink{stitch-antiderivative}{antiderivatives} of $\beta$}\} \\
        &= \sum_{\beta \in \glslink{calS-lparallel-tau}{\calS^\lparallel(\glslink{stitch-type-derivative}{\partial\tau})}} \#\{\text{\glslink{left-leaning-stitch}{left-leaning} \glslink{stitch-extension}{extensions} of $\beta$}\} \\
        &= \sum_{\beta \in \glslink{calS-lparallel-tau}{\calS^\lparallel(\glslink{stitch-type-derivative}{\partial\tau})}} \sum_{\omega \in \glslink{calT-n-lparallel}{\calT_k^\lparallel}} \#\{\text{\glslink{left-leaning-stitch}{left-leaning} \glslink{stitch-extension}{extensions} of $\beta$ with \glslink{stitch-type}{type} $\omega$}\} \\
        &= \sum_{\omega \in \glslink{calT-n-lparallel}{\calT_k^\lparallel}} \sum_{\beta \in \glslink{calS-lparallel-tau}{\calS^\lparallel(\glslink{stitch-type-derivative}{\partial\tau})}} \#\{\beta' \in \glslink{calS-lparallel-tau}{\calS^\lparallel(\omega)} : \glslink{stitch-subseteq}{\beta \subseteq \beta'}\} \\
        &= \sum_{\omega \in \glslink{calT-n-lparallel}{\calT_k^\lparallel}} \sum_{\beta' \in \glslink{calS-lparallel-tau}{\calS^\lparallel(\omega)}} \#\{\beta \in \glslink{calS-lparallel-tau}{\calS^\lparallel(\glslink{stitch-type-derivative}{\partial\tau})} : \glslink{stitch-subseteq}{\beta \subseteq \beta'}\}.
    \end{align*} Since a \glslink{substitch}{substitch} of a \glslink{left-leaning-stitch}{left-leaning} \glslink{stitch}{stitch} is still \glslink{left-leaning-stitch}{left-leaning}, this in turn becomes \begin{align*}
        \lvert\glslink{calS-lparallel-tau}{\calS^\lparallel(\tau)}\rvert &= \sum_{\omega \in \glslink{calT-n-lparallel}{\calT_k^\lparallel}} \sum_{\beta' \in \glslink{calS-lparallel-tau}{\calS^\lparallel(\omega)}} \#\{\beta \in \glslink{calS-tau}{\calS(\glslink{stitch-type-derivative}{\partial\tau})} : \glslink{stitch-subseteq}{\beta \subseteq \beta'}\} = \sum_{\omega \in \glslink{calT-n-lparallel}{\calT_k^\lparallel}} \glslink{stitch-binom}{\binom{\omega}{\glslink{stitch-type-derivative}{\partial\tau}}} \cdot \lvert\glslink{calS-lparallel-tau}{\calS^\lparallel(\omega)}\rvert.
    \end{align*}
\end{proof}

While the recursive procedure for computing $\lvert\glslink{calS-lparallel-tau}{\calS^\lparallel(\tau)}\rvert$ given in Lemma~\ref{recursion-lem} already gives us a more efficient algorithm for computing $\lvert\glslink{calS-lparallel-tau}{\calS^\lparallel(\tau)}\rvert$ than brute-force counting, it turns out we can do much better: we can obtain a surprisingly simple explicit formula for $\lvert\glslink{calS-lparallel-tau}{\calS^\lparallel(\tau)}\rvert$.

\begin{theorem} \label{left-leaning-count-thm}
    Let $\tau$ be any \glslink{cycle}{cycle}-free \glslink{stitch-type}{stitch type} with $n$ \glslink{column}{columns} and $k$ \glslink{wire}{wires}. Then \begin{equation*}
        \lvert\glslink{calS-lparallel-tau}{\calS^\lparallel(\tau)}\rvert = \frac{1}{\glslink{z-tau}{z_\tau}} \cdot \frac{n!}{(k+1)!}.
    \end{equation*} Equivalently, the probability $\frac{\lvert\glslink{calS-lparallel-tau}{\calS^\lparallel(\tau)}\rvert}{\lvert\glslink{calS-tau}{\calS(\tau)}\rvert}$ that a random \glslink{stitch}{stitch} of \glslink{stitch-type}{type} $\tau$ is \glslink{left-leaning-stitch}{left-leaning} is given by $\frac{1}{(k+1)!}$, independently of $\tau$.
\end{theorem}

\noindent
Before giving the proof of Theorem~\ref{left-leaning-count-thm}, we require one additional lemma.

\begin{lemma} \label{brackbinom-sum-lem}
    Let $\tau$ be a \glslink{cycle}{cycle}-free \glslink{stitch-type}{stitch type} with $n$ \glslink{column}{columns} and $k$ \glslink{wire}{wires}, and let $k \le k' \le n$. Then we have \begin{equation*}
        \sum_{\tau' \in \glslink{calT-n-k-lparallel}{\calT_{n,k'}^\lparallel}} \glslink{stitch-brackbinom}{\brackbinom{\tau'}{\tau}} = \frac{(n-k)!}{(n-k')!} \cdot \glslink{binomial}{\binom{n-k-1}{k'-k}},
    \end{equation*} where $\glslink{calT-n-k-lparallel}{\calT_{n,k'}^\lparallel}$ denotes the set of \glslink{stitch-type}{stitch types} with $n$ \glslink{column}{columns} and $k'$ \glslink{wire}{wires}.
\end{lemma}

\begin{proof}
    By the definition of the coefficient $\glslink{stitch-brackbinom}{\brackbinom{\tau'}{\tau}}$, we have \begin{equation}
        \sum_{\tau' \in \glslink{calT-n-k-lparallel}{\calT_{n,k'}^\lparallel}} \glslink{stitch-brackbinom}{\brackbinom{\tau'}{\tau}} = \#\{\alpha' \in \glslink{calS-n-k}{\calS_{n,k'}} : \text{$\glslink{stitch-subseteq}{\alpha_0 \subseteq \alpha'}$, $\alpha'$ \glslink{cycle}{cycle}-free}\}, \label{brackbinom-sum-eq}
    \end{equation} where $\alpha_0 \in \glslink{calS-tau}{\calS(\tau)}$ is an arbitrary representative. The fact that $\alpha_0$ has $n$ \glslink{column}{columns} and $k$ \glslink{wire}{wires} implies that it consists of $n-k$ \glslink{chain}{chains}. Similarly, any \glslink{stitch}{stitch} $\alpha' \in \glslink{calS-n-k}{\calS_{n,k'}}$ which is \glslink{cycle}{cycle}-free consists of $n-k'$ \glslink{chain}{chains}. Thus, the right hand side of \eqref{brackbinom-sum-eq} counts the number of ways one can ``glue'' the $n-k$ \glslink{chain}{chains} of $\alpha_0$ together to form $n-k'$ \glslink{chain}{chains}. Now, we can represent such a ``gluing'' of the \glslink{chain}{chains} of $\alpha_0$ as a way of listing these $n-k$ \glslink{chain}{chains} out in a particular order, then choosing $(n-k)-(n-k') = k'-k$ of the $n-k-1$ spaces between the \glslink{chain}{chains} to glue. The number of ways to order the \glslink{chain}{chains} is then $(n-k)!$, and the number of ways to choose which spaces to glue is then $\glslink{binomial}{\binom{n-k-1}{k'-k}}$. Using this system, each possible gluing of the \glslink{chain}{chains} of $\alpha_0$ can be represented in $(n-k')!$ ways, based on how the final glued \glslink{chain}{chains} are ordered. 
\end{proof}

\begin{proof}[Proof of Theorem~\ref{left-leaning-count-thm}]
    We proceed by induction on the number of \glslink{column}{columns} of $\tau$. If $\tau = (\varnothing,\varnothing)$, then the theorem holds trivially. Meanwhile, if $\tau$ is non-\glslink{reduced}{reduced} with $n$ \glslink{column}{columns} and $m$ \glslink{gap}{gaps}, then combining Lemma~\ref{recursion-lem}(ii) with the induction hypothesis yields \begin{equation*}
        \lvert\glslink{calS-lparallel-tau}{\calS^\lparallel(\tau)}\rvert = \glslink{binomial}{\binom{n}{m}} \cdot \lvert\glslink{calS-lparallel-tau}{\calS^\lparallel(\glslink{bar-tau}{\bar{\tau}})}\rvert = \glslink{binomial}{\binom{n}{m}} \cdot \frac{1}{\glslink{z-tau}{z_{\glslink{bar-tau}{\bar{\tau}}}}} \cdot \frac{(n-m)!}{(k+1)!} = \frac{1}{\glslink{z-tau}{z_\tau}} \cdot \frac{n!}{(k+1)!},
    \end{equation*} using the fact that $z_\tau = m!z_{\bar{\tau}}$. Finally, suppose $\tau \ne (\varnothing,\varnothing)$ is \glslink{reduced}{reduced} with $n$ \glslink{column}{columns} and $k$ \glslink{wire}{wires}, noting that this implies $\glslink{stitch-type-derivative}{\partial\tau}$ must have $s := 2k-n$ \glslink{wire}{wires}. Use Lemma~\ref{recursion-lem}(iii) together with the induction hypothesis to write \begin{align*}
        \lvert\glslink{calS-lparallel-tau}{\calS^\lparallel(\tau)}\rvert &= \sum_{\omega \in \glslink{calT-n-lparallel}{\calT_k^\lparallel}} \glslink{stitch-binom}{\binom{\omega}{\glslink{stitch-type-derivative}{\partial\tau}}} \cdot \lvert\glslink{calS-lparallel-tau}{\calS^\lparallel(\omega)}\rvert = \sum_{s \le s' \le k} \sum_{\omega \in \glslink{calT-n-k-lparallel}{\calT_{k,s'}^\lparallel}} \glslink{stitch-binom}{\binom{\omega}{\glslink{stitch-type-derivative}{\partial\tau}}} \cdot \frac{1}{\glslink{z-tau}{z_\omega}} \cdot \frac{k!}{(s'+1)!}.
    \end{align*} Applying Lemma~\ref{binom-relation-lem} (noting that the fact that $\tau$ is \glslink{reduced}{reduced} implies $\glslink{z-tau}{z_\tau} = \glslink{z-tau}{z_{\glslink{stitch-type-derivative}{\partial\tau}}}$) along with Lemma~\ref{brackbinom-sum-lem}, we then obtain \begin{align*}
        \glslink{z-tau}{z_\tau}\lvert\glslink{calS-lparallel-tau}{\calS^\lparallel(\tau)}\rvert &= \sum_{s \le s' \le k} \frac{k!}{(s'+1)!} \cdot \sum_{\omega \in \glslink{calT-n-k-lparallel}{\calT_{k,s'}^\lparallel}} \glslink{stitch-brackbinom}{\brackbinom{\omega}{\glslink{stitch-type-derivative}{\partial\tau}}} = \sum_{s \le s' \le k} \frac{k!}{(s'+1)!} \cdot \frac{(k-s)!}{(k-s')!} \cdot \glslink{binomial}{\binom{k-s-1}{s'-s}} \\
        &= \frac{k!(k-s)!}{(k+1)!} \cdot \sum_{t=0}^{k-s} \glslink{binomial}{\binom{k+1}{k-s-t}} \glslink{binomial}{\binom{k-s-1}{t}} \\
        &= \frac{k!(k-s)!}{(k+1)!} \glslink{binomial}{\binom{2k-s}{k-s}} = \frac{n!}{(k+1)!},
    \end{align*} where the second-to-last step follows by Vandermonde's identity.
\end{proof}

\subsection{Counting \texorpdfstring{\glslink{left-leaning-stitch}{left-leaning}}{left-leaning} \texorpdfstring{\glslink{stitch}{stitches}}{stitches} with a fixed number of \texorpdfstring{\glslink{irreducible-stitch}{irreducible}}{irreducible} components}

\noindent
With an explicit formula for $\lvert\glslink{calS-lparallel-tau}{\calS^\lparallel(\tau)}\rvert$ in hand, we now turn to the task of counting \glslink{left-leaning-stitch}{left-leaning} \glslink{stitch}{stitches} of a given \glslink{stitch-type}{type} with a fixed number of \glslink{irreducible-stitch}{irreducible} components. Towards this end, let $\tau$ be a \glslink{cycle}{cycle}-free \glslink{stitch-type}{stitch type}, and for $\nu \ge 1$, let $\glslink{tilde-calS-nu-lparallel-tau}{\tilde{\calS}_\nu^\lparallel(\tau)}$ denote the set of \glslink{left-leaning-stitch}{left-leaning} \glslink{stitch}{stitches} of \glslink{stitch-type}{type} $\tau$ containing exactly $\nu$ \glslink{irreducible-stitch}{irreducible} components. While our ultimate goal is to compute the size of $\glslink{tilde-calS-nu-lparallel-tau}{\tilde{\calS}_\nu^\lparallel(\tau)}$, it helps to start with the closely related set $\glslink{hat-calS-nu-lparallel-tau}{\hat{\calS}_\nu^\lparallel(\tau)}$ of tuples $(\alpha_1,\dots,\alpha_\nu)$, where each $\alpha_i$ is a \glslink{empty-stitch}{nonempty} \glslink{left-leaning-stitch}{left-leaning} \glslink{stitch}{stitch} and $\glslink{oplus}{\alpha_1 \oplus \dots \oplus \alpha_\nu}$ has \glslink{stitch-type}{type} $\tau$. Intuitively, instead of counting \glslink{stitch}{stitches} with exactly $\nu$ parts in their \glslink{irreducible-stitch}{irreducible} decomposition, we begin by counting \glslink{stitch}{stitches} equipped with any chosen decomposition into $\nu$ (\glslink{empty-stitch}{nonempty}) parts.

When $\nu = 1$, $\glslink{hat-calS-nu-lparallel-tau}{\hat{\calS}_\nu^\lparallel(\tau)}$ just counts the number of \glslink{left-leaning-stitch}{left-leaning} \glslink{stitch}{stitches} of \glslink{stitch-type}{type} $\tau$, which we have already computed in Theorem~\ref{left-leaning-count-thm}. More generally, we will use Theorem~\ref{left-leaning-count-thm} together with the following lemma to derive an explicit formula for $\lvert\glslink{hat-calS-nu-lparallel-tau}{\hat{\calS}_\nu^\lparallel(\tau)}\rvert$.

\begin{lemma} \label{cayley-lem}
    Let $n_1, \dots, n_r \ge 1$ and $1 \le d \le r$, and for any subset $R \subseteq \glslink{bracket-n}{[r]}$, let $n_R := \sum_{i \in R} n_i$. Adopt the notation $\glslink{falling-factorial}{(n)_k}$ for the falling factorial $\frac{n!}{(n-k)!}$. Then \begin{equation}
        \sum_{\substack{R_1 \sqcup \cdots \sqcup R_d = \glslink{bracket-n}{[r]} \\ R_1,\dots,R_d \ne \varnothing}} \glslink{falling-factorial}{(n_{R_1})_{\lvert R_1\rvert-1}} \cdots \glslink{falling-factorial}{(n_{R_d})_{\lvert R_d\rvert-1}} = d \cdot \frac{(r-1)!}{(r-d)!} \cdot \glslink{falling-factorial}{(n_{\glslink{bracket-n}{[r]}})_{r-d}}. \label{n-R-sum-eq}
    \end{equation}
\end{lemma}

\begin{proof}
    Interpret each side of the equation as a polynomial in $n_1,\dots,n_r$, and define an invertible linear map $\phi \colon \Q[n_1,\dots,n_r] \to \Q[x_1,\dots,x_r]$ sending $\glslink{falling-factorial}{(n_1)_{a_1}} \cdots \glslink{falling-factorial}{(n_r)_{a_r}}$ to $x_1^{a_1} \cdots x_r^{a_r}$. Then since falling factorials form a polynomial sequence of binomial type, $\phi$ sends $\glslink{falling-factorial}{(n_R)_a}$ to $x_R^a$ for any subset $R \subseteq \glslink{bracket-n}{[r]}$ and any $a \ge 0$, where $x_R := \sum_{i \in R} x_i$. Applying $\phi$ to both sides of \eqref{n-R-sum-eq}, we thus see that it suffices to prove \begin{equation*}
        \sum_{\substack{R_1 \sqcup \cdots \sqcup R_d = \glslink{bracket-n}{[r]} \\ R_1,\dots,R_d \ne \varnothing}} x_{R_1}^{\lvert R_1\rvert-1} \cdots x_{R_d}^{\lvert R_d\rvert-1} = d \cdot \frac{(r-1)!}{(r-d)!} \cdot x_{\glslink{bracket-n}{[r]}}^{r-d}.
    \end{equation*} Multiplying by $x_{\glslink{bracket-n}{[r]}}^{d-2} \prod_{i \in \glslink{bracket-n}{[r]}} x_i$ and rearranging, we may equivalently state this as \begin{equation}
        \sum_{\substack{R_1 \sqcup \cdots \sqcup R_d = \glslink{bracket-n}{[r]} \\ R_1,\dots,R_d \ne \varnothing}} \left(x_{\glslink{bracket-n}{[r]}}^{d-2}\prod_{j \in \glslink{bracket-n}{[d]}} x_{R_j}\right) \prod_{j \in \glslink{bracket-n}{[d]}} \left(x_{R_j}^{\lvert R_j\rvert-2}\prod_{i \in R_j} x_i\right) = d \cdot \frac{(r-1)!}{(r-d)!} \cdot \left(x_{\glslink{bracket-n}{[r]}}^{r-2}\prod_{i \in \glslink{bracket-n}{[r]}} x_i\right). \label{x-R-sum-eq-cayley}
    \end{equation}
    
    To prove equation \eqref{x-R-sum-eq-cayley}, recall the following generalization of Cayley's formula:

    \begin{proposition}[Weighted form of Cayley's formula] \label{cayley-prop}
        Let $y_1, \dots, y_s$ be formal variables. Then \begin{equation*}
            (y_1+\dots+y_s)^{s-2} \cdot y_1 \cdots y_s = \sum_{\textnormal{tree $T$ on $\glslink{bracket-n}{[s]}$}} \prod_{i \in \glslink{bracket-n}{[s]}} y_i^{\glslink{graph-deg}{\deg_T(i)}}.
        \end{equation*}
    \end{proposition}

    \noindent
    By Proposition~\ref{cayley-prop} (applied first with $s = d$ and $y_j = x_{R_j}$, then, for each $j \in \glslink{bracket-n}{[d]}$, with $s = \lvert R_j\rvert$ and $\{y_1,\dots,y_s\} = \{x_i : i \in R_j\}$), the left hand side of \eqref{x-R-sum-eq-cayley} is equal to \begin{align}
        \sum_{\substack{R_1, \dots, R_d \\ U, T_1, \dots, T_d}} \left(\prod_{j \in \glslink{bracket-n}{[d]}} x_{R_j}^{\glslink{graph-deg}{\deg_U(j)}}\right) \prod_{j \in \glslink{bracket-n}{[d]}} \left(\prod_{i \in R_j} x_i^{\glslink{graph-deg}{\deg_{T_j}(i)}}\right), \label{R-U-T-formula}
    \end{align} where the sum is taken over (i) nonempty subsets $R_1,\dots,R_d$ satisfying $R_1 \sqcup \dots \sqcup R_d = \glslink{bracket-n}{[r]}$, together with (ii) trees $U$ on $\glslink{bracket-n}{[d]}$ and (iii) trees $T_j$ on $R_j$ for each $j$. Next, observe that for fixed $R_1, \dots, R_d$ and $U$, we may write $\prod_{j \in [d]} x_{R_j}^{\glslink{graph-deg}{\deg_U(j)}} = \sum_{U'} \prod_{i \in [r]} x_i^{\deg_{U'}(i)}$, where $U'$ varies over all graphs on $[r]$ which become the tree $U$ when we collapse each set $R_j$ to a single vertex. Thus, defining $\glslink{graph-deg}{\deg_{T_j}(i)} := 0$ whenever $i \notin R_j$, we can rewrite \eqref{R-U-T-formula} as \begin{equation}
        \sum_{\substack{R_1, \dots, R_d \\ U', T_1, \dots, T_d}} \prod_{i \in \glslink{bracket-n}{[r]}} x_i^{\glslink{graph-deg}{\deg_{U'}(i)} + \glslink{graph-deg}{\deg_{T_1}(i)} + \dots + \glslink{graph-deg}{\deg_{T_d}(i)}}, \label{split-tree-lhs}
    \end{equation} where we have replaced the summation over trees $U$ on $\glslink{bracket-n}{[d]}$ with a summation over graphs $U'$ on $\glslink{bracket-n}{[r]}$ which become trees on $\glslink{bracket-n}{[d]}$ when we collapse each set $R_j$ to a single vertex. Now observe that the data of $(R_1,\dots,R_d, U', T_1,\dots,T_d)$ is equivalent to the data of a single tree $T$ on $\glslink{bracket-n}{[r]}$, together with (i) a way of splitting $T$ into $d$ subtrees by removing $d-1$ edges, and (ii) a way of ordering these $d$ subtrees. Moreover, under this correspondence, the quantity $\glslink{graph-deg}{\deg_{U'}(i)} + \glslink{graph-deg}{\deg_{T_1}(i)} + \dots + \glslink{graph-deg}{\deg_{T_d}(i)}$ becomes simply $\glslink{graph-deg}{\deg_T(i)}$. For any fixed tree $T$ on $\glslink{bracket-n}{[r]}$, there are $\glslink{binomial}{\binom{r-1}{d-1}}$ choices of edges to remove to split it into $d$ subtrees, followed by $d!$ ways to order the resulting subtrees, so we conclude that \eqref{split-tree-lhs} is equal to \begin{equation*}
        d! \glslink{binomial}{\binom{r-1}{d-1}} \sum_{\text{tree $T$ on $\glslink{bracket-n}{[r]}$}} \prod_{i \in \glslink{bracket-n}{[r]}} x_i^{\glslink{graph-deg}{\deg_T(i)}} = d \cdot \frac{(r-1)!}{(r-d)!} \cdot \left(x_{\glslink{bracket-n}{[r]}}^{r-2} \prod_{i \in \glslink{bracket-n}{[r]}} x_i\right),
    \end{equation*} as desired.
\end{proof}

\begin{theorem} \label{left-leaning-proper-decomp-count-thm}
    Let $\tau$ be a \glslink{cycle}{cycle}-free \glslink{stitch-type}{stitch type} with $n$ \glslink{column}{columns}, $k$ \glslink{wire}{wires} and $j = n-k$ \glslink{chain}{chains}, and let $1 \le \nu \le j$. Then \begin{equation*}
        \lvert\glslink{hat-calS-nu-lparallel-tau}{\hat{\calS}_\nu^\lparallel(\tau)}\rvert = \frac{1}{\glslink{z-tau}{z_\tau}} \cdot \nu \cdot \frac{(j-1)!}{(j-\nu)!} \cdot \frac{n!}{(k+\nu)!}.
    \end{equation*}
\end{theorem}

\begin{proof}
    Write $\tau = (\varnothing, (n_1,\dots,n_j))$ for some $n_1,\dots,n_j \ge 1$, and observe that for any element $(\alpha_1,\dots,\alpha_\nu) \in \glslink{hat-calS-nu-lparallel-tau}{\hat{\calS}_\nu^\lparallel(\tau)}$, $\glslink{z-tau}{z_\tau}$ is equal to the number of ways of assigning a \glslink{stitch-coloring}{$j$-coloring} to the \glslink{stitch}{stitch} $\alpha := \glslink{oplus}{\alpha_1 \oplus \dots \oplus \alpha_\nu}$ such that for all $i \in \glslink{bracket-n}{[j]}$, the restriction $\glslink{alpha-rvert-i}{\alpha\rvert_i}$ is a \glslink{chain}{chain} of size $n_i$. Moreover, for each $i \in \glslink{bracket-n}{[j]}$, the \glslink{chain}{chain} of color $i$ must lie entirely in the ``component'' of $\alpha$ corresponding to one of the $\alpha_\ell$. Letting $R_\ell \subseteq \glslink{bracket-n}{[j]}$ denote the set of colors $i$ such that the \glslink{chain}{chain} of color $i$ lies in the components corresponding to $\alpha_\ell$, we thus obtain a partition $R_1 \sqcup \dots \sqcup R_\nu$ of $\glslink{bracket-n}{[j]}$ into nonempty subsets. If we fix such a partition $\glslink{bracket-n}{[j]} = R_1 \sqcup \dots \sqcup R_\nu$, then an element $(\alpha_1,\dots,\alpha_\nu) \in \glslink{hat-calS-nu-lparallel-tau}{\hat{\calS}_\nu^\lparallel(\tau)}$ with \glslink{stitch-coloring}{$j$-coloring} as described above gives rise to the partition $R_1 \sqcup \dots \sqcup R_\nu$ if and only if for all $\ell \in [\nu]$, $\alpha_\ell$ consists of exactly one \glslink{chain}{chain} of color $i$ and length $n_i$ for each $i \in R_\ell$. For a fixed $\ell \in [\nu]$, the number of \glslink{left-leaning-stitch}{left-leaning} \glslink{stitch}{stitches} with colored \glslink{chain}{chains} of lengths $(n_i)_{i \in R_\ell}$ is just \begin{equation*}
        \frac{(\sum_{i \in R_\ell} n_i)!}{(\sum_{i \in R_\ell}(n_i-1)+1)!}
    \end{equation*} by Theorem~\ref{left-leaning-count-thm}. Adopting the notation from Lemma~\ref{cayley-lem}, we deduce that \begin{equation*}
        \glslink{z-tau}{z_\tau}\lvert\glslink{hat-calS-nu-lparallel-tau}{\hat{\calS}_\nu^\lparallel(\tau)}\rvert = \sum_{\substack{R_1 \sqcup \cdots \sqcup R_\nu = \glslink{bracket-n}{[j]} \\ R_1,\dots,R_\nu \ne \varnothing}} \glslink{falling-factorial}{(n_{R_1})_{\lvert R_1\rvert-1}} \cdots \glslink{falling-factorial}{(n_{R_\nu})_{\lvert R_\nu\rvert-1}},
    \end{equation*} at which point Lemma~\ref{cayley-lem} gives the desired formula for $\glslink{z-tau}{z_\tau}\lvert\glslink{hat-calS-nu-lparallel-tau}{\hat{\calS}_\nu^\lparallel(\tau)}\rvert$.
\end{proof}

Having proven Theorem~\ref{left-leaning-proper-decomp-count-thm}, we quickly arrive at our desired formula for the number of \glslink{left-leaning-stitch}{left-leaning} \glslink{stitch}{stitches} of a fixed \glslink{stitch-type}{type} with a given number of \glslink{irreducible-stitch}{irreducible} components.

\begin{corollary} \label{left-leaning-irred-count-thm}
    Let $\tau$ be a \glslink{cycle}{cycle}-free \glslink{stitch-type}{stitch type} with $n$ \glslink{column}{columns}, $k \ge 1$ \glslink{wire}{wires}, and $j = n-k$ \glslink{chain}{chains}, let $1 \le \nu \le j$, and let $\glslink{tilde-calS-nu-lparallel-tau}{\tilde{\calS}_\nu^\lparallel(\tau)}$ denote the set of \glslink{left-leaning-stitch}{left-leaning} \glslink{stitch}{stitches} of \glslink{stitch-type}{type} $\tau$ containing exactly $\nu$ \glslink{irreducible-stitch}{irreducible} components. Then \begin{equation*}
        \lvert\glslink{tilde-calS-nu-lparallel-tau}{\tilde{\calS}_\nu^\lparallel(\tau)}\rvert = \frac{1}{\glslink{z-tau}{z_\tau}} \cdot \nu \cdot \frac{(j-1)!}{(j-\nu)!} \cdot \frac{(n-\nu-1)!}{(k-1)!}.
    \end{equation*}
\end{corollary}

\begin{proof}
    Observe that for any $c, \nu' \ge 1$, a \glslink{left-leaning-stitch}{left-leaning} \glslink{stitch}{stitch} $\alpha$ with $\nu'$ \glslink{irreducible-stitch}{irreducible} components can be decomposed as $\glslink{oplus}{\alpha_1 \oplus \dots \oplus \alpha_c}$ for some \glslink{empty-stitch}{nonempty} \glslink{stitch}{stitches} $\alpha_1,\dots,\alpha_c$ in exactly $\glslink{binomial}{\binom{\nu'-1}{c-1}}$ ways, allowing us to write \begin{equation*}
        \lvert\glslink{hat-calS-nu-lparallel-tau}{\hat{\calS}_c(\tau)}\rvert = \sum_{\nu' \ge c} \glslink{binomial}{\binom{\nu'-1}{c-1}} \cdot \lvert\glslink{tilde-calS-nu-lparallel-tau}{\tilde{\calS}_{\nu'}(\tau)}\rvert.
    \end{equation*} For any $\nu \ge 1$, we thus obtain \begin{align*}
        \sum_{c \ge \nu} (-1)^{c-\nu} \glslink{binomial}{\binom{c-1}{\nu-1}} \cdot \lvert\glslink{hat-calS-nu-lparallel-tau}{\hat{\calS}_c(\tau)}\rvert &= \sum_{\nu' \ge c \ge \nu} (-1)^{c-\nu} \glslink{binomial}{\binom{c-1}{\nu-1}} \glslink{binomial}{\binom{\nu'-1}{c-1}} \cdot \lvert\glslink{tilde-calS-nu-lparallel-tau}{\tilde{\calS}_{\nu'}(\tau)}\rvert \\
        &= \sum_{\nu' \ge \nu} \lvert\glslink{tilde-calS-nu-lparallel-tau}{\tilde{\calS}_{\nu'}(\tau)}\rvert \glslink{binomial}{\binom{\nu'-1}{\nu-1}}\sum_{\nu \le c \le \nu'} (-1)^{c-\nu} \glslink{binomial}{\binom{\nu'-\nu}{c-\nu}} \\
        &= \lvert\glslink{tilde-calS-nu-lparallel-tau}{\tilde{\calS}_\nu(\tau)}\rvert.
    \end{align*} From here, Theorem~\ref{left-leaning-proper-decomp-count-thm} gives \begin{align*}
        \glslink{z-tau}{z_\tau} \lvert\glslink{tilde-calS-nu-lparallel-tau}{\tilde{\calS}_\nu(\tau)}\rvert &= \sum_{\nu \le c \le j} (-1)^{c-\nu} \glslink{binomial}{\binom{c-1}{\nu-1}} \cdot c \cdot \frac{(j-1)!}{(j-c)!} \cdot \frac{n!}{(k+c)!} \\
        &= \sum_{0 \le i \le j-\nu} (-1)^i \frac{(\nu+i)!(j-1)!n!}{(\nu-1)!i!(j-\nu-i)!(k+\nu+i)!} \\
        &= \nu(j-1)! \sum_{0 \le i \le j-\nu} \glslink{binomial}{\binom{-\nu-1}{i}} \glslink{binomial}{\binom{n}{j-\nu-i}} \\
        &= \nu(j-1)! \cdot \glslink{binomial}{\binom{n-\nu-1}{j-\nu}} = \nu \cdot \frac{(j-1)!}{(j-\nu)!} \cdot \frac{(n-\nu-1)!}{(k-1)!},
    \end{align*} where the second-to-last step follows by Vandermonde's identity.
\end{proof}

We will also find the following slight variant of Theorem~\ref{left-leaning-proper-decomp-count-thm} with a similar proof to that of Corollary~\ref{left-leaning-irred-count-thm} to be useful in Section~\ref{generating-functions-for-colored-stitches-section}.

\begin{corollary} \label{left-leaning-decomp-count-thm}
    Let $\tau$ be a \glslink{cycle}{cycle}-free \glslink{stitch-type}{stitch type} with $n$ \glslink{column}{columns}, $k$ \glslink{wire}{wires}, and $j = n-k$ \glslink{chain}{chains}, let $\nu \ge 1$, and let $\glslink{calS-nu-lparallel-tau}{\calS_\nu^\lparallel(\tau)}$ denote the set of tuples $(\alpha_1,\dots,\alpha_\nu)$, where each $\alpha_i$ is a (possibly \glslink{empty-stitch}{empty}) \glslink{left-leaning-stitch}{left-leaning} \glslink{stitch}{stitch} and $\glslink{oplus}{\alpha_1 \oplus \dots \oplus \alpha_\nu}$ has \glslink{stitch-type}{type} $\tau$. Then \begin{equation*}
        \lvert\glslink{calS-nu-lparallel-tau}{\calS_\nu^\lparallel(\tau)}\rvert = \frac{1}{\glslink{z-tau}{z_\tau}} \cdot \frac{\nu(n+\nu-1)!}{(k+\nu)!}.
    \end{equation*}
\end{corollary}

\begin{proof}
    If $j = 0$, then $\tau$ must equal $(\varnothing, \varnothing)$, and the result holds trivially. Otherwise, observe that we have a bijection $\glslink{calS-nu-lparallel-tau}{\calS_\nu^\lparallel(\tau)} \simto \bigsqcup_{R\subseteq [\nu]} \glslink{hat-calS-nu-lparallel-tau}{\hat{\calS}_{\lvert R\rvert}^\lparallel(\tau)}$ sending a tuple $(\alpha_1, \dots, \alpha_\nu) \in \glslink{calS-nu-lparallel-tau}{\calS_\nu^\lparallel(\tau)}$ to the pair $(R, (\alpha_{i_1},\dots,\alpha_{i_{\lvert R\rvert}}))$, where $R$ is the set of indices $i \in [\nu]$ such that $\alpha_i$ is \glslink{empty-stitch}{nonempty} and $i_1,\dots,i_{\lvert R\rvert}$ are the elements of $R$ listed in increasing order. Consequently, we may write \begin{align*}
        \lvert\glslink{calS-nu-lparallel-tau}{\calS_\nu^\lparallel(\tau)}\rvert &= \sum_{c=1}^{\min(\nu,j)} \glslink{binomial}{\binom{\nu}{c}} \lvert\glslink{hat-calS-nu-lparallel-tau}{\hat{\calS}_c^\lparallel(\tau)}\rvert \\
        &= \frac{1}{\glslink{z-tau}{z_\tau}} \sum_{c=1}^{\min(\nu,j)} \glslink{binomial}{\binom{\nu}{c}} \cdot c \cdot \frac{(j-1)!}{(j-c)!} \cdot \frac{n!}{(k+c)!} = \frac{1}{\glslink{z-tau}{z_\tau}} \nu(j-1)! \cdot \sum_{c=1}^{\min(\nu,j)} \glslink{binomial}{\binom{\nu-1}{c-1}} \cdot \glslink{binomial}{\binom{n}{j-c}} \\
        &= \frac{1}{\glslink{z-tau}{z_\tau}} \nu(j-1)! \cdot \sum_{c=0}^{\min(\nu,j)-1} \glslink{binomial}{\binom{\nu-1}{c}} \cdot \glslink{binomial}{\binom{n}{j-1-c}} \\
        &= \frac{1}{\glslink{z-tau}{z_\tau}} \nu(j-1)! \cdot \glslink{binomial}{\binom{n+\nu-1}{j-1}} = \frac{1}{\glslink{z-tau}{z_\tau}} \cdot \frac{\nu(n+\nu-1)!}{(k+\nu)!},
    \end{align*} where the second line follows by Theorem~\ref{left-leaning-proper-decomp-count-thm} and the second-to-last step follows by Vandermonde's identity.
\end{proof}

\subsection{Counting \texorpdfstring{\glslink{parallel-stitch}{parallel}}{parallel} \texorpdfstring{\glslink{stitch}{stitches}}{stitches}}

With Corollary~\ref{left-leaning-irred-count-thm}, we are finally able to accomplish the stated goal of this section of computing $\lvert\glslink{calS-parallel-tau}{\calS^\parallel(\tau)}\rvert$ for a fixed \glslink{stitch-type}{stitch type} $\tau$.

\begin{proof}[Proof of Theorem~\ref{parallel-count-thm}]
    Suppose $\tau$ has $n$ \glslink{column}{columns} and $m$ \glslink{gap}{gaps}. Then by Lemma~\ref{reduction-lem}, we may write \begin{equation*}
        \frac{\lvert\glslink{calS-parallel-tau}{\calS^\parallel(\tau)}\rvert}{\lvert\glslink{calS-tau}{\calS(\tau)}\rvert} = \frac{\glslink{z-tau}{z_\tau}}{n!} \cdot \glslink{binomial}{\binom{n}{m}} \lvert\glslink{calS-parallel-tau}{\calS^\parallel(\glslink{bar-tau}{\bar{\tau}})}\rvert = \frac{\glslink{z-tau}{z_{\glslink{bar-tau}{\bar{\tau}}}}}{(n-m)!} \cdot \lvert\glslink{calS-parallel-tau}{\calS^\parallel(\glslink{bar-tau}{\bar{\tau}})}\rvert = \frac{\lvert\glslink{calS-parallel-tau}{\calS^\parallel(\glslink{bar-tau}{\bar{\tau}})}\rvert}{\lvert\glslink{calS-tau}{\calS(\glslink{bar-tau}{\bar{\tau}})}\rvert}.
    \end{equation*} Thus, we may assume without loss of generality that $\tau$ is \glslink{reduced}{reduced}. Moreover, if $\tau$ is \glslink{reduced}{reduced} with no \glslink{proper-chain}{proper chains}, then $\tau = (1^\nu, \varnothing)$, and there is a unique stitch of type $\tau$ consisting of $\nu$ vertical wires; thus, we have $\frac{\lvert\glslink{calS-parallel-tau}{\calS^\parallel(\tau)}\rvert}{\lvert\glslink{calS-tau}{\calS(\tau)}\rvert} = \frac{1}{1} = 1$ as claimed. We may therefore assume for the remainder of the proof that the number of \glslink{proper-chain}{proper chains} $\ell$ is at least $1$.
    
    Now, let $\alpha$ be a \glslink{parallel-stitch}{parallel} \glslink{stitch}{stitch} of \glslink{stitch-type}{type} $\tau$, and decompose $\alpha$ as $\glslink{oplus}{\alpha_1 \oplus \dots \oplus \alpha_c}$ for some \glslink{irreducible-stitch}{irreducible} \glslink{stitch}{stitches} $\alpha_1, \dots, \alpha_c$. Clearly, $\nu$ of the $\alpha_i$ must be \glslink{pole}{poles}, and by Proposition~\ref{irreducible-parallel-classification-prop}, the rest of the $\alpha_i$ are either \glslink{left-leaning-stitch}{left-leaning} or \glslink{left-leaning-stitch}{right-leaning} (and not both). Letting $\glslink{tilde-calS-d-parallel-tau}{\tilde{\calS}_d^\parallel(\tau)}$ denote the set of \glslink{parallel-stitch}{parallel} \glslink{stitch}{stitches} of \glslink{stitch-type}{type} $\tau$ with exactly $d$ non-\glslink{pole}{poles} in their \glslink{irreducible-stitch}{irreducible} decompositions and letting $\tau'$ denote the \glslink{stitch-type}{stitch type} obtained by removing all \glslink{pole}{poles} from $\tau$, we conclude that \begin{equation*}
        \lvert\glslink{tilde-calS-d-parallel-tau}{\tilde{\calS}_d^\parallel(\tau)}\rvert = 2^d \glslink{binomial}{\binom{\nu+d}{\nu}} \cdot \#\left\{(\alpha_1,\dots,\alpha_d)\,\middle|\, \begin{gathered}
            \text{$\alpha_i$ \glslink{left-leaning-stitch}{left-leaning} and \glslink{irreducible-stitch}{irreducible}} \\
            \text{$\glslink{oplus}{\alpha_1 \oplus \dots \oplus \alpha_d}$ has \glslink{stitch-type}{type} $\tau'$}
        \end{gathered}\right\} = 2^d \glslink{binomial}{\binom{\nu+d}{\nu}} \cdot \lvert\glslink{tilde-calS-nu-lparallel-tau}{\tilde{\calS}_d^\lparallel(\tau')}\rvert,
    \end{equation*} where the factor $2^d$ comes from the fact that each non-\glslink{pole}{pole} can be either \glslink{left-leaning-stitch}{left-leaning} or \glslink{left-leaning-stitch}{right-leaning}, and the factor $\glslink{binomial}{\binom{\nu+d}{\nu}}$ comes from the different ways of interlacing \glslink{pole}{pole} and non-\glslink{pole}{pole} \glslink{irreducible-stitch}{irreducible} components. Using the fact that $\tau'$ has $n-\nu$ \glslink{column}{columns}, $k-\nu$ \glslink{wire}{wires}, and $\ell$ (proper) \glslink{chain}{chains}, we may thus apply Corollary~\ref{left-leaning-irred-count-thm} to deduce that \begin{align*}
        \lvert\glslink{tilde-calS-d-parallel-tau}{\tilde{\calS}_d^\parallel(\tau)}\rvert &= 2^d \glslink{binomial}{\binom{\nu+d}{\nu}} \cdot \frac{1}{\glslink{z-tau}{z_{\tau'}}} \cdot d \cdot \frac{(\ell-1)!}{(\ell-d)!} \cdot \frac{(n-\nu-d-1)!}{(k-\nu-1)!} \\
        &= \frac{1}{\glslink{z-tau}{z_\tau}} \cdot 2^d \glslink{binomial}{\binom{\ell-1}{d-1}} \frac{(\nu+d)!(n-\nu-d-1)!}{(k-\nu-1)!}.
    \end{align*} Summing over the possible values of $d$, we conclude that \begin{align*}
        \lvert\glslink{calS-parallel-tau}{\calS^\parallel(\tau)}\rvert &= \frac{1}{\glslink{z-tau}{z_\tau}} \cdot \sum_{1 \le d \le \ell} 2^d \glslink{binomial}{\binom{\ell-1}{d-1}} \frac{(\nu+d)!(n-\nu-d-1)!}{(k-\nu-1)!} \\
        &= \frac{1}{\glslink{z-tau}{z_\tau}} \cdot \sum_{0 \le i \le \ell-1} 2^{\ell-i} \glslink{binomial}{\binom{\ell-1}{i}} \frac{(\nu+\ell-i)!(n-\nu-\ell+i-1)!}{(k-\nu-1)!} \\
        &= \frac{2^\ell(\nu+\ell)!}{\glslink{z-tau}{z_\tau}} \cdot \sum_{0 \le i \le \ell-1} \frac{(\ell-1)!}{(\ell-1-i)!} \cdot \frac{(k-\nu-1+i)!}{(k-\nu-1)!} \cdot \frac{(\nu+\ell-i)!}{(\nu+\ell)!} \cdot \frac{2^{-i}}{i!} \\
        &= \frac{2^\ell(\nu+\ell)!}{\glslink{z-tau}{z_\tau}} \cdot \glslink{2F1}{{}_2F_1\left(1-\ell, k-\nu; -\nu-\ell, \tfrac{1}{2}\right)},
    \end{align*} where $\glslink{2F1}{{}_2F_1(a,b;c;z)}$ denotes Gauss's hypergeometric function. (Here, $1-\ell$ is a nonpositive integer, so the hypergeometric series terminates before any vanishing denominator arising from the parameter $-\nu-\ell$.) From here, applying the Pfaff transformation identity \begin{equation*}
        \glslink{2F1}{{}_2F_1\left(a,b;c;\tfrac{1}{2}\right)} = 2^a \glslink{2F1}{{}_2F_1(a,c-b;c;-1)}
    \end{equation*} (the first equation of \cite[\href{http://dlmf.nist.gov/15.8.E1}{(15.8.E1)}]{dlmf} applied with $z = \frac{1}{2}$) gives \begin{align*}
        \glslink{z-tau}{z_\tau}\lvert\glslink{calS-parallel-tau}{\calS^\parallel(\tau)}\rvert &= 2^\ell(\nu+\ell)! \cdot 2^{1-\ell}\glslink{2F1}{{}_2F_1(1-\ell,-n;-\nu-\ell;-1)} \\
        &= 2(\nu+\ell)! \sum_{i=0}^{\ell-1} \glslink{binomial}{\binom{\ell-1}{i}} \cdot \frac{n!}{(n-i)!} \cdot \frac{(\nu+\ell-i)!}{(\nu+\ell)!} \\
        &= 2n!\sum_{i=0}^{\ell-1} \glslink{binomial}{\binom{\ell-1}{i}} \frac{(\nu+\ell-i)!}{(k+\ell-i)!} \\
        &= 2n!\sum_{i=0}^{\ell-1} \glslink{binomial}{\binom{\ell-1}{i}} \frac{(\nu+1+i)!}{(k+1+i)!},
    \end{align*} which is equivalent to the formula in the theorem statement.
\end{proof}

\subsection{A simplified formula for \texorpdfstring{$\glslink{calE-f}{\calE_f}$}{Ef}}

While the explicit formula for $\lvert\glslink{calS-parallel-tau}{\calS^\parallel(\tau)}\rvert$ given in Theorem~\ref{parallel-count-thm} does not yet give a clear path to computing $\glslink{h-tau-j-k}{h_\tau(j, k)}$ or the generating functions $\glslink{H-tau}{H_\tau}$ from Lemma~\ref{E-from-H-tau-lem}, it brings us one step closer with the following lemma.

\begin{lemma} \label{h-comparison-lem}
    Let $\tau_1$ and $\tau_2$ be two \glslink{parallelizable-stitch-type}{parallelizable} \glslink{stitch-type}{stitch types} with the same numbers of \glslink{column}{columns}, \glslink{wire}{wires}, \glslink{gap}{gaps}, and \glslink{pole}{poles}. Then for all $j, k \ge 0$, we have \begin{equation*}
        \glslink{z-tau}{z_{\tau_1}}\glslink{h-tau-j-k}{h_{\tau_1}(j, k)} = \glslink{z-tau}{z_{\tau_2}}\glslink{h-tau-j-k}{h_{\tau_2}(j, k)} \qquad \text{and} \qquad \glslink{z-tau}{z_{\tau_1}} \glslink{H-tau}{H_{\tau_1}} = \glslink{z-tau}{z_{\tau_2}} \glslink{H-tau}{H_{\tau_2}}.
    \end{equation*}
\end{lemma}

\begin{proof}
    Let $\tau$ be another \glslink{parallelizable-stitch-type}{parallelizable} \glslink{stitch-type}{stitch type}, and observe that the number of \glslink{stitch-coloring}{$2$-colored} \glslink{parallel-stitch}{parallel} \glslink{stitch}{stitches} $\alpha$ with $\glslink{alpha-rvert-i}{\alpha\rvert_1} \in \glslink{calS-parallel-tau}{\calS^\parallel(\tau_i)}$ and $\glslink{alpha-rvert-i}{\alpha\rvert_2} \in \glslink{calS-parallel-tau}{\calS^\parallel(\tau)}$ is given by $\frac{\glslink{z-tau}{z_{\glslink{oplus}{\tau_i\oplus\tau}}}}{\glslink{z-tau}{z_{\tau_i}}\glslink{z-tau}{z_\tau}} \cdot \lvert\glslink{calS-parallel-tau}{\calS^\parallel(\glslink{oplus}{\tau_i\oplus\tau})}\rvert$: we can construct such a stitch by first choosing any \glslink{parallel-stitch}{parallel} \glslink{stitch}{stitch} of \glslink{stitch-type}{type} $\glslink{oplus}{\tau_i \oplus \tau}$, then choosing which of the \glslink{pole}{poles} and \glslink{chain}{chains} to assign each color. Consequently, we may write \begin{equation*}
        \glslink{z-tau}{z_{\tau_i}}\glslink{h-tau-j-k}{h_{\tau_i}(j,k)} = \sum_{\tau \in \glslink{calT-n-k-parallel}{\calT_{j+k,k}^\parallel}} \frac{\glslink{z-tau}{z_{\glslink{oplus}{\tau_i\oplus\tau}}}}{\glslink{z-tau}{z_\tau}} \cdot \lvert\glslink{calS-parallel-tau}{\calS^\parallel(\glslink{oplus}{\tau_i\oplus\tau})}\rvert = \sum_{\tau \in \glslink{calT-n-k-parallel}{\calT_{j+k,k}^\parallel}} \frac{(n_0+j+k)!}{\glslink{z-tau}{z_\tau}} \cdot \frac{\lvert\glslink{calS-parallel-tau}{\calS^\parallel(\glslink{oplus}{\tau_i\oplus\tau})}\rvert}{\lvert\glslink{calS-tau}{\calS(\glslink{oplus}{\tau_i \oplus \tau})}\rvert},
    \end{equation*} where $n_0$ is the number of \glslink{column}{columns} of $\tau_1$ (and hence of $\tau_2$). By Theorem~\ref{parallel-count-thm}, for each $\tau \in \glslink{calT-n-k-parallel}{\calT_{j+k,k}^\parallel}$, the ratios $\frac{\lvert\glslink{calS-parallel-tau}{\calS^\parallel(\glslink{oplus}{\tau_1\oplus\tau})}\rvert}{\lvert\glslink{calS-tau}{\calS(\glslink{oplus}{\tau_1 \oplus \tau})}\rvert}$ and $\frac{\lvert\glslink{calS-parallel-tau}{\calS^\parallel(\glslink{oplus}{\tau_2\oplus\tau})}\rvert}{\lvert\glslink{calS-tau}{\calS(\glslink{oplus}{\tau_2 \oplus \tau})}\rvert}$ agree, as $\glslink{oplus}{\tau_1\oplus\tau}$ and $\glslink{oplus}{\tau_2\oplus\tau}$ have the same numbers of \glslink{wire}{wires}, \glslink{pole}{poles}, and \glslink{proper-chain}{proper chains}. Thus, $\glslink{z-tau}{z_{\tau_1}}\glslink{h-tau-j-k}{h_{\tau_1}(j,k)} = \glslink{z-tau}{z_{\tau_2}}\glslink{h-tau-j-k}{h_{\tau_2}(j,k)}$ as claimed. The second equation then follows immediately from the definition of the $\glslink{H-tau}{H_\tau}$.
\end{proof}

As we alluded to in Section~\ref{proof-overview-section}, Lemma~\ref{h-comparison-lem} is useful because it allows us to eliminate the need to sum over \glslink{stitch-type}{stitch types} in the formula from Lemma~\ref{E-from-H-tau-lem}. With the help of the following lemma, we thus obtain the simplified formula for the functions $\glslink{calE-f}{\calE_f}$ given in Corollary~\ref{E-from-H-ell-m-t-nu-cor}.

\begin{lemma} \label{H-sum-lem}
    Let $\ell,m,t,\nu \ge 0$ with either $\ell \ge 1$ or $t = 0$, and let $\glslink{calT-ell-m-t-nu-parallel}{\calT_{\ell,m,t,\nu}^\parallel}$ denote the set of \glslink{parallelizable-stitch-type}{parallelizable} \glslink{stitch-type}{stitch types} with $\ell$ \glslink{proper-chain}{proper chains}, $m$ \glslink{gap}{gaps}, $t$ \glslink{jump}{jumps}, and $\nu$ \glslink{pole}{poles}. Then \begin{equation*}
        \sum_{\tau \in \glslink{calT-ell-m-t-nu-parallel}{\calT_{\ell,m,t,\nu}^\parallel}} \frac{1}{\glslink{z-tau}{z_\tau}} = \frac{1}{\ell!m!\nu!} \glslink{binomial}{\binom{\ell+t-1}{t}}
    \end{equation*}
\end{lemma}

\begin{proof}
    Write \begin{align*}
        \sum_{\tau \in \glslink{calT-ell-m-t-nu-parallel}{\calT_{\ell,m,t,\nu}^\parallel}} \frac{1}{\glslink{z-tau}{z_\tau}} &= \frac{1}{\nu!} \sum_{\substack{\mu \vdash 2\ell+m+t \\ \text{with $\ell$ parts $\ge2$} \\ \text{and $m$ parts $=1$}}} \frac{1}{\glslink{z-tau}{z_{(\varnothing,\mu)}}} = \frac{1}{m!\nu!} \sum_{\substack{\mu \vdash 2\ell+t \\ \text{with $\ell$ parts $\ge2$} \\ \text{and $0$ parts $=1$}}} \frac{1}{\glslink{z-tau}{z_{(\varnothing,\mu)}}} = \frac{1}{m!\nu!} \sum_{\substack{\mu \vdash \ell+t \\ \text{with $\ell$ parts}}} \frac{1}{\glslink{z-tau}{z_{(\varnothing,\mu)}}} \\
        &= \frac{1}{\ell!m!\nu!} \sum_{\substack{\mu \vdash \ell+t \\ \text{with $\ell$ parts}}} \frac{\ell!}{\glslink{z-tau}{z_{(\varnothing,\mu)}}}.
    \end{align*} For each partition $\mu \vdash \ell+t$ with $\ell$ parts, $\frac{\ell!}{\glslink{z-tau}{z_{(\varnothing,\mu)}}}$ counts the number of integer compositions of $\ell+t$ with $\ell$ parts giving rise to the partition $\mu$. Thus, the sum in the above expression is equal to the total number of integer compositions of $\ell+t$ with $\ell$ parts, which is just $\glslink{binomial}{\binom{\ell+t-1}{t}}$.
\end{proof}

\begin{proof}[Proof of Corollaries~\ref{H-tau-from-H-ell-m-t-nu-cor} and \ref{E-from-H-ell-m-t-nu-cor}]
    Let $\ell,m,t,\nu \ge 0$ with either $\ell \ge 1$ or $t = 0$, and let $\glslink{calT-ell-m-t-nu-parallel}{\calT_{\ell,m,t,\nu}^\parallel}$ denote the set of \glslink{parallelizable-stitch-type}{parallelizable} \glslink{stitch-type}{stitch types} with $\ell$ \glslink{proper-chain}{proper chains}, $m$ \glslink{gap}{gaps}, $t$ \glslink{jump}{jumps}, and $\nu$ \glslink{pole}{poles}. Then by Lemma~\ref{h-comparison-lem}, the functions $\glslink{z-tau}{z_\tau} \glslink{H-tau}{H_\tau}$ for $\tau \in \glslink{calT-ell-m-t-nu-parallel}{\calT_{\ell,m,t,\nu}^\parallel}$ are all given by some shared power series $\glslink{tilde-H-ell-m-t-nu}{\tilde{H}_{\ell,m,t,\nu}}$. Moreover, by the definitions of $\glslink{h-ell-m-t-nu-j-k}{h_{\ell,m,t,\nu}(j,k)}$ and $\glslink{h-tau-j-k}{h_\tau(j,k)}$, we have \begin{equation*}
        \glslink{H-ell-m-t-nu}{H_{\ell,m,t,\nu}} = \sum_{\tau \in \glslink{calT-ell-m-t-nu-parallel}{\calT_{\ell,m,t,\nu}^\parallel}} \glslink{H-tau}{H_\tau} = \left(\sum_{\tau \in \glslink{calT-ell-m-t-nu-parallel}{\calT_{\ell,m,t,\nu}^\parallel}} \frac{1}{\glslink{z-tau}{z_\tau}}\right)\glslink{tilde-H-ell-m-t-nu}{\tilde{H}_{\ell,m,t,\nu}} = \frac{1}{\ell!m!\nu!} \glslink{binomial}{\binom{\ell+t-1}{t}} \cdot \glslink{tilde-H-ell-m-t-nu}{\tilde{H}_{\ell,m,t,\nu}},
    \end{equation*} where the last step follows by Lemma~\ref{H-sum-lem}. Consequently, we see that \begin{align*}
        \sum_{\tau \in \glslink{calT-ell-m-t-nu-parallel}{\calT_{\ell,m,t,\nu}^\parallel}} \glslink{stitch-binom}{\binom{\rho}{\tau}} \cdot \glslink{z-tau}{z_\tau} \glslink{H-tau}{H_\tau} &= \sum_{\tau \in \glslink{calT-ell-m-t-nu-parallel}{\calT_{\ell,m,t,\nu}^\parallel}} \glslink{stitch-binom}{\binom{\rho}{\tau}} \cdot \glslink{tilde-H-ell-m-t-nu}{\tilde{H}_{\ell,m,t,\nu}} \\
        &= \left(\sum_{\tau \in \glslink{calT-ell-m-t-nu-parallel}{\calT_{\ell,m,t,\nu}^\parallel}} \glslink{stitch-binom}{\binom{\rho}{\tau}}\right)\frac{\ell!m!\nu!}{\glslink{binomial}{\binom{\ell+t-1}{t}}} \cdot \glslink{H-ell-m-t-nu}{H_{\ell,m,t,\nu}} \\
        &= \glslink{stitch-ell-m-t-nu-binom}{\binom{\rho}{\ell, m, t, \nu}} \frac{\ell!m!\nu!}{\glslink{binomial}{\binom{\ell+t-1}{t}}} \cdot \glslink{H-ell-m-t-nu}{H_{\ell,m,t,\nu}}.
    \end{align*} The statement of the corollary then follows immediately by rewriting the formula in Lemma~\ref{E-from-H-tau-lem} as a sum over $\ell$, $m$, $t$, and $\nu$, then applying the above identity.
\end{proof}

\section{Generating functions for colored \texorpdfstring{\glslink{parallel-stitch}{parallel}}{parallel} \texorpdfstring{\glslink{stitch}{stitches}}{stitches}} \label{generating-functions-for-colored-stitches-section}

We now turn our attention to generating functions for certain families of \glslink{stitch-coloring}{colored} \glslink{parallel-stitch}{parallel} \glslink{stitch}{stitches}. As we alluded to in Section~\ref{proof-overview-section}, we will ultimately be able to specialize our results to the case of \glslink{stitch-coloring}{$2$-colored} \glslink{stitch}{stitches} in order to obtain formulas for the generating functions $\glslink{H-ell-m-t-nu}{H_{\ell,m,t,\nu}}$. However, we choose to state and prove our main results pertaining to these generating functions for an arbitrary number of colors $r$, as they may be of independent interest for use in other applications.

\subsection{Algebraic equations for generating functions of \texorpdfstring{\glslink{stitch-coloring}{colored}}{colored} \texorpdfstring{\glslink{left-leaning-stitch}{left-leaning}}{left-leaning} \texorpdfstring{\glslink{stitch}{stitches}}{stitches}}

We will begin by studying the following collection of quantities and associated generating function.

\begin{definition}
    For each $j_1,\dots,j_r,k_1,\dots,k_r\ge0$, we let $\glslink{f-lparallel}{f_\lparallel(j_1,k_1;\dots;j_r,k_r)}$ denote the number of \glslink{stitch-coloring}{$r$-colored} \glslink{left-leaning-stitch}{left-leaning} \glslink{stitch}{stitches} $\alpha$ such that for each $1 \le i \le r$, the $i$th restriction $\glslink{alpha-rvert-i}{\alpha\rvert_i}$ has exactly $j_i$ \glslink{chain}{chains} and $k_i$ \glslink{wire}{wires}. In addition, we define a generating function \begin{equation*}
        \glslink{F-lparallel}{F_\lparallel} := \sum_{\substack{j_1,\dots,j_r \ge 0 \\ k_1,\dots,k_r \ge 0}} \glslink{f-lparallel}{f_\lparallel(j_1,k_1;\dots;j_r,k_r)} \cdot x_1^{j_1} y_1^{k_1} \cdots x_r^{j_r} y_r^{k_r} \in \Z[[x_1,y_1; \dots; x_r,y_r]].
    \end{equation*}
\end{definition}

A priori, it is not at all clear that the quantities $\glslink{f-lparallel}{f_\lparallel(j_1,k_1;\dots;j_r,k_r)}$ should admit a simple closed formula. However, using the results of Section~\ref{counting-parallel-stitches-section}, it turns out that we can derive explicit formulas not only for $\glslink{f-lparallel}{f_\lparallel(j_1,k_1;\dots;j_r,k_r)}$, but also for the series coefficients of $\glslink{F-lparallel}{F_\lparallel^\nu}$ for arbitrary $\nu \ge 1$.

\begin{lemma} \label{colored-left-leaning-decomp-count-lem}
    Let $j_1,\dots,j_r,k_1,\dots,k_r \ge 0$ and $\nu \ge 1$, and let $\glslink{f-lparallel-nu}{f_\lparallel^{(\nu)}(j_1,k_1;\dots;j_r,k_r)}$ denote the number of tuples of (possibly \glslink{empty-stitch}{empty}) \glslink{stitch-coloring}{$r$-colored} \glslink{left-leaning-stitch}{left-leaning} \glslink{stitch}{stitches} $(\alpha_1,\dots,\alpha_\nu)$ such that for each $1 \le i \le r$, the $i$th restriction $\glslink{alpha-rvert-i}{(\glslink{oplus}{\alpha_1 \oplus \dots \oplus \alpha_\nu})\rvert_i}$ has exactly $j_i$ \glslink{chain}{chains} and $k_i$ \glslink{wire}{wires}. Then \begin{equation*}
        \glslink{f-lparallel-nu}{f_\lparallel^{(\nu)}(j_1,k_1;\dots;j_r,k_r)} = \frac{\nu(j+k+\nu-1)!}{(k+\nu)!} \prod_{i \in \glslink{bracket-n}{[r]}} \frac{1}{j_i!} \glslink{binomial}{\binom{j_i+k_i-1}{k_i}},
    \end{equation*} where $j := j_1+\dots+j_r$ and $k := k_1+\dots+k_r$.
\end{lemma}

\begin{proof}
    To start, let $T$ denote the set of \glslink{stitch-coloring}{$r$-colored} \glslink{cycle}{cycle}-free (but not necessarily \glslink{parallel-stitch}{parallel}) \glslink{stitch}{stitches} $\alpha$ such that for each $1 \le i \le r$, the $i$th restriction $\glslink{alpha-rvert-i}{\alpha\rvert_i}$ has $j_i$ \glslink{chain}{chains} and $k_i$ \glslink{wire}{wires}. In addition, let $\tilde{T}$ denote the set of pairs $\tilde{\alpha} = (\alpha, (c_1,\dots,c_j))$, where $\alpha \in T$ and $c_1,\dots,c_j$ is an ordering of the $j$ \glslink{chain}{chains} of $\alpha$ such that the first $j_1$ \glslink{chain}{chains} have color $1$, the next $j_2$ \glslink{chain}{chains} have color $2$, and so on; it is clear that the fibers of the map $\tilde{T} \to T$ sending $(\alpha, (c_1,\dots,c_j))$ to $\alpha$ all have size $j_1! \cdots j_r!$, so $\lvert T\rvert = \frac{1}{j_1! \cdots j_r!} \lvert\tilde{T}\rvert$. Now, observe that given an element $\tilde{\alpha} = (\alpha, (c_1,\dots,c_j)) \in \tilde{T}$, we can obtain a permutation $\pi(\tilde{\alpha}) \in \glslink{symmetric-group}{S_{j+k}}$ by listing out the vertices of the \glslink{chain}{chain} $c_1$ in order, followed by the vertices of the \glslink{chain}{chain} $c_2$, and so on. In addition, for each $i \in \glslink{bracket-n}{[r]}$, we can obtain a composition $K_i(\tilde{\alpha})$ of $j_i+k_i$ as a sum of $j_i$ positive integers by writing \begin{equation*}
        j_i+k_i = \#\{\text{vertices of $c_{j_1+\dots+j_{i-1}+1}$}\} + \dots + \#\{\text{vertices of $c_{j_1+\dots+j_{i-1}+j_i}$}\},
    \end{equation*} since the total number of vertices of color $i$ is given by $j_i+k_i$. Moreover, it is not hard to see that, given any permutation $\pi \in \glslink{symmetric-group}{S_{j+k}}$ and compositions $K_i$ of $j_i+k_i$ as sums of $j_i$ positive integers for each $i \in \glslink{bracket-n}{[r]}$, one can recover a unique $\tilde{\alpha} \in \tilde{T}$ satisfying $\pi(\tilde{\alpha}) = \pi$ and $K_i(\tilde{\alpha}) = K_i$ for all $i$. We conclude that $\lvert\tilde{T}\rvert = (j+k)! \prod_{i \in \glslink{bracket-n}{[r]}} \glslink{binomial}{\binom{j_i+k_i-1}{k_i}}$, hence that \begin{equation*}
        \lvert T\rvert = (j+k)! \prod_{i \in \glslink{bracket-n}{[r]}} \frac{1}{j_i!} \glslink{binomial}{\binom{j_i+k_i-1}{k_i}}.
    \end{equation*}

    Returning now to the computation of $\glslink{f-lparallel-nu}{f_\lparallel^{(\nu)}(j_1,k_1;\dots;j_r,k_r)}$, observe that \begin{equation*}
        \frac{\glslink{f-lparallel-nu}{f_\lparallel^{(\nu)}(j_1,k_1;\dots;j_r,k_r)}}{\lvert T\rvert}
    \end{equation*} represents the expected number of ways that a uniformly random element of $T$ can be expressed as $\glslink{oplus}{\alpha_1 \oplus \dots \oplus \alpha_\nu}$ for (possibly \glslink{empty-stitch}{empty}) \glslink{stitch-coloring}{$r$-colored} \glslink{left-leaning-stitch}{left-leaning} \glslink{stitch}{stitches} $\alpha_1,\dots,\alpha_\nu$. Now, by the definition of $T$, each \glslink{stitch}{stitch} $\alpha \in T$ has exactly $j+k$ \glslink{column}{columns} and $k$ \glslink{wire}{wires}. Moreover, for each \glslink{stitch-type}{stitch type} $\tau \in \glslink{calT-n-k-lparallel}{\calT_{j+k,k}^\lparallel}$, Corollary~\ref{left-leaning-decomp-count-thm} implies that the expected number of such decompositions $\glslink{oplus}{\alpha_1 \oplus \dots \oplus \alpha_\nu}$ for a uniformly random element of $T \cap \glslink{calS-tau}{\calS(\tau)}$ is \begin{equation*}
        \frac{1}{(j+k)!} \cdot \frac{\nu(j+k+\nu-1)!}{(k+\nu)!},
    \end{equation*} independent of $\tau$. It thus follows that \begin{equation*}
        \frac{\glslink{f-lparallel-nu}{f_\lparallel^{(\nu)}(j_1,k_1;\dots;j_r,k_r)}}{\lvert T\rvert} = \frac{1}{(j+k)!} \cdot \frac{\nu(j+k+\nu-1)!}{(k+\nu)!},
    \end{equation*} hence that \begin{align*}
        \glslink{f-lparallel-nu}{f_\lparallel^{(\nu)}(j_1,k_1;\dots;j_r,k_r)} &= \frac{1}{(j+k)!} \cdot \frac{\nu(j+k+\nu-1)!}{(k+\nu)!} \cdot (j+k)! \prod_{i \in \glslink{bracket-n}{[r]}} \frac{1}{j_i!} \glslink{binomial}{\binom{j_i+k_i-1}{k_i}} \\
        &= \frac{\nu(j+k+\nu-1)!}{(k+\nu)!} \prod_{i \in \glslink{bracket-n}{[r]}} \frac{1}{j_i!} \glslink{binomial}{\binom{j_i+k_i-1}{k_i}}.
    \end{align*}
\end{proof}

One might hope that the generating function $\glslink{F-lparallel}{F_\lparallel}$ itself also has a simple closed form. Unfortunately, it turns out that this is not the case whenever $r \ge 2$ (for an appropriate definition of ``simple''). However, using Lemma~\ref{colored-left-leaning-decomp-count-lem}, we can nonetheless obtain a relatively simple algebraic characterization of $\glslink{F-lparallel}{F_\lparallel}$.

\begin{lemma} \label{F-lparallel-algebraicity-lem}
    The generating function $\glslink{F-lparallel}{F_\lparallel}$ satisfies the degree $r+1$ algebraic equation \begin{equation*}
        \glslink{F-lparallel}{F_\lparallel} = 1 + \sum_{i \in \glslink{bracket-n}{[r]}} \frac{x_i\glslink{F-lparallel}{F_\lparallel}}{1-y_i\glslink{F-lparallel}{F_\lparallel}}.
    \end{equation*}
\end{lemma}

\begin{proof}
    Observe that by the definition of the functions $\glslink{f-lparallel-nu}{f_\lparallel^{(\nu)}(j_1,k_1;\dots;j_r,k_r)}$ from Lemma~\ref{colored-left-leaning-decomp-count-lem}, we have \begin{equation*}
        \glslink{F-lparallel}{F_\lparallel^\nu} = \sum_{\substack{j_1,\dots,j_r \ge 0 \\ k_1,\dots,k_r \ge 0}} \glslink{f-lparallel-nu}{f_\lparallel^{(\nu)}(j_1,k_1;\dots;j_r,k_r)} \cdot x_1^{j_1} y_1^{k_1} \cdots x_r^{j_r} y_r^{k_r}
    \end{equation*} for all $\nu \ge 1$. Consequently, for each $i \in \glslink{bracket-n}{[r]}$, we may write \begin{align}
        &\frac{x_i \glslink{F-lparallel}{F_\lparallel}}{1-y_i \glslink{F-lparallel}{F_\lparallel}} = \sum_{\nu \ge 0} x_i y_i^\nu \glslink{F-lparallel}{F_\lparallel^{\nu+1}} \notag \\
        &= \sum_{\substack{j_1,\dots,j_r \ge 0 \\ k_1,\dots,k_r \ge 0 \\ j_i \ge 1}} x_1^{j_1} y_1^{k_1} \cdots x_r^{j_r} y_r^{k_r} \sum_{\nu=0}^{k_i} \frac{(\nu+1)(j+k-1)!}{(k+1)!(j_i-1)!} \glslink{binomial}{\binom{j_i+k_i-\nu-2}{k_i-\nu}} \prod_{i' \in \glslink{bracket-n}{[r]} \setminus \{i\}} \frac{1}{j_{i'}!} \glslink{binomial}{\binom{j_{i'}+k_{i'}-1}{k_{i'}}} \label{F-lparallel-fraction-initial-rhs}
    \end{align} where $j$ and $k$ are shorthand for $j_1+\dots+j_r$ and $k_1+\dots+k_r$, respectively, and the second line follows by Lemma~\ref{colored-left-leaning-decomp-count-lem}. Using standard techniques, it is not difficult to verify that \begin{align*}
        \sum_{\nu=0}^{k_i} \frac{(\nu+1)(j+k-1)!}{(k+1)!(j_i-1)!} \glslink{binomial}{\binom{j_i+k_i-\nu-2}{k_i-\nu}} &= \frac{(j+k-1)!(j_i+k_i)!}{(k+1)!j_i!(j_i-1)!k_i!} \\
        &= (j_i+k_i) \cdot \frac{(j+k-1)!}{(k+1)!} \cdot \frac{1}{j_i!} \glslink{binomial}{\binom{j_i+k_i-1}{k_i}},
    \end{align*} and so equation \eqref{F-lparallel-fraction-initial-rhs} reduces to \begin{equation}
        \sum_{\substack{j_1,\dots,j_r \ge 0 \\ k_1,\dots,k_r \ge 0 \\ j_i \ge 1}} \left((j_i+k_i) \cdot \frac{(j+k-1)!}{(k+1)!} \prod_{i' \in \glslink{bracket-n}{[r]}} \frac{1}{j_{i'}!} \glslink{binomial}{\binom{j_{i'}+k_{i'}-1}{k_{i'}}}\right) x_1^{j_1} y_1^{k_1} \cdots x_r^{j_r} y_r^{k_r}. \label{F-lparallel-fraction-simplified-rhs}
    \end{equation} Moreover, when $j_i = 0$ and $k_i \ge 1$, the expression in parentheses in equation~\eqref{F-lparallel-fraction-simplified-rhs} evaluates to zero because $\frac{1}{j_i!} \glslink{binomial}{\binom{j_i+k_i-1}{k_i}} = \glslink{binomial}{\binom{k_i-1}{k_i}} = 0$, and when $j_i = k_i = 0$ and $j+k \ge 1$, the expression evaluates to zero because $j_i+k_i = 0$. Consequently, we may re-express \eqref{F-lparallel-fraction-simplified-rhs} as \begin{equation*}
        \sum_{\substack{j_1,\dots,j_r \ge 0 \\ k_1,\dots,k_r \ge 0 \\ j+k \ge 1}} \left((j_i+k_i) \cdot \frac{(j+k-1)!}{(k+1)!} \prod_{i' \in \glslink{bracket-n}{[r]}} \frac{1}{j_{i'}!} \glslink{binomial}{\binom{j_{i'}+k_{i'}-1}{k_{i'}}}\right) x_1^{j_1} y_1^{k_1} \cdots x_r^{j_r} y_r^{k_r}.
    \end{equation*} Summing the contributions from each expression $\frac{x_i\glslink{F-lparallel}{F_\lparallel}}{1-y_i\glslink{F-lparallel}{F_\lparallel}}$ as $i$ varies, we thus obtain \begin{align*}
        \sum_{i \in \glslink{bracket-n}{[r]}} \frac{x_i\glslink{F-lparallel}{F_\lparallel}}{1-y_i\glslink{F-lparallel}{F_\lparallel}} &= \sum_{\substack{j_1,\dots,j_r \ge 0 \\ k_1,\dots,k_r \ge 0 \\ j+k \ge 1}} \left((j+k) \cdot \frac{(j+k-1)!}{(k+1)!} \prod_{i' \in \glslink{bracket-n}{[r]}} \frac{1}{j_{i'}!} \glslink{binomial}{\binom{j_{i'}+k_{i'}-1}{k_{i'}}}\right) x_1^{j_1} y_1^{k_1} \cdots x_r^{j_r} y_r^{k_r} \\
        &= \sum_{\substack{j_1,\dots,j_r \ge 0 \\ k_1,\dots,k_r \ge 0 \\ j+k \ge 1}} \left(\frac{(j+k)!}{(k+1)!} \prod_{i' \in \glslink{bracket-n}{[r]}} \frac{1}{j_{i'}!} \glslink{binomial}{\binom{j_{i'}+k_{i'}-1}{k_{i'}}}\right) x_1^{j_1} y_1^{k_1} \cdots x_r^{j_r} y_r^{k_r} \\
        &= \glslink{F-lparallel}{F_\lparallel} - 1,
    \end{align*} which is equivalent to the claimed formula.
\end{proof}

\subsection{Generating functions of \texorpdfstring{\glslink{stitch-coloring}{colored}}{colored} \texorpdfstring{\glslink{integral}{integral}}{integral} \texorpdfstring{\glslink{left-leaning-stitch}{left-leaning}}{left-leaning} \texorpdfstring{\glslink{stitch}{stitches}}{stitches}}

Although the generating function $\glslink{F-lparallel}{F_\lparallel}$ proved particularly accessible via the results of Section~\ref{counting-parallel-stitches-section}, most of the other generating functions we might be interested in are more naturally expressed using the power series $\glslink{F-tilde-lparallel}{\tilde{F}_\lparallel}$ defined in Theorem~\ref{F-tilde-lparallel-algebraicity-thm}. Intuitively, this is because more complicated types of \glslink{parallel-stitch}{parallel} \glslink{stitch}{stitches} can be expressed in terms of \glslink{integral}{integral} \glslink{stitch}{stitches} using a combination of the \glslink{reduction}{reduction} operation and Proposition~\ref{stitch-unique-factorization-prop}. We will thus find the following lemma relating the two functions $\glslink{F-lparallel}{F_\lparallel}$ and $\glslink{F-tilde-lparallel}{\tilde{F}_\lparallel}$ to be useful.

\begin{lemma} \label{F-parallel-from-F-tilde-parallel-lem}
    Let $\xi := \sum_{i\in\glslink{bracket-n}{[r]}} x_i$. Then the power series $\glslink{F-tilde-lparallel}{\tilde{F}_\lparallel}$ defined in Theorem~\ref{F-tilde-lparallel-algebraicity-thm} is related to the power series $\glslink{F-lparallel}{F_\lparallel}$ from Lemma~\ref{F-lparallel-algebraicity-lem} by the formula \begin{equation*}
        \glslink{F-lparallel}{F_\lparallel} = \frac{1}{(1-\xi)\left(1-\glslink{F-tilde-lparallel}{\tilde{F}_\lparallel}\left(\frac{x_1y_1}{(1-\xi)^2}, \frac{y_1}{1-\xi}; \dots; \frac{x_ry_r}{(1-\xi)^2}, \frac{y_r}{1-\xi}\right)\right)}.
    \end{equation*}
\end{lemma}

\begin{proof}
    Intuitively, the key idea is that any \glslink{left-leaning-stitch}{left-leaning} \glslink{stitch}{stitch} can be uniquely reconstructed by first taking a direct sum of \glslink{integral}{integral} \glslink{left-leaning-stitch}{left-leaning} \glslink{stitch}{stitches}, then adding in some \glslink{gap}{gaps} between the \glslink{column}{columns}. To make this intuition concrete, observe first that since the number of \glslink{wire}{wires} in a \glslink{cycle}{cycle}-free \glslink{reduced}{reduced} \glslink{stitch}{stitch} is equal to the number of \glslink{jump}{jumps} plus the number of \glslink{chain}{chains}, the coefficient of $x_1^{j_1} y_1^{k_1} \cdots x_r^{j_r} y_r^{k_r}$ in the function \begin{equation}
        G_\lparallel := \glslink{F-tilde-lparallel}{\tilde{F}_\lparallel}(x_1y_1,y_1;\dots;x_ry_r,y_r) \label{G-from-F-tilde}
    \end{equation} counts the number of \glslink{stitch-coloring}{$r$-colored} \glslink{integral}{integral} \glslink{left-leaning-stitch}{left-leaning} \glslink{stitch}{stitches} $\alpha$ such that the $i$th restriction has $j_i$ \glslink{chain}{chains} and $k_i$ \glslink{wire}{wires}. Next, note that since a \glslink{reduced}{reduced} \glslink{left-leaning-stitch}{left-leaning} \glslink{stitch}{stitch} decomposes uniquely as a direct sum of \glslink{integral}{integral} \glslink{left-leaning-stitch}{left-leaning} \glslink{stitch}{stitches}, the coefficient of $x_1^{j_1} y_1^{k_1} \cdots x_r^{j_r} y_r^{k_r}$ in the function \begin{equation}
        \hat{G}_\lparallel := \frac{1}{1-G_\lparallel} = \sum_{\nu\ge0} G_\lparallel^\nu \label{G-prime-from-G}
    \end{equation} counts the number of \glslink{stitch-coloring}{$r$-colored} \glslink{reduced}{reduced} \glslink{left-leaning-stitch}{left-leaning} \glslink{stitch}{stitches} $\alpha$ such that the $i$th restriction has $j_i$ \glslink{chain}{chains} and $k_i$ \glslink{wire}{wires}. Finally, observe that a general \glslink{stitch-coloring}{$r$-colored} \glslink{left-leaning-stitch}{left-leaning} \glslink{stitch}{stitch} $\alpha$ is determined by its \glslink{reduction}{reduction} $\glslink{bar-alpha}{\bar{\alpha}}$ together with its set of \glslink{gap}{gaps} and their colorings. Conversely, if we start from a \glslink{reduced}{reduced} \glslink{left-leaning-stitch}{left-leaning} \glslink{stitch}{stitch} $\beta$, we can obtain any \glslink{left-leaning-stitch}{left-leaning} \glslink{stitch}{stitch} with \glslink{reduction}{reduction} $\beta$ by inserting \glslink{gap}{gaps} of various colors in between the \glslink{column}{columns} of $\beta$. If $\beta$ has $j$ \glslink{chain}{chains} and $k$ \glslink{wire}{wires}, then it has $j+k$ \glslink{column}{columns}, and so there are a total of $j+k+1$ places between the \glslink{column}{columns} (including at the beginning and the end) where \glslink{gap}{gaps} can be placed, and the power series $\frac{1}{1-\xi}$ counts the numbers of ways to insert a given number of \glslink{gap}{gaps} of each color into one of these locations. Consequently, we deduce that $\glslink{F-lparallel}{F_\lparallel}$ is given in terms of $\hat{G}_\lparallel$ by the formula \begin{equation}
        \glslink{F-lparallel}{F_\lparallel} = \frac{1}{1-\xi} \cdot \hat{G}_\lparallel\left(\frac{x_1}{1-\xi},\frac{y_1}{1-\xi};\dots;\frac{x_r}{1-\xi},\frac{y_r}{1-\xi}\right). \label{F-from-G-prime}
    \end{equation} Combining equations \eqref{G-from-F-tilde}, \eqref{G-prime-from-G}, and \eqref{F-from-G-prime} gives the claimed equation.
\end{proof}

\noindent
In light of Lemma~\ref{F-parallel-from-F-tilde-parallel-lem}, it is short work to prove the algebraic relation for $\glslink{F-tilde-lparallel}{\tilde{F}_\lparallel}$ stated in Section~\ref{proof-overview-section}.

\begin{proof}[Proof of Theorem~\ref{F-tilde-lparallel-algebraicity-thm}]
    This follows straightforwardly by substituting the formula from Lemma~\ref{F-parallel-from-F-tilde-parallel-lem} into the equation in Lemma~\ref{F-lparallel-algebraicity-lem}, applying the substitutions $x_i \mapsto \frac{x_i/y_i}{1+\frac{x_1}{y_1}+\dots+\frac{x_r}{y_r}}$ and $y_i \mapsto \frac{y_i}{1+\frac{x_1}{y_1}+\dots+\frac{x_r}{y_r}}$, and rearranging.
\end{proof}

\subsection{The generating function \texorpdfstring{$\glslink{H}{H}$}{H}}

We now turn our attention to the case where $r = 2$ and the generating function $\glslink{H}{H}$ defined in Lemma~\ref{H-from-F-lparallel-tilde-lem}. The proof of the claimed expression for $\glslink{H}{H}$ in terms of $\glslink{F-tilde-lparallel}{\tilde{F}_\lparallel}$ is essentially an exercise in translating certain bijections into corresponding generating function identities.

\begin{proof}[Proof of Lemma~\ref{H-from-F-lparallel-tilde-lem}]
    Observe that the coefficient of $L^\ell T^t N^\nu x^j y^k$ in $\glslink{H}{H}$ is equal to the number $\glslink{h-ell-m-t-nu-j-k}{h_{\ell,0,t,\nu}(j, k)}$ of \glslink{stitch-coloring}{$2$-colored} \glslink{parallel-stitch}{parallel} \glslink{stitch}{stitches} $\alpha$ such that $\glslink{alpha-rvert-i}{\alpha\rvert_1}$ has $\ell$ \glslink{proper-chain}{proper chains}, $0$ \glslink{gap}{gaps}, $t$ \glslink{jump}{jumps}, and $\nu$ \glslink{pole}{poles} and $\glslink{alpha-rvert-i}{\alpha\rvert_2}$ has $j$ \glslink{chain}{chains} and $k$ \glslink{wire}{wires}. Now, given such a \glslink{stitch}{stitch} $\alpha$, note that the \glslink{irreducible-stitch}{irreducible} components of $\alpha$ come in the following five forms: \begin{itemize}
        \item[(i)] \glslink{pole}{poles} with color $1$;
        \item[(ii)] \glslink{gap}{gaps} with color $2$;
        \item[(iii)] \glslink{pole}{poles} with color $2$;
        \item[(iv)] \glslink{irreducible-stitch}{irreducible} \glslink{stitch-coloring}{$2$-colored} \glslink{left-leaning-stitch}{left-leaning} \glslink{stitch}{stitches} with at least one \glslink{wire}{wire}; and
        \item[(v)] \glslink{irreducible-stitch}{irreducible} \glslink{stitch-coloring}{$2$-colored} \glslink{left-leaning-stitch}{right-leaning} \glslink{stitch}{stitches} with at least one \glslink{wire}{wire}.
    \end{itemize} Consequently, we may write \begin{equation}
        \glslink{H}{H} = \frac{1}{1-x-y-N-2G}, \label{H-from-G-fraction}
    \end{equation} where $G \in \Z[[L,T;x,y]]$ is the generating function defined by the property that the coefficient of $L^\ell T^t x^j y^k$ is the number of \glslink{irreducible-stitch}{irreducible} \glslink{stitch-coloring}{$2$-colored} \glslink{left-leaning-stitch}{left-leaning} \glslink{stitch}{stitches} $\alpha$ with at least one \glslink{wire}{wire} such that $\glslink{alpha-rvert-i}{\alpha\rvert_1}$ has $\ell$ \glslink{proper-chain}{proper chains}, $0$ \glslink{gap}{gaps}, and $t$ \glslink{jump}{jumps} and $\glslink{alpha-rvert-i}{\alpha\rvert_2}$ has $j$ \glslink{chain}{chains} and $k$ \glslink{wire}{wires}.

    Now, define $\hat{G} \in \Z[[L,T;x,y]]$ by the property that the coefficient of $L^\ell T^t x^j y^k$ is the number of \glslink{integral}{\textit{integral}} \glslink{stitch-coloring}{$2$-colored} \glslink{left-leaning-stitch}{left-leaning} \glslink{stitch}{stitches} $\alpha$ such that $\glslink{alpha-rvert-i}{\alpha\rvert_1}$ has $\ell$ \glslink{proper-chain}{proper chains} and $t$ \glslink{jump}{jumps} and $\glslink{alpha-rvert-i}{\alpha\rvert_2}$ has $j$ \glslink{chain}{chains} and $k$ \glslink{wire}{wires}. Then it is easy to see from the definition of the function $\glslink{F-tilde-lparallel}{\tilde{F}_\lparallel}$ that \begin{equation}
        \hat{G} = \glslink{F-tilde-lparallel}{\tilde{F}_\lparallel}\left(L, T; xy, y\right). \label{G-prime-from-F-lparallel}
    \end{equation} Meanwhile, using reasoning similar to what we used to derive equation~\eqref{F-from-G-prime}, we find that \begin{equation}
        G = (1-x) \cdot \hat{G}\left(\frac{L}{(1-x)^2}, \frac{T}{1-x}; \frac{x}{1-x}, \frac{y}{1-x}\right); \label{G-from-G-prime}
    \end{equation} this follows because, given an \glslink{integral}{integral} \glslink{stitch-coloring}{$2$-colored} \glslink{left-leaning-stitch}{left-leaning} \glslink{stitch}{stitch} $\alpha$ such that $\glslink{alpha-rvert-i}{\alpha\rvert_1}$ has $\ell$ \glslink{proper-chain}{proper chains} and $t$ \glslink{jump}{jumps} and $\glslink{alpha-rvert-i}{\alpha\rvert_2}$ has $j$ \glslink{chain}{chains} and $k$ \glslink{wire}{wires}, there are $2\ell+t+j+k-1$ places between \glslink{column}{columns} of $\alpha$ where \glslink{gap}{gaps} of color $2$ can be added while maintaining irreducibility of $\alpha$. From here, combining equations~\eqref{H-from-G-fraction}, \eqref{G-from-G-prime}, and \eqref{G-prime-from-F-lparallel} gives the desired result.
\end{proof}

\noindent
At this point, by combining Lemma~\ref{H-from-F-lparallel-tilde-lem} and Theorem~\ref{F-tilde-lparallel-algebraicity-thm}, we are able to prove the main takeaway of this section---the claimed degree-$3$ algebraic equation for $\glslink{H}{H}$.

\begin{proof}[Proof of Corollary~\ref{H-algebraicity-cor}]
    This follows tediously but straightforwardly by making the substitutions $x_1 \mapsto \frac{L}{(1-x)^2}$, $y_1 \mapsto \frac{T}{1-x}$, $x_2 \mapsto \frac{xy}{(1-x)^2}$, and $y_2 \mapsto \frac{y}{1-x}$ in the $r=2$ case of Theorem~\ref{F-tilde-lparallel-algebraicity-thm}, then substituting the formula for $\glslink{F-tilde-lparallel}{\tilde{F}_\lparallel}\left(\frac{L}{(1-x)^2}, \frac{T}{1-x}; \frac{xy}{(1-x)^2}, \frac{y}{1-x}\right)$ in terms of $\glslink{H}{H}$ that is implied by Lemma~\ref{H-from-F-lparallel-tilde-lem} and rearranging.
\end{proof}

\section{Explicit formulas for \texorpdfstring{$\glslink{H-ell-m-t-nu}{H_{\ell,m,t,\nu}}$}{Hℓ,m,t,ν}} \label{explicit-formulas-for-H-section}

The sole purpose of this section will be to use Corollary~\ref{H-algebraicity-cor} to prove Theorem~\ref{H-from-p-thm} giving explicit formulas for the $\glslink{H-ell-m-t-nu}{H_{\ell,m,t,\nu}}$. We break up the proof of the theorem into a number of smaller lemmas that each deal with different possible values for $\ell$, $m$, $t$, $\nu$, and $i$.

We start with two lemmas that will act as ``base cases'' for an inductive argument on the variables $m$ and $i$. The first of these lemmas addresses the case where $\ell = m = t = i = 0$, whereas the second addresses the case where $\ell \ge 1$ and $m = i = 0$; together, the two lemmas address all possible cases where $m = i = 0$.

\begin{lemma} \label{H-formula-ell-equals-zero-lem}
    For all $\nu \ge 0$, we have \begin{equation*}
        \glslink{H-ell-m-t-nu}{H_{0,0,0,\nu}} = \glslink{Q}{Q}^{-1-\nu}.
    \end{equation*}
\end{lemma}

\begin{proof}
    Adopt the notation $\glslink{Q}{Q}$, $\glslink{U}{U}$, $\glslink{V}{V}$, and $\glslink{W}{W}$ from Section~\ref{proof-overview-section}, and observe that Corollary~\ref{H-algebraicity-cor} may be restated as \begin{equation}
        L = -\frac{(1+(2\glslink{V}{V}\glslink{Q}{Q}-\glslink{Q}{Q}-2T+N)\glslink{H}{H})(1+(\glslink{Q}{Q}+N)\glslink{H}{H})(1-(\glslink{Q}{Q}-N)\glslink{H}{H})}{4\glslink{H}{H}^2(1+(2W\glslink{Q}{Q}-\glslink{Q}{Q}+N)\glslink{H}{H})}. \label{L-from-H}
    \end{equation} In particular, since the ring of formal power series $\Q[[L,T,N;x,y]]$ has no zero divisors, this implies that $\glslink{H}{H}\rvert_{L=0}$ is either $-\frac{1}{1-x+y-2T+N}$, $-\frac{1}{\glslink{Q}{Q}+N}$, or $\frac{1}{\glslink{Q}{Q}-N}$. But the first two power series have constant coefficient $-1$ whereas $\glslink{H}{H}\rvert_{L=0}$ has constant coefficient $1$, so in fact we must have $\glslink{H}{H}\rvert_{L=0} = \frac{1}{\glslink{Q}{Q}-N}$. Since $\glslink{H-ell-m-t-nu}{H_{0,0,0,\nu}}$ is the coefficient of $L^0T^0N^\nu$ in $\glslink{H}{H}$, the statement of the lemma follows immediately.
\end{proof}

\begin{lemma} \label{H-formula-m-equals-zero-lem}
    For all $\ell \ge 1$ and $t, \nu \ge 0$, we have \begin{equation*}
        \glslink{H-ell-m-t-nu}{H_{\ell,0,t,\nu}} = \glslink{W}{W}^\ell \glslink{V}{V}^{-\ell-t} \glslink{Q}{Q}^{-1-2\ell-t-\nu} \glslink{p-ell-m-t-nu-i}{p_{\ell,0,t,\nu}^{(0)}}\left(-\frac{1}{\glslink{W}{W}}, \frac{1}{\glslink{V}{V}}\right),
    \end{equation*} where $\glslink{p-ell-m-t-nu-i}{p_{\ell,0,t,\nu}^{(0)}} \in \Q[z, w]$ is a polynomial given explicitly by \begin{equation*}
        \glslink{p-ell-m-t-nu-i}{p_{\ell,0,t,\nu}^{(0)}} := \frac{\nu+1}{\ell} \sum_{a+b+c+d=\ell-1} \glslink{binomial}{\binom{\ell}{a}} \binom{b+t+\ell-1}{b,\,t,\,\ell-1} \glslink{binomial}{\binom{c+\ell-1}{c}} 2^{d+1} \glslink{binomial}{\binom{\nu+1+d}{d}} z^a w^b.
    \end{equation*}
\end{lemma}

\begin{proof}
    For each $\ell \ge 1$ and $t, \nu \ge 0$, define a power series \begin{equation}
        P_{\ell,0,t,\nu} := \glslink{W}{W}^{-\ell} \glslink{V}{V}^{\ell+t} \glslink{Q}{Q}^{1+2\ell+t+\nu} \glslink{H-ell-m-t-nu}{H_{\ell,0,t,\nu}} \in \Q[[x, y]], \label{P-ell-t-nu-def}
    \end{equation} and let \begin{equation*}
        \tilde{P} := (1-N)\sum_{\ell \ge 1,\, t,\nu \ge 0} P_{\ell,0,t,\nu} \cdot L^\ell T^t N^\nu \in \Q[[L,T,N; x,y]].
    \end{equation*} Then we may write \begin{align}
        \tilde{P} &= (1-N)\sum_{\ell \ge 1,\, t,\nu \ge 0} \glslink{W}{W}^{-\ell} \glslink{V}{V}^{\ell+t} \glslink{Q}{Q}^{1+2\ell+t+\nu} \glslink{H-ell-m-t-nu}{H_{\ell,0,t,\nu}} L^\ell T^t N^\nu \notag \\
        &= \glslink{Q}{Q}(1-N) \sum_{\ell \ge 1,\, t,\nu \ge 0} \glslink{H-ell-m-t-nu}{H_{\ell,0,t,\nu}} \cdot \left(\frac{\glslink{V}{V}\glslink{Q}{Q}^2L}{\glslink{W}{W}}\right)^\ell (\glslink{V}{V}\glslink{Q}{Q}T)^t (\glslink{Q}{Q}N)^\nu \notag \\
        &= \glslink{Q}{Q}(1-N)\left[\glslink{H}{H} - \frac{1}{\glslink{Q}{Q}-N}\right]_{L=\frac{VQ^2L}{\glslink{W}{W}},\,T=\glslink{V}{V}\glslink{Q}{Q}T\,N=\glslink{Q}{Q}N} \notag \\
        &= \glslink{Q}{Q}(1-N) \glslink{H}{H}\left(\frac{\glslink{V}{V}\glslink{Q}{Q}^2L}{\glslink{W}{W}},\, \glslink{V}{V}\glslink{Q}{Q}T,\, \glslink{Q}{Q}N;\, x,\, y\right) - 1, \label{P-hat-from-H}
    \end{align} where the third step uses our observation from Lemma~\ref{H-formula-ell-equals-zero-lem} that $\glslink{H}{H}\rvert_{L=0} = \frac{1}{Q-N}$. If we make the substitutions $L \mapsto \frac{VQ^2L}{W}$, $T \mapsto VQT$, and $N \mapsto QN$ in equation \eqref{L-from-H} and then apply equation \eqref{P-hat-from-H}, then, after some manipulation, we find that \begin{equation*}
        L = \frac{(1-N)\big(2(1-T)(1+\tilde{P})-(1-N)\tilde{P} \cdot \frac{1}{\glslink{V}{V}}\big)\big(2+(1+N)\tilde{P}\big)\tilde{P}}{4(1+\tilde{P})^2\big(2(1+\tilde{P})-(1-N)\tilde{P} \cdot \frac{1}{\glslink{W}{W}}\big)}.
    \end{equation*} Equivalently, adopting the slight abuse of notation $z := -\frac{1}{\glslink{W}{W}} \in \Q[[x,y]]$ and $w := \frac{1}{\glslink{V}{V}} \in \Q[[x,y]]$, we have $L = g(\tilde{P})$, where \begin{equation*}
        g(p) := \frac{(1-N)\big(2(1-T)(1+p)-(1-N)p \cdot w\big)\big(2+(1+N)p\big)p}{4(1+p)^2\big(2(1+p)+(1-N)p \cdot z\big)}.
    \end{equation*}

    Next, write \begin{equation*}
        \tilde{P} =: \sum_{\ell\ge1} \tilde{P}_\ell \cdot L^\ell
    \end{equation*} for some formal power series $\tilde{P}_\ell \in \Q[[T,N; x,y]]$, and observe that by the Lagrange inversion theorem, for all $\ell \ge 1$, $\tilde{P}_\ell$ is equal to the coefficient of $p^{\ell-1}$ in the formal power series \begin{equation}
        \frac{1}{\ell}\left(\frac{p}{g(p)}\right)^\ell = \frac{4^\ell (1+p)^{2\ell}\big(2(1+p)+(1-N)pz\big)^\ell}{\ell(1-N)^\ell\big(2(1-T)(1+p)-(1-N)pw\big)^\ell\big(2+(1+N)p\big)^\ell} \in \Q[[T, N; x, y]][[p]]. \label{lagrange-inversion-p-series}
    \end{equation} Now, if we expand the expression in \eqref{lagrange-inversion-p-series} as a power series in $z$ and $w$ using the binomial theorem, then for each $a, b \ge 0$, we find that the coefficient of $z^a w^b$ is given by \begin{align*}
        &\frac{1}{\ell} \cdot \glslink{binomial}{\binom{\ell}{a}} \glslink{binomial}{\binom{b+\ell-1}{b}} \frac{2^{2\ell-a-b} (1-T)^{-\ell-b} (1-N)^{a+b-\ell} (1+p)^{2\ell-a-b}}{\big(2(1+p)-(1-N)p\big)^\ell} \cdot p^{a+b} \\
        &= \frac{1}{\ell} \sum_{c \ge 0} \glslink{binomial}{\binom{\ell}{a}} \glslink{binomial}{\binom{b+\ell-1}{b}} \glslink{binomial}{\binom{c+\ell-1}{c}} 2^{\ell-a-b-c} (1-T)^{-\ell-b} (1-N)^{a+b+c-\ell} (1+p)^{\ell-a-b-c} \cdot p^{a+b+c}.
    \end{align*} At this point, another application of the binomial theorem to the factor $(1+p)^{\ell-a-b-c}$ shows that the coefficient of $z^a w^b p^{\ell-1}$ in \eqref{lagrange-inversion-p-series} is \begin{equation*}
        \frac{1}{\ell} \sum_{c+d=\ell-1-a-b} \glslink{binomial}{\binom{\ell}{a}} \glslink{binomial}{\binom{b+\ell-1}{b}} \glslink{binomial}{\binom{c+\ell-1}{c}} 2^{d+1} (d+1) \cdot (1-T)^{-\ell-b} (1-N)^{-1-d}.
    \end{equation*} It therefore follows that for each $\ell \ge 1$, $\tilde{P}_\ell$ is given by the expression \begin{equation*}
        \frac{1}{\ell} \sum_{a+b+c+d=\ell-1} \glslink{binomial}{\binom{\ell}{a}} \glslink{binomial}{\binom{b+\ell-1}{b}} \glslink{binomial}{\binom{c+\ell-1}{c}} 2^{d+1} (d+1) \cdot (1-T)^{-\ell-b} (1-N)^{-1-d} \cdot z^a w^b.
    \end{equation*}

    Finally, observe that by the definitions of $\tilde{P}$ and $\tilde{P}_\ell$, for each $\ell \ge 1$, we have \begin{equation*}
        \frac{\tilde{P}_\ell}{1-N} = \sum_{t,\nu\ge0} P_{\ell,0,t,\nu} \cdot T^t N^\nu.
    \end{equation*} In other words, $P_{\ell,0,t,\nu}$ is equal to the coefficient of $T^tN^\nu$ in the expression \begin{equation*}
        \frac{1}{\ell} \sum_{a+b+c+d=\ell-1} \glslink{binomial}{\binom{\ell}{a}} \glslink{binomial}{\binom{b+\ell-1}{b}} \glslink{binomial}{\binom{c+\ell-1}{c}} 2^{d+1} (d+1) \cdot (1-T)^{-\ell-b} (1-N)^{-2-d} \cdot z^a w^b.
    \end{equation*} From here, applying the binomial theorem to the factors $(1-T)^{-\ell-b}$ and $(1-N)^{-2-d}$ and simplifying gives \begin{equation*}
        P_{\ell,0,t,\nu} = \glslink{p-ell-m-t-nu-i}{p_{\ell,0,t,\nu}^{(0)}}(z, w),
    \end{equation*} and applying \eqref{P-ell-t-nu-def} gives the desired formula for $\glslink{H-ell-m-t-nu}{H_{\ell,0,t,\nu}}$.
\end{proof}

Having addressed the main ``base cases,'' we next prove two lemmas that show first how to generalize to all $m \ge 0$, and then how to generalize to all $i \ge 0$.

\begin{lemma} \label{H-formula-i-equals-zero-lem}
    Let $\ell, m, t, \nu \ge 0$, and suppose that there exists a polynomial $\glslink{p-ell-m-t-nu-i}{p_{\ell,m,t,\nu}^{(0)}} \in \Q[z, w]$ such that \begin{equation*}
        \glslink{H-ell-m-t-nu}{H_{\ell,m,t,\nu}} = \glslink{W}{W}^{\ell+m} \glslink{V}{V}^{-\ell-t} Q^{-1-2\ell-m-t-\nu} \glslink{p-ell-m-t-nu-i}{p_{\ell,m,t,\nu}^{(0)}}\left(-\frac{1}{\glslink{W}{W}}, \frac{1}{\glslink{V}{V}}\right).
    \end{equation*} Then we have \begin{equation*}
        \glslink{H-ell-m-t-nu}{H_{\ell,m+1,t,\nu}} = \glslink{W}{W}^{\ell+(m+1)} \glslink{V}{V}^{-\ell-t} Q^{-1-2\ell-(m+1)-t-\nu} \glslink{p-ell-m-t-nu-i}{p_{\ell,m+1,t,\nu}^{(0)}}\left(-\frac{1}{\glslink{W}{W}}, \frac{1}{\glslink{V}{V}}\right),
    \end{equation*} where $\glslink{p-ell-m-t-nu-i}{p_{\ell,m+1,t,\nu}^{(0)}} \in \Q[z, w]$ is defined by \begin{multline*}
        \glslink{p-ell-m-t-nu-i}{p_{\ell,m+1,t,\nu}^{(0)}} := \frac{1}{m+1}\Big((3\ell+3m+\nu+1)z+(\ell+t)w+(4\ell+4m+2\nu+2) \\
        -2z(z+1)\partial_z-w(z-w+2)\partial_w\Big) \glslink{p-ell-m-t-nu-i}{p_{\ell,m,t,\nu}^{(0)}}.
    \end{multline*}
\end{lemma}

\begin{proof}
    Observe that for all $\ell, m, t, \nu \ge 0$ and $j, k \ge 0$, we have \begin{equation}
        \glslink{h-ell-m-t-nu-j-k}{h_{\ell,m,t,\nu}(j, k)} = \glslink{binomial}{\binom{2\ell+m+t+\nu+j+k}{m}} \glslink{h-ell-m-t-nu-j-k}{h_{\ell,0,t,\nu}(j, k)}. \label{h-reduction}
    \end{equation} Indeed, this follows by similar reasoning to Lemma~\ref{reduction-lem}: a \glslink{stitch-coloring}{$2$-colored} \glslink{stitch}{stitch} $\alpha$ such that $\glslink{alpha-rvert-i}{\alpha\rvert_1}$ has $m$ \glslink{gap}{gaps} is determined from the \glslink{stitch}{stitch} $\alpha'$ obtained by removing these \glslink{gap}{gaps} by specifying where to ``add the $m$ \glslink{gap}{gaps} back in.'' The quantity on the top of the binomial then comes from the fact that a \glslink{stitch-coloring}{$2$-colored} \glslink{stitch}{stitch} $\alpha$ such that $\glslink{alpha-rvert-i}{\alpha\rvert_1}$ has $\ell$ \glslink{proper-chain}{proper chains}, $m$ \glslink{gap}{gaps}, $t$ \glslink{jump}{jumps}, and $\nu$ \glslink{pole}{poles} and $\glslink{alpha-rvert-i}{\alpha\rvert_2}$ has $j$ \glslink{chain}{chains} and $k$ \glslink{wire}{wires} has a total of $2\ell+m+t+\nu+j+k$ \glslink{column}{columns}. By combining \eqref{h-reduction} with the same formula with $m+1$ substituted for $m$, we thus obtain the recursive formula \begin{equation*}
        \glslink{h-ell-m-t-nu-j-k}{h_{\ell,m+1,t,\nu}(j, k)} = \frac{2\ell+m+1+t+\nu+j+k}{m+1} \glslink{h-ell-m-t-nu-j-k}{h_{\ell,m,t,\nu}(j, k)}.
    \end{equation*} In terms of the generating functions $\glslink{H-ell-m-t-nu}{H_{\ell,m,t,\nu}}$, this translates to the equation \begin{equation}
        \glslink{H-ell-m-t-nu}{H_{\ell,m+1,t,\nu}} = \frac{1}{m+1} (2\ell+m+1+t+\nu+x\partial_x+y\partial_y) \glslink{H-ell-m-t-nu}{H_{\ell,m,t,\nu}}. \label{H-m-recursion}
    \end{equation} From here, the statement of the lemma is a tedious but routine exercise in substituting the formula for $\glslink{H-ell-m-t-nu}{H_{\ell,m,t,\nu}}$ in terms of $\glslink{p-ell-m-t-nu-i}{p_{\ell,m,t,\nu}^{(0)}}$ on the right hand side of \eqref{H-m-recursion} and applying the chain rule.
\end{proof}

\begin{lemma} \label{H-formula-general-i-lem}
    Let $\ell, m, t, \nu, i \ge 0$, and suppose that there exists a polynomial $\glslink{p-ell-m-t-nu-i}{p_{\ell,m,t,\nu}^{(i)}} \in \Q[z, w]$ such that \begin{equation*}
        \glslink{partial-x-i}{\partial_x^i} \glslink{H-ell-m-t-nu}{H_{\ell,m,t,\nu}} = \glslink{W}{W}^{\ell+m} \glslink{V}{V}^{i-\ell-t} \glslink{Q}{Q}^{-1-2\ell-m-t-\nu-i} \glslink{p-ell-m-t-nu-i}{p_{\ell,m,t,\nu}^{(i)}}\left(-\frac{1}{\glslink{W}{W}}, \frac{1}{\glslink{V}{V}}\right).
    \end{equation*} Then we have \begin{equation*}
        \glslink{partial-x-i}{\partial_x^{i+1}} \glslink{H-ell-m-t-nu}{H_{\ell,m,t,\nu}} = \glslink{W}{W}^{\ell+m} \glslink{V}{V}^{(i+1)-\ell-t} \glslink{Q}{Q}^{-1-2\ell-m-t-\nu-(i+1)} \glslink{p-ell-m-t-nu-i}{p_{\ell,m,t,\nu}^{(i+1)}}\left(-\frac{1}{\glslink{W}{W}}, \frac{1}{\glslink{V}{V}}\right),
    \end{equation*} where $\glslink{p-ell-m-t-nu-i}{p_{\ell,m,t,\nu}^{(i+1)}} \in \Q[z,w]$ is defined by \begin{multline*}
        \glslink{p-ell-m-t-nu-i}{p_{\ell,m,t,\nu}^{(i+1)}} := \Big((\ell+m)z-(\ell+2m+3i-t+\nu+1)w + (4\ell+4m+4i+2\nu+2) \\
        -z(z-w+2)\partial_z+2w(w-1)\partial_w\Big) \glslink{p-ell-m-t-nu-i}{p_{\ell,m,t,\nu}^{(i)}}(z,w).
    \end{multline*}
\end{lemma}

\begin{proof}
    Similarly to Lemma~\ref{H-formula-i-equals-zero-lem}, this follows straightforwardly by starting with the equation \begin{equation*}
        \glslink{partial-x-i}{\partial_x^{i+1}} \glslink{H-ell-m-t-nu}{H_{\ell,m,t,\nu}} = \partial_x(\glslink{partial-x-i}{\partial_x^i} \glslink{H-ell-m-t-nu}{H_{\ell,m,t,\nu}}),
    \end{equation*} substituting the formula for $\glslink{partial-x-i}{\partial_x^i} \glslink{H-ell-m-t-nu}{H_{\ell,m,t,\nu}}$ in terms of $\glslink{p-ell-m-t-nu-i}{p_{\ell,m,t,\nu}^{(i)}}$, and applying the chain rule.
\end{proof}

The previous four lemmas address all cases covered by Theorem~\ref{H-from-p-thm} which satisfy $i \ge 0$. To address the case where $i = -1$, we require two additional lemmas.

\begin{lemma} \label{H-formula-i-negative-one-base-case-lem}
    We have \begin{equation*}
        \glslink{partial-x-i}{\partial_x^{-1}} \glslink{H-ell-m-t-nu}{H_{0,1,0,0}} = \glslink{W}{W} \glslink{V}{V}^{-1} \glslink{Q}{Q}^{-1} - 1.
    \end{equation*}
\end{lemma}

\begin{proof}
    By Lemma~\ref{H-formula-ell-equals-zero-lem}, we have $\glslink{H-ell-m-t-nu}{H_{0,0,0,0}} = \glslink{Q}{Q}^{-1} \cdot 1$. Setting $\glslink{p-ell-m-t-nu-i}{p_{0,0,0,0}^{(0)}} := 1$ and applying Lemma~\ref{H-formula-i-equals-zero-lem} then produces \begin{equation*}
        \glslink{p-ell-m-t-nu-i}{p_{0,1,0,0}^{(0)}} = \left(z+2 - 2z(z+1)\partial_z - w(z-w+2)\partial_w\right) \cdot 1 = z+2
    \end{equation*} and hence \begin{equation*}
        \glslink{H-ell-m-t-nu}{H_{0,1,0,0}} = \glslink{W}{W} \glslink{Q}{Q}^{-2} \left(2-\frac{1}{\glslink{W}{W}}\right) = (2\glslink{W}{W}-1)\glslink{Q}{Q}^{-2}.
    \end{equation*} A quick computation shows that \begin{equation*}
        \partial_x(\glslink{W}{W} \glslink{V}{V}^{-1} \glslink{Q}{Q}^{-1} - 1) = (2\glslink{W}{W}-1)\glslink{Q}{Q}^{-2} = \glslink{H-ell-m-t-nu}{H_{0,1,0,0}}
    \end{equation*} and \begin{equation*}
        \left[\glslink{W}{W} \glslink{V}{V}^{-1} \glslink{Q}{Q}^{-1} - 1\right]_{x=0} = 0.
    \end{equation*} The result then follows by the definition of the operator $\glslink{partial-x-i}{\partial_x^{-1}}$.
\end{proof}

\begin{lemma} \label{H-formula-i-negative-one-inductive-step-lem}
    Let $m \ge 1$, and suppose that there exists a polynomial $\glslink{p-ell-m-t-nu-i}{p_{0,m,0,0}^{(-1)}} \in \Q[z, w]$ such that \begin{equation*}
        \glslink{partial-x-i}{\partial_x^{-1}} \glslink{H-ell-m-t-nu}{H_{0,m,0,0}} = \glslink{W}{W}^m \glslink{V}{V}^{-1} \glslink{Q}{Q}^{-m} \glslink{p-ell-m-t-nu-i}{p_{0,m,0,0}^{(-1)}}\left(-\frac{1}{\glslink{W}{W}}, \frac{1}{V}\right) - \frac{1}{m}.
    \end{equation*} Then we have \begin{equation*}
        \glslink{partial-x-i}{\partial_x^{-1}} \glslink{H-ell-m-t-nu}{H_{0,m+1,0,0}} = \glslink{W}{W}^{m+1} \glslink{V}{V}^{-1} \glslink{Q}{Q}^{-(m+1)} \glslink{p-ell-m-t-nu-i}{p_{0,m+1,0,0}^{(-1)}}\left(-\frac{1}{\glslink{W}{W}}, \frac{1}{\glslink{V}{V}}\right) - \frac{1}{m+1},
    \end{equation*} where $\glslink{p-ell-m-t-nu-i}{p_{0,m+1,0,0}^{(-1)}} \in \Q[z, w]$ is defined by \begin{equation*}
        \glslink{p-ell-m-t-nu-i}{p_{0,m+1,0,0}^{(-1)}} := \frac{1}{m+1} \big((3m-1)z + w + (4m-2) - 2z(z+1)\partial_z - w(z-w+2)\partial_w\big) \glslink{p-ell-m-t-nu-i}{p_{0,m,0,0}^{(-1)}}.
    \end{equation*}
\end{lemma}

\begin{proof}
    Observe that in the case where $\ell = t = \nu = 0$, equation \eqref{H-m-recursion} gives \begin{equation}
        \glslink{H-ell-m-t-nu}{H_{0,m+1,0,0}} = \frac{1}{m+1}(m+1+x\partial_x+y\partial_y) \glslink{H-ell-m-t-nu}{H_{0,m,0,0}}. \label{H-m-recursion-special-case}
    \end{equation} It is then not difficult to obtain the relation \begin{equation}
        \glslink{partial-x-i}{\partial_x^{-1}} \glslink{H-ell-m-t-nu}{H_{0,m+1,0,0}} = \frac{1}{m+1}(m+x\partial_x+y\partial_y) \glslink{partial-x-i}{\partial_x^{-1}} \glslink{H-ell-m-t-nu}{H_{0,m,0,0}}: \label{H-m-integral-recursion}
    \end{equation} we readily verify that the right hand side of \eqref{H-m-integral-recursion} is an antiderivative of the right hand side of \eqref{H-m-recursion-special-case} and evaluates to $0$ when $x = 0$. From here, we can once again substitute the formula for $\glslink{partial-x-i}{\partial_x^{-1}}\glslink{H-ell-m-t-nu}{H_{0,m,0,0}}$ in terms of $\glslink{p-ell-m-t-nu-i}{p_{0,m,0,0}^{(-1)}}$ and apply the chain rule to obtain the desired result.
\end{proof}

Finally, we combine the six lemmas of this section to give a proof of Theorem~\ref{H-from-p-thm}.

\begin{proof}[Proof of Theorem~\ref{H-from-p-thm}]
    To start, assume that $i \ge 0$, and proceed by induction, first on $m$ and then on $i$. If $m = i = 0$, then either $\ell = t = 0$, in which case the result follows with $\glslink{p-ell-m-t-nu-i}{p_{0,0,0,\nu}^{(0)}} = 1$ by Lemma~\ref{H-formula-ell-equals-zero-lem}; or $\ell \ge 1$, and the result follows by Lemma~\ref{H-formula-m-equals-zero-lem}. If $i = 0$ and $m \ge 1$ and the result holds with $m$ replaced by $m-1$, then the claim follows by Lemma~\ref{H-formula-i-equals-zero-lem}. Finally, if $i \ge 1$ and the result holds with $i$ replaced by $i-1$, then the claim follows by Lemma~\ref{H-formula-general-i-lem}. At each step, we can easily verify the claimed degree constraints on the polynomials $\glslink{p-ell-m-t-nu-i}{p_{\ell,m,t,\nu}^{(i)}}$.

    Next, suppose that $i = -1$, so $\ell = t = \nu = 0$. Then the case $m = 1$ follows by Lemma~\ref{H-formula-i-negative-one-base-case-lem}, and Lemma~\ref{H-formula-i-negative-one-inductive-step-lem} shows that if the claim holds with $m$ replaced by $m-1$, then it holds for $m$. Once again, the degree constraints are easy to verify at every step.

    Lastly, to prove claim (iv), suppose that $\ell = m = t = 0$, and proceed by induction on $i$. When $i = 0$, $\glslink{p-ell-m-t-nu-i}{p_{0,0,0,\nu}^{(i)}} = 1$, and the claim is immediate. For the inductive step, note first that the degree constraint we proved in part (i) of the theorem tells us that $\glslink{deg-x}{\deg_z(\glslink{p-ell-m-t-nu-i}{p_{0,0,0,\nu}^{(i)}})} = 0$ for all $i \ge 0$. Consequently, the recursive definition of $\glslink{p-ell-m-t-nu-i}{p_{0,0,0,\nu}^{(i+1)}}$ from Lemma~\ref{H-formula-general-i-lem} simplifies to \begin{equation*}
        \glslink{p-ell-m-t-nu-i}{p_{0,0,0,\nu}^{(i+1)}}(z, w) = (-(3i+\nu+1)w + (4i+2\nu+2) + 2w(w-1)\partial_w) \glslink{p-ell-m-t-nu-i}{p_{0,0,0,\nu}^{(i)}}(z, w),
    \end{equation*} which in particular implies that \begin{equation*}
        \glslink{p-ell-m-t-nu-i}{p_{0,0,0,\nu}^{(i+1)}}(z, 0) = (4i+2\nu+2)\glslink{p-ell-m-t-nu-i}{p_{0,0,0,\nu}^{(i)}}(z, 0)
    \end{equation*} and \begin{equation*}
        (\partial_w \glslink{p-ell-m-t-nu-i}{p_{0,0,0,\nu}^{(i+1)}})(z, 0) = -(3i+\nu+1)\glslink{p-ell-m-t-nu-i}{p_{0,0,0,\nu}^{(i)}}(z, 0) + (4i+2\nu) (\partial_w \glslink{p-ell-m-t-nu-i}{p_{0,0,0,\nu}^{(i)}})(z, 0).
    \end{equation*} If we assume the claimed formulas for $\glslink{p-ell-m-t-nu-i}{p_{0,0,0,\nu}^{(i)}}$ and $(\partial_w \glslink{p-ell-m-t-nu-i}{p_{0,0,0,\nu}^{(i)}})(z, 0)$, it is then routine to verify that the claimed formulas also hold for $\glslink{p-ell-m-t-nu-i}{p_{0,0,0,\nu}^{(i+1)}}$ and $(\partial_w \glslink{p-ell-m-t-nu-i}{p_{0,0,0,\nu}^{(i+1)}})(z, 0)$.
\end{proof}

\section{Explicit characterizations of \texorpdfstring{$\glslink{calE-f}{\calE_f}$}{Ef}} \label{explicit-formulas-for-E-section}

The goal of this section is to combine Corollary~\ref{E-from-H-ell-m-t-nu-cor} with the explicit formulas given by Theorem~\ref{H-from-p-thm} to obtain various explicit characterizations of the generating functions $\glslink{calE-f}{\calE_f}$, which will in turn give rise to proofs of our main theorems---Theorems~\ref{main-thm-1}, \ref{main-thm-2}, and \ref{main-thm-3}. More precisely, we will start by considering the restricted generating functions $\left.\glslink{partial-x-i}{\partial_x^j} \glslink{calE-f}{\calE_f}\right\rvert_{x=0} \in \Q[[y]]$ for fixed $j \ge 0$, which encode the information of the values $\glslink{E-f-n-k}{E_f(j+k, k)}$ for $k \ge 0$. By obtaining explicit formulas for these restricted generating functions, we will obtain proofs of Theorems~\ref{main-thm-2} and \ref{main-thm-3}(ii). We will then turn to computing closed forms for $\glslink{partial-x-i}{\partial_x^i} \glslink{calE-f}{\calE_f}$ for sufficiently large $i$, which will allow us to prove Theorems~\ref{main-thm-1} and \ref{main-thm-3}(i). In the latter case, we will see that the formulas for $\glslink{partial-x-i}{\partial_x^i}\glslink{calE-f}{\calE_f}$ are particularly simple when $f$ is \glslink{pure-character-polynomial}{pure} in the sense of Section~\ref{background-on-character-polynomials-section}.

\subsection{Explicit formulas for \texorpdfstring{$\left.\glslink{partial-x-i}{\partial_x^j}\glslink{calE-f}{\calE_{\glslink{I-rho}{I_\rho}}}\right\rvert_{x=0}$}{∂jxEIρ|x=0}}

On the subject of the restricted functions $\left.\glslink{partial-x-i}{\partial_x^j} \glslink{calE-f}{\calE_f}\right\rvert_{x=0}$, our key result is the following theorem, which shows that these functions are always given by Laurent polynomials in $\frac{1}{1-y}$. To start, we focus on the case where $f = \glslink{I-rho}{I_\rho}$, where slightly more can be said about $\left.\glslink{partial-x-i}{\partial_x^j} \glslink{calE-f}{\calE_f}\right\rvert_{x=0}$.

\begin{theorem} \label{E-I-rho-from-delta-thm}
    Let $\rho$ be any partition, and let $j \ge 0$. Then there exists a Laurent polynomial $\glslink{delta-f-j}{\delta_{\glslink{I-rho}{I_\rho}}^{(j)}} \in \Q[v, \frac{1}{v}]$ such that \begin{equation*}
        \left.\glslink{partial-x-i}{\partial_x^j}\glslink{calE-f}{\calE_{\glslink{I-rho}{I_\rho}}}\right|_{x=0} = \glslink{delta-f-j}{\delta_{\glslink{I-rho}{I_\rho}}^{(j)}}\left(\frac{1}{1-y}\right).
    \end{equation*} In particular, if $0 \le j < \glslink{len}{\len(\rho)} - \glslink{m-j-rho}{m_1(\rho)}$, then $\glslink{delta-f-j}{\delta_{\glslink{I-rho}{I_\rho}}^{(j)}} = 0$, and otherwise, the following hold: \begin{itemize}
        \item[(i)] The smallest power of $v$ appearing in $\glslink{delta-f-j}{\delta_{\glslink{I-rho}{I_\rho}}^{(j)}}$ is at least \begin{equation}
            \begin{cases}
                1+j & \text{if $\rho = 1^r$ for some $r \ge 0$}, \\
                2-\glslink{abs-rho}{\lvert\rho\rvert}+\max(\glslink{len}{\len(\rho)},j) & \text{otherwise}.
            \end{cases} \label{delta-smallest-v-power}
        \end{equation}
        \item[(ii)] The largest power of $v$ appearing in $\glslink{delta-f-j}{\delta_{\glslink{I-rho}{I_\rho}}^{(j)}}$ is at most $1+2j+2\glslink{m-j-rho}{m_1(\rho)}-\glslink{len}{\len(\rho)}$.
        \item[(iii)] If $\rho = 1^r$ for some $r \ge 0$, then the coefficients of $v^{1+2j+r}$ and $v^{2j+r}$ in $\glslink{delta-f-j}{\delta_{\glslink{I-rho}{I_\rho}}^{(j)}}$ are \begin{equation*}
            \frac{2^{2j}(\frac{r-1}{2}+j)!}{(\frac{r-1}{2})!} \qquad \text{and} \qquad -\frac{2^{2j-1}(j+2r)(\frac{r-1}{2}+j)!}{(\frac{r-1}{2})!}
        \end{equation*} respectively.
    \end{itemize}
\end{theorem}

\begin{proof}
    To start, define \begin{equation}
        \glslink{hat-p-ell-m-t-nu-i}{\hat{p}_{\ell,m,t,\nu}^{(i)}} := \frac{\ell!m!\nu!}{\glslink{binomial}{\binom{\ell+t-1}{t}}} \glslink{p-ell-m-t-nu-i}{p_{\ell,m,t,\nu}^{(i)}} \label{p-hat-def}
    \end{equation} for all $\ell,m,t,\nu,i \ge 0$, and use Corollary~\ref{E-from-H-ell-m-t-nu-cor} and Theorem~\ref{H-from-p-thm} to write \begin{align*}
        &\left.\glslink{z-rho}{z_\rho}\glslink{partial-x-i}{\partial_x^j}\glslink{calE-f}{\calE_{\glslink{I-rho}{I_\rho}}}\right\rvert_{x=0} = \sum_{\substack{2\ell+m+t+\nu=\glslink{abs-rho}{\lvert\rho\rvert} \\ \text{$\ell \ge 1$ or $t = 0$} \\ \ell+m \le j}} \glslink{stitch-ell-m-t-nu-binom}{\binom{\rho}{\ell,m,t,\nu}} \cdot \frac{\ell!m!\nu!}{\glslink{binomial}{\binom{\ell+t-1}{t}}} \cdot y^{\ell+t+\nu} \left[\glslink{partial-x-i}{\partial_x^{j-\ell-m}} \glslink{H-ell-m-t-nu}{H_{\ell,m,t,\nu}}\right]_{x=0} \\
        &= \sum_{\substack{2\ell+m+t+\nu=\glslink{abs-rho}{\lvert\rho\rvert} \\ \text{$\ell \ge 1$ or $t = 0$} \\ \ell+m \le j}} \glslink{stitch-ell-m-t-nu-binom}{\binom{\rho}{\ell,m,t,\nu}} \cdot y^{\ell+t+\nu} \left[\glslink{W}{W}^{\ell+m} \glslink{V}{V}^{j-2\ell-m-t} \glslink{Q}{Q}^{-1-\ell-t-\nu-j} \glslink{hat-p-ell-m-t-nu-i}{\hat{p}_{\ell,m,t,\nu}^{(j-\ell-m)}}\left(-\frac{1}{\glslink{W}{W}}, \frac{1}{\glslink{V}{V}}\right)\right]_{x=0} \\
        &= \sum_{\substack{2\ell+m+t+\nu=\glslink{abs-rho}{\lvert\rho\rvert} \\ \text{$\ell \ge 1$ or $t = 0$} \\ \ell+m \le j}} \glslink{stitch-ell-m-t-nu-binom}{\binom{\rho}{\ell,m,t,\nu}} \cdot y^{\ell+t+\nu} (1-y)^{-1+\ell+m-\nu-2j} \glslink{hat-p-ell-m-t-nu-i}{\hat{p}_{\ell,m,t,\nu}^{(j-\ell-m)}}(-1, 1-y),
    \end{align*} where in the third line we use the fact that $\glslink{W}{W}\rvert_{x=0} = 1$, $\glslink{V}{V}\rvert_{x=0} = \frac{1}{1-y}$, and $\glslink{Q}{Q}\rvert_{x=0} = 1-y$. If we set $v := \frac{1}{1-y}$, then this becomes \begin{equation}
        \left.\glslink{z-rho}{z_\rho}\glslink{partial-x-i}{\partial_x^j}\glslink{calE-f}{\calE_{\glslink{I-rho}{I_\rho}}}\right\rvert_{x=0} = \sum_{\substack{2\ell+m+t+\nu=\glslink{abs-rho}{\lvert\rho\rvert} \\ \text{$\ell \ge 1$ or $t = 0$} \\ \ell+m \le j}} \glslink{stitch-ell-m-t-nu-binom}{\binom{\rho}{\ell,m,t,\nu}} \cdot \left(1-\frac{1}{v}\right)^{\ell+t+\nu} v^{1-\ell-m+\nu+2j} \glslink{hat-p-ell-m-t-nu-i}{\hat{p}_{\ell,m,t,\nu}^{(j-\ell-m)}}\left(-1, \frac{1}{v}\right), \label{E-from-p-hat-formula}
    \end{equation} which immediately proves the existence of the $\glslink{delta-f-j}{\delta_{\glslink{I-rho}{I_\rho}}^{(j)}}$. Observe that in order for the coefficient $\glslink{stitch-ell-m-t-nu-binom}{\binom{\rho}{\ell,m,t,\nu}}$ to be nonzero, we must have $\nu \le \glslink{m-j-rho}{m_1(\rho)}$ and $\ell+m+\nu \ge \glslink{len}{\len(\rho)}$: the first inequality follows because \glslink{pole}{poles} are precisely $1$-\glslink{cycle}{cycles}, and the second inequality follows because removing \glslink{wire}{wires} from a stitch can never decrease its total number of \glslink{cycle}{cycles} and \glslink{chain}{chains}. In particular, this tells us that the coefficient $\glslink{stitch-ell-m-t-nu-binom}{\binom{\rho}{\ell,m,t,\nu}}$ can only be nonzero if $\ell+m \ge \glslink{len}{\len(\rho)}-\glslink{m-j-rho}{m_1(\rho)}$, which in turn implies $j \ge \glslink{len}{\len(\rho)}-\glslink{m-j-rho}{m_1(\rho)}$; that is, the right hand side of \eqref{E-from-p-hat-formula} is zero whenever $j < \glslink{len}{\len(\rho)}-\glslink{m-j-rho}{m_1(\rho)}$.

    Next, suppose that $j \ge \glslink{len}{\len(\rho)}-\glslink{m-j-rho}{m_1(\rho)}$, and consider the Laurent polynomial \begin{equation}
        \left(1-\frac{1}{v}\right)^{\ell+t+\nu} v^{1-\ell-m+\nu+2j} \glslink{hat-p-ell-m-t-nu-i}{\hat{p}_{\ell,m,t,\nu}^{(j-\ell-m)}}\left(-1, \frac{1}{v}\right) \label{v-term}
    \end{equation} for some fixed $\ell,m,t,\nu \ge 0$ with $\ell+m \le j$ and $\glslink{stitch-ell-m-t-nu-binom}{\binom{\rho}{\ell,m,t,\nu}} > 0$. We start by examining the smallest power of $v$ that could appear in \eqref{v-term}. If $\ell = t = 0$, then by the $w$-degree constraint in Theorem~\ref{H-from-p-thm}(i), the lowest power of $v$ in \eqref{v-term} is at least \begin{equation*}
        (-\nu) + (1-m+\nu+2j) - (j-m) = 1+j,
    \end{equation*} which we verify is greater than or equal to \eqref{delta-smallest-v-power} in either case. On the other hand, if $\ell \ge 1$, then $\rho$ must have at least one part of size greater than $1$ (a \glslink{proper-chain}{proper chain} cannot be obtained by removing \glslink{wire}{wires} from a $1$-\glslink{cycle}{cycle}), so we are in the second case of \eqref{delta-smallest-v-power}. By the total degree constraint in Theorem~\ref{H-from-p-thm}(ii), the lowest possible power of $v$ in \eqref{v-term} is \begin{align}
        (-\ell-t-\nu) + (1-\ell-m+\nu+2j) - (j-1) &= 2-2\ell-m-t+j \notag \\
        &= 2-\glslink{abs-rho}{\lvert\rho\rvert}+\nu+j. \label{ell-nonzero-delta-v-degree-lower-bound}
    \end{align} Clearly, this is greater than or equal to $2-\glslink{abs-rho}{\lvert\rho\rvert}+j$. Moreover, by combining the relations $\ell+m+\nu \ge \glslink{len}{\len(\rho)}$ and $\ell+m \le j$, we find that $\nu \ge \glslink{len}{\len(\rho)}-j$, and so the expression in \eqref{ell-nonzero-delta-v-degree-lower-bound} is also greater than or equal to $2-\glslink{abs-rho}{\lvert\rho\rvert}+\glslink{len}{\len(\rho)}$. Thus, in all possible cases, the constraint in part (i) of the theorem statement is satisfied.

    Now, observe that the largest possible power of $v$ that can appear in \eqref{v-term} is \begin{align*}
        1-\ell-m+\nu+2j &\le 1 + 2\nu - \glslink{len}{\len(\rho)} + 2j \\
        &\le 1+2\glslink{m-j-rho}{m_1(\rho)}-\glslink{len}{\len(\rho)}+2j,
    \end{align*} where the first step uses the inequality $\ell+m+\nu \ge \glslink{len}{\len(\rho)}$ and the second step uses the inequality $\nu \le \glslink{m-j-rho}{m_1(\rho)}$. This implies part (ii) of the theorem.
    
    Finally, suppose that $\rho = 1^r$ for some $r \ge 0$, and suppose that $1-\ell-m+\nu+2j \ge 2j+r$. Then $1-\ell-m+\nu \ge r = 2\ell+m+t+\nu \ge 2\ell+m+\nu$, so $3\ell+2m \le 1$, implying that $\ell = m = 0$ and hence that $t = 0$ and $\nu = r$. In other words, the case where $\ell = m = t = 0$ and $\nu = r$ is the only possible case that can contribute to the coefficients of $v^{1+2j+r}$ and $v^{2j+r}$. Moreover, in this case, it is easy to verify that $\glslink{stitch-ell-m-t-nu-binom}{\binom{\rho}{\ell,m,t,\nu}} = \glslink{stitch-ell-m-t-nu-binom}{\binom{\rho}{0,0,0,r}} = 1$. Meanwhile, Theorem~\ref{H-from-p-thm}(iv) tells us that \begin{equation*}
        \glslink{p-ell-m-t-nu-i}{p_{0,0,0,r}^{(j)}}(-1, 0) = \frac{2^{2j}(\frac{r-1}{2}+j)!}{(\frac{r-1}{2})!} \qquad \text{and} \qquad (\partial_w \glslink{p-ell-m-t-nu-i}{p_{0,0,0,r}^{(j)})}(-1, 0) = -\frac{2^{2j-1}j(\frac{r-1}{2}+j)!}{(\frac{r-1}{2})!},
    \end{equation*} hence that \begin{equation*}
        \glslink{hat-p-ell-m-t-nu-i}{\hat{p}_{0,0,0,r}^{(j)}}(-1, 0) = \frac{2^{2j}r!(\frac{r-1}{2}+j)!}{(\frac{r-1}{2})!} \qquad \text{and} \qquad (\partial_w \glslink{hat-p-ell-m-t-nu-i}{\hat{p}_{0,0,0,r}^{(j)}})(-1, 0) = -\frac{2^{2j-1}r!j(\frac{r-1}{2}+j)!}{(\frac{r-1}{2})!}.
    \end{equation*} We thus compute that \begin{align*}
        &\left(1-\frac{1}{v}\right)^{\ell+t+\nu} v^{1-\ell-m+\nu+2j} \glslink{hat-p-ell-m-t-nu-i}{\hat{p}_{\ell,m,t,\nu}^{(j-\ell-m)}}\left(-1, \frac{1}{v}\right) \\
        &= \frac{2^{2j}r!(\frac{r-1}{2}+j)!}{(\frac{r-1}{2})!} \cdot \left(1-\frac{1}{v}\right)^r v^{1+r+2j} \left(1 - \frac{j}{2v} + O\left(\frac{1}{v^2}\right)\right) \\
        &= \frac{2^{2j}r!(\frac{r-1}{2}+j)!}{(\frac{r-1}{2})!} \cdot v^{1+r+2j} \left(1-\frac{r}{v} + O\left(\frac{1}{v^2}\right)\right) \left(1 - \frac{j}{2v} + O\left(\frac{1}{v^2}\right)\right) \\
        &= \frac{2^{2j}r!(\frac{r-1}{2}+j)!}{(\frac{r-1}{2})!} \cdot v^{1+r+2j} \left(1 - \left(\frac{j}{2}+r\right)\frac{1}{v} + O\left(\frac{1}{v^2}\right)\right) \\
        &= \frac{2^{2j}r!(\frac{r-1}{2}+j)!}{(\frac{r-1}{2})!} \cdot v^{1+r+2j} - \frac{2^{2j-1}r!(j+2r)(\frac{r-1}{2}+j)!}{(\frac{r-1}{2})!} \cdot v^{r+2j} + O\left(v^{r+2j-1}\right).
    \end{align*} From here, dividing both sides of \eqref{E-from-p-hat-formula} by $\glslink{z-rho}{z_\rho} = r!$ gives the claimed coefficients.
\end{proof}

Having proven our list of key results on $\left.\glslink{partial-x-i}{\partial_x^j} \glslink{calE-f}{\calE_f}\right\rvert_{x=0}$ for $f$ belonging to the basis $\glslink{I-rho}{I_\rho}$, it is straightforward to prove the general-$f$ analogue of Theorem~\ref{E-I-rho-from-delta-thm} stated earlier as Theorem~\ref{E-from-delta-thm}.

\begin{proof}[Proof of Theorem~\ref{E-from-delta-thm}]
    The existence of the Laurent polynomials $\glslink{delta-f-j}{\delta_f^{(j)}}$ follows immediately from Theorem~\ref{E-I-rho-from-delta-thm} together with the fact that $\glslink{calE-f}{\calE_f}$ is linear in $f$, and the degree constraints in parts (i) and (ii) of the theorem statement follow by considering the smallest possible value of the quantity in \eqref{delta-smallest-v-power} and the largest possible value of $1+2j+2\glslink{m-j-rho}{m_1(\rho)}-\glslink{len}{\len(\rho)}$ as $\rho$ ranges through all partitions of size at most $\glslink{character-polynomial-deg}{\deg(f)}$. For part (iii) of the theorem statement, we note that in order to have $1+2j+2\glslink{m-j-rho}{m_1(\rho)}-\glslink{len}{\len(\rho)} \ge 2j+\glslink{character-polynomial-deg}{\deg(f)}$ for a partition $\rho$ of size at most $\glslink{character-polynomial-deg}{\deg(f)}$, we must have $\glslink{m-j-rho}{m_1(\rho)}-(\glslink{len}{\len(\rho)}-\glslink{m-j-rho}{m_1(\rho)})+1 \ge \glslink{character-polynomial-deg}{\deg(f)}$, which is impossible unless $\rho = 1^{\glslink{character-polynomial-deg}{\deg(f)}}$ or $\rho = 1^{\glslink{character-polynomial-deg}{\deg(f)}-1}$. The claim then follows immediately by Theorem~\ref{E-I-rho-from-delta-thm}(iii).
\end{proof}

\noindent
From here, Theorems~\ref{main-thm-2} and \ref{main-thm-3}(ii) quickly follow.

\begin{proof}[Proof of Theorem~\ref{main-thm-2}]
    Let $j \ge 0$, and use Theorem~\ref{E-from-delta-thm} to write \begin{equation}
        \sum_{n\ge j} n! \glslink{E-f-n-k}{E_{\glslink{X-lam}{X^\lambda}}(n, n-j)} \cdot y^{n-j} = \left.\glslink{partial-x-i}{\partial_x^j}\glslink{calE-f}{\calE_{\glslink{X-lam}{X^\lambda}}}\right|_{x=0} = \glslink{delta-f-j}{\delta_{\glslink{X-lam}{X^\lambda}}^{(j)}}\left(\frac{1}{1-y}\right). \label{E-polynomial-from-delta}
    \end{equation} By Theorem~\ref{E-from-delta-thm}(i-ii), the smallest power of $v$ appearing in $\glslink{delta-f-j}{\delta_{\glslink{X-lam}{X^\lambda}}^{(j)}}$ is at least $2+j-\glslink{abs-rho}{\lvert\lambda\rvert}$ (using the assumption $\glslink{abs-rho}{\lvert\lambda\rvert} \ge 1$), while the largest power of $v$ appearing in $\glslink{delta-f-j}{\delta_{\glslink{X-lam}{X^\lambda}}^{(j)}}$ is at most $1+2j+\glslink{abs-rho}{\lvert\lambda\rvert}$. The first statement tells us that the coefficient of $y^k$ in $\glslink{delta-f-j}{\delta_{\glslink{X-lam}{X^\lambda}}^{(j)}}\left(\frac{1}{1-y}\right)$ agrees with a polynomial for all $k \ge \glslink{abs-rho}{\lvert\lambda\rvert}-j-1$, whereas the second statement tells us that the degree of this polynomial is at most $2j+\glslink{abs-rho}{\lvert\lambda\rvert}$. Substituting $k = n-j$ then gives the claim.
\end{proof}

\begin{proof}[Proof of Theorem~\ref{main-thm-3}(ii)]
    Recall the definition of the \glslink{character-polynomial}{character polynomials} $\glslink{X-lam}{X^\lambda}$ from Theorem~\ref{macdonald-thm}. Based on this definition, it is clear that the coefficient of $\glslink{I-rho}{I_{1^{\glslink{abs-rho}{\lvert\lambda\rvert}}}}$ in the $\glslink{I-rho}{I_\rho}$-expansion of $\glslink{X-lam}{X^\lambda}$ is $\chi_{1^{\glslink{abs-rho}{\lvert\lambda\rvert}}}^\lambda = \glslink{chi-lambda}{\chi^\lambda}(\glslink{id-k}{\id_{\glslink{abs-rho}{\lvert\lambda\rvert}}})$, whereas the coefficient of $\glslink{I-rho}{I_{1^{\glslink{abs-rho}{\lvert\lambda\rvert-1}}}}$ is given by \begin{equation*}
        -\sum_\mu \glslink{chi-rho-lambda}{\chi_{1^{\glslink{abs-rho}{\lvert\lambda\rvert-1}}}^\mu} = -\sum_\mu \glslink{chi-lambda}{\chi^\mu}(\glslink{id-k}{\id_{\glslink{abs-rho}{\lvert\lambda\rvert}-1}}) = -\glslink{chi-lambda}{\chi^\lambda}(\glslink{id-k}{\id_{\glslink{abs-rho}{\lvert\lambda\rvert}}}),
    \end{equation*} where the sum is taken over partitions $\mu$ obtained by removing a single box from $\lambda$ and the second equality follows by Young's branching rule. From this, we deduce by Theorem~\ref{E-from-delta-thm}(iii) that the coefficients of $v^{1+2j+\glslink{abs-rho}{\lvert\lambda\rvert}}$ and $v^{2j+\glslink{abs-rho}{\lvert\lambda\rvert}}$ in $\glslink{delta-f-j}{\delta_{\glslink{X-lam}{X^\lambda}}^{(j)}}$ are given by \begin{equation*}
        \glslink{chi-lambda}{\chi^\lambda}(\glslink{id-k}{\id_{\glslink{abs-rho}{\lvert\lambda\rvert}}}) \cdot \frac{2^{2j}(\frac{\glslink{abs-rho}{\lvert\lambda\rvert}-1}{2}+j)!}{(\frac{\glslink{abs-rho}{\lvert\lambda\rvert}-1}{2})!} \qquad \text{and} \qquad -\glslink{chi-lambda}{\chi^\lambda}(\glslink{id-k}{\id_{\glslink{abs-rho}{\lvert\lambda\rvert}}}) \left(\frac{2^{2j}(\frac{\glslink{abs-rho}{\lvert\lambda\rvert}-2}{2}+j)!}{(\frac{\glslink{abs-rho}{\lvert\lambda\rvert}-2}{2})!} + \frac{2^{2j-1}(j+2\glslink{abs-rho}{\lvert\lambda\rvert})(\frac{\glslink{abs-rho}{\lvert\lambda\rvert}-1}{2}+j)!}{(\frac{\glslink{abs-rho}{\lvert\lambda\rvert}-1}{2})!}\right),
    \end{equation*} respectively. If we refer to the above expressions as $C_0$ and $C_1$, respectively, then it follows from \eqref{E-polynomial-from-delta} that $n!\glslink{E-f-n-k}{E_{\glslink{X-lam}{X^\lambda}}(n, n-j)}$ is given by \begin{align*}
        C_0 \cdot \glslink{binomial}{\binom{n+j+\glslink{abs-rho}{\lvert\lambda\rvert}}{2j+\glslink{abs-rho}{\lvert\lambda\rvert}}} + C_1 \cdot \glslink{binomial}{\binom{n+j+\glslink{abs-rho}{\lvert\lambda\rvert}-1}{2j+\glslink{abs-rho}{\lvert\lambda\rvert}-1}} + O(n^{2j+\glslink{abs-rho}{\lvert\lambda\rvert}-2}).
    \end{align*} Applying the well-known approximation \begin{align*}
        \glslink{binomial}{\binom{n+a}{b}} = \frac{n^b}{b!} + \frac{b(1+2a-b)n^{b-1}}{2b!} + O\left(n^{b-2}\right)
    \end{align*} (see \cite[\href{https://dlmf.nist.gov/5.11.E13}{(5.11.13)} and \href{https://dlmf.nist.gov/5.11.E15}{(5.11.15)}]{dlmf}) and simplifying then gives the claimed formulas.
\end{proof}

\subsection{Explicit formulas for \texorpdfstring{$\glslink{partial-x-i}{\partial_x^{\glslink{abs-rho}{\lvert\rho\rvert}-1}}\glslink{calE-f}{\calE_{\glslink{J-rho}{J_\rho}}}$}{∂|ρ|−1xEJρ}}

As explained in Section~\ref{proof-overview-section}, our first step towards computing the functions $\glslink{partial-x-i}{\partial_x^i}\glslink{calE-f}{\calE_f}$ for sufficiently large $i$ will be to prove that these functions are given by polynomials in $\glslink{U}{U}$ and $\glslink{V}{V}$.

\begin{proof}[Proof of Lemma~\ref{E-from-epsilon-hat-lem}]
    Let $\rho$ be any partition, and suppose that $i \ge \max(0, \glslink{abs-rho}{\lvert\rho\rvert}-1)$. Then by Corollary~\ref{E-from-H-ell-m-t-nu-cor}, we may write \begin{equation*}
        \glslink{z-rho}{z_\rho}\glslink{partial-x-i}{\partial_x^i}\glslink{calE-f}{\calE_{\glslink{I-rho}{I_\rho}}} = \sum_{\substack{2\ell+m+t+\nu=\glslink{abs-rho}{\lvert\rho\rvert} \\ \text{$\ell \ge 1$ or $t = 0$}}} \glslink{stitch-ell-m-t-nu-binom}{\binom{\rho}{\ell,m,t,\nu}} \cdot \frac{\ell!m!\nu!}{\glslink{binomial}{\binom{\ell+t-1}{t}}} \cdot y^{\ell+t+\nu} \glslink{partial-x-i}{\partial_x^{i-\ell-m}} \glslink{H-ell-m-t-nu}{H_{\ell,m,t,\nu}}.
    \end{equation*} Now, the smallest possible value of $i-\ell-m$ in the sum is $-1$, achieved only when $i = \glslink{abs-rho}{\lvert\rho\rvert}-1$, $\ell = t = \nu = 0$, and $m = \glslink{abs-rho}{\lvert\rho\rvert}$. We readily verify that $\glslink{stitch-ell-m-t-nu-binom}{\binom{\rho}{0,\glslink{abs-rho}{\lvert\rho\rvert},0,0}} = 1$, so applying Theorem~\ref{H-from-p-thm} thus gives \begin{align*}
        \glslink{z-rho}{z_\rho}\glslink{partial-x-i}{\partial_x^i}\glslink{calE-f}{\calE_{\glslink{I-rho}{I_\rho}}} &= \sum_{\substack{2\ell+m+t+\nu=\glslink{abs-rho}{\lvert\rho\rvert} \\ \text{$\ell \ge 1$ or $t = 0$}}} \glslink{stitch-ell-m-t-nu-binom}{\binom{\rho}{\ell,m,t,\nu}} \cdot y^{\ell+t+\nu} \glslink{W}{W}^{\ell+m} \glslink{V}{V}^{i-2\ell-m-t} \glslink{Q}{Q}^{-1-\ell-t-\nu-i} \glslink{hat-p-ell-m-t-nu-i}{\hat{p}_{\ell,m,t,\nu}^{(i-\ell-m)}}\left(-\frac{1}{\glslink{W}{W}}, \frac{1}{\glslink{V}{V}}\right) \\
        &\quad - \glslink{W}{W}^{\glslink{abs-rho}{\lvert\rho\rvert}} \glslink{V}{V}^{-1} \glslink{Q}{Q}^{-\glslink{abs-rho}{\lvert\rho\rvert}} \cdot (\glslink{abs-rho}{\lvert\rho\rvert}-1)! \cdot \glslink{indicator}{\mathbbm{1}_{i=\glslink{abs-rho}{\lvert\rho\rvert}-1}},
    \end{align*} where $\glslink{hat-p-ell-m-t-nu-i}{\hat{p}_{\ell,m,t,\nu}^{(i)}}$ is defined as in \eqref{p-hat-def} and $\glslink{indicator}{\mathbbm{1}_{i=\glslink{abs-rho}{\lvert\rho\rvert}-1}}$ is $1$ if $i = \glslink{abs-rho}{\lvert\rho\rvert}-1$ and $0$ otherwise. (If $\glslink{abs-rho}{\lvert\rho\rvert} = 0$, we simply omit the last term.) If we apply the relations \begin{equation*}
        \glslink{W}{W} = \frac{\glslink{U}{U}\glslink{V}{V}}{\glslink{U}{U}+\glslink{V}{V}-1}, \qquad \glslink{Q}{Q} = \frac{1}{\glslink{U}{U}+\glslink{V}{V}-1}, \qquad y = \frac{\glslink{V}{V}(\glslink{V}{V}-1)}{(\glslink{U}{U}+\glslink{V}{V}-1)^2}
    \end{equation*} and simplify, then this in turn becomes \begin{align*}
        &(\glslink{U}{U}+\glslink{V}{V}-1)^{1-\glslink{abs-rho}{\lvert\rho\rvert}+i} \sum_{\substack{2\ell+m+t+\nu=\glslink{abs-rho}{\lvert\rho\rvert} \\ \text{$\ell \ge 1$ or $t = 0$}}} \glslink{stitch-ell-m-t-nu-binom}{\binom{\rho}{\ell,m,t,\nu}} \cdot \glslink{U}{U}^{\ell+m} \glslink{V}{V}^{i+\nu} (\glslink{V}{V}-1)^{\ell+t+\nu} \glslink{hat-p-ell-m-t-nu-i}{\hat{p}_{\ell,m,t,\nu}^{(i-\ell-m)}}\left(\frac{1-\glslink{U}{U}-\glslink{V}{V}}{\glslink{U}{U}\glslink{V}{V}}, \frac{1}{\glslink{V}{V}}\right) \\
        &\quad - (\glslink{abs-rho}{\lvert\rho\rvert}-1)! \cdot \glslink{U}{U}^{\glslink{abs-rho}{\lvert\rho\rvert}} \glslink{V}{V}^{\glslink{abs-rho}{\lvert\rho\rvert}-1} \cdot \glslink{indicator}{\mathbbm{1}_{i=\glslink{abs-rho}{\lvert\rho\rvert}-1}}.
    \end{align*} Now fix $\ell, m, t, \nu \ge 0$ such that $2\ell+m+t+\nu=\glslink{abs-rho}{\lvert\rho\rvert}$ and either $\ell \ge 1$ or $t = 0$, and consider the expression \begin{equation}
        \glslink{U}{U}^{\ell+m} \glslink{V}{V}^{i+\nu} (\glslink{V}{V}-1)^{\ell+t+\nu} \glslink{hat-p-ell-m-t-nu-i}{\hat{p}_{\ell,m,t,\nu}^{(i-\ell-m)}}\left(\frac{1-\glslink{U}{U}-\glslink{V}{V}}{\glslink{U}{U}\glslink{V}{V}}, \frac{1}{\glslink{V}{V}}\right). \label{eps-hat-term}
    \end{equation} If $i-\ell-m \ge 0$, then Theorem~\ref{H-from-p-thm}(i-ii) imply that \eqref{eps-hat-term} is a polynomial in $\glslink{U}{U}$ and $\glslink{V}{V}$. Otherwise, we have $\ell = t = \nu = 0$ and $m = \glslink{abs-rho}{\lvert\rho\rvert}$, at which point Theorem~\ref{H-from-p-thm}(iii) allows us to again conclude that \eqref{eps-hat-term} is a polynomial in $\glslink{U}{U}$ and $\glslink{V}{V}$. In other words, the expression \begin{align}
        \glslink{hat-varepsilon-f-i}{\hat{\varepsilon}_{\glslink{I-rho}{I_\rho}}^{(i)}}(u, v) &:= \frac{1}{\glslink{z-rho}{z_\rho}} \cdot (u+v-1)^{1-\glslink{abs-rho}{\lvert\rho\rvert}+i} \!\!\!\!\!\!\!\!\sum_{\substack{2\ell+m+t+\nu=\glslink{abs-rho}{\lvert\rho\rvert} \\ \text{$\ell \ge 1$ or $t = 0$}}} \glslink{stitch-ell-m-t-nu-binom}{\binom{\rho}{\ell,m,t,\nu}} u^{\ell+m} v^{i+\nu} (v-1)^{\ell+t+\nu} \glslink{hat-p-ell-m-t-nu-i}{\hat{p}_{\ell,m,t,\nu}^{(i-\ell-m)}}\left(\tfrac{1-u-v}{uv}, \tfrac{1}{v}\right) \notag \\
        &\quad - \frac{(\glslink{abs-rho}{\lvert\rho\rvert}-1)!}{\glslink{z-rho}{z_\rho}} \cdot u^{\glslink{abs-rho}{\lvert\rho\rvert}} v^{\glslink{abs-rho}{\lvert\rho\rvert}-1} \cdot \glslink{indicator}{\mathbbm{1}_{i=\glslink{abs-rho}{\lvert\rho\rvert}-1}}. \label{epsilon-hat-formula}
    \end{align} is a polynomial in $u$ and $v$ such that \begin{equation*}
        \glslink{partial-x-i}{\partial_x^i} \glslink{calE-f}{\calE_{\glslink{I-rho}{I_\rho}}} = \glslink{hat-varepsilon-f-i}{\hat{\varepsilon}_{\glslink{I-rho}{I_\rho}}^{(i)}}(\glslink{U}{U}, \glslink{V}{V}).
    \end{equation*} From here, the result for general $f \in \glslink{Delta-Q}{\Delta_\Q}$ follows by linearity.
\end{proof}

At this point, in order to relate what Lemma~\ref{E-from-epsilon-hat-lem} implies about the values $\glslink{E-f-n-k}{E_f(n, k)}$ to what we already know about these values in the case where $f = \glslink{X-lam}{X^\lambda}$, it is useful to have an explicit formula for the series coefficients of the monomials $\glslink{U}{U}^a \glslink{V}{V}^b$ appearing in the expressions $\glslink{hat-varepsilon-f-i}{\hat{\varepsilon}_f^{(i)}}(\glslink{U}{U}, \glslink{V}{V})$. As explained in Section~\ref{proof-overview-section}, such explicit formulas are provided by Theorem~\ref{q-a-b-r-formula-thm}.

\begin{proof}[Proof of Theorem~\ref{q-a-b-r-formula-thm}]
    By repeatedly differentiating with respect to $x$ and $y$, it is easy to see that the coefficient of $x^j y^k$ in $\glslink{U}{U}^a \glslink{V}{V}^b \glslink{Q}{Q}^{-1-r}$ is a polynomial in $r$ of degree at most $j+k$. Thus, throughout the proof, we need not concern ourselves with the subtlety of how to interpret factorials of negative integers---as long as an expression for $\glslink{q-a-b-r-j-k}{q_{a,b}^{(r)}(j, k)}$ is well-defined for $r$ outside of some discrete subset of $\R$, we can determine the correct value of $\glslink{q-a-b-r-j-k}{q_{a,b}^{(r_0)}(j, k)}$ for any $r_0$ where the expression does not make sense by taking the limit as $r \to r_0$.

    We proceed by induction on $a + b$. In the case $a = b = 0$, the claimed identity reduces to \begin{equation*}
        \glslink{q-a-b-r-j-k}{q_{0,0}^{(r)}(j, k)} = \frac{(r+j+k)!(\frac{r}{2}-1)!(\frac{r}{2}+j+k)!}{2j!k!(r-1)!(\frac{r}{2}+j)!(\frac{r}{2}+k)!} = \frac{(j+k+r)!}{j!} \cdot \frac{2^k}{(r-1)!!(2k+r)!!} \glslink{binomial}{\binom{j+k+\frac{r}{2}}{k}}.
    \end{equation*} This is equivalent to Lemma~2.6 from \cite{gaetz-pierson}. Indeed, if $F = (1-2x^2(1-y^2)+x^4(1-y^2)^2)^{-\frac{1}{2}}$ is the function defined in \cite{gaetz-pierson}, our function $\glslink{Q}{Q}$ is related to $F$ by $F(x, y) = \glslink{Q}{Q}(x^2, x^2y^2)^{-1}$. Meanwhile, our coefficients $\glslink{q-a-b-r-j-k}{q_{0,0}^{(r)}(j, k)}$ are related to Gaetz and Pierson's quantities $E(j, k, r)$ by $\glslink{q-a-b-r-j-k}{q_{0,0}^{(r)}(j, k)} = \frac{(j+k+r)!}{j!} E(j+k+r, k+r, r)$.

    Now suppose that $a + b \ge 1$; by symmetry, we may assume without loss of generality that $b \ge 1$. For each $r \in \R \setminus \Z$, define \begin{equation*}
        P_{\alpha,\beta}^{(r)} := \sum_{j, k \ge 0} \frac{(r+j+k)!(\frac{r+\alpha+\beta}{2}+j+k)!}{2j!k!(\frac{r-\alpha+\beta}{2}+j)!(\frac{r+\alpha-\beta}{2}+k)!} \cdot x^j y^k \in \R[[x, y]]
    \end{equation*} for all $\alpha, \beta \in \Z$ and \begin{equation*}
        c_{\alpha,\beta}^{a,b}(r) := \glslink{binomial}{\binom{a}{\alpha}} \glslink{binomial}{\binom{b}{\beta}} \frac{(\frac{r-\alpha+\beta}{2}+a-1)!(\frac{r+\alpha-\beta}{2}+b-1)!}{(r+a+b-1)!(\frac{r+\alpha+\beta}{2}-1)!} \in \R
    \end{equation*} for all $a, b, \alpha, \beta \in \Z$. It suffices to prove that \begin{equation*}
        \glslink{U}{U}^a \glslink{V}{V}^b \glslink{Q}{Q}^{-1-r} = \sum_{\alpha, \beta \in \Z} c_{\alpha,\beta}^{a,b}(r) \cdot P_{\alpha,\beta}^{(r)}
    \end{equation*} for all $j, k \ge 0$ and $r \in \R \setminus \Z$. For this purpose, start by noting that \begin{equation*}
        \partial_x(\glslink{U}{U}^a \glslink{V}{V}^{b-1} \glslink{Q}{Q}^{-r}) = (2(a+b+r-1) \glslink{U}{U}^a \glslink{V}{V}^b - (a+2b+r-2) \glslink{U}{U}^a \glslink{V}{V}^{b-1} - a \glslink{U}{U}^{a-1} \glslink{V}{V}^b + a \glslink{U}{U}^{a-1} \glslink{V}{V}^{b-1}) \glslink{Q}{Q}^{-1-r},
    \end{equation*} or equivalently, \begin{align*}
        2(a+b+r-1) \glslink{U}{U}^a \glslink{V}{V}^b \glslink{Q}{Q}^{-1-r} &= \partial_x(\glslink{U}{U}^a \glslink{V}{V}^{b-1} \glslink{Q}{Q}^{-r}) + (a+2b+r-2) \glslink{U}{U}^a \glslink{V}{V}^{b-1} \glslink{Q}{Q}^{-1-r} \\
        &\quad + a \glslink{U}{U}^{a-1} \glslink{V}{V}^b \glslink{Q}{Q}^{-1-r} - a \glslink{U}{U}^{a-1} \glslink{V}{V}^{b-1} \glslink{Q}{Q}^{-1-r}.
    \end{align*} Now, by the inductive hypothesis, we may write \begin{alignat*}{2}
        \glslink{U}{U}^a \glslink{V}{V}^{b-1} \glslink{Q}{Q}^{-r} &= \sum_{\alpha, \beta \in \Z} c_{\alpha,\beta}^{a,b-1}(r-1) \cdot P_{\alpha,\beta}^{(r-1)}, & \qquad \glslink{U}{U}^a \glslink{V}{V}^{b-1} \glslink{Q}{Q}^{-1-r} &= \sum_{\alpha, \beta \in \Z} c_{\alpha,\beta}^{a,b-1}(r) \cdot P_{\alpha,\beta}^{(r)}, \\
        a \glslink{U}{U}^{a-1} \glslink{V}{V}^b \glslink{Q}{Q}^{-1-r} &= a\sum_{\alpha, \beta \in \Z} c_{\alpha,\beta}^{a-1,b}(r) \cdot P_{\alpha,\beta}^{(r)}, &\qquad a U^{a-1} \glslink{V}{V}^{b-1} \glslink{Q}{Q}^{-1-r} &= a\sum_{\alpha, \beta \in \Z} c_{\alpha,\beta}^{a-1,b-1}(r) \cdot P_{\alpha,\beta}^{(r)},
    \end{alignat*} where we note that the bottom two identities are vacuously true if $a = 0$. Moreover, observing that \begin{equation*}
        \partial_x P_{\alpha,\beta}^{(r)} = \sum_{j, k \ge 0} \frac{(r+j+k+1)!(\frac{r+\alpha+\beta}{2}+j+k+1)!}{2j!k!(\frac{r-\alpha+\beta}{2}+j+1)!(\frac{r+\alpha-\beta}{2}+k)!} \cdot x^j y^k = P_{\alpha,\beta+1}^{(r+1)}
    \end{equation*} for all $\alpha, \beta \in \Z$, we may use the first of the four identities to obtain \begin{equation*}
        \partial_x(\glslink{U}{U}^a \glslink{V}{V}^{b-1} \glslink{Q}{Q}^{-r}) = \sum_{\alpha, \beta \in \Z} c_{\alpha,\beta}^{a,b-1}(r-1) \cdot P_{\alpha,\beta+1}^{(r)} = \sum_{\alpha, \beta \in \Z} c_{\alpha,\beta-1}^{a,b-1}(r-1) \cdot P_{\alpha,\beta}^{(r)}.
    \end{equation*} Putting everything together, we thus find that \begin{align*}
        &2(a+b+r-1) \glslink{U}{U}^a \glslink{V}{V}^b \glslink{Q}{Q}^{-1-r} \\
        &= \sum_{\alpha, \beta \in \Z} \Big(c_{\alpha,\beta-1}^{a,b-1}(r-1) + (a+2b+r-1)c_{\alpha,\beta}^{a,b-1}(r) + ac_{\alpha,\beta}^{a-1,b}(r) - ac_{\alpha,\beta}^{a-1,b-1}(r)\Big) P_{\alpha,\beta}^{(r)}.
    \end{align*} It is then a tedious but routine exercise to show that \begin{equation*}
        c_{\alpha,\beta-1}^{a,b-1}(r-1) + (a+2b+r-1)c_{\alpha,\beta}^{a,b-1}(r) + ac_{\alpha,\beta}^{a-1,b}(r) - ac_{\alpha,\beta}^{a-1,b-1}(r) = 2(a+b+r-1) c_{\alpha,\beta}^{a,b}(r),
    \end{equation*} at which point the claimed formula for $\glslink{U}{U}^a \glslink{V}{V}^b \glslink{Q}{Q}^{-1-r}$ follows.
\end{proof}

Using Theorem~\ref{q-a-b-r-formula-thm}, we can then show that any monomials of the form $\glslink{U}{U}^a \glslink{V}{V}^b$ with $a \ge 1$ appearing in $\glslink{hat-varepsilon-f-i}{\hat{\varepsilon}_{\glslink{X-lam}{X^\lambda}}}(\glslink{U}{U}, \glslink{V}{V})$ would force the degrees of the polynomials $\glslink{a-sigma-lambda}{a_{\glslink{id-k}{\id_k}}^\lambda}(n)$ to be larger than required by Theorem~\ref{gp-poly-thm}.

\begin{proof}[Proof of Lemma~\ref{E-from-epsilon-lem}]
    Let $\lambda$ be any partition of size at least $1$, and use Lemma~\ref{E-from-epsilon-hat-lem} to write \begin{equation}
        \glslink{partial-x-i}{\partial_x^{\glslink{abs-rho}{\lvert\lambda\rvert}-1}}\glslink{calE-f}{\calE_{\glslink{X-lam}{X^\lambda}}} = \sum_{a, b \ge 0} \epsilon_{a,b} \cdot U^a \glslink{V}{V}^b \label{E-from-U-V-powers}
    \end{equation} for some coefficients $\epsilon_{a,b} \in \Q$ finitely many of which are nonzero. For each $\alpha, \beta, j, k \ge 0$ with $j \ge \frac{1+\alpha-\beta}{2}$ and $k \ge \frac{1+\beta-\alpha}{2}$, let \begin{equation*}
        g_{\alpha,\beta}(j, k) := \frac{(j+k-1)!(j+k+\frac{\alpha+\beta-1}{2})!}{2j!k!(j+\frac{\beta-\alpha-1}{2})!(k+\frac{\alpha-\beta-1}{2})!} \in \R.
    \end{equation*} Then by Theorem~\ref{q-a-b-r-formula-thm}, the coefficient of $x^j y^k$ in $\glslink{U}{U}^a \glslink{V}{V}^b$ is given by \begin{equation}
        \sum_{\substack{0 \le \alpha \le a \\ 0 \le \beta \le b}} c_{\alpha,\beta}^{a,b} \cdot g_{\alpha,\beta}(j,k) \label{g-alpha-beta-sum}
    \end{equation} for all $j \ge \frac{1+\alpha-\beta}{2}$ and $k \ge \frac{1+\beta-\alpha}{2}$, where \begin{equation*}
        c_{\alpha,\beta}^{a,b} := \glslink{binomial}{\binom{a}{\alpha}} \glslink{binomial}{\binom{b}{\beta}} \lim_{r \to -1}\left(\frac{(\frac{r-\alpha+\beta}{2}+a-1)!(\frac{r+\alpha-\beta}{2}+b-1)!}{(r+a+b-1)!(\frac{r+\alpha+\beta}{2}-1)!}\right) \in \R.
    \end{equation*} In particular, for $a+b \ge 1$, we note that the coefficient \begin{equation*}
        c_{a,b}^{a,b} = \lim_{r \to -1} \left(\frac{(\frac{r+a+b}{2}-1)!}{(r+a+b-1)!}\right) = \lim_{r \to -1} \left(\frac{(-\frac{1}{2})!}{2^{a+b+r-1}(\frac{a+b+r-1}{2})!}\right) = \frac{(-\frac{1}{2})!}{2^{a+b-2}(\frac{a+b}{2}-1)!}
    \end{equation*} is nonzero. Combining \eqref{E-from-U-V-powers} and \eqref{g-alpha-beta-sum}, we deduce that \begin{equation*}
        \frac{(j+k+\glslink{abs-rho}{\lvert\lambda\rvert}-1)!}{j!} \cdot \glslink{E-f-n-k}{E_{\glslink{X-lam}{X^\lambda}}(j+k+\glslink{abs-rho}{\lvert\lambda\rvert}-1, k)} = \sum_{\alpha, \beta \ge 0} C_{\alpha,\beta} \cdot g_{\alpha,\beta}(j, k)
    \end{equation*} for all sufficiently large $j$ and $k$, where \begin{equation}
        C_{\alpha,\beta} := \sum_{\substack{a \ge \alpha \\ b \ge \beta}} \epsilon_{a,b} \cdot c_{\alpha,\beta}^{a,b}. \label{C-from-eps}
    \end{equation} This in turn implies \begin{equation}
        \frac{n!}{(n-k-\glslink{abs-rho}{\lvert\lambda\rvert}+1)!} \cdot \glslink{a-sigma-lambda}{a_{\glslink{id-k}{\id_k}}^\lambda}(n) = \sum_{\alpha, \beta \ge 0} C_{\alpha,\beta} \cdot g_{\alpha,\beta}(n-k-\glslink{abs-rho}{\lvert\lambda\rvert}+1, k) \label{a-id-k-from-g}
    \end{equation} as polynomials in $n$ for all sufficiently large $k$. Now, observe that by Theorem~\ref{gp-poly-thm}, the left hand side of \eqref{a-id-k-from-g} has degree at most $k-\glslink{abs-rho}{\lvert\lambda\rvert}+k+\glslink{abs-rho}{\lvert\lambda\rvert}-1 = 2k-1$. Meanwhile, observe that for any fixed $\alpha, \beta, k \ge 0$, $g_{\alpha,\beta}(j, k)$ is a polynomial in $j$ of degree $2k-1+\alpha$ with leading coefficient $\frac{1}{2k!(k+\frac{\alpha-\beta-1}{2})!}$. Since \begin{equation*}
        \lim_{k\to\infty}\frac{(k+\frac{\alpha-\beta'-1}{2})!}{(k+\frac{\alpha-\beta-1}{2})!} = \lim_{k\to\infty} O(k^{\frac{\beta-\beta'}{2}}) = 0
    \end{equation*} whenever $\beta' > \beta$, it follows that if $\alpha_0$ is the largest integer such that $C_{\alpha_0, \beta}$ is nonzero for some $\beta \ge 0$ and $\beta_0$ is the largest integer such that $C_{\alpha_0, \beta_0}$ is nonzero, then the right hand side of \eqref{a-id-k-from-g} is a polynomial in $n$ of degree $2k-1+\alpha_0$ for all sufficiently large $k$. By comparing with the degree of the left hand side, we thus conclude that $\alpha_0 = 0$, hence that $C_{\alpha,\beta} = 0$ for all $\alpha \ge 1$ and $\beta \ge 0$. This in turn implies $\epsilon_{a,b} = 0$ for all $a \ge 1$ and $b \ge 0$---otherwise, we could choose an $a_0 \ge 1$ and $b_0 \ge 0$ such that $\epsilon_{a_0,b_0} \ne 0$ and $\epsilon_{a,b} = 0$ for all $a \ge a_0$ and $b \ge b_0$ with $(a, b) \ne (a_0, b_0)$, at which point \eqref{C-from-eps} would imply $C_{a_0,b_0} = \epsilon_{a_0,b_0} \cdot c_{a_0,b_0}^{a_0,b_0} \ne 0$. In other words, $\glslink{partial-x-i}{\partial_x^{\glslink{abs-rho}{\lvert\lambda\rvert}-1}} \glslink{calE-f}{\calE_{\glslink{X-lam}{X^\lambda}}}$ is given by a polynomial in $\glslink{V}{V}$.
\end{proof}

Having proven that $\glslink{partial-x-i}{\partial_x^{\glslink{character-polynomial-deg}{\deg(f)}-1}}\glslink{calE-f}{\calE_f}$ is given by a polynomial $\glslink{varepsilon-f}{\varepsilon_f}(\glslink{V}{V})$ in $\glslink{V}{V}$ whenever $f$ is \glslink{pure-character-polynomial}{pure}, we can now reuse Theorem~\ref{E-I-rho-from-delta-thm} to deduce a variety of properties of the polynomials $\glslink{varepsilon-f}{\varepsilon_{\glslink{J-rho}{J_\rho}}}$.

\begin{theorem} \label{eps-J-rho-properties-thm}
    Let $\rho$ be a partition of size at least $1$. Then the polynomial $\glslink{varepsilon-f}{\varepsilon_{\glslink{J-rho}{J_\rho}}}$ from Lemma~\ref{E-from-epsilon-lem} satisfies the following properties: \begin{itemize}
        \item[(i)] $\glslink{total-deg}{\deg(\glslink{varepsilon-f}{\varepsilon_{\glslink{J-rho}{J_\rho}}})} \le 2\glslink{abs-rho}{\lvert\rho\rvert}+\glslink{m-j-rho}{m_1(\rho)}-1$;
        \item[(ii)] $\glslink{varepsilon-f}{\varepsilon_{\glslink{J-rho}{J_\rho}}}(0) = 0$;
        \item[(iii)] $(v-1)^{\glslink{abs-rho}{\lvert\rho\rvert}-1} \mid \glslink{varepsilon-f}{\varepsilon_{\glslink{J-rho}{J_\rho}}}'$.
        \item[(iv)] The coefficient of $v^{2\glslink{abs-rho}{\lvert\rho\rvert}+\glslink{m-j-rho}{m_1(\rho)}-1}$ in $\glslink{varepsilon-f}{\varepsilon_{\glslink{J-rho}{J_\rho}}}$ is given by \begin{equation*}
            (-1)^{\glslink{len}{\len(\rho)}-\glslink{m-j-rho}{m_1(\rho)}} \cdot \frac{2^{2\glslink{abs-rho}{\lvert\rho\rvert}-2}\glslink{m-j-rho}{m_1(\rho)}!(\frac{2\glslink{abs-rho}{\lvert\rho\rvert}+\glslink{m-j-rho}{m_1(\rho)}-3}{2})!}{\glslink{z-rho}{z_\rho}(\frac{\glslink{m-j-rho}{m_1(\rho)}-1}{2})!}.
        \end{equation*} Moreover, if $\rho = 1^r$ for some $r \ge 1$, then the coefficient of $v^{2\glslink{abs-rho}{\lvert\rho\rvert}+\glslink{m-j-rho}{m_1(\rho)}-2} = v^{3r-2}$ in $\glslink{varepsilon-f}{\varepsilon_{\glslink{J-rho}{J_\rho}}}$ is given by \begin{equation*}
            -2^{2r-2}\left(\frac{(\frac{3r-1}{2})!}{(\frac{r-1}{2})!} + \frac{(\frac{3r-4}{2})!}{(\frac{r-2}{2})!}\right).
        \end{equation*}
    \end{itemize}
\end{theorem}

\begin{proof}
    Observe that by combining Theorem~\ref{E-from-delta-thm} and Lemma~\ref{E-from-epsilon-lem}, we may write \begin{equation*}
        \glslink{delta-f-j}{\delta_{\glslink{J-rho}{J_\rho}}^{(\glslink{abs-rho}{\lvert\rho\rvert}-1)}}\left(\frac{1}{1-y}\right) = \left.\glslink{partial-x-i}{\partial_x^{\glslink{abs-rho}{\lvert\rho\rvert}-1}}\glslink{calE-f}{\calE_{\glslink{J-rho}{J_\rho}}}\right\rvert_{x=0} = \left.\glslink{varepsilon-f}{\varepsilon_{\glslink{J-rho}{J_\rho}}}\left(\glslink{V}{V}\right)\right\rvert_{x=0} = \glslink{varepsilon-f}{\varepsilon_{\glslink{J-rho}{J_\rho}}}\left(\frac{1}{1-y}\right).
    \end{equation*} Applying Proposition~\ref{J-prop}(iii), we thus obtain \begin{equation*}
        \glslink{varepsilon-f}{\varepsilon_{\glslink{J-rho}{J_\rho}}} = \glslink{delta-f-j}{\delta_{\glslink{J-rho}{J_\rho}}^{(\glslink{abs-rho}{\lvert\rho\rvert}-1)}} = \sum_{\rho' \sqcup \mu = \rho} \frac{(-1)^{\glslink{len}{\len(\mu)}}}{\glslink{z-rho}{z_\mu}} \cdot \glslink{delta-f-j}{\delta_{\glslink{I-rho}{I_{\rho'}}}^{(\glslink{abs-rho}{\lvert\rho\rvert}-1)}}.
    \end{equation*} Now, by Theorem~\ref{E-I-rho-from-delta-thm}(ii), for each $\rho'$ appearing in the sum, the largest power of $v$ appearing in $\glslink{delta-f-j}{\delta_{\glslink{I-rho}{I_{\rho'}}}^{(\glslink{abs-rho}{\lvert\rho\rvert}-1)}}$ is at most \begin{equation}
        2\glslink{abs-rho}{\lvert\rho\rvert} + 2\glslink{m-j-rho}{m_1(\rho')} - \glslink{len}{\len(\rho')} - 1 = 2\glslink{abs-rho}{\lvert\rho\rvert} + \glslink{m-j-rho}{m_1(\rho')} - (\glslink{len}{\len(\rho')} - \glslink{m-j-rho}{m_1(\rho')}) - 1 \le 2\glslink{abs-rho}{\lvert\rho\rvert} + \glslink{m-j-rho}{m_1(\rho)} - 1. \label{delta-I-rho-prime-degree-bound}
    \end{equation} This proves claim (i). In order for the left and right hand sides of \eqref{delta-I-rho-prime-degree-bound} to be equal, we must have $\glslink{m-j-rho}{m_1(\rho')}-(\glslink{len}{\len(\rho')}-\glslink{m-j-rho}{m_1(\rho')}) \ge \glslink{m-j-rho}{m_1(\rho)}$, implying that $\glslink{m-j-rho}{m_1(\rho')} = \glslink{m-j-rho}{m_1(\rho)}$ and $\glslink{len}{\len(\rho')} = \glslink{m-j-rho}{m_1(\rho')}$ and hence that $\rho' = 1^{\glslink{m-j-rho}{m_1(\rho)}}$. In this case, Theorem~\ref{E-I-rho-from-delta-thm}(iii) tells us that the leading coefficient of $\glslink{delta-f-j}{\delta_{\glslink{I-rho}{I_{\rho'}}}^{(\glslink{abs-rho}{\lvert\rho\rvert}-1)}}$ is \begin{equation*}
        \frac{2^{2\glslink{abs-rho}{\lvert\rho\rvert}-2}(\frac{2\glslink{abs-rho}{\lvert\rho\rvert}+\glslink{m-j-rho}{m_1(\rho)}-3}{2})!}{(\frac{\glslink{m-j-rho}{m_1(\rho)}-1}{2})!},
    \end{equation*} at which point multiplying by $\frac{(-1)^{\glslink{len}{\len(\mu)}}}{\glslink{z-rho}{z_\mu}} = \frac{(-1)^{\glslink{len}{\len(\rho)}-\glslink{m-j-rho}{m_1(\rho)}}\glslink{m-j-rho}{m_1(\rho)}!}{\glslink{z-rho}{z_\rho}}$ gives the first of the two formulas in claim (iv). If $\rho = 1^r$ for some $r \ge 1$, then the next highest power of $v$ must come only from the terms $\rho' = 1^r$ and $\rho' = 1^{r-1}$, and applying Theorem~\ref{E-I-rho-from-delta-thm}(iii) and simplifying gives the second claimed formula in (iv). Meanwhile, by Theorem~\ref{E-I-rho-from-delta-thm}(i), it is clear that the smallest power of $v$ appearing in $\glslink{delta-f-j}{\delta_{\glslink{I-rho}{I_{\rho'}}}^{(\glslink{abs-rho}{\lvert\rho\rvert}-1)}}$ is at least $1$---if $\rho' = 1^{r'}$ for some $r'$, then this is immediate, and otherwise, it follows because $2-\glslink{abs-rho}{\lvert\rho'\rvert}+(\glslink{abs-rho}{\lvert\rho\rvert}-1) \ge 1$. This proves claim (ii).
    
    All that remains is claim (iii). To see this, use the fact that $\partial_x \glslink{V}{V} = 2\glslink{V}{V}(\glslink{V}{V}-1)\glslink{Q}{Q}^{-1}$ to write \begin{equation*}
        \glslink{partial-x-i}{\partial_x^{\glslink{abs-rho}{\lvert\rho\rvert}}} \glslink{calE-f}{\calE_{\glslink{J-rho}{J_\rho}}} = \partial_x\glslink{varepsilon-f}{\varepsilon_{\glslink{J-rho}{J_\rho}}}(\glslink{V}{V}) = 2\glslink{V}{V}(\glslink{V}{V}-1)\glslink{Q}{Q}^{-1} \cdot \glslink{varepsilon-f}{\varepsilon_{\glslink{J-rho}{J_\rho}}}'(\glslink{V}{V}),
    \end{equation*} then set $x$ equal to zero to obtain \begin{equation*}
        \sum_{k \ge 0} (k+\glslink{abs-rho}{\lvert\rho\rvert})! \cdot \glslink{E-f-n-k}{E_{\glslink{J-rho}{J_\rho}}(k+\glslink{abs-rho}{\lvert\rho\rvert}, k)} \cdot y^k = -\frac{2y}{(1-y)^3} \cdot \glslink{varepsilon-f}{\varepsilon_{\glslink{J-rho}{J_\rho}}}'\left(\frac{1}{1-y}\right).
    \end{equation*} Then by Theorem~\ref{gp-poly-thm}, the coefficient of $y^k$ on the left hand side is zero for $0 \le k < \glslink{abs-rho}{\lvert\rho\rvert}$, implying that $y^{\glslink{abs-rho}{\lvert\rho\vert}-1}$ divides $\glslink{varepsilon-f}{\varepsilon_{\glslink{J-rho}{J_\rho}}}'\left(\frac{1}{1-y}\right)$. Since $\frac{1}{1-y} - 1 = \frac{y}{1-y}$, this is equivalent to the statement that $(v-1)^{\glslink{abs-rho}{\lvert\rho\rvert}-1} \mid \glslink{varepsilon-f}{\varepsilon_{\glslink{J-rho}{J_\rho}}}'$.
\end{proof}

\noindent
As before, the above specialized result about the polynomials $\glslink{varepsilon-f}{\varepsilon_{\glslink{J-rho}{J_\rho}}}$ implies the corresponding result from Section~\ref{proof-overview-section} stated for general $f \in \glslink{Delta-Q-d}{\Delta_\Q^d}$.

\begin{proof}[Proof of Theorem~\ref{eps-properties-thm}]
    Since the $\glslink{J-rho}{J_\rho}$ with $\glslink{abs-rho}{\lvert\rho\rvert} = d$ form a basis for the space $\glslink{Delta-Q-d}{\Delta_\Q^d}$ of \glslink{character-polynomial}{character polynomials} of \glslink{pure-character-polynomial}{pure degree $d$}, we can write $\glslink{varepsilon-f}{\varepsilon_f} = \sum_{\rho \vdash d} c_\rho \cdot \glslink{varepsilon-f}{\varepsilon_{\glslink{J-rho}{J_\rho}}}$ for some coefficients $c_\rho \in \Q$. Claims (i), (ii), and (iii) then follow immediately from parts (i), (ii), and (iii) of Theorem~\ref{eps-J-rho-properties-thm}, respectively. Claim (iv) follows similarly from Theorem~\ref{eps-J-rho-properties-thm}(iv) once we note that the only partition $\rho \vdash d$ satisfying $2\glslink{abs-rho}{\lvert\rho\rvert}+\glslink{m-j-rho}{m_1(\rho)}-1 \ge 3d-2$ is $\rho = 1^d$.
\end{proof}

In order to prove Theorems~\ref{main-thm-1} and \ref{main-thm-3}(i), we start by proving the analogous result for the basis $\glslink{J-rho}{J_\rho}$. While we could jump directly to proving Theorems~\ref{main-thm-1} and \ref{main-thm-3}(i), we note that the following theorem is strictly stronger than what is implied by Theorems~\ref{main-thm-1} and \ref{main-thm-3}(i) in the case where $\rho$ is not of the form $1^r$. This stronger statement will turn out to be useful later for proving Theorem~\ref{conditional-expectation-thm}.

\begin{theorem} \label{E-J-rho-main-thm}
    Let $\rho$ be a partition of size at least $1$. Then: \begin{itemize}
        \item[(i)] There exist rational coefficients $\glslink{theta-f-i}{\theta_{\glslink{J-rho}{J_\rho}}^{(1)}}, \dots, \glslink{theta-f-i}{\theta_{\glslink{J-rho}{J_\rho}}^{(2\glslink{abs-rho}{\lvert\rho\rvert}+\glslink{m-j-rho}{m_1(\rho)})}}$ such that for all $n \ge k \ge 0$ with $n \ge k+\glslink{abs-rho}{\lvert\rho\rvert}-1$ and $(n, k) \ne (\glslink{abs-rho}{\lvert\rho\rvert}-1, 0)$, we have \begin{equation*}
            \frac{n!}{(n-\glslink{abs-rho}{\lvert\rho\rvert})!} \cdot \glslink{E-f-n-k}{E_{\glslink{J-rho}{J_\rho}}(n, k)} = \sum_{1 \le i \le 2\glslink{abs-rho}{\lvert\rho\rvert}+\glslink{m-j-rho}{m_1(\rho)}} \glslink{theta-f-i}{\theta_{\glslink{J-rho}{J_\rho}}^{(i)}} \cdot \frac{(\frac{i}{2})!}{(k-\frac{i}{2})!} \glslink{binomial}{\binom{n-\glslink{abs-rho}{\lvert\rho\rvert}+\frac{i}{2}}{k}}.
        \end{equation*}
        \item[(ii)] The last of these coefficients is given by \begin{equation*}
            \glslink{theta-f-i}{\theta_{\glslink{J-rho}{J_\rho}}^{(2\glslink{abs-rho}{\lvert\rho\rvert}+\glslink{m-j-rho}{m_1(\rho)})}} = \frac{(-1)^{\glslink{len}{\len(\rho)}-\glslink{m-j-rho}{m_1(\rho)}}}{\glslink{z-rho}{z_\rho}} \cdot \frac{(\frac{\glslink{m-j-rho}{m_1(\rho)}}{2})!}{(\glslink{abs-rho}{\lvert\rho\rvert}+\frac{\glslink{m-j-rho}{m_1(\rho)}}{2})!}.
        \end{equation*}
        \item[(iii)] If $\rho = 1^r$ for some $r \ge 1$, then the second-to-last coefficient is given by \begin{equation*}
            \glslink{theta-f-i}{\theta_{\glslink{J-rho}{J_\rho}}^{(3r-1)}} = -\frac{(\frac{r-1}{2})!}{(r-1)!(\frac{3r-1}{2})!}.
        \end{equation*}
    \end{itemize}
\end{theorem}

\begin{figure}[h]
    \renewcommand{\arraystretch}{1.5}
    \setlength{\tabcolsep}{2pt}
    \begin{tabular}{|c|c|c|c|c|c|c|c|c|c|c|c|c|}
        \hline
        $\rho$ & $\glslink{theta-f-i}{\theta_{\glslink{J-rho}{J_\rho}}^{(1)}}$ & $\glslink{theta-f-i}{\theta_{\glslink{J-rho}{J_\rho}}^{(2)}}$ & $\glslink{theta-f-i}{\theta_{\glslink{J-rho}{J_\rho}}^{(3)}}$ & $\glslink{theta-f-i}{\theta_{\glslink{J-rho}{J_\rho}}^{(4)}}$ & $\glslink{theta-f-i}{\theta_{\glslink{J-rho}{J_\rho}}^{(5)}}$ & $\glslink{theta-f-i}{\theta_{\glslink{J-rho}{J_\rho}}^{(6)}}$ & $\glslink{theta-f-i}{\theta_{\glslink{J-rho}{J_\rho}}^{(7)}}$ & $\glslink{theta-f-i}{\theta_{\glslink{J-rho}{J_\rho}}^{(8)}}$ & $\glslink{theta-f-i}{\theta_{\glslink{J-rho}{J_\rho}}^{(9)}}$ & $\glslink{theta-f-i}{\theta_{\glslink{J-rho}{J_\rho}}^{(10)}}$ & $\glslink{theta-f-i}{\theta_{\glslink{J-rho}{J_\rho}}^{(11)}}$ & $\glslink{theta-f-i}{\theta_{\glslink{J-rho}{J_\rho}}^{(12)}}$ \\ \hline
        $(1)$ & $\phantom{+}\frac{1}{2}$ & $-1$ & $\phantom{+}\frac{2}{3}$ & & & & & & & & & \\ \hline
        $(2)$ & $\phantom{+}\frac{1}{2}$ & $-1$ & $\phantom{+}\frac{2}{3}$ & $-\frac{1}{4}$ & & & & & & & & \\ \hline
        $(1, 1)$ & $\phantom{+}\frac{1}{4}$ & $\phantom{+}0$ & $-\frac{2}{3}$ & $\phantom{+}\frac{3}{4}$ & $-\frac{4}{15}$ & $\phantom{+}\frac{1}{12}$ & & & & & & \\ \hline
        $(3)$ & $\phantom{+}\frac{3}{4}$ & $-2$ & $\phantom{+}2$ & $-1$ & $\phantom{+}\frac{4}{15}$ & $-\frac{1}{18}$ & & & & & & \\ \hline
        $(2, 1)$ & $\phantom{+}\frac{15}{16}$ & $-1$ & $-\frac{3}{4}$ & $\phantom{+}\frac{3}{2}$ & $-\frac{13}{15}$ & $\phantom{+}\frac{1}{4}$ & $-\frac{4}{105}$ & & & & & \\ \hline
        $(1, 1, 1)$ & $\phantom{+}\frac{3}{16}$ & $\phantom{+}0$ & $-\frac{5}{16}$ & $\phantom{+}0$ & $\phantom{+}\frac{7}{20}$ & $-\frac{5}{18}$ & $\phantom{+}\frac{13}{105}$ & $-\frac{1}{48}$ & $\phantom{+}\frac{4}{945}$ & & & \\ \hline
        $(4)$ & $\phantom{+}\frac{15}{8}$ & $-5$ & $\phantom{+}\frac{11}{2}$ & $-\frac{7}{2}$ & $\phantom{+}\frac{22}{15}$ & $-\frac{5}{12}$ & $\phantom{+}\frac{8}{105}$ & $-\frac{1}{96}$ & & & & \\ \hline
        $(3, 1)$ & $\phantom{+}\frac{35}{16}$ & $-4$ & $\phantom{+}\frac{2}{3}$ & $\phantom{+}3$ & $-\frac{44}{15}$ & $\phantom{+}\frac{4}{3}$ & $-\frac{16}{45}$ & $\phantom{+}\frac{1}{18}$ & $-\frac{16}{2835}$ & & & \\ \hline
        $(2, 2)$ & $\phantom{+}\frac{15}{16}$ & $-2$ & $\phantom{+}\frac{5}{4}$ & $\phantom{+}0$ & $-\frac{1}{3}$ & $\phantom{+}\frac{1}{6}$ & $-\frac{4}{105}$ & $\phantom{+}\frac{1}{192}$ & & & & \\ \hline
        $(2, 1, 1)$ & $\phantom{+}\frac{45}{32}$ & $-1$ & $-\frac{27}{16}$ & $\phantom{+}\frac{3}{2}$ & $\phantom{+}\frac{23}{60}$ & $-\frac{11}{12}$ & $\phantom{+}\frac{17}{35}$ & $-\frac{13}{96}$ & $\phantom{+}\frac{4}{189}$ & $-\frac{1}{480}$ & & \\ \hline
        $(1, 1, 1, 1)$ & $\phantom{+}\frac{5}{32}$ & $\phantom{+}0$ & $-\frac{25}{96}$ & $\phantom{+}0$ & $\phantom{+}\frac{1}{6}$ & $\phantom{+}0$ & $-\frac{1}{9}$ & $\phantom{+}\frac{49}{576}$ & $-\frac{16}{567}$ & $\phantom{+}\frac{11}{1440}$ & $-\frac{8}{10395}$ & $\phantom{+}\frac{1}{8640}$ \\ \hline
    \end{tabular}
    \caption{The values of the coefficients $\glslink{theta-f-i}{\theta_{\glslink{J-rho}{J_\rho}}^{(i)}}$ from Theorem~\ref{E-J-rho-main-thm} for all partitions $\rho$ of size at most $4$.} \label{theta-J-figure}
\end{figure}

\begin{proof}
    Use Lemma~\ref{E-from-epsilon-lem} to write \begin{equation}
        \sum_{j,k\ge0} \frac{(j+k+\glslink{abs-rho}{\lvert\rho\rvert}-1)!}{j!} \cdot \glslink{E-f-n-k}{E_{J_\rho}(j+k+\glslink{abs-rho}{\lvert\rho\rvert}-1, k)} \cdot x^j y^k = \glslink{partial-x-i}{\partial_x^{\glslink{abs-rho}{\lvert\rho\rvert}-1}} \glslink{calE-f}{\calE_{J_\rho}} = \glslink{varepsilon-f}{\varepsilon_{\glslink{J-rho}{J_\rho}}}(V). \label{E-J-lam-from-eps-J-lam}
    \end{equation} Then by Theorem~\ref{eps-J-rho-properties-thm}, $\glslink{varepsilon-f}{\varepsilon_{\glslink{J-rho}{J_\rho}}}$ has degree at most $2\glslink{abs-rho}{\lvert\rho\rvert}+\glslink{m-j-rho}{m_1(\rho)}-1$ and no constant coefficient. Now, in the special case where $a = 0$ and $b \ge 1$, the formula from Theorem~\ref{q-a-b-r-formula-thm} for the series coefficients of $\glslink{V}{V}^b \glslink{Q}{Q}^{-r}$ reduces to \begin{align}
        \glslink{q-a-b-r-j-k}{q_{0,b}^{(r-1)}(j, k)} &= \sum_{0 \le \beta \le b} \glslink{binomial}{\binom{b}{\beta}} \frac{(j+k-1+r)!(b-\frac{\beta+3}{2}+\frac{r}{2})!(j+k+\frac{\beta-1}{2}+\frac{r}{2})!}{2j!k!(b-2+r)!(j+\frac{\beta-1}{2}+\frac{r}{2})!(k-\frac{\beta+1}{2}+\frac{r}{2})!} \notag \\
        &= \sum_{0 \le \beta \le b} \glslink{binomial}{\binom{b}{\beta}} \frac{(b-\frac{\beta+3}{2}+\frac{r}{2})!}{2(b-2+r)!(\frac{\beta+1}{2})!} \cdot \frac{(j+k-1+r)!(\frac{\beta+1}{2})!}{j!(k-\frac{\beta+1}{2}+\frac{r}{2})!} \glslink{binomial}{\binom{j+k+\frac{\beta-1}{2}+r}{k}} \label{q-formula-a-equals-0}
    \end{align} as polynomials in $r$. In our case, we are interested in the limit of the right hand side as $r \to 0$. Note that for any fixed $0 \le \beta \le b$, the quantity $\glslink{binomial}{\binom{b}{\beta}} \frac{(b-\frac{\beta+3}{2}+\frac{r}{2})!}{2(b-2+r)!(\frac{\beta+1}{2})!}$ approaches a finite rational limit as $r \to 0$: if $b \ge 2$, then the expression is already a well-defined rational number when $r = 0$, whereas if $b = 1$, then a quick computation shows that the limit is $0$ if $\beta = 0$ and $1$ if $\beta = 1$. Meanwhile, the other factor in each summand of \eqref{q-formula-a-equals-0} gives a well-defined value for $r = 0$ whenever $j+k \ge 1$, as long as we interpret the value to be zero whenever $k-\frac{\beta+1}{2}$ is a negative integer. Consequently, we deduce that there exist values $c_0, \dots, c_{2\glslink{abs-rho}{\lvert\rho\rvert}+\glslink{m-j-rho}{m_1(\rho)}-1} \in \Q$ such that for all $j, k \ge 0$ with $j+k \ge 1$, the coefficient of $x^j y^k$ in the right hand side of \eqref{E-J-lam-from-eps-J-lam} is given by \begin{align*}
        &\frac{(j+k-1)!}{j!} \cdot \sum_{\beta=0}^{2\glslink{abs-rho}{\lvert\rho\rvert}+\glslink{m-j-rho}{m_1(\rho)}-1} c_\beta \cdot \frac{(\frac{\beta+1}{2})!}{(k-\frac{\beta+1}{2})!} \glslink{binomial}{\binom{j+k+\frac{\beta-1}{2}}{k}},
    \end{align*} where we interpret a term in the sum as zero whenever $k-\frac{\beta+1}{2}$ is a negative integer. Applying equation \eqref{E-J-lam-from-eps-J-lam}, we deduce that \begin{equation*}
        \frac{(j+k+\glslink{abs-rho}{\lvert\rho\rvert}-1)!}{(j+k-1)!} \cdot \glslink{E-f-n-k}{E_{\glslink{J-rho}{J_\rho}}(j+k+\glslink{abs-rho}{\lvert\rho\rvert}-1, k)} = \sum_{i=1}^{2\glslink{abs-rho}{\lvert\rho\rvert}+\glslink{m-j-rho}{m_1(\rho)}} c_{i-1} \cdot \frac{(\frac{i}{2})!}{(k-\frac{i}{2})!} \glslink{binomial}{\binom{j+k-1+\frac{i}{2}}{k}}
    \end{equation*} for all $j, k \ge 0$ with $j+k \ge 1$. Equivalently, upon substituting $j = n-k-\glslink{abs-rho}{\lvert\rho\rvert}+1$, we obtain \begin{equation*}
        \frac{n!}{(n-\glslink{abs-rho}{\lvert\rho\rvert})!} \cdot \glslink{E-f-n-k}{E_{\glslink{J-rho}{J_\rho}}(n, k)} = \sum_{i=1}^{2\glslink{abs-rho}{\lvert\rho\rvert}+\glslink{m-j-rho}{m_1(\rho)}} c_{i-1} \cdot \frac{(\frac{i}{2})!}{(k-\frac{i}{2})!} \glslink{binomial}{\binom{n-\glslink{abs-rho}{\lvert\rho\rvert}+\frac{i}{2}}{k}}
    \end{equation*} for all $n \ge k \ge 0$ with $n \ge k+\glslink{abs-rho}{\lvert\rho\rvert}-1$ and $(n, k) \ne (\glslink{abs-rho}{\lvert\rho\rvert}-1, 0)$. Setting $\glslink{theta-f-i}{\theta_{\glslink{J-rho}{J_\rho}}^{(i)}} := c_{i-1}$ for all $i$, claim (i) follows immediately.

    For claim (ii), note that the above computation makes it clear that the coefficient $\glslink{theta-f-i}{\theta_{\glslink{J-rho}{J_\rho}}^{(2\glslink{abs-rho}{\lvert\rho\rvert}+\glslink{m-j-rho}{m_1(\rho)})}}$ comes entirely from the coefficient of $v^{2\glslink{abs-rho}{\lvert\rho\rvert}+\glslink{m-j-rho}{m_1(\rho)}-1}$ in $\glslink{varepsilon-f}{\varepsilon_{\glslink{J-rho}{J_\rho}}}$, which is given explicitly by Theorem~\ref{eps-J-rho-properties-thm}(iv). By \eqref{q-formula-a-equals-0}, the coefficient $\glslink{theta-f-i}{\theta_{\glslink{J-rho}{J_\rho}}^{(2\glslink{abs-rho}{\lvert\rho\rvert}+\glslink{m-j-rho}{m_1(\rho)})}} = c_{2\glslink{abs-rho}{\lvert\rho\rvert}+\glslink{m-j-rho}{m_1(\rho)}-1}$ is thus given explicitly by \begin{multline*}
        (-1)^{\glslink{len}{\len(\rho)}-\glslink{m-j-rho}{m_1(\rho)}} \cdot \frac{2^{2\glslink{abs-rho}{\lvert\rho\rvert}-2}\glslink{m-j-rho}{m_1(\rho)}!(\frac{2\glslink{abs-rho}{\lvert\rho\rvert}+\glslink{m-j-rho}{m_1(\rho)}-3}{2})!}{\glslink{z-rho}{z_\rho}(\frac{\glslink{m-j-rho}{m_1(\rho)}-1}{2})!} \cdot \frac{(2\glslink{abs-rho}{\lvert\rho\rvert}+\glslink{m-j-rho}{m_1(\rho)}-1-\frac{2\glslink{abs-rho}{\lvert\rho\rvert}+\glslink{m-j-rho}{m_1(\rho)}+2}{2})!}{2(2\glslink{abs-rho}{\lvert\rho\rvert}+\glslink{m-j-rho}{m_1(\rho)}-3)!(\frac{2\glslink{abs-rho}{\lvert\rho\rvert}+\glslink{m-j-rho}{m_1(\rho)}}{2})!} \\
        = \frac{(-1)^{\glslink{len}{\len(\rho)}-\glslink{m-j-rho}{m_1(\rho)}}}{\glslink{z-rho}{z_\rho}} \cdot \frac{(\frac{\glslink{m-j-rho}{m_1(\rho)}}{2})!}{(\glslink{abs-rho}{\lvert\rho\rvert}+\frac{\glslink{m-j-rho}{m_1(\rho)}}{2})!}.
    \end{multline*}

    For claim (iii), suppose that $\rho = 1^r$ for some $r \ge 1$, and note that, similarly to before, the coefficient $\glslink{theta-f-i}{\theta_{\glslink{J-rho}{J_\rho}}^{(3r-1)}}$ comes entirely from the coefficients of $v^{3r-1}$ and $v^{3r-2}$ in $\glslink{varepsilon-f}{\varepsilon_{\glslink{J-rho}{J_\rho}}}$, which are both given explicitly by Theorem~\ref{eps-J-rho-properties-thm}(iv). We thus find that the coefficient $\glslink{theta-f-i}{\theta_{\glslink{J-rho}{J_\rho}}^{(3r-1)}} = c_{3r-2}$ is given by \begin{multline*}
        \frac{2^{2r-2}(\frac{3r-3}{2})!}{(\frac{r-1}{2})!} \cdot \frac{(3r-1) \cdot (3r-1-\frac{3r+1}{2})!}{2(3r-3)!(\frac{3r-1}{2})!} \\
        \quad - 2^{2r-2} \left(\frac{(\frac{3r-1}{2})!}{(\frac{r-1}{2})!} + \frac{(\frac{3r-4}{2})!}{(\frac{r-2}{2})!}\right) \cdot \frac{(3r-2-\frac{3r+1}{2})!}{2(3r-4)!(\frac{3r-1}{2})!} = -\frac{(\frac{r-1}{2})!}{(r-1)!(\frac{3r-1}{2})!},
    \end{multline*} as claimed.
\end{proof}

\noindent
Having proven Theorem~\ref{E-J-rho-main-thm}, Theorems~\ref{main-thm-1} and \ref{main-thm-3}(i) become immediate corollaries.

\begin{proof}[Proofs of Theorems~\ref{main-thm-1} and \ref{main-thm-3}(i)]
    These results follow immediately by expressing $\glslink{X-lam}{X^\lambda}$ as a linear combination of $\glslink{J-rho}{J_\rho}$ for $\rho \vdash \glslink{abs-rho}{\lvert\lambda\rvert}$. In particular, the computation of the coefficients in Theorem~\ref{main-thm-3}(i) makes use of Proposition~\ref{J-prop}(ii), which tells us that the coefficient of $\glslink{J-rho}{J_{1^{\glslink{abs-rho}{\lvert\lambda\rvert}}}}$ in the expansion of $\glslink{X-lam}{X^\lambda}$ in the $\glslink{J-rho}{J_\rho}$-basis is $\glslink{chi-lambda}{\chi^\lambda}(\glslink{id-k}{\id_{\glslink{abs-rho}{\lvert\lambda\rvert}}})$.
\end{proof}

\section{Applications} \label{applications-section}

In this section, we tie up some loose ends by showing how the various corollaries and applications stated in the introduction follow from our main results.

\subsection{Near real-rootedness of the \texorpdfstring{$\glslink{a-sigma-lambda}{a_{\glslink{id-k}{\id_k}}^\lambda}(n)$}{aλidk(n)}}

Our first application is a proof of Corollary~\ref{real-roots-cor}, which represents a slight weakening of Conjecture~\ref{root-conj}. To prove this corollary, we make use of the following two lemmas.

\begin{lemma} \label{alternating-signs-root-lem}
    Let $f \in \R[x]$ be a nonzero polynomial, and suppose that there exists a sequence $x_0 < \dots < x_m$ of real numbers such that for all $0 \le i \le m$, we have $(-1)^i f(x_i) \ge 0$. Then $f$ has at least $m$ real roots in the interval $[x_0, x_m]$, counting multiplicity.
\end{lemma}

\begin{proof}
    If $f(x_i) \ne 0$ for all $0 \le i \le m$, the result follows immediately by the intermediate value theorem. Otherwise, proceed by induction on $m$. If $m = 0$, there is nothing to prove. Otherwise, suppose that $f(x_j) = 0$ for some $j$. Then writing $f = (x_j-x)g$ for some $g \in \R[x]$, we see that $(-1)^ig(x_i) \ge 0$ for $0 \le i \le j-1$ and $(-1)^ig(x_{i+1}) \ge 0$ for $j \le i \le m-1$. Applying the inductive hypothesis to the values $x_0 < \dots < x_{j-1} < x_{j+1} < \dots < x_m$, we conclude that $g$ has at least $m-1$ roots on the interval $[x_0, x_m]$ counting multiplicity, and so $f$ has at least $m$ roots in $[x_0, x_m]$ counting multiplicity.
\end{proof}

\begin{lemma} \label{alternating-signs-subseq-lem}
    Let $m \ge 0$ and $0 \le d \le m$, and let $x_0 < \dots < x_m$ be real numbers and $g \in \R[x]$ be a polynomial of degree at most $d$. Then there exists a subsequence $x_{i_0} < \dots < x_{i_{m-d}}$ of the $x_i$ and a sign $\sigma \in \{\pm1\}$ such that $\sigma (-1)^j\cdot (-1)^{i_j}g(x_{i_j}) \ge 0$ for all $0 \le j \le m-d$.
\end{lemma}

\begin{proof}
    If $g$ is zero, take any subsequence of the required length and either sign $\sigma$. Otherwise, proceed by induction on the number of real roots of $g$. If $g$ is nonzero and has no real roots, then take the subsequence $x_0 < \dots < x_{m-d}$, and take $\sigma$ to be the sign of $g(x)$ for all $x \in \R$. Next, suppose that $g = (x-X)\hat{g}$ for some $X \in \R$ and $\hat{g} \in \R[x]$, and use the inductive hypothesis to obtain a subsequence $x_{i_0'} < \dots < x_{i_{m-d+1}'}$ and a sign $\sigma'$ such that $\sigma'(-1)^j \cdot (-1)^{i_j}\hat{g}(x_{i_j'}) \ge 0$ for all $0 \le j \le m-d+1$. If $X \notin [x_{i_0'}, x_{i_{m-d+1}'}]$, then we may use the subsequence $x_{i_0'} < \dots < x_{i_{m-d}'}$ and take $\sigma$ to be $\sigma'$ times the sign of $x_{i_0'}-X$. Otherwise, we may choose a $0 \le j \le m-d$ such that $X \in [x_{i_j'}, x_{i_{j+1}'}]$. Then we may take $\sigma := -\sigma'$ and take $x_{i_0} < \dots < x_{i_{m-d}}$ to be the sequence $x_{i_0'} < \dots < x_{i_{m-d+1}'}$ with $x_{i_j'}$ removed.
\end{proof}

With these lemmas in hand, Corollary~\ref{real-roots-cor} follows by noting that in the formula for $\glslink{a-sigma-lambda}{a_{\glslink{id-k}{\id_k}}^\lambda}(n)$ from Theorem~\ref{main-thm-1}, it can be fruitful to separate terms corresponding to even $i$ and those corresponding to odd $i$.

\begin{proof}[Proof of Corollary~\ref{real-roots-cor}]
    If $k \le 3\glslink{abs-rho}{\lvert\lambda\rvert}-1$, then $\glslink{total-deg}{\deg(\glslink{a-sigma-lambda}{a_{\glslink{id-k}{\id_k}}^\lambda}(n))} \le k-\glslink{abs-rho}{\lvert\lambda\rvert} \le 2\glslink{abs-rho}{\lvert\lambda\rvert}-1$ and the result holds trivially; thus, suppose that $k \ge 3\glslink{abs-rho}{\lvert\lambda\rvert}$. Use Theorem~\ref{main-thm-1} to write \begin{equation*}
        \frac{n!}{(n-\glslink{abs-rho}{\lvert\lambda\rvert})!} \cdot \glslink{a-sigma-lambda}{a_{\glslink{id-k}{\id_k}}^\lambda}(n) = g_0(n) + g_1(n),
    \end{equation*} where \begin{equation*}
        g_\epsilon := \sum_{\substack{1 \le i \le 3\glslink{abs-rho}{\lvert\lambda\rvert} \\ i \equiv \epsilon \!\!\!\!\pmod{2}}} \glslink{theta-f-i}{\theta_{\glslink{X-lam}{X^\lambda}}^{(i)}} \cdot \frac{(\frac{i}{2})!}{(k-\frac{i}{2})!} \cdot \glslink{binomial}{\binom{n-\glslink{abs-rho}{\lvert\lambda\rvert}+\frac{i}{2}}{k}} \in \Q[n] \qquad \text{for $\epsilon \in \{0, 1\}$.}
    \end{equation*} Then since $\glslink{binomial}{\binom{x}{k}} = 0$ for all integers $-1 < x < k$, we find that $g_0(n) = 0$ for $n \in \Z$ such that $\glslink{abs-rho}{\lvert\lambda\rvert}-\frac{3}{2} < n < k-\frac{\glslink{abs-rho}{\lvert\lambda\rvert}}{2}$ (a total of $\lceil k-\frac{3\glslink{abs-rho}{\lvert\lambda\rvert}}{2}+1\rceil$ values), whereas $g_1(n) = 0$ for $n \in \frac{1}{2}+\Z$ such that $\glslink{abs-rho}{\lvert\lambda\rvert}-\frac{3}{2} < n < k-\frac{\glslink{abs-rho}{\lvert\lambda\rvert}}{2}$ (a total of $\lfloor k-\frac{3\glslink{abs-rho}{\lvert\lambda\rvert}}{2}+1\rfloor$ values). In other words, letting \begin{equation*}
        P_\epsilon := \prod_{\substack{2\glslink{abs-rho}{\lvert\lambda\rvert}-3 < i < 2k-\glslink{abs-rho}{\lvert\lambda\rvert} \\ i\equiv\epsilon\!\!\!\!\pmod{2}}} \left(n-\frac{i}{2}\right) \in \Q[n] \qquad \text{for $\epsilon \in \{0, 1\}$},
    \end{equation*} we have $P_0 \mid g_0$ and $P_1 \mid g_1$. For each root $n$ of $P_1(n)$, we may thus write \begin{equation*}
        \frac{n!}{(n-\glslink{abs-rho}{\lvert\lambda\rvert})!} \cdot \glslink{a-sigma-lambda}{a_{\glslink{id-k}{\id_k}}^\lambda}(n) = P_0(n) \cdot \frac{g_0}{P_0}(n).
    \end{equation*} Now, as $n$ ranges through these $\lfloor k-\frac{3\glslink{abs-rho}{\lvert\lambda\rvert}}{2}+1\rfloor$ roots, it is clear that the sign of $P_0(n)$ alternates between positive and negative. Moreover, it is clear that the polynomial $\frac{g_0}{P_0}$ has degree at most $k - \lceil k-\frac{3\glslink{abs-rho}{\lvert\lambda\rvert}}{2}+1\rceil = \lfloor\frac{3\glslink{abs-rho}{\lvert\lambda\rvert}}{2}\rfloor-1$. Applying Lemma~\ref{alternating-signs-subseq-lem} and noting that $\lfloor k-\frac{3\glslink{abs-rho}{\lvert\lambda\rvert}}{2}+1\rfloor-1-\lfloor\frac{3\glslink{abs-rho}{\lvert\lambda\rvert}}{2}-1\rfloor = k-3\glslink{abs-rho}{\lvert\lambda\rvert}+1$, it follows that there is some subsequence $n_0 < \dots < n_{k-3\glslink{abs-rho}{\lvert\lambda\rvert}+1}$ of the roots of $P_1(n)$ and some sign $\sigma \in \{\pm1\}$ such that $\sigma(-1)^i \cdot P_0(n_i) \cdot \frac{g_0}{P_0}(n_i) \ge 0$ for all $i$. Since $\frac{n!}{(n-\glslink{abs-rho}{\lvert\lambda\rvert})!}$ is positive for $n \in (\glslink{abs-rho}{\lvert\lambda\rvert}-1, k-\frac{\glslink{abs-rho}{\lvert\lambda\rvert}}{2})$, this in turn implies that \begin{equation*}
        \sigma(-1)^i \cdot \glslink{a-sigma-lambda}{a_{\glslink{id-k}{\id_k}}^\lambda}(n_i) = \frac{(n_i-\glslink{abs-rho}{\lvert\lambda\rvert})!}{n_i!} \cdot P_0(n_i) \cdot \frac{g_0}{P_0}(n_i) \ge 0
    \end{equation*} for all $0 \le i \le k-3\glslink{abs-rho}{\lvert\lambda\rvert}+1$. Applying Lemma~\ref{alternating-signs-root-lem}, we conclude that $\glslink{a-sigma-lambda}{a_{\glslink{id-k}{\id_k}}^\lambda}(n)$ has at least $k-3\glslink{abs-rho}{\lvert\lambda\rvert}+1$ roots in the interval $(\glslink{abs-rho}{\lvert\lambda\rvert}-1, k-\frac{\glslink{abs-rho}{\lvert\lambda\rvert}}{2})$, counting multiplicity. Since the total degree of $\glslink{a-sigma-lambda}{a_{\glslink{id-k}{\id_k}}^\lambda}(n)$ is at most $k-\glslink{abs-rho}{\lvert\lambda\rvert}$, this means it has at most $2\glslink{abs-rho}{\lvert\lambda\rvert}-1$ non-real roots, hence at most $\glslink{abs-rho}{\lvert\lambda\rvert}-1$ conjugate pairs of non-real roots.
\end{proof}

\subsection{Positivity of \texorpdfstring{$\glslink{A-sigma-lambda}{A_{\glslink{id-k}{\id_k}}^{\glslink{lambda-bracket-n}{\lambda[n]}}}$}{Aλ[n]idk} for sufficiently large \texorpdfstring{$k$}{k}}

Next, we use Theorems~\ref{main-thm-1} and \ref{main-thm-3}(i) to give a short proof that the values $\glslink{A-sigma-lambda}{A_{\glslink{id-k}{\id_k}}^{\glslink{lambda-bracket-n}{\lambda[n]}}}$ are nonnegative whenever $k$ is sufficiently large.

\begin{proof}[Proof of Corollary~\ref{eventually-pos-cor}]
    For $\lambda = \varnothing$, it is already known that $\glslink{A-sigma-lambda}{A_{\glslink{id-k}{\id_k}}^{\glslink{lambda-bracket-n}{\lambda[n]}}} \ge 0$ for all $n \ge k \ge 0$. Thus, let $\lambda$ be any partition of size at least $1$. Then we may explicitly compute the polynomials $\glslink{b-f-j}{b_{\glslink{X-lam}{X^\lambda}}^{(0)}}, \dots, \glslink{b-f-j}{b_{\glslink{X-lam}{X^\lambda}}^{(\glslink{abs-rho}{\lvert\lambda\rvert}-1)}}$ from Theorem~\ref{main-thm-2} as well as the coefficients $\glslink{theta-f-i}{\theta_{\glslink{X-lam}{X^\lambda}}^{(1)}}, \dots, \glslink{theta-f-i}{\theta_{\glslink{X-lam}{X^\lambda}}^{(3\glslink{abs-rho}{\lvert\lambda\rvert})}}$ from Theorem~\ref{main-thm-1}. By Theorem~\ref{main-thm-3}(ii), the leading coefficients of the polynomials $\glslink{b-f-j}{b_{\glslink{X-lam}{X^\lambda}}^{(0)}}, \dots, \glslink{b-f-j}{b_{\glslink{X-lam}{X^\lambda}}^{(\glslink{abs-rho}{\lvert\lambda\rvert}-1)}}$ are all positive, so there exists a $t_b \ge 0$ such that $\glslink{A-sigma-lambda}{A_{\glslink{id-k}{\id_k}}^{\glslink{lambda-bracket-n}{\lambda[j+k]}}} = \frac{\glslink{b-f-j}{b_{\glslink{X-lam}{X^\lambda}}^{(j)}}(j+k)}{(j+k)!} \ge 0$ for all $0 \le j \le \glslink{abs-rho}{\lvert\lambda\rvert}-1$ and $k \ge t_b$.
    
    Meanwhile, by Theorem~\ref{main-thm-1}, for all $j \ge \glslink{abs-rho}{\lvert\lambda\rvert}$ and $k \ge \frac{3\glslink{abs-rho}{\lvert\lambda\rvert}}{2}$, we may write \begin{equation}
        \glslink{A-sigma-lambda}{A_{\glslink{id-k}{\id_k}}^{\glslink{lambda-bracket-n}{\lambda[j+k]}}} = \frac{(j+k-\glslink{abs-rho}{\lvert\lambda\rvert})!}{(j+k)!} \sum_{1 \le i \le 3\glslink{abs-rho}{\lvert\lambda\rvert}} \glslink{theta-f-i}{\theta_{\glslink{X-lam}{X^\lambda}}^{(i)}} \cdot \frac{(\frac{i}{2})!}{(k-\frac{i}{2})!} \cdot \glslink{binomial}{\binom{j+k-\glslink{abs-rho}{\lvert\lambda\rvert}+\frac{i}{2}}{k}}. \label{A-from-thetas}
    \end{equation} Now, by Theorem~\ref{main-thm-3}(i), the coefficient $\glslink{theta-f-i}{\theta_{\glslink{X-lam}{X^\lambda}}^{(3\glslink{abs-rho}{\lvert\lambda\rvert})}}$ is positive, and by the definition of the generalized binomial coefficient, we may write \begin{equation*}
        \glslink{binomial}{\binom{j+k-\glslink{abs-rho}{\lvert\lambda\rvert}+\frac{i}{2}}{k}} = \frac{1}{k!} \prod_{\ell=0}^{k-1} \left(j+k-\glslink{abs-rho}{\lvert\lambda\rvert}+\frac{i}{2}-\ell\right) \le \frac{1}{k!} \prod_{\ell=0}^{k-1} \left(j+k+\frac{\glslink{abs-rho}{\lvert\lambda\rvert}}{2}-\ell\right) = \glslink{binomial}{\binom{j+k+\frac{\glslink{abs-rho}{\lvert\lambda\rvert}}{2}}{k}}
    \end{equation*} for all $1 \le i \le 3\glslink{abs-rho}{\lvert\lambda\rvert}$. Consequently, we may use \eqref{A-from-thetas} to deduce that \begin{align*}
        \glslink{A-sigma-lambda}{A_{\glslink{id-k}{\id_k}}^{\glslink{lambda-bracket-n}{\lambda[j+k]}}} \ge \frac{(j+k-\glslink{abs-rho}{\lvert\lambda\rvert})!}{(j+k)!(k-\frac{3\glslink{abs-rho}{\lvert\lambda\rvert}}{2})!} \glslink{binomial}{\binom{j+k+\frac{\glslink{abs-rho}{\lvert\lambda\rvert}}{2}}{k}} \!\! \left((\tfrac{3\glslink{abs-rho}{\lvert\lambda\rvert}}{2})!\glslink{theta-f-i}{\theta_{\glslink{X-lam}{X^\lambda}}^{(3\glslink{abs-rho}{\lvert\lambda\rvert})}} - \!\!\!\!\!\!\sum_{1 \le i \le 3\glslink{abs-rho}{\lvert\lambda\rvert}-1}\!\!\!\! (\tfrac{i}{2})!\max\left(-\glslink{theta-f-i}{\theta_{\glslink{X-lam}{X^\lambda}}^{(i)}}, 0\right) \frac{(k-\frac{3\glslink{abs-rho}{\lvert\lambda\rvert}}{2})!}{(k-\frac{i}{2})!}\right).
    \end{align*} For each $1 \le i \le 3\glslink{abs-rho}{\lvert\lambda\rvert}-1$, the ratio $\frac{(k-\frac{3\glslink{abs-rho}{\lvert\lambda\rvert}}{2})!}{(k-\frac{i}{2})!}$ is decreasing in $k$, since for all $k \ge \frac{3\glslink{abs-rho}{\lvert\lambda\rvert}}{2}$, we have \begin{equation*}
        \frac{(k+1-\frac{3\glslink{abs-rho}{\lvert\lambda\rvert}}{2})!}{(k+1-\frac{i}{2})!} = \frac{k+1-\frac{3\glslink{abs-rho}{\lvert\lambda\rvert}}{2}}{k+1-\frac{i}{2}} \cdot \frac{(k-\frac{3\glslink{abs-rho}{\lvert\lambda\rvert}}{2})!}{(k-\frac{i}{2})!}
    \end{equation*} with $\frac{k+1-\frac{3\glslink{abs-rho}{\lvert\lambda\rvert}}{2}}{k+1-\frac{i}{2}} < 1$. Thus, in order to find a $t_a \ge 0$ such that $\glslink{A-sigma-lambda}{A_{\glslink{id-k}{\id_k}}^{\glslink{lambda-bracket-n}{\lambda[j+k]}}} \ge 0$ for all $j \ge \glslink{abs-rho}{\lvert\lambda\rvert}$ and $k \ge t_a$, it suffices to choose any $t_a \ge \frac{3\glslink{abs-rho}{\lvert\lambda\rvert}}{2}$ such that \begin{equation*}
        \sum_{1 \le i \le 3\glslink{abs-rho}{\lvert\lambda\rvert}-1} (\tfrac{i}{2})!\max\left(-\glslink{theta-f-i}{\theta_{\glslink{X-lam}{X^\lambda}}^{(i)}}, 0\right) \cdot \frac{(t_a-\frac{3\glslink{abs-rho}{\lvert\lambda\rvert}}{2})!}{(t_a-\frac{i}{2})!} \le (\tfrac{3\glslink{abs-rho}{\lvert\lambda\rvert}}{2})!\glslink{theta-f-i}{\theta_{\glslink{X-lam}{X^\lambda}}^{(3\glslink{abs-rho}{\lvert\lambda\rvert})}}.
    \end{equation*} The fact that some such $t_a$ exists is immediate from the asymptotic $\frac{(x-a)!}{(x-b)!} \sim x^{b-a}$, which holds by \cite[\href{https://dlmf.nist.gov/5.11.E12}{(5.11.12)}]{dlmf}. The claim therefore follows with $t_\lambda := \max(t_a, t_b)$.
\end{proof}

For completeness, we also give a brief explanation of the algorithm used to verify Proposition~\ref{nonnegativity-prop}.

\begin{proof}[Algorithm used for verifying Proposition~\ref{nonnegativity-prop}]
    Fix any partition $\lambda$ of size at least $1$. Start by finding the smallest integer $t \ge \frac{3\glslink{abs-rho}{\lvert\lambda\rvert}}{2}$ such that \begin{equation*}
        \sum_{1 \le i \le 3\glslink{abs-rho}{\lvert\lambda\rvert}-1} \max\left(-\glslink{theta-f-i}{\theta_{\glslink{X-lam}{X^\lambda}}^{(i)}}, 0\right) \cdot \frac{(\tfrac{i}{2})!(t-\frac{3\glslink{abs-rho}{\lvert\lambda\rvert}}{2})!}{(\tfrac{3\glslink{abs-rho}{\lvert\lambda\rvert}}{2})!(t-\frac{i}{2})!} \le \glslink{theta-f-i}{\theta_{\glslink{X-lam}{X^\lambda}}^{(3\glslink{abs-rho}{\lvert\lambda\rvert})}}.
    \end{equation*} Then by the reasoning from the proof of Corollary~\ref{eventually-pos-cor}, we are guaranteed to have $\glslink{A-sigma-lambda}{A_{\glslink{id-k}{\id_k}}^{\glslink{lambda-bracket-n}{\lambda[j+k]}}} \ge 0$ for all $j \ge \glslink{abs-rho}{\lvert\lambda\rvert}$ and $k \ge t$. From here, for each fixed $0 \le j < \glslink{abs-rho}{\lvert\lambda\rvert}$, we can compute the polynomial $\glslink{b-f-j}{b_{\glslink{X-lam}{X^\lambda}}^{(j)}}(n)$ and verify that it takes nonnegative values for all $n \ge \glslink{abs-rho}{\lvert\lambda\rvert}+\glslink{lambda-i}{\lambda_1}$. Similarly, for each fixed $0 \le k < t$, we can compute the polynomial $\glslink{a-f-k}{a_{\glslink{X-lam}{X^\lambda}}^{(k)}}(n) = \glslink{a-sigma-lambda}{a_{\glslink{id-k}{\id_k}}^\lambda}(n)$ and verify that it takes nonnegative values for all $n \ge \max(k, \glslink{abs-rho}{\lvert\lambda\rvert}+\glslink{lambda-i}{\lambda_1})$. Throughout, use Proposition~\ref{stitch-binom-prop} to speed up the computation of the coefficients $\glslink{stitch-ell-m-t-nu-binom}{\binom{\rho}{\ell,m,t,\nu}}$ wherever they are used. To be concrete, in order to algorithmically show that a polynomial $g \in \Q[n]$ takes nonnegative values for all $n \ge n_0$, we can first find an $n_1 \ge n_0$ such that the coefficients of $g(n+n_1)$ are all nonnegative, then manually check that the values $g(n)$ are nonnegative for each $n_0 \le n < n_1$.
\end{proof}

\subsection{A local rule for the polynomials \texorpdfstring{$\glslink{a-sigma-lambda}{a_{\glslink{id-k}{\id_k}}^\lambda}(n)$}{aλidk(n)}}

Theorem~\ref{main-thm-1} also implies that the polynomials $\glslink{a-sigma-lambda}{a_{\glslink{id-k}{\id_k}}^\lambda}(n)$ satisfy a certain recurrence relation, whose existence is essentially made possible by the fact that the functions $\glslink{partial-x-i}{\partial_x^{\glslink{character-polynomial-deg}{\deg(f)}-1}}\glslink{calE-f}{\calE_f}$ can be expressed as polynomials in $\glslink{V}{V}$ alone when $f$ is \glslink{pure-character-polynomial}{pure}. It would be an interesting topic for future study to investigate whether this recurrence admits a more direct combinatorial proof.

\begin{proof}[Proof of Corollary~\ref{local-rule-cor}]
    By multiplying both sides by $\frac{n!}{(n-\glslink{abs-rho}{\lvert\lambda\rvert})!}$, we see that the claimed equation is equivalent to \begin{equation*}
        (n+1-\glslink{abs-rho}{\lvert\lambda\rvert}) \cdot \frac{(n+1)!}{(n+1-\glslink{abs-rho}{\lvert\lambda\rvert})!} \glslink{a-sigma-lambda}{a_{\glslink{id-k}{\id_{k+1}}}^\lambda}(n+1) = \frac{n!}{(n-\glslink{abs-rho}{\lvert\lambda\rvert})!} \glslink{a-sigma-lambda}{a_{\glslink{id-k}{\id_k}}^\lambda}(n) + (n+k+2-\glslink{abs-rho}{\lvert\lambda\rvert}) \cdot \frac{n!}{(n-\glslink{abs-rho}{\lvert\lambda\rvert})!} \glslink{a-sigma-lambda}{a_{\glslink{id-k}{\id_{k+1}}}^\lambda}(n).
    \end{equation*} We routinely verify that \begin{multline*}
        (n+1-\glslink{abs-rho}{\lvert\lambda\rvert}) \cdot \frac{(\frac{i}{2})!}{(k+1-\frac{i}{2})!} \glslink{binomial}{\binom{n+1-\glslink{abs-rho}{\lvert\lambda\rvert}+\frac{i}{2}}{k+1}} \\
        = \frac{(\frac{i}{2})!}{(k-\frac{i}{2})!} \glslink{binomial}{\binom{n-\glslink{abs-rho}{\lvert\lambda\rvert}+\frac{i}{2}}{k}} + (n+k+2-\glslink{abs-rho}{\lvert\lambda\rvert}) \cdot \frac{(\frac{i}{2})!}{(k+1-\frac{i}{2})!} \glslink{binomial}{\binom{n-\glslink{abs-rho}{\lvert\lambda\rvert}+\frac{i}{2}}{k+1}}
    \end{multline*} as polynomials in $n$ for all $k \ge 0$ and $i \ge 1$, and so the claim follows by Theorem~\ref{main-thm-1}.
\end{proof}

\subsection{Existence of \texorpdfstring{\glslink{character-polynomial}{character polynomials}}{character polynomials} whose covariances with \texorpdfstring{$\glslink{N-sigma}{N_{\glslink{id-k}{\id_k}}}$}{Nidk} are always zero}

Our results also give rise to the somewhat counterintuitive consequence that the map $\glslink{Delta-Q}{\Delta_\Q} \to \Q[[x, y]]$ given by $f \mapsto \glslink{calE-f}{\calE_f}$ is highly non-injective, despite the fact that $\glslink{calE-f}{\calE_f}$ carries information about the infinite array of values $\glslink{E-f-n-k}{E_f(n, k)}$ for all $n \ge k \ge 0$. The task of characterizing the kernel of this map in more detail could be an interesting topic for future study.

\begin{proof}[Proof of Corollary~\ref{kernel-cor}]
    Let $d \ge 1$, and define linear maps \begin{align*}
        \phi_j \colon \glslink{Delta-Q-d}{\Delta_\Q^d} &\to \Q[[y]] \\
        f &\mapsto \left.\glslink{partial-x-i}{\partial_x^j} \glslink{calE-f}{\calE_f}\right|_{x=0}
    \end{align*} for all $0 \le j \le d-2$ as well as \begin{align*}
        \psi \colon \glslink{Delta-Q-d}{\Delta_\Q^d} &\to \Q[[x, y]] \\
        f &\mapsto \glslink{partial-x-i}{\partial_x^{d-1}} \glslink{calE-f}{\calE_f}.
    \end{align*} Then by the definition of the generating function $\glslink{calE-f}{\calE_f}$, a \glslink{character-polynomial}{character polynomial} $f \in \glslink{Delta-Q-d}{\Delta_\Q^d}$ satisfies $\glslink{E-f-n-k}{E_f(n, k)} = 0$ for all $n \ge k \ge 0$ if and only if $f \in \ker(\psi) \cap \bigcap_{j=0}^{d-2} \ker(\phi_j)$. Now, by Theorem~\ref{E-from-delta-thm}, for each $0 \le j \le d-2$, the image of $\phi_j$ has dimension at most $(1+2j+d) - (2+j-d) + 1 = 2d+j$. Meanwhile, by Lemma~\ref{E-from-epsilon-lem} together with Theorem~\ref{eps-properties-thm}(i--iii), the image of $\psi$ has dimension at most $2d$. Consequently, the total dimension of the image of the combined map $\psi \oplus \bigoplus_{j=0}^{d-2} \phi_j$ is at most \begin{equation*}
        \sum_{j=0}^{d-2} (2d+j) + 2d = \frac{1}{2}\left(5d^2-3d+2\right).
    \end{equation*} Since $\glslink{Delta-Q-d}{\Delta_\Q^d}$ has dimension equal to $p(d)$, the number of partitions of size $d$, it follows that the dimension of $\ker(\psi) \cap \bigcap_{j=0}^{d-2} \ker(\phi_j)$ is at least \begin{equation*}
        p(d) - \frac{1}{2}\left(5d^2-3d+2\right).
    \end{equation*} We verify that this lower bound is positive for $d \ge 24$.
\end{proof}

For completeness, we also give a brief justification for Example~\ref{E-f-n-k-zero-example} from the introduction.

\begin{proof}[Proof of Example~\ref{E-f-n-k-zero-example}]
    By combining the relations between the $\glslink{m-j}{m_j}$ and $\glslink{I-rho}{I_\rho}$ bases and those between the $\glslink{I-rho}{I_\rho}$ and $\glslink{J-rho}{J_\rho}$ bases, one can compute that $f = 6\glslink{J-rho}{J_{(6,2)}} - 12\glslink{J-rho}{J_{(4,4)}} - 48\glslink{J-rho}{J_{(2,2,2,2)}}$. By Theorem~\ref{E-J-rho-main-thm} and linearity of $\glslink{E-f-n-k}{E_f(n, k)}$ in $f$, it follows that there exist rational coefficients $\glslink{theta-f-i}{\theta_f^{(1)}}, \dots, \glslink{theta-f-i}{\theta_f^{(16)}}$ such that \begin{equation*}
        \frac{n!}{(n-8)!} \cdot \glslink{E-f-n-k}{E_f(n, k)} = \sum_{1 \le i \le 16} \glslink{theta-f-i}{\theta_f^{(i)}} \cdot \frac{\left(\frac{i}{2}\right)!}{\left(k-\frac{i}{2}\right)!} \glslink{binomial}{\binom{n-8+\frac{i}{2}}{k}}
    \end{equation*} whenever $n-k \ge 7$ and $(n, k) \ne (7, 0)$. Explicitly, using our computer implementation of the algorithms from the paper, we find that $\glslink{theta-f-i}{\theta_f^{(i)}} = 0$ for all $i$, and so $\glslink{E-f-n-k}{E_f(n, k)} = 0$ for all $n \ge k \ge 0$ with $n-k \ge 7$ and $(n, k) \ne (7, 0)$.
    
    To address the remaining quantity $\glslink{E-f-n-k}{E_f(7, 0)}$, use Theorem~\ref{main-thm-combo}(iii) to write $\glslink{E-f-n-k}{E_f(n, n-7)} = \frac{\glslink{b-f-j}{b_f^{(7)}}(n)}{n!}$ for all $n \ge 7$. Since $\glslink{E-f-n-k}{E_f(n, n-7)} = 0$ for all $n \ge 8$, it follows that $\glslink{b-f-j}{b_f^{(7)}}$ is the zero polynomial, at which point we may conclude that $\glslink{E-f-n-k}{E_f(7, 0)} = \frac{\glslink{b-f-j}{b_f^{(7)}}(7)}{7!} = 0$ as well.
\end{proof}

\subsection{Conditional expected values of \texorpdfstring{$\glslink{N-sigma}{N_{\glslink{id-k}{\id_k}}(\pi)}$}{Nidk(π)}}

Lastly, we prove the claimed asymptotic formula from Theorem~\ref{conditional-expectation-thm} for the conditional expectation of $\glslink{N-sigma}{N_{\glslink{id-k}{\id_k}}(\pi)}$ given an observed sub-permutation of $\pi$ with cycle type $\rho$. Before proceeding to the proof, we recall a standard approximation formula for ratios of gamma functions.

\begin{proposition}[Weakened form of {\cite[\href{https://dlmf.nist.gov/5.11.E13}{(5.11.13)}]{dlmf}}] \label{gamma-ratio-prop}
    Fix real numbers $a, b \in \R$. Then for real $z$, as $z \to +\infty$, we have \begin{equation*}
        \frac{\Gamma(z+a)}{\Gamma(z+b)} = z^{a-b}\left(1+O\left(\frac{1}{z}\right)\right).
    \end{equation*}
\end{proposition}

\noindent
In our case, we will use the following slightly more precise formulation of this result.

\begin{corollary} \label{gamma-ratio-cor}
    Fix $a, b \in \R$. Then we have \begin{equation*}
        \frac{(x+a)!}{(x+b)!} = x^{a-b}\left(1+O\left(\frac{1}{x}\right)\right)
    \end{equation*} for all $x \ge \max(1, -a)$, where we interpret the left hand side as $0$ whenever $x+b$ is a negative integer.
\end{corollary}

\begin{proof}
    By Proposition~\ref{gamma-ratio-prop}, we have \begin{equation}
        x\left(\frac{(x+a)!}{(x+b)!x^{a-b}} - 1\right) = O(1) \label{gamma-ratio-o-1}
    \end{equation} for sufficiently large $x$, say, greater than some $x_0$. Moreover, with our specified interpretation of $\frac{(x+a)!}{(x+b)!}$, the left hand side of \eqref{gamma-ratio-o-1} is continuous and hence bounded for $x \in [\max(1, -a), x_0]$. Together, these two facts imply that the left hand side of \eqref{gamma-ratio-o-1} is bounded for $x \in [\max(1, -a), \infty)$.
\end{proof}

With the above approximation in mind, we start with the following lemma, which uses Theorem~\ref{E-J-rho-main-thm} to give a simpler approximate formula for the values $\glslink{E-f-n-k}{E_{\glslink{J-rho}{J_\rho}}(n, k)}$.

\begin{lemma} \label{E-J-rho-approx-lem}
    Fix a partition $\rho$ and a real number $0 < \epsilon < \frac{1}{2}$. Then for all $n \ge 1$ and $\epsilon n \le k \le (1-\epsilon)n$, setting $t := \frac{k}{n}$, we have \begin{equation*}
        \frac{\glslink{E-f-n-k}{E_{\glslink{J-rho}{J_\rho}}(n, k)}}{\frac{1}{k!}\glslink{binomial}{\binom{n}{k}}} = \frac{(-1)^{\glslink{len}{\len(\rho)}-\glslink{m-j-rho}{m_1(\rho)}}(\tfrac{\glslink{m-j-rho}{m_1(\rho)}}{2})!}{\glslink{z-rho}{z_\rho}} \cdot \left(\frac{nt}{1-t}\right)^{\frac{\glslink{m-j-rho}{m_1(\rho)}}{2}} t^{\glslink{abs-rho}{\lvert\rho\rvert}} + O\left(n^{\frac{m_1(\rho)-1}{2}}\right),
    \end{equation*} where the implicit constant of the big $O$ depends only on $\rho$ and $\epsilon$.
\end{lemma}

\begin{proof}
    To start, note that in the case $\rho = \varnothing$, we have $\frac{E_{J_\rho}(n, k)}{\frac{1}{k!}\glslink{binomial}{\binom{n}{k}}} = 1$, and the result is immediate. Thus, we assume $\lvert\rho\rvert \ge 1$ for the remainder of the proof.

    Next, observe that only finitely many pairs of values $(n, k)$ subject to the given constraints satisfy $n-k \le \lvert\rho\rvert$---this condition implies $\varepsilon n \le n-k \le \lvert\rho\rvert$ and hence $n \le \frac{\lvert\rho\rvert}{\epsilon}$, and for each value of $n$, there are only finitely many possible values of $k$. Thus, we may prove the result assuming that $n-k \ge \lvert\rho\rvert$, then adjust the implicit constant after the fact to accommodate these finitely many additional cases.

    Suppose then that $n-k \ge \glslink{abs-rho}{\lvert\rho\rvert}$, and use Theorem~\ref{E-J-rho-main-thm} to write \begin{align*}
        &\frac{\glslink{E-f-n-k}{E_{\glslink{J-rho}{J_\rho}}(n, k)}}{\frac{1}{k!}\glslink{binomial}{\binom{n}{k}}} = \frac{(n-\glslink{abs-rho}{\lvert\rho\rvert})!}{n!} \cdot \frac{k!k!(n-k)!}{n!} \sum_{1 \le i \le 2\glslink{abs-rho}{\lvert\rho\rvert}+\glslink{m-j-rho}{m_1(\rho)}} \glslink{theta-f-i}{\theta_{\glslink{J-rho}{J_\rho}}^{(i)}} \cdot \frac{(\frac{i}{2})!}{(k-\frac{i}{2})!} \glslink{binomial}{\binom{n-\glslink{abs-rho}{\lvert\rho\rvert}+\frac{i}{2}}{k}} \\
        &= \sum_{1 \le i \le 2\glslink{abs-rho}{\lvert\rho\rvert}+\glslink{m-j-rho}{m_1(\rho)}} (\tfrac{i}{2})!\glslink{theta-f-i}{\theta_{\glslink{J-rho}{J_\rho}}^{(i)}} \cdot \frac{k!}{(k-\frac{i}{2})!} \cdot \frac{(n-\glslink{abs-rho}{\lvert\rho\rvert})!}{n!} \frac{(n-\glslink{abs-rho}{\lvert\rho\rvert}+\frac{i}{2})!}{n!} \cdot \frac{(n-k)!}{(n-k-\glslink{abs-rho}{\lvert\rho\rvert}+\frac{i}{2})!}.
    \end{align*} Then by Corollary~\ref{gamma-ratio-cor}, we may rewrite this as \begin{align*}
        &\sum_{1 \le i \le 2\glslink{abs-rho}{\lvert\rho\rvert}+\glslink{m-j-rho}{m_1(\rho)}} (\tfrac{i}{2})!\glslink{theta-f-i}{\theta_{\glslink{J-rho}{J_\rho}}^{(i)}} \cdot k^{\frac{i}{2}} \left(1+O\left(\tfrac{1}{k}\right)\right) \cdot n^{\frac{i}{2}-2\glslink{abs-rho}{\lvert\rho\rvert}} \left(1+O\left(\tfrac{1}{n}\right)\right) \cdot (n-k)^{\glslink{abs-rho}{\lvert\rho\rvert}-\frac{i}{2}} \left(1+O\left(\tfrac{1}{n-k}\right)\right) \\
        &= \left(\frac{n-k}{n^2}\right)^{\glslink{abs-rho}{\lvert\rho\rvert}}\sum_{1 \le i \le 2\glslink{abs-rho}{\lvert\rho\rvert}+\glslink{m-j-rho}{m_1(\rho)}} (\tfrac{i}{2})!\glslink{theta-f-i}{\theta_{\glslink{J-rho}{J_\rho}}^{(i)}} \cdot \left(\frac{kn}{n-k}\right)^{\frac{i}{2}} \left(1+O\left(\frac{1}{n}\right)\right),
    \end{align*} where in the last line we use the fact that $k, n-k \ge \epsilon n$ to combine the $O(\frac{1}{k})$ and $O(\frac{1}{n-k})$ into a single $O(\frac{1}{n})$ term. We may then use the fact that $\frac{kn}{n-k} \ge \frac{\epsilon n^2}{n} = \epsilon n$ to obtain \begin{align*}
        \frac{\glslink{E-f-n-k}{E_{\glslink{J-rho}{J_\rho}}(n, k)}}{\frac{1}{k!}\glslink{binomial}{\binom{n}{k}}} &= (\glslink{abs-rho}{\lvert\rho\rvert}+\tfrac{\glslink{m-j-rho}{m_1(\rho)}}{2})!\glslink{theta-f-i}{\theta_{\glslink{J-rho}{J_\rho}}^{(2\glslink{abs-rho}{\lvert\rho\rvert}+\glslink{m-j-rho}{m_1(\rho)})}} \cdot \left(\frac{n-k}{n^2}\right)^{\glslink{abs-rho}{\lvert\rho\rvert}}\left(\frac{kn}{n-k}\right)^{\glslink{abs-rho}{\lvert\rho\rvert}+\tfrac{\glslink{m-j-rho}{m_1(\rho)}}{2}} \left(1+O\left(\frac{1}{\sqrt{n}}\right)\right) \\
        &= \frac{(-1)^{\len(\rho)-\glslink{m-j-rho}{m_1(\rho)}}(\frac{\glslink{m-j-rho}{m_1(\rho)}}{2})!}{\glslink{z-rho}{z_\rho}} \cdot \left(\frac{nt}{1-t}\right)^{\frac{\glslink{m-j-rho}{m_1(\rho)}}{2}} t^{\glslink{abs-rho}{\lvert\rho\rvert}} \cdot \left(1 + O\left(\frac{1}{\sqrt{n}}\right)\right) \\
        &= \frac{(-1)^{\len(\rho)-\glslink{m-j-rho}{m_1(\rho)}}(\frac{\glslink{m-j-rho}{m_1(\rho)}}{2})!}{\glslink{z-rho}{z_\rho}} \cdot \left(\frac{nt}{1-t}\right)^{\frac{\glslink{m-j-rho}{m_1(\rho)}}{2}} t^{\glslink{abs-rho}{\lvert\rho\rvert}} + \left(\frac{t}{1-t}\right)^{\frac{\glslink{m-j-rho}{m_1(\rho)}}{2}} t^{\glslink{abs-rho}{\lvert\rho\rvert}} O\left(n^{\frac{\glslink{m-j-rho}{m_1(\rho)}-1}{2}}\right)
    \end{align*} Since $\epsilon \le t \le 1-\epsilon$ and $\epsilon \le 1-t \le 1-\epsilon$, this is equivalent to the statement in the lemma.
\end{proof}

\noindent
Using the above lemma, Theorem~\ref{conditional-expectation-thm} essentially reduces to expanding $\glslink{I-rho}{I_\rho}$ in the $\glslink{J-rho}{J_\rho}$-basis.

\begin{proof}[Proof of Theorem~\ref{conditional-expectation-thm}]
    Observe that for fixed $\pi$, the number of subsets $R \subseteq [n]$ such that $\pi(R) = R$ and $\pi\rvert_R$ has cycle type $\rho$ is given by $\glslink{I-rho}{I_\rho}(\pi)$. Consequently, we see that \begin{equation*}
        \E[\glslink{N-sigma}{N_{\glslink{id-k}{\id_k}}(\pi)} \mid \text{$\pi(R) = R$ and $\pi\rvert_R$ has cycle type $\rho$}] = \frac{\sum_{\pi \in \glslink{symmetric-group}{S_n}} \glslink{I-rho}{I_\rho}(\pi) \cdot \glslink{N-sigma}{N_{\glslink{id-k}{\id_k}}(\pi)}}{\sum_{\pi \in \glslink{symmetric-group}{S_n}} \glslink{I-rho}{I_\rho}(\pi)} = \glslink{z-rho}{z_\rho} \cdot \glslink{E-f-n-k}{E_{\glslink{I-rho}{I_\rho}}(n, k)},
    \end{equation*} where the last step follows by the definition of $\glslink{E-f-n-k}{E_f(n, k)}$, together with Proposition~\ref{I-rho-expectation-prop}. Meanwhile, it is a standard result that \begin{equation*}
        \E[\glslink{N-sigma}{N_{\glslink{id-k}{\id_k}}(\pi)}] = \frac{1}{k!} \glslink{binomial}{\binom{n}{k}}.
    \end{equation*} Applying Proposition~\ref{J-prop}(iii), we may thus write \begin{equation*}
        \frac{\E[\glslink{N-sigma}{N_{\glslink{id-k}{\id_k}}(\pi)} \mid \text{$\pi(R) = R$ and $\pi\rvert_R$ has cycle type $\rho$}]}{\E[\glslink{N-sigma}{N_{\glslink{id-k}{\id_k}}(\pi)}]} = \frac{\glslink{z-rho}{z_\rho} \cdot \glslink{E-f-n-k}{E_{\glslink{I-rho}{I_\rho}}(n, k)}}{\frac{1}{k!} \glslink{binomial}{\binom{n}{k}}} = \sum_{\rho' \sqcup \mu = \rho} \frac{\glslink{z-rho}{z_\rho}}{\glslink{z-rho}{z_\mu}} \cdot \frac{\glslink{E-f-n-k}{E_{\glslink{J-rho}{J_{\rho'}}}(n, k)}}{\frac{1}{k!} \glslink{binomial}{\binom{n}{k}}}.
    \end{equation*} Applying Lemma~\ref{E-J-rho-approx-lem}, the right hand side becomes \begin{align}
        &\sum_{\rho' \sqcup \mu = \rho} \frac{\glslink{z-rho}{z_\rho}}{\glslink{z-rho}{z_\mu}} \left(\frac{(-1)^{\glslink{len}{\len(\rho')}-\glslink{m-j-rho}{m_1(\rho')}}(\tfrac{\glslink{m-j-rho}{m_1(\rho')}}{2})!}{\glslink{z-rho}{z_{\rho'}}} \cdot \left(\frac{nt}{1-t}\right)^{\frac{\glslink{m-j-rho}{m_1(\rho')}}{2}} t^{\glslink{abs-rho}{\lvert\rho'\rvert}} + O\left(n^{\frac{m_1(\rho')-1}{2}}\right)\right) \notag \\
        &= \sum_{\substack{\rho' \sqcup \mu = \rho \\ \glslink{m-j-rho}{m_1(\rho')} = \glslink{m-j-rho}{m_1(\rho)}}} \frac{\glslink{z-rho}{z_\rho}}{\glslink{z-rho}{z_\mu}\glslink{z-rho}{z_{\rho'}}} \left((-1)^{\glslink{len}{\len(\rho')}-\glslink{m-j-rho}{m_1(\rho)}} \left(\frac{\glslink{m-j-rho}{m_1(\rho)}}{2}\right)! \cdot \left(\frac{nt}{1-t}\right)^{\frac{\glslink{m-j-rho}{m_1(\rho)}}{2}} t^{\glslink{abs-rho}{\lvert\rho'\rvert}}\right) + O\left(n^{\frac{\glslink{m-j-rho}{m_1(\rho)}-1}{2}}\right), \label{conditional-expectation-sum}
    \end{align} where the second line uses the fact that $t, 1-t \in [\epsilon, 1-\epsilon]$. Now, observe that any way of decomposing $\rho$ as $\rho' \sqcup \mu$ with $\glslink{m-j-rho}{m_1(\rho')} = \glslink{m-j-rho}{m_1(\rho)}$ is determined by values $0 \le m_j' \le \glslink{m-j-rho}{m_j(\rho)}$ for all $j \ge 2$, where $m_j'$ represents the number of parts of size $j$ in $\rho'$. Consequently, we can rewrite the sum in \eqref{conditional-expectation-sum} as a product \begin{multline*}
        \left(\frac{\glslink{m-j-rho}{m_1(\rho)}}{2}\right)! \cdot \left(\frac{nt^3}{1-t}\right)^{\frac{\glslink{m-j-rho}{m_1(\rho)}}{2}} \left(\prod_{j\ge2} \sum_{m_j'=0}^{\glslink{m-j-rho}{m_j(\rho)}} \glslink{binomial}{\binom{\glslink{m-j-rho}{m_j(\rho)}}{m_j'}} (-1)^{m_j'} t^{jm_j'}\right) \\
        = \left(\frac{\glslink{m-j-rho}{m_1(\rho)}}{2}\right)! \cdot \left(\frac{nt^3}{1-t}\right)^{\frac{\glslink{m-j-rho}{m_1(\rho)}}{2}} \left(\prod_{j\ge2} (1-t^j)^{\glslink{m-j-rho}{m_j(\rho)}}\right).
    \end{multline*} Together with the error term in \eqref{conditional-expectation-sum}, this is precisely the claimed formula.
\end{proof}

\appendix

\section{Properties of the \texorpdfstring{$\glslink{J-rho}{J_\rho}$}{Jρ}} \label{J-appendix}

\begin{proof}[Proof of Proposition~\ref{J-prop}]
    We start with claim (ii). Let $\glslink{Y-lam}{Y^\lambda} := \sum_{\rho\vdash\glslink{abs-rho}{\lvert\lambda\rvert}} \glslink{chi-rho-lambda}{\chi_\rho^\lambda} \cdot \glslink{I-rho}{I_\rho}$ be as in the statement of Theorem~\ref{macdonald-thm}, and given two partitions $\mu$ and $\lambda$, write $\glslink{prec}{\mu \prec \lambda}$ (resp.\@ $\glslink{prec-prime}{\mu \prec' \lambda}$) to mean that $\mu$ may be obtained from $\lambda$ by removing at most one box in each \glslink{column}{column} (resp.\@ row). Then by Theorem~\ref{macdonald-thm}, we may write \begin{equation}
        \sum_{\glslink{prec}{\mu \prec \lambda}} \glslink{X-lam}{X^\mu} = \sum_{\glslink{prec}{\mu \prec \lambda}} \sum_{\glslink{prec-prime}{\lambda' \prec' \mu}} (-1)^{\glslink{abs-rho}{\lvert\mu\rvert}-\glslink{abs-rho}{\lvert\lambda'\rvert}} \cdot \glslink{Y-lam}{Y^{\lambda'}} = \sum_{\lambda'} \glslink{Y-lam}{Y^{\lambda'}} \sum_{\lambda' \glslink{prec-prime}{\prec'} \mu \glslink{prec}{\prec} \lambda} (-1)^{\glslink{abs-rho}{\lvert\mu\rvert}-\glslink{abs-rho}{\lvert\lambda'\rvert}}. \label{X-mu-double-sum}
    \end{equation} Now, for a fixed partition $\lambda'$, consider the sum \begin{equation}
        \sum_{\lambda' \glslink{prec-prime}{\prec'} \mu \glslink{prec}{\prec} \lambda} (-1)^{\glslink{abs-rho}{\lvert\mu\rvert}-\glslink{abs-rho}{\lvert\lambda'\rvert}}. \label{mu-sum}
    \end{equation} For any $\mu$ appearing in the sum, the condition $\glslink{prec-prime}{\glslink{prec-prime}{\lambda' \prec' \mu}}$ is equivalent to the statement that $\glslink{lambda-i}{\mu_i} \in \{\glslink{lambda-i}{\lambda_i'}, \glslink{lambda-i}{\lambda_i'}+1\}$ for all $i$, while the condition $\glslink{prec}{\mu \prec \lambda}$ is equivalent to the condition that $\glslink{lambda-i}{\lambda_{i+1}} \le \glslink{lambda-i}{\mu_i} \le \glslink{lambda-i}{\lambda_i}$ for all $i$. Moreover, the latter condition alone is enough to ensure that any list of value $\glslink{lambda-i}{\mu_1},\glslink{lambda-i}{\mu_2},\dots$ is weakly decreasing, as it implies $\glslink{lambda-i}{\mu_{i+1}} \le \glslink{lambda-i}{\lambda_{i+1}} \le \glslink{lambda-i}{\mu_i}$ for all $i$. Consequently, we may factor the sum in \eqref{mu-sum} as \begin{equation*}
        \prod_{i\ge1} \sum_{\glslink{lambda-i}{\mu_i} \in \{\glslink{lambda-i}{\lambda_i'}, \glslink{lambda-i}{\lambda_i'}+1\} \cap [\glslink{lambda-i}{\lambda_{i+1}}, \glslink{lambda-i}{\lambda_i}]} (-1)^{\glslink{lambda-i}{\mu_i}-\glslink{lambda-i}{\lambda_i'}}.
    \end{equation*} Now, suppose this product is nonzero. Then for all $i \ge 1$, we must have either $\glslink{lambda-i}{\lambda_i'} = \glslink{lambda-i}{\lambda_i}$ or $\glslink{lambda-i}{\lambda_i'}+1 = \glslink{lambda-i}{\lambda_{i+1}}$. Consequently, if we have $\glslink{lambda-i}{\lambda_i'} < \glslink{lambda-i}{\lambda_i}$ for some $i$, then we must have $\glslink{lambda-i}{\lambda_i'}+1 = \glslink{lambda-i}{\lambda_{i+1}}$ and hence $\glslink{lambda-i}{\lambda_{i+1}'} \le \glslink{lambda-i}{\lambda_i'} = \glslink{lambda-i}{\lambda_{i+1}}-1 < \glslink{lambda-i}{\lambda_{i+1}}$. By induction, it then follows that $\glslink{lambda-i}{\lambda_{i'}'} < \glslink{lambda-i}{\lambda_{i'}}$ for all $i' \ge i$, contradicting the fact that the $\glslink{lambda-i}{\lambda_i}$ and $\glslink{lambda-i}{\lambda_i'}$ are both eventually zero. We conclude that \eqref{mu-sum} is only nonzero if $\lambda' = \lambda$, in which case it is equal to $1$, and so \eqref{X-mu-double-sum} becomes \begin{equation*}
        \glslink{Y-lam}{Y^\lambda} = \sum_{\glslink{prec}{\mu \prec \lambda}} \glslink{X-lam}{X^\mu}.
    \end{equation*} From here, observe that projecting each side of the equation $\glslink{Y-lam}{Y^\lambda} = \sum_{\rho \vdash \glslink{abs-rho}{\lvert\lambda\rvert}} \glslink{chi-rho-lambda}{\chi_\rho^\lambda} \cdot \glslink{I-rho}{I_\rho}$ onto $\glslink{Delta-Q-d}{\Delta_\Q^{\glslink{abs-rho}{\lvert\lambda\rvert}}}$ produces \begin{equation*}
        \glslink{X-lam}{X^\lambda} = \sum_{\rho \vdash \glslink{abs-rho}{\lvert\lambda\rvert}} \glslink{chi-rho-lambda}{\chi_\rho^\lambda} \cdot \glslink{J-rho}{J_\rho},
    \end{equation*} which is the first of the two claimed relations from (ii).
    The other relation $\glslink{J-rho}{J_\rho} = \frac{1}{\glslink{z-rho}{z_\rho}}\sum_{\lambda \vdash \glslink{abs-rho}{\lvert\rho\rvert}} \glslink{chi-rho-lambda}{\chi_\rho^\lambda} \cdot \glslink{X-lam}{X^\lambda}$ then follows immediately from orthogonality of symmetric group characters.

    From here, to prove claim (i), for any two partitions $\rho$ and $\rho'$ use the first of the two relations from (ii) along with the orthonormality of the $\glslink{X-lam}{X^\lambda}$ to write \begin{align*}
        \glslink{inner-product}{\langle \glslink{J-rho}{J_{\rho}}, \glslink{J-rho}{J_{\rho'}} \rangle} = \begin{cases}
            \frac{1}{\glslink{z-rho}{z_\rho} \glslink{z-rho}{z_{\rho'}}} \sum_{\lambda \vdash \glslink{abs-rho}{\lvert\rho\rvert}} \glslink{chi-rho-lambda}{\chi_\rho^\lambda} \cdot \glslink{chi-rho-lambda}{\chi_{\rho'}^\lambda} & \text{if $\glslink{abs-rho}{\lvert\rho\rvert} = \glslink{abs-rho}{\lvert\rho'\rvert}$,} \\
            0 & \text{otherwise.}
        \end{cases}
    \end{align*} If $\glslink{abs-rho}{\lvert\rho\vert} = \glslink{abs-rho}{\lvert\rho'\rvert}$, then orthogonality of irreducible symmetric group characters tells us that $\sum_{\lambda \vdash \glslink{abs-rho}{\lvert\rho\rvert}} \glslink{chi-rho-lambda}{\chi_\rho^\lambda} \cdot \glslink{chi-rho-lambda}{\chi_{\rho'}^\lambda}$ is equal to $\glslink{z-rho}{z_\rho}$ if $\rho' = \rho$ and $0$ otherwise.

    Finally, we address claim (iii). Fix a partition $\rho$, and start by using claim (ii) along with the definition of $\glslink{X-lam}{X^\lambda}$ to write \begin{align}
        \glslink{J-rho}{J_\rho} &= \frac{1}{\glslink{z-rho}{z_\rho}} \sum_{\lambda \vdash \glslink{abs-rho}{\lvert\rho\rvert}} \glslink{chi-rho-lambda}{\chi_\rho^\lambda} \cdot \glslink{X-lam}{X^\lambda} \notag \\
        &= \frac{1}{\glslink{z-rho}{z_\rho}} \sum_{\lambda \vdash \glslink{abs-rho}{\lvert\rho\rvert}} \glslink{chi-rho-lambda}{\chi_\rho^\lambda} \sum_{\glslink{prec-prime}{\mu \prec' \lambda}} (-1)^{\glslink{abs-rho}{\lvert\lambda\rvert}-\glslink{abs-rho}{\lvert\mu\rvert}} \sum_{\rho' \vdash \glslink{abs-rho}{\lvert\mu\rvert}} \glslink{chi-rho-lambda}{\chi_{\rho'}^\mu} \cdot \glslink{I-rho}{I_{\rho'}} \notag \\
        &= \frac{1}{\glslink{z-rho}{z_\rho}} \sum_{\substack{\rho',\,\glslink{abs-rho}{\lvert\rho'\rvert} \le \glslink{abs-rho}{\lvert\rho\rvert}}} \glslink{I-rho}{I_{\rho'}} \cdot (-1)^{\glslink{abs-rho}{\lvert\rho\rvert}-\glslink{abs-rho}{\lvert\rho'\rvert}}\sum_{\mu \vdash \glslink{abs-rho}{\lvert\rho'\rvert}} \glslink{chi-rho-lambda}{\chi_{\rho'}^\mu} \sum_{\substack{\lambda \vdash \glslink{abs-rho}{\lvert\rho\rvert} \\ \glslink{prec-prime}{\mu \prec' \lambda}}} \glslink{chi-rho-lambda}{\chi_\rho^\lambda}. \label{J-from-I-rho-prime}
    \end{align} Next, fix a partition $\rho'$ of size at most $\glslink{abs-rho}{\lvert\rho\rvert}$, and consider the quantity \begin{equation*}
        (-1)^{\glslink{abs-rho}{\lvert\rho\rvert}-\glslink{abs-rho}{\lvert\rho'\rvert}}\sum_{\mu \vdash \glslink{abs-rho}{\lvert\rho'\rvert}} \glslink{chi-rho-lambda}{\chi_{\rho'}^\mu} \sum_{\substack{\lambda \vdash \glslink{abs-rho}{\lvert\rho\rvert} \\ \glslink{prec-prime}{\mu \prec' \lambda}}} \glslink{chi-rho-lambda}{\chi_\rho^\lambda}.
    \end{equation*} Let $e_k$, $p_\rho$, and $s_\lambda$ denote the usual elementary, power sum, and Schur symmetric functions, and recall Pieri's rule for the $e_k$, which states that \begin{equation*}
        s_\mu e_k = \sum_{\substack{\lambda \vdash \glslink{abs-rho}{\lvert\mu\rvert}+k \\ \glslink{prec-prime}{\mu \prec' \lambda}}} s_\lambda
    \end{equation*} for all partitions $\mu$ and integers $k \ge 0$. Applying the relation $\glslink{hall-inner-product}{\langle s_\lambda, p_\rho \rangle} = \glslink{chi-rho-lambda}{\chi_\rho^\lambda}$, it follows that \begin{equation*}
        \glslink{hall-inner-product}{\langle s_\mu e_k, p_\rho \rangle} = \sum_{\substack{\lambda \vdash \glslink{abs-rho}{\lvert\mu\rvert}+k \\ \glslink{prec-prime}{\mu \prec' \lambda}}} \glslink{chi-rho-lambda}{\chi_\rho^\lambda}
    \end{equation*} for all partitions $\mu$ and all $k \ge 0$, and so we may write \begin{align*}
        (-1)^{\glslink{abs-rho}{\lvert\rho\rvert}-\glslink{abs-rho}{\lvert\rho'\rvert}}\sum_{\mu \vdash \glslink{abs-rho}{\lvert\rho'\rvert}} \glslink{chi-rho-lambda}{\chi_{\rho'}^\mu} \sum_{\substack{\lambda \vdash \glslink{abs-rho}{\lvert\rho\rvert} \\ \glslink{prec-prime}{\mu \prec' \lambda}}} \chi_\rho^\lambda &= (-1)^{\glslink{abs-rho}{\lvert\rho\rvert}-\glslink{abs-rho}{\lvert\rho'\rvert}}\sum_{\mu \vdash \glslink{abs-rho}{\lvert\rho'\rvert}} \glslink{chi-rho-lambda}{\chi_{\rho'}^\mu} \glslink{hall-inner-product}{\langle s_\mu e_{\glslink{abs-rho}{\lvert\rho\rvert}-\glslink{abs-rho}{\lvert\rho'\rvert}},\, p_\rho \rangle} \\
        &= (-1)^{\glslink{abs-rho}{\lvert\rho\rvert}-\glslink{abs-rho}{\lvert\rho'\rvert}} \glslink{hall-inner-product}{\left\langle \left(\sum_{\mu \vdash \glslink{abs-rho}{\lvert\rho'\rvert}} \glslink{chi-rho-lambda}{\chi_{\rho'}^\mu} s_\mu\right) e_{\glslink{abs-rho}{\lvert\rho\rvert}-\glslink{abs-rho}{\lvert\rho'\rvert}},\, p_\rho \right\rangle} \\
        &= (-1)^{\glslink{abs-rho}{\lvert\rho\rvert}-\glslink{abs-rho}{\lvert\rho'\rvert}} \glslink{hall-inner-product}{\left\langle p_{\rho'} e_{\glslink{abs-rho}{\lvert\rho\rvert}-\glslink{abs-rho}{\lvert\rho'\rvert}},\, p_\rho \right\rangle}.
    \end{align*} Now, let $p_{\rho'}^\perp$ denote the adjoint operator (with respect to the Hall inner product) to multiplication by $p_{\rho'}$, and recall the standard result that $p_k^\perp$ acts as $k \cdot \partial_{p_k}$ in the power sum basis. Consequently, we may write \begin{align*}
        (-1)^{\glslink{abs-rho}{\lvert\rho\rvert}-\glslink{abs-rho}{\lvert\rho'\rvert}} \glslink{hall-inner-product}{\left\langle p_{\rho'} e_{\glslink{abs-rho}{\lvert\rho\rvert}-\glslink{abs-rho}{\lvert\rho'\rvert}},\, p_\rho \right\rangle} &= (-1)^{\glslink{abs-rho}{\lvert\rho\rvert}-\glslink{abs-rho}{\lvert\rho'\rvert}} \glslink{hall-inner-product}{\left\langle e_{\glslink{abs-rho}{\lvert\rho\rvert}-\glslink{abs-rho}{\lvert\rho'\rvert}},\, p_{\rho'}^\perp p_\rho \right\rangle} \\
        &= (-1)^{\glslink{abs-rho}{\lvert\rho\rvert}-\glslink{abs-rho}{\lvert\rho'\rvert}} \glslink{hall-inner-product}{\left\langle e_{\glslink{abs-rho}{\lvert\rho\rvert}-\glslink{abs-rho}{\lvert\rho'\rvert}},\, \left(\prod_{j\ge1} (p_j^\perp)^{\glslink{m-j-rho}{m_j(\rho')}}\right) p_\rho \right\rangle} \\
        &= (-1)^{\glslink{abs-rho}{\lvert\rho\rvert}-\glslink{abs-rho}{\lvert\rho'\rvert}} \glslink{hall-inner-product}{\left\langle e_{\glslink{abs-rho}{\lvert\rho\rvert}-\glslink{abs-rho}{\lvert\rho'\rvert}},\, \left(\prod_{j\ge1} j^{\glslink{m-j-rho}{m_j(\rho')}} \partial_{p_j}^{\glslink{m-j-rho}{m_j(\rho')}}\right) p_\rho \right\rangle} \\
        &= \begin{cases}
            (-1)^{\glslink{abs-rho}{\lvert\mu\rvert}} \frac{\glslink{z-rho}{z_\rho}}{\glslink{z-rho}{z_\mu}} \glslink{hall-inner-product}{\left\langle e_{\glslink{abs-rho}{\lvert\mu\rvert}}, p_\mu \right\rangle} & \text{if $\rho = \rho' \sqcup \mu$ for some partition $\mu$,} \\
            0 & \text{otherwise.}
        \end{cases}
    \end{align*} If $\rho = \rho' \sqcup \mu$ for some $\mu$, we may then use the standard relation $e_k = \sum_{\lambda \vdash k} \frac{(-1)^{k-\glslink{len}{\len(\lambda)}}}{\glslink{z-rho}{z_\lambda}} \cdot p_\lambda$ to write \begin{equation*}
        (-1)^{\glslink{abs-rho}{\lvert\mu\rvert}} \frac{\glslink{z-rho}{z_\rho}}{\glslink{z-rho}{z_\mu}} \glslink{hall-inner-product}{\left\langle e_{\glslink{abs-rho}{\lvert\mu\rvert}}, p_\mu \right\rangle} = (-1)^{\glslink{abs-rho}{\lvert\mu\rvert}} \frac{\glslink{z-rho}{z_\rho}}{\glslink{z-rho}{z_\mu}} \cdot (-1)^{\glslink{abs-rho}{\lvert\mu\rvert}-\glslink{len}{\len(\mu)}} = (-1)^{\glslink{len}{\len(\mu)}} \cdot \frac{\glslink{z-rho}{z_\rho}}{\glslink{z-rho}{z_\mu}}.
    \end{equation*} Substituting this back into \eqref{J-from-I-rho-prime}, we conclude that \begin{equation*}
        \glslink{J-rho}{J_\rho} = \frac{1}{\glslink{z-rho}{z_\rho}} \sum_{\mu \sqcup \rho' = \rho} (-1)^{\glslink{len}{\len(\mu)}} \cdot \frac{\glslink{z-rho}{z_\rho}}{\glslink{z-rho}{z_\mu}} \cdot \glslink{I-rho}{I_{\rho'}} = \sum_{\mu \sqcup \rho' = \rho} \frac{(-1)^{\glslink{len}{\len(\mu)}}}{\glslink{z-rho}{z_\mu}} \cdot \glslink{I-rho}{I_{\rho'}}
    \end{equation*} as claimed. From here, the other claimed relation can be verified by first writing \begin{equation*}
        \sum_{\mu \sqcup \rho' = \rho} \frac{1}{\glslink{z-rho}{z_\mu}} \cdot \glslink{J-rho}{J_{\rho'}} = \sum_{\mu \sqcup \rho' = \rho} \frac{1}{\glslink{z-rho}{z_\mu}} \sum_{\substack{\mu' \sqcup \rho'' = \rho'}} \frac{(-1)^{\glslink{len}{\len(\mu')}}}{\glslink{z-rho}{z_{\mu'}}} \cdot \glslink{I-rho}{I_{\rho''}} = \sum_{\lambda \sqcup \rho' = \rho} \glslink{I-rho}{I_{\rho'}} \sum_{\mu \sqcup \mu' = \lambda} \frac{(-1)^{\glslink{len}{\len(\mu')}}}{\glslink{z-rho}{z_\mu} \glslink{z-rho}{z_{\mu'}}},
    \end{equation*} then noting that \begin{equation*}
        \sum_{\mu \sqcup \mu' = \lambda} \frac{(-1)^{\glslink{len}{\len(\mu')}}}{\glslink{z-rho}{z_\mu} \glslink{z-rho}{z_{\mu'}}} = \prod_{j \ge 1} \sum_{m_j=0}^{\glslink{m-j-rho}{m_j(\lambda)}} \frac{(-1)^{m_j}}{j^{\glslink{m-j-rho}{m_j(\lambda)}}m_j!(\glslink{m-j-rho}{m_j(\lambda)}-m_j)!} = \prod_{j\ge1} \frac{(1-1)^{\glslink{m-j-rho}{m_j(\lambda)}}}{j^{\glslink{m-j-rho}{m_j(\lambda)}}\glslink{m-j-rho}{m_j(\lambda)}!} = \begin{cases}
            1 & \text{if $\lambda = \varnothing$}, \\
            0 & \text{otherwise}.
        \end{cases}
    \end{equation*}
\end{proof}

\section{Efficient computation of \texorpdfstring{$\glslink{stitch-ell-m-t-nu-binom}{\binom{\rho}{\ell, m, t, \nu}}$}{(ρℓ,m,t,ν)}} \label{stitch-binom-appendix}

A key step in our computer implementation of the formulas underlying Theorem~\ref{E-I-rho-from-delta-thm} is the computation of the coefficients $\glslink{stitch-ell-m-t-nu-binom}{\binom{\rho}{\ell,m,t,\nu}}$. When $\rho$ is large, computing $\glslink{stitch-ell-m-t-nu-binom}{\binom{\rho}{\ell,m,t,\nu}}$ directly from the definition quickly becomes computationally infeasible. To remedy this, we instead make use of the following alternative method for computing the coefficients $\glslink{stitch-ell-m-t-nu-binom}{\binom{\rho}{\ell,m,t,\nu}}$ in our code.

\begin{proposition} \label{stitch-binom-prop}
    Let $\rho$ be a partition, and for each $j \ge 1$, define a polynomial $g_j \in \Z[L, M, T, N]$ by \begin{equation*}
        g_j := \begin{cases}
            M + N & \text{if $j = 1$,} \\
            \displaystyle \sum_{\substack{2\ell+m+t = j \\ \text{$\ell \ge 1$ or $t = 0$}}} \left(\glslink{binomial}{\binom{\ell+m}{m}}\glslink{binomial}{\binom{\ell+t-1}{t}} + \glslink{binomial}{\binom{\ell+m-1}{m}}\glslink{binomial}{\binom{\ell+t}{t}}\right) L^\ell M^m T^t & \text{otherwise.}
        \end{cases}
    \end{equation*} Then we have \begin{equation*}
        \sum_{2\ell+m+t+\nu=\glslink{abs-rho}{\lvert\rho\rvert}} \glslink{stitch-ell-m-t-nu-binom}{\binom{\rho}{\ell,m,t,\nu}} L^\ell M^m T^t N^\nu = \prod_{j\ge1} g_j^{\glslink{m-j-rho}{m_j(\rho)}}.
    \end{equation*}
\end{proposition}

\begin{proof}
    Fix a permutation $\pi_0$ with cycle type $\rho$, and observe that by the definition of $\glslink{stitch-ell-m-t-nu-binom}{\binom{\rho}{\ell,m,t,\nu}}$ we have \begin{equation*}
        \sum_{2\ell+m+t+\nu=\glslink{abs-rho}{\lvert\rho\rvert}} \glslink{stitch-ell-m-t-nu-binom}{\binom{\rho}{\ell,m,t,\nu}} L^\ell M^m T^t N^\nu = \sum_{\substack{\beta \in \glslink{calS-n}{\calS_{\glslink{abs-rho}{\lvert\rho\rvert}}},\, \beta \subseteq \pi_0 \\ \text{$\beta$ has no \glslink{proper-cycle}{proper cycles}}}} L^{\ell(\beta)} M^{m(\beta)} T^{t(\beta)} N^{\nu(\beta)},
    \end{equation*} where $\ell(\beta)$, $m(\beta)$, $t(\beta)$, and $\nu(\beta)$ denote the numbers of \glslink{proper-chain}{proper chains}, \glslink{gap}{gaps}, \glslink{jump}{jumps}, and \glslink{pole}{poles} of $\beta$, respectively. Moreover, observe that specifying a \glslink{substitch}{substitch} $\beta \subseteq \pi_0$ with no \glslink{proper-cycle}{proper cycles} is equivalent to specifying one \glslink{substitch}{substitch} $\beta_c$ with no \glslink{proper-cycle}{proper cycles} for each \glslink{cycle}{cycle} $c$ of $\pi_0$, in which case we have $\ell(\beta) = \sum_c \ell(\beta_c)$, $m(\beta) = \sum_c m(\beta_c)$, $t(\beta) = \sum_c t(\beta_c)$, and $\nu(\beta) = \sum_c \nu(\beta_c)$. Consequently, we may rewrite the right hand side as $\prod_{j\ge 1} \tilde{g}_j^{\glslink{m-j-rho}{m_j(\rho)}}$, where for each $j \ge 1$, we let $c_j \in S_j$ denote the $j$-cycle given by \begin{equation*}
        c_j(i) := \begin{cases}
            i+1 & \text{if $1 \le i \le j-1$,} \\
            1 & \text{if $i = j$}
        \end{cases}
    \end{equation*} and we define \begin{equation}
        \tilde{g}_j := \sum_{\substack{\beta \in \glslink{calS-n}{\calS_j},\, \beta \subseteq c_j \\ \text{$\beta$ has no \glslink{proper-cycle}{proper cycles}}}} L^{\ell(\beta)} M^{m(\beta)} T^{t(\beta)} N^{\nu(\beta)} \label{tilde-g-def}
    \end{equation}
    
    Now, when $j = 1$, there are exactly two possible values of $\beta$ appearing in the sum on the right hand side of \eqref{tilde-g-def}: $\beta = (\glslink{bracket-n}{[1]}, \varnothing)$ and $\beta = (\glslink{bracket-n}{[1]}, \{(1,1)\})$. The first case contributes a term $M$ and the second case contributes a term $N$, and so $\tilde{g}_1 = g_1$.
    
    Next, suppose that $j \ge 2$. We claim that $\tilde{g}_j = g_j$. Since $c_j$ has no fixed points, no \glslink{substitch}{substitch} $\beta \subseteq c_j$ can have a \glslink{pole}{pole}. Moreover, since a \glslink{jump}{jump} can only appear inside a \glslink{proper-chain}{proper chain} or a \glslink{proper-cycle}{proper cycle}, any $\beta$ appearing in the sum on the right hand side of \eqref{tilde-g-def} must satisfy either $\ell(\beta) \ge 1$ or $t(\beta) = 0$. Since the number of \glslink{column}{columns} of $\beta$ is equal to $2\ell(\beta)+m(\beta)+t(\beta)+\nu(\beta)$, we conclude that any $\beta$ appearing in the sum must satisfy $\nu(\beta) = 0$, $2\ell(\beta)+m(\beta)+t(\beta) = j$, and either $\ell(\beta) \ge 1$ or $t(\beta) = 0$. Fixing values $\ell, m, t \ge 0$ with $2\ell+m+t=j$ and either $\ell \ge 1$ or $t = 0$, it thus suffices to show that \begin{equation*}
        \#\left\{\beta \subsetneq c_j \,\middle|\, \begin{gathered}
            \text{$\beta$ has $\ell$ proper chains,} \\
            \text{$m$ gaps, $t$ jumps}
        \end{gathered}\right\} = \glslink{binomial}{\binom{\ell+m}{m}}\glslink{binomial}{\binom{\ell+t-1}{t}} + \glslink{binomial}{\binom{\ell+m-1}{m}}\glslink{binomial}{\binom{\ell+t}{t}}.
    \end{equation*}

    To count \glslink{substitch}{substitches} $\beta \subsetneq c_j$ with $\ell$ \glslink{proper-chain}{proper chains}, $m$ \glslink{gap}{gaps}, and $t$ \glslink{jump}{jumps}, we split into cases based on whether $\beta$ contains the \glslink{wire}{wire} $(j, 1)$. If $\beta$ does not contain the \glslink{wire}{wire} $(j, 1)$, then $\beta$ corresponds to a way of partitioning the ordered list $1, \dots, j$ into $\ell+m$ consecutive nonempty blocks, exactly $m$ of which are singletons. There are $\glslink{binomial}{\binom{\ell+m}{m}}$ possible ways to choose which $m$ blocks out of the total $\ell+m$ blocks are singletons. Once such a choice is fixed, there are $\glslink{binomial}{\binom{\ell+t-1}{t}}$ ways to choose the sizes of the remaining $\ell$ blocks while ensuring that each of these sizes is at least $2$ (equivalently, after assigning two \glslink{column}{columns} to each such block, distributing the remaining $t$ \glslink{column}{columns} among the $\ell$ blocks). We thus conclude that \begin{equation}
        \#\left\{\glslink{stitch-subseteq}{\beta \subsetneq c_j} \,\middle|\, \begin{gathered}
            \text{$\beta$ has $\ell$ \glslink{proper-chain}{proper chains}, $m$ \glslink{gap}{gaps}, $t$ \glslink{jump}{jumps}} \\
            \text{$\beta$ does not contain the \glslink{wire}{wire} $(j, 1)$}
        \end{gathered}\right\} = \glslink{binomial}{\binom{\ell+m}{m}} \glslink{binomial}{\binom{\ell+t-1}{t}} \label{cycle-substitches-formula-1}
    \end{equation} If instead $\beta$ does contain the \glslink{wire}{wire} $(j, 1)$, then this \glslink{wire}{wire} must be part of a \glslink{proper-chain}{proper chain} that ``wraps around'' from $j$ to $1$. If we remove the \glslink{wire}{wire} $(j, 1)$ from this \glslink{proper-chain}{proper chain}, splitting it into an initial chunk containing \glslink{column}{column} $1$ and a terminal chunk containing \glslink{column}{column} $j$, then we find that $\beta$ corresponds to a way of partitioning the ordered list $1, \dots, j$ into $\ell+m+1$ consecutive nonempty blocks so that, apart from the first and last blocks, exactly $m$ of the remaining blocks are singletons. There are $\glslink{binomial}{\binom{\ell+m-1}{m}}$ possible ways to choose which $m$ blocks out of the total $\ell+m-1$ ``middle'' blocks are singletons. Once such a choice is fixed, there are $\glslink{binomial}{\binom{\ell+t}{t}}$ ways to choose the sizes of the first block, the last block, and the remaining $\ell-1$ ``middle'' blocks while ensuring that the remaining middle blocks all have size at least $2$ (equivalently, after assigning one \glslink{column}{column} to each of the first and last blocks and two \glslink{column}{columns} to each of the other $\ell-1$ nonsingleton blocks, distributing the remaining $t$ \glslink{column}{columns} among these $\ell+1$ blocks). We thus conclude that \begin{equation}
        \#\left\{\beta \subsetneq c_j \,\middle|\, \begin{gathered}
            \text{$\beta$ has $\ell$ \glslink{proper-chain}{proper chains}, $m$ \glslink{gap}{gaps}, $t$ \glslink{jump}{jumps}} \\
            \text{$\beta$ contains the \glslink{wire}{wire} $(j, 1)$}
        \end{gathered}\right\} = \glslink{binomial}{\binom{\ell+m-1}{m}} \glslink{binomial}{\binom{\ell+t}{t}}. \label{cycle-substitches-formula-2}
    \end{equation} Combining \eqref{cycle-substitches-formula-1} and \eqref{cycle-substitches-formula-2} then gives the claimed formula for the total number of \glslink{substitch}{substitches} $\glslink{stitch-subseteq}{\beta \subsetneq c_j}$ with $\ell$ \glslink{proper-chain}{proper chains}, $m$ \glslink{gap}{gaps}, and $t$ \glslink{jump}{jumps}.
\end{proof}

\glsaddallunused
\setlength{\glsdescwidth}{0.6\linewidth}
\setlength{\glspagelistwidth}{0.2\linewidth}
\printglossary[type=symbols]
\setlength{\glsdescwidth}{0.4\linewidth}
\setlength{\glspagelistwidth}{0.2\linewidth}
\printglossary[type=terms]

\printbibliography

\end{document}